\documentclass[a4paper, 10pt]{book}
\usepackage[font=small,labelfont=bf]{caption}
\usepackage{a4wide}  
\newcommand{\comment}[1]{}

\usepackage{amsmath, amsthm, amssymb, amsfonts}
\usepackage[comma, numbers]{natbib}
\usepackage{tikz} 
\usepackage{enumitem}  \setenumerate[2]{label=(\Roman*)} \setenumerate[1]{label=(\roman*)} 
\usepackage{url}			
\usepackage{diagbox}		
\usepackage{longtable}		
\usepackage{cancel}			
\usepackage{rotating}		
\usepackage{comment}		
\usepackage{mathrsfs, yfonts}	
\usepackage{xfrac}
\usepackage{datetime}
\usepackage[normalem]{ulem}
\usepackage{bbm}
\usepackage[toc,page]{appendix}
\usepackage{float}
\usepackage{listings}
\usepackage{color}
\usepackage{setspace}
\usepackage{fancyhdr}
\usepackage{lscape}
\usepackage{tocloft}
\usepackage{multicol}
\usepackage{mathdots}

\usepackage[sf,sl,outermarks]{titlesec}

\definecolor{darkgray}{gray}{0.375}
\titleformat{\chapter}[display]
	{\normalfont\Large\filcenter\normalfont\bfseries}
	{\color{darkgray}\titlerule[1.25pt]\vspace{1.5pt}\titlerule[1.25pt]\vspace{-0.5pc}\normalfont\bfseries\Huge\sffamily Chapter \thechapter \\ \vspace{1pc}\titlerule[1.25pt]\vspace{1.5pt}\titlerule[1.25pt]}
	{-4.5pc}
	{\Huge\sffamily}
	[\vspace{-1ex}]
	
\titleformat{name=\chapter, numberless}[display]
	{\normalfont\Large\filcenter\normalfont\bfseries}
	{}
	{-10pc}
	{\Huge\sffamily}
	
\titleformat{\section}[display]
	{\normalfont\Large\normalfont\bfseries}
	{}
	{-3pc}
	{\Large\normalfont\bfseries$\blacksquare$\ \thesection\ }
	[\vspace{-1ex}]
	
\titleformat{name=\section, numberless}[display]
	{\normalfont\Large\normalfont\bfseries}
	{}
	{-3pc}
	{\Large\normalfont\bfseries }
	
\titleformat{\subsection}[display]
	{}
	{}
	{-3pc}
	{\large\normalfont\bfseries\S\ \thesubsection\ }

\titleformat{\subsubsection}[display]
	{}
	{}
	{-3pc}
	{\normalfont\bfseries\thesubsubsection\ }

\definecolor{dkgreen}{rgb}{0,0.6,0}
\definecolor{gray}{rgb}{0.5,0.5,0.5}
\definecolor{mauve}{rgb}{0.58,0,0.82}

\newdateformat{monthyeardate}{\monthname[\THEMONTH], \THEYEAR}

\newcommand{\verteq}{\rotatebox{90}{$\,=$}}
\newcommand{\nullset}{\varnothing}							
\newcommand{\setdiff}[2]{#1-#2}								
\newcommand{\powerset}[1]{\mathscr{P}(#1)}					
\newcommand{\set}[1]{\left\{#1\right\}}						
\newcommand{\tuple}[1]{\left(#1\right)}						
\newcommand{\gens}[1]{\left\langle#1\right\rangle}			
\newcommand{\card}[1]{\left|#1\right|}						

\newcommand{\ints}{\mathbb{Z}}								
\newcommand{\nonnegints}{\mathbb{Z}_{\geq 0}}				
\newcommand{\posints}{\mathbb{Z}_{> 0}}						
\newcommand{\intsge}[1]{\mathbb{Z}_{\geq #1}}				
\newcommand{\intsg}[1]{\mathbb{Z}_{> #1}}					
\newcommand{\intsmod}[1]{\mathbb{Z}/#1\mathbb{Z}}			

\newcommand{\bellnos}[1]{\mathscr{B}_{#1}} 		
\newcommand{\fibno}[1]{\mathscr{F}(#1)}			
\newcommand{\intparts}[1]{p(#1)} 				
\newcommand{\ointparts}[1]{op(#1)}				

\newcommand{\floor}[1]{\left\lfloor#1\right\rfloor}			
\newcommand{\ceil}[1]{\left\lceil#1\right\rceil}			

\newcommand{\congmod}[3]{#1 \equiv #2\ (\text{mod } #3)}			
\newcommand{\ncongmod}[3]{#1 \not\equiv #2\ (\text{mod } #3)}		
\newcommand{\tcongmod}[4]{#1 \equiv #2 \equiv #3\ (\text{mod } #4)}	

\newcommand{\congrel}[1]{#1^{\dagger}}

\newcommand{\planar}[1]{\mathbbmss{P}#1}			
\newcommand{\crossed}[1]{\mathbbmss{X}#1}			

\newcommand{\freesemigrp}[1]{#1^{+}} 			
\newcommand{\freemon}[1]{#1^{\star}} 				

\newcommand{\symgrpsymbol}{\mathcal{S}}			
\newcommand{\symgrp}[1]{\symgrpsymbol_{#1}}		

\newcommand{\pttnsymbol}{\mathcal{P}}
\newcommand{\pttnmon}[1]{\pttnsymbol_{#1}}			
\newcommand{\ppttnmon}[1]{\planar{\pttnmon{#1}}}	

\newcommand{\syminvmonsymbol}{\mathcal{IS}}
\newcommand{\syminvmon}[1]{\syminvmonsymbol_{#1}}	

\newcommand{\psyminvmon}[1]{\planar{\syminvmon{#1}}}	

\newcommand{\blockbijmonsymbol}{\mathcal{I}}
\newcommand{\blockbijmon}[1]{\blockbijmonsymbol_{#1}^*}	

\newcommand{\uniblockbijmonsymbol}{\mathfrak{F}}
\newcommand{\uniblockbijmon}[1]{\uniblockbijmonsymbol_{#1}}	
\newcommand{\puniblockbijmon}[1]{\planar{\uniblockbijmon{#1}}}

\newcommand{\brauersymbol}{\mathcal{B}}
\newcommand{\brauermon}[1]{\brauersymbol_{#1}}		
\newcommand{\jonessymbol}{\mathcal{J}}
\newcommand{\jonesmon}[1]{\jonessymbol_{#1}}			

\newcommand{\modmonsymbol}{\textgoth{M}}
\newcommand{\modmon}[2]{\modmonsymbol^{#1}_{#2}}		
\newcommand{\pmodmon}[2]{\planar{\modmon{#1}{#2}}}		

\newcommand{\apsissymbol}{\textgoth{A}}
\newcommand{\apsismon}[2]{\apsissymbol^{#1}_{#2}}		
\newcommand{\uniapmon}[1]{\apsismon{1}{#1}}				
\newcommand{\diapmon}[1]{\apsismon{2}{#1}}				
\newcommand{\triapmon}[1]{\apsismon{3}{#1}}				

\newcommand{\crossedmon}[1]{\mathbbmss{X}#1}				
\newcommand{\capsismon}[2]{\crossedmon{\apsismon{#1}{#2}}}	

\newcommand{\overbar}[1]{\mkern 1.5mu\overline{\mkern-1.5mu#1\mkern-1.5mu}\mkern 1.5mu}
\newcommand{\apsisbound}[2]{\overbar{\apsissymbol}^{#1}_{#2}}	
\newcommand{\capsisbound}[2]{\uline{\crossed{\apsissymbol}}^{#1}_{#2}} 
\newcommand{\pmodbound}[2]{\uline{\planar{\modmonsymbol}}^{#1}_{#2}} 

\newcommand{\symgen}[1]{\sigma_{#1}}					
\newcommand{\apgen}[2]{\textswab{a}^{#1}_{#2}}			
\newcommand{\transapgen}[1]{\textswab{t}_{#1}}		
\newcommand{\monapgen}[1]{\apgen{1}{#1}}				
\newcommand{\diapgen}[1]{\apgen{2}{#1}}					
\newcommand{\hookgen}[1]{\textswab{d}_{#1}}					
\newcommand{\triapgen}[1]{\apgen{3}{#1}}				
\newcommand{\ftransgen}[1]{\textswab{f}_{#1}}
\newcommand{\btransgen}[1]{\textswab{b}_{#1}}

\newcommand{\productgraph}[1]{\Gamma\big(#1\big)}

\newcommand{\diagclass}[1]{\textswab{D}_{#1}}
\newcommand{\diagrams}[1]{\textswab{D}_{#1}}
\newcommand{\pdiagrams}[1]{\planar{\diagrams{#1}}}
\newcommand{\diagequiv}{\sim_{\textswab{D}}}

\newcommand{\id}[1]{\text{id}_{#1}}					
\newcommand{\rank}[1]{\text{rank}(#1)}				

\newcommand{\upverts}[1]{U\tuple{#1}}			%
\newcommand{\lowverts}[1]{L\tuple{#1}}			%
\newcommand{\noupverts}[1]{u(#1)}					
\newcommand{\nolowverts}[1]{l(#1)}					
\newcommand{\trans}[1]{T(#1)}						

\newcommand{\uppat}[1]{U(#1)}						
\newcommand{\uptrans}[1]{UT(#1)}					
\newcommand{\upnontrans}[1]{UN(#1)}					
\newcommand{\lowpat}[1]{L(#1)}						
\newcommand{\lowtrans}[1]{LT(#1)}					
\newcommand{\lownontrans}[1]{LN(#1)}				

\newcommand{\Run}[2]{\mathfrak{R}^{#1}_{#2}}				
\newcommand{\run}[2]{\mathfrak{r}^{#1}_{#2}}				
\newcommand{\apmorph}[2]{\theta^{#1}_{#2}}					
\newcommand{\apmorphs}[2]{\Theta^{#1}_{#2}}					
\newcommand{\aptrans}[2]{\textswab{T}^{#1}_{#2}}			
\newcommand{\pmodtrans}[2]{\textswab{T}^{#1}_{#2}}			
\newcommand{\apupnontrans}[2]{\textswab{U}^{#1}_{#2}}		
\newcommand{\aplownontrans}[2]{\textswab{L}^{#1}_{#2}}		
\newcommand{\pmodupnontrans}[2]{\textswab{U}^{#1}_{#2}}		
\newcommand{\pmodlownontrans}[2]{\textswab{L}^{#1}_{#2}}	

\newcommand{\pfeastrans}[3]{\planar{T}(#1, #2, #3)}		
\newcommand{\pfeasnontrans}[1]{\planar{N}(#1)}		
\newcommand{\pfeasapsisnontrans}[1]{\planar{\overline{N}}(#1)} 
\newcommand{\pfeasmodminusapsisnontrans}[1]{\tuple{\planar{N} - \planar{\overline{N}}}(#1)} 
\newcommand{\feastrans}[3]{\crossedmon{T}(#1, #2, #3)}		
\newcommand{\feasnontrans}[1]{\crossedmon{N}(#1)}		
\newcommand{\feasapsisnontrans}[1]{\crossedmon{\overline{N}}(#1)} 

\newcommand{\tomega}[2]{\omega^{#1}_{#2}}			
\newcommand{\tupsilon}[2]{\upsilon^{#1}_{#2}}			
\newcommand{\tomegabar}[2]{\overbar{\omega}^{#1}_{#2}}	

\newcommand{\diagalg}[2][\xi]{\mathbb{C}^{#1}[#2]}			
\newcommand{\grpalg}[1]{\mathbb{C}[#1]}
\newcommand{\GL}[1]{\text{GL}_{#1}(\mathbb{C})}
\newcommand{\orthgrp}[1]{\mathcal{O}_{#1}(\mathbb{C})}

\newtheoremstyle{thesis}
{3pt}
{3pt}
{}
{}
{\bfseries}
{:}
{.5em}
{}

\theoremstyle{thesis}
\swapnumbers
\newtheorem{theorem}{Theorem}[section]
\newtheorem{definition}[theorem]{Definition}

\newtheorem{conjecture}[theorem]{Conjecture}
\newtheorem{example}[theorem]{Example}
\newtheorem{proposition}[theorem]{Proposition}
\newtheorem{corollary}[theorem]{Corollary}

\newtheorem{definition-proposition}[theorem]{Definition-Proposition}

\begin{document}
	\frontmatter
	\pagestyle{fancy}
    \fancyhf{}
    \renewcommand{\headrulewidth}{0pt}
    \renewcommand{\sectionmark}[1]{\markright{\textbf{Sec \thesection:}\ #1}} 
    \renewcommand{\chaptermark}[1]{\markboth{\textbf{\chaptername\ \thechapter:}\ #1}{}}
    \fancyhead[RO,LE]{}
    \fancyhead[RO]{\small\rightmark}
    \fancyhead[LE]{\small\leftmark}
    \fancyfoot[RO,LE]{\vspace{0.25cm}\textbf{\thepage}}
	
	\fancypagestyle{plain}{%
	\fancyhf{} 
	\fancyfoot[RO,LE]{\vspace{0.25cm}\textbf{\thepage}} 
	\renewcommand{\headrulewidth}{0pt}
	\renewcommand{\footrulewidth}{0pt}}
	
	\begin{titlepage}
	\begin{center}
		\vspace*{-0.25cm}
		
		{\huge\sffamily\bfseries The Planar Modular Partition Monoid}
		
		\vspace{0.7cm}
		
		{\large by}
		
		\vspace{0.1cm}
		
		{\Large Nicholas Charles Ham}
		
		\vspace{0.1cm}
				
		{\large BEc, BSc (Hons)}
		
		\vspace{0.7cm}
		
  		{\Large Supervised by Des FitzGerald \& Peter Jarvis}
  
  		\vspace{2.2cm}
  		\includegraphics[scale = 1.1]{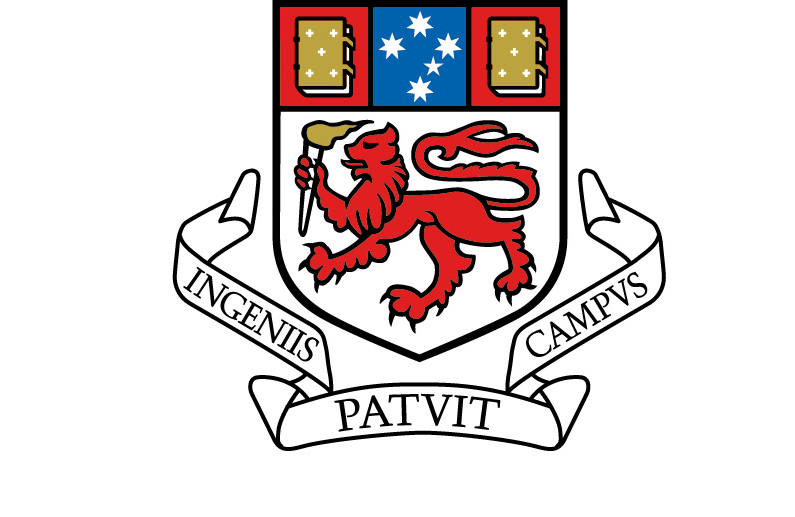}
  		
  		\vspace{1.8cm}
      
		{\large A thesis submitted in fulfilment of the requirements}
		
		\vspace{0.25cm}
		
		{\large for the Degree of Doctor of Philosophy}
		
		\vspace{0.7cm}
		
		{\large School of Physical Sciences}
		
		\vspace{0.25cm}
		
		{\large University of Tasmania}
		
		\vspace{0.7cm}
		
		{\large December, 2015}
		
		\vspace{1.5cm}
	\end{center}
\end{titlepage}
	
	\clearpage{\cleardoublepage}
\section*{Abstract}
\addcontentsline{toc}{chapter}{Abstract}
	The primary contribution of this thesis is to introduce and examine the planar modular partition monoid for parameters $m, k \in \posints$, which has simultaneously and independently generated interest from other researchers as outlined within.
	
	Our collective understanding of related monoids, in particular the Jones, Brauer and partition monoids, along with the algebras they generate, has heavily influenced the direction of research by a significant number of mathematicians and physicists. Examples include Schur-Weyl type dualities in representation theory along with Potts, ice-type and Andrew-Baxter-Forrester models from statistical mechanics, giving strong motivation for the planar modular partition monoid to be examined.
	
	The original results contained within this thesis relating to the planar modular partition monoid are: the establishment of generators; recurrence relations for the cardinality of the monoid; recurrence relations for the cardinality of Green's $\mathcal{R}$, $\mathcal{L}$ and $\mathcal{D}$ relations; and a conjecture on relations that appear to present the planar modular partition monoid when $m=2$. For diagram semigroups that are closed under vertical reflections, characterisations of Green's $\mathcal{R}$, $\mathcal{L}$ and $\mathcal{H}$ relations have previously been established using the upper and lower patterns of bipartitions. We give a characterisation of Green's $\mathcal{D}$ relation with a similar flavour for diagram semigroups that are closed under vertical reflections. 
	
	We also give a number of analogous results for the modular partition monoid, the monoid generated by replacing diapses with $m$-apses in the generators of the Jones monoid, later referred to as the $m$-apsis monoid, and the join of the $m$-apsis monoid with the symmetric group. 
	
	A further contribution of this thesis is a reasonably comprehensive exposition of the fundamentals of diagram semigroups, which have traditionally been approached from the representation theory side and have since blossomed into a thriving area of research in their own right. 
	\clearpage{\cleardoublepage}
\section*{Acknowledgements}
\addcontentsline{toc}{chapter}{Acknowledgements}
	First and foremost I wish to express my deepest gratitude towards my supervisors, Des FitzGerald and Peter Jarvis, who have provided insightful and encouraging feedback on my progress throughout my candidacy. Furthermore, following an idea put forward to Peter by Bertfried Fauser, they originally suggested that examining the monoid generated when replacing diapses with triapses in the generators of the Jones and Brauer monoids may be a suitable topic for a doctoral thesis, and indeed it revealed itself to be both a fruitful and enjoyable topic to get my hands dirty with. It also cannot be overstated how dramatically my ability to expose my own mathematical thoughts in a clear, fluent, unambiguous and grammatically correct manner has improved under the guidance of Des, of which I am especially thankful. 
	
	I am immensely grateful to James East who provided a number of insightful suggestions, references and counterexamples at various stages, and who has been beyond generous with his own time during the final year of my candidacy. I was further included and felt like I was treated as a respected colleague within some wider research circles throughout my candidacy - in regards to which I would particularly like to mention Igor Dolinka, James Mitchell, Attila Egri-Nagy, Athanasios Evangelou and Nicholas Loughlin.
	
	My time at the University of Tasmania, throughout undergrad and candidacy, was made all the more enjoyable and constructive both by opportunities afforded to me and by interactions with a number of people. In that regard I give thanks to: Jet Holloway and JJ Harrison for their help catching up to speed with mathematics classes after having turned up to university somewhat behind and only having decided to pick mathematics up in my second year; Mardi Dungey, Sarah Jennings and Simon Wotherspoon who first gave me the opportunity to obtain some research experience during my undergraduate years; Barry Gardner who very kindly took three units on top of his usual workload during our honours year; Karen Bradford, Michael Brideson, Kumudini Dharmadasa, Tracy Kostiuk and Jing Tian for the opportunity to gain teaching experience along with their assistance with doing so; Arwin Kahlon for putting up with my friendly banter about physicists; and Jeremy Sumner who took the time to provide a number of useful comments on my honours research despite it being somewhat outside his area of expertise.
	
	Identifying a number of the original results contained within this thesis was made substantially easier by the online encyclopaedia of integer sequences \cite{man:OEIS} and the free software package GAP \cite{man:GAP4} for computational discrete algebra. I would especially like to mention James Mitchell for the GAP semigroups package \cite{man:GAPsemigrps} which was not only useful for obtaining computational results, but was also used to generate the majority of the figures.
	
	For the first three and a half years of my candidacy I was supported financially by a Tasmanian Graduate Research Scholarship (186), and each year of my candidacy the Victorian Algebra Group provided me with financial assistance for travelling to the annual Victorian Algebra Conference. On that note, both Marcel Jackson and Brian Davey have been notably helpful towards ensuring all Australian post graduate algebra students feel both included and valued within the Australian algebra community.
	
	Finally, I would like to express my eternal gratitude towards my family, especially my parents and my brother, in particular for both their support and willingness to tolerate my own stubbornness, which would be almost impossible to understate at times though has also contributed significantly to many of my most cherished achievements.
	
	\vfill
	\line(1,0){454} 
	\vspace{-0.3cm}
	{\small To the best knowledge and belief of the author this thesis contains no material which has been accepted for a degree or diploma by the University of Tasmania or any other institution, no material previously published or written by another person, except by way of background information where duly acknowledged, and no material that infringes copyright.
	
	This thesis may be made available for loan, copying and communication in
	accordance with the Australian \textit{Copyright Act 1968}.
	
	\hspace{0.5cm} \begin{tikzpicture}[scale=0.6]
	\draw(12.5, 1.1) node {\includegraphics[scale=0.6]{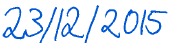}};
	
	\draw(0,0.5)--(7,0.5);
	\draw(9,0.5)--(16,0.5);
	
	\draw(1.25, 0) node {Signature};
	\draw(9.75, 0) node {Date};
\end{tikzpicture}

	}
		
	\vspace{0.3cm} \line(1,0){454} 
	
	\cleardoublepage
	\singlespacing
    \tableofcontents
    \doublespacing
    
    \mainmatter  
    
    \pagestyle{fancy}
    \fancyhf{}
    \renewcommand{\headrulewidth}{0.5pt}
    \renewcommand{\sectionmark}[1]{\markright{\textbf{Sec \thesection:}\ #1}} 
    \renewcommand{\chaptermark}[1]{\markboth{\textbf{\chaptername\ \thechapter:}\ #1}{}}
    \fancyhead[RO,LE]{}
    \fancyhead[RO]{\small\rightmark}
    \fancyhead[LE]{\small\leftmark}
    \fancyfoot[RO,LE]{\vspace{0.25cm}\textbf{\thepage}}
    	
    \clearpage{\pagestyle{empty}\cleardoublepage}
\chapter{Introduction}
	As with any increasingly vast body of human knowledge, to give an exhaustive historical account would quickly turn into its own body of work. Nevertheless the author has attempted to provide the reader with a historical account of where, how and why the diagram semigroups relevant to this thesis arose and have since blossomed into a thriving area of research in their own right.
	
	In 1896, Eliakim Moore \cite{art:Moore:SymGrpPresentation} established a presentation by way of generators and relations for the symmetric group $\symgrp{k}$. Later in 1927, Issai Schur \cite{art:Schur:Rationalen} identified that a duality, now famously known as Schur-Weyl duality, holds between the symmetric group algebra $\grpalg{\symgrp{k}}$ and the general linear group $\GL{n}$ (invertible $n\times{n}$ matrices over the complex field $\mathbb{C}$ with a fixed ordered basis). Ten years later Richard Brauer \cite{art:Brauer:AlgebrasConnectedWithSemisimpleContGrps} identified that an analogous duality holds between what is now-known as the Brauer algebra $\diagalg{\brauermon{k}}$ and the orthogonal group $\orthgrp{n}$, it would be almost impossible to overstate the influence that the identification of these dualities has since had on the direction of research by a significant number of mathematicians and physicists.  

	Then in 1971, Neville Temperley and Elliott Lieb \cite{proc:TempLieb:RltnsBetweenPercAndColProblem} identified what is now known as the Temperley-Lieb algebra, from which certain transfer matrices may be built. An understanding of the Temperley-Lieb algebra has played an important role in statistical mechanics, in particular, as noted by Ridout and Saint-Aubin \cite{art:Ridout:StandardModules}, with Potts models \cite{art:Martin:Potts}, ice-type models \cite{book:Baxter:ExactlySolvedModelsStatMech} and Andrew-Baxter-Forrester models \cite{art:Andrews:EightVertSOSModel}. The Temperley-Lieb algebra was later rediscovered independently by Jones \cite{art:Jones:IndexForSubfactors} in 1983, identifying what is now famously known as the Jones polynomial in knot theory. The planar diagrams in the Brauer monoid, typically referred to as either the Jones or Temperley-Lieb monoid, form a basis of the Temperley-Lieb algebra.
	
	The Partition monoid was then introduced independently by Jones \cite{proc:Jones:PottsModelAndSymGrp} and Martin \cite{art:Martin:RepresentationsGraphTLAlgebras} after both considered generalisations of the Temperley-Lieb algebra and Potts models from statistical mechanics. Jones \cite{proc:Jones:PottsModelAndSymGrp} did so when considering the centraliser of the tensor representation of the symmetric group $\symgrp{n}$ (when treated as the group of all $n\times{n}$ permutation matrices) in the endomorphism ring $\text{End}(V^{\otimes k})$, where $V$ is an $n$-dimensional vector space on which $\symgrp{n}$ acts diagonally on $V^{\otimes k}$. 
	
	Graham and Lehrer \cite{art:Graham:CellularAlgebras} introduced and provided a unified framework for obtaining a great deal of information about the representation theory of what are referred to as \textit{cellular algebras}. East \cite{art:East:CellularAlgebrasInvSemigrps} then showed that under certain compatibility assumptions, the semigroup algebra of an inverse semigroup is cellular if the group algebras of its maximal subgroups are cellular, while Guo and Xi \cite{art:Guo:CellularityTwistedSemigrpAlgebras} and Wilcox \cite{art:Wilcox:CullularityDiagramAlgebras} later examined cellularity for twisted semigroup algebras, allowing the cellularity of such algebras to be established using the theory on cellular algebras introduced in \cite{art:Graham:CellularAlgebras}.
	
	Motivated to generalise the aforementioned duality results from Jones \cite{proc:Jones:PottsModelAndSymGrp}, Tanabe \cite{art:Tanabe:CentralizerAlgebraOfUniRefGrp} examined the centraliser algebra in the endomorphism ring of $V ^{\otimes k}$ on which the unitary reflection groups of type $G(m, p, n)$ act diagonally, where $G(m, p, n)$ is an index-$p$ subgroup of $G(m, 1, n)$, and $G(m, 1, n)$ is a group of $n\times{n}$ matrices whose non-zero entries are $m$th roots of unity. Note that the unitary reflection groups were introduced by Shephard and Todd \cite{art:Shephard:FiniteUnitaryRefGrps}. A year later, FitzGerald and Leech \cite{art:Fitzgerald:DualSymInvMonsAndRepThry} independently examined the structure of both the monoid of block bijections and monoid of uniform block bijections $\uniblockbijmon{k}$, which has previously been described as the largest factorisable inverse submonoid of the dual symmetric inverse semigroup $\blockbijmon{k}$ (see \cite{art:Fitzgerald:PresentationMonUniBlockBij}). Kosuda \cite{art:Kosuda:CharacterizationPartAlgebras} then identified that in the case where $m > k$, the centraliser algebra considered by Tanabe corresponds to the monoid of uniform block bijections $\uniblockbijmon{k}$.
	
	Fitzgerald \cite{art:Fitzgerald:PresentationMonUniBlockBij} and Kosuda \cite{art:Kosuda:CharacterizationPartAlgebras} both independently gave the same presentation by way of generators and relations for the monoid of uniform block bijections $\uniblockbijmon{k}$. While the aforementioned presentation of $\uniblockbijmon{k}$ was economical in terms of the number of generators used, the inherent symmetry possessed by the monoid of uniform block bijections $\uniblockbijmon{k}$ was not reflected in the defining relations. With the addition of more generators, Kosuda \cite{art:Kosuda:CharacterizationPartAlgebras} and Kudryavtseva and Mazorchuk \cite{art:Kudryavtseva:PresentationsBrauerTypeMonoids} also independently gave an equivalent set of relations that better reflected the symmetry possessed, and using the same set of generators East \cite {art:East:FactorizableBraidMon} also arrived at an equivalent set of defining relations. Kosuda \cite{art:Kosuda:IrrepsofPartyAlgebra} further constructed a complete set of representatives of the irreducible representations of the algebra of uniform block bijections $\diagalg{\uniblockbijmon{k}}$.
	
	Kosuda \cite{art:Kosuda:StructurePartyAlgebraTypeB} later identified that in the case where $m \leq k$ then the centraliser algebra considered by Tanabe corresponds to the modular partition monoid $\modmon{m}{k}$. Furthermore, in \cite{art:Kosuda:StdExpForPartyAlgebra} Kosuda presented a candidate of the standard expression of the modular partition monoid $\modmon{m}{k}$, and in \cite{art:Kosuda:CharacterizationModularPartyAlgebra} established a presentation by way of generators and relations for the modular partition monoid $\modmon{m}{k}$ for all $1 \leq m \leq k$.
	
	A number of presentations by way of generators and relations have recently been given for diagram semigroups: Halverson and Ram \cite{art:Halverson:PartitionAlgebras} gave presentations for the planar partition monoid $\ppttnmon{k}$ and the partition monoid $\pttnmon{k}$, with their exposition since becoming reasonably famous as a survey-style treatment of the partition algebras; Kudryavtseva and Mazorchuk \cite{art:Kudryavtseva:PresentationsBrauerTypeMonoids} gave presentations for the Brauer monoid $\brauermon{k}$, the partial Brauer monoid and, as previously mentioned, the monoid of uniform block bijections $\uniblockbijmon{k}$; Posner, Hatch and Ly \cite{art:Posner:PresentationMotzkinMon} have suggested a presentation of the Motzkin monoid, though the paper appears to remain in preprint; and Easdown, East and FitzGerald \cite{art:Easdown:PresDualSymInvMonoid} gave a presentation for the dual symmetric inverse monoid, which may also be described as the monoid of block bijections.
	
	Motivated by Fauser, Jarvis and King's work on symmetric functions and generalised universal characters \cite{art:Fauser:NewBranchingRulesInducedByPlethysm, art:Fauser:PlethysmsReplicatedSchurFnctnsAndSeries, art:Fauser:HopfAlgebraStructureOfCharRingsOfClassicalGrps,  art:Fauser:RibbonHopfAlgebrasFromGroupCharRings}, Fauser put the idea forward to Jarvis that it may be fruitful to consider the consequences of replacing the diapses in the generators of the Jones and Brauer monoids with triapses. When the author sought suggestions from potential supervisors on possibly suitable doctoral research topics, Jarvis shared this idea from Fauser. The author felt like the topic had more than enough potential to be fun, which it most certainly has been, while also appearing to have the potential for applications which seems increasingly promising.
	
	During a visit to Leeds in July of 2015, James East had a discussion with Chwas Ahmed, who is currently one of Paul Martin's doctoral students. It turned out Ahmed has also been examining the planar modular partition monoid during his candidacy. In a preprint on the arXiv, Ahmed, Martin and Mazorchuk \cite{art:Ahmed:OnTheNoOfPrincipalIdealsInDTonalPartMons} study the number of principal ideals, along with the number of principal ideals generated by an element of fixed rank, of the (planar) modular partition monoids. 

	The remainder of the thesis is organised as follows. In Chapter \ref{chap:background} we: outline general notation, terminology and results that will be useful at various stages throughout; give a reasonably comprehensive construction of the fundamentals of diagram semigroups from scratch; and review known presentations for a number of contextually relevant diagram semigroups. In Chapter \ref{chap:characterisations} we: characterise the monoid generated by $m$-apsis generators, which we refer to as the $m$-apsis generated diagram monoid $\apsismon{m}{k}$; characterise the join of the $m$-apsis generated diagram monoid $\apsismon{m}{k}$ with the symmetric group $\symgrp{k}$, which we refer to as the crossed $m$-apsis generated diagram monoid $\capsismon{m}{k}$; and establish a minimal generating set for the planar modular partition monoid $\pmodmon{m}{k}$. In Chapter \ref{chap:cardinalities} we give recurrence relations for the cardinalities of the (crossed) $m$-apsises generated diagram monoids $\card{\apsismon{m}{k}}$ and $\card{\capsismon{m}{k}}$, along with the (planar) modular partition monoids $\card{\pmodmon{m}{k}}$ and $\card{\modmon{m}{k}}$. In Chapter \ref{chap:greensrltns} we begin by giving a characterisation of Green's $\mathcal{D}$ relation for diagram semigroups closed under the vertical flip involution $^*$, then count the number of Green's $\mathcal{D}$, $\mathcal{R}$ and $\mathcal{L}$ relations for the modular and planar modular partition monoids. In Chapter \ref{chap:presentations} we conjecture a presentation of the planar modular-$2$ partition monoid, establish a bound on reduced $\pmodmon{2}{k}$-words and conjecture a number of further results on the quest to identify candidates for $\pmodmon{2}{k}$-words in normal form.

    \clearpage{\pagestyle{plain}\cleardoublepage}
\chapter{Background} \label{chap:background}
    \section{Preliminary notation and terminology}
	The reader may find it of use to have precise definitions for the terminology and notation that will be used consistently throughout this thesis. Whenever it has been possible to do so the author has attempted to select terminology and notation that, at least in the author's opinion, is preferably descriptive, succinct, unambiguous and grammatically correct. Furthermore, the author has also attempted to select terminology and notation that is commonly used, at least in the author's experience, however the terminology and notation that the author has selected comes from a range of different sources rather than being identical to any particular reference on related literature.
	
	\subsection{Sets}
		\begin{definition}
			We denote by:
			\begin{enumerate}
				\item $\mathbb{R}$ the set of real numbers; and
				\item $\ints$ the set of integers, that is $\ints = \{\ldots, -2, -1, 0, 1, 2, \ldots\}$.
			\end{enumerate}
			
			Furthermore for each integer $z \in \ints$ and subset of the integers $A \subseteq \ints$, we denote by:
			\begin{enumerate}
				\item $\intsge{z}$ the set of integers greater than or equal to $z$, for example $\intsge{0} = \{0, 1, 2, \ldots\}$ is the set of non-negative integers;
				\item $\intsg{z}$ the set of integers greater than $z$, for example $\intsg{0} = \{1, 2, 3, \ldots\}$ is the set of positive integers; and
				\item $zA$ the set $\{za: a \in A\}$, for example $z\ints = \{\ldots, -2z, -z, 0, z, 2z, \ldots\}$, $z\posints = \{z, 2z, 3z, \ldots\}$ and $z\nonnegints = \{0, z, 2z, \ldots\}$.
			\end{enumerate}
		\end{definition}
		
		\begin{definition}
			For each real number $x \in \mathbb{R}$, we denote by $\floor{x}$ the greatest integer less than $x$, which is typically referred to as \textit{the floor of $x$}.
		\end{definition}
		
		\begin{definition}
			Given sets $A$ and $B$, we denote by:
			\begin{enumerate}
				\item $\setdiff{A}{B}$ the set containing elements of $A$ that are not elements of $B$ as $\setdiff{A}{B}$, that is $\setdiff{A}{B} = \{a \in A: a \notin B\}$;
				\item $A\times{B}$ the Cartesian product of $A$ and $B$, that is the set $\{(a, b): a \in A, b \in B\}$.
			\end{enumerate}
			
			Furthermore given $n \in \intsge{2}$ sets $A_1, \ldots, A_n$, we denote by:
			\begin{enumerate}
				\item $\bigcup_{i=1}^nA_i$ the union of $A_1, \ldots, A_n$; and
				\item $\bigcap_{i=1}^nA_i$ the intersection of $A_1, \ldots, A_n$.
			\end{enumerate}
		\end{definition}
		
	\subsection{Families of subsets}		
		\begin{definition} 
			Let $X$ be a set. A \textit{family of subsets of $X$}, which we will refer to more succinctly as a \textit{family of $X$} when we may do so without any contextual ambiguity, is a set $\mathcal{A}$ such that each element of $\mathcal{A}$ is a subset of $X$, that is $A \subseteq X$ for all $A \in \mathcal{A}$. The \textit{power set of $X$}, which we denote as $\powerset{X}$, is the family of all subsets of $X$.
		\end{definition}
		
		The reader should note it follows by definition that given a set $X$, $\powerset{\powerset{X}}$ is the family of families of subsets of $X$.
		
		\begin{definition}
			Let $X$ be a set and $\mathcal{A}$ be a family of subsets of $X$. We denote by:
			\begin{enumerate}
				\item $\bigcup_{A \in \mathcal{A}}A$ the union of the elements from $\mathcal{A}$; and
				\item $\bigcap_{A \in \mathcal{A}}A$ the intersection of the elements from $\mathcal{A}$.
			\end{enumerate}
			
			Furthermore let $Y$ be a set, $\mathcal{B}$ be a family of subsets of $Y$ and $f:\mathcal{A}\rightarrow\mathcal{B}$ be a function. We denote by:
			\begin{enumerate}
				\item $\bigcup_{A \in \mathcal{A}} f(A)$ the union of the elements from $\mathcal{A}$ mapped under $f$; and
				\item $\bigcap_{A \in \mathcal{A}} f(A)$ the intersection of the elements from $\mathcal{A}$ mapped under $f$.
			\end{enumerate} 
		\end{definition}
		
	\subsection{Binary relations}
	    \begin{definition}
	        Let $X$ be a set. A \textit{binary relation on $X$} is a subset $\mathcal{R}$ of $X\times X$.
	    \end{definition}
	
	    \begin{definition}
	        A binary relation $\mathcal{R}$ on a set $X$ is said to be:
	        \begin{enumerate}
	            \item \textit{reflexive} if $(x, x) \in \mathcal{R}$ for all $x \in X$;
	            \item \textit{symmetric} if $(x,y) \in \mathcal{R}$ implies $(y, x) \in \mathcal{R}$ for all $x, y \in X$;
	            \item \textit{antisymmetric} if $(x,y) \in \mathcal{R}$ implies $(y, x)\ \cancel{\in}\ \mathcal{R}$ for all distinct $x, y \in X$;
	            \item \textit{transitive} if $(x, y), (y, z) \in \mathcal{R}$ implies $(x, z) \in \mathcal{R}$ for all $x, y, z \in X$; and
	            \item \textit{total} if $(x, y) \in \mathcal{R}$ or $(y, x) \in \mathcal{R}$ for all $x , y \in X$.
	        \end{enumerate}
	        We refer to reflexivity, symmetry, antisymmetry, transitivity and totality as \textit{properties of binary relations}.
	    \end{definition}
    
	\subsection{Order relations}
	    \begin{definition}
	     A binary relation $\mathcal{R}$ on a set $X$ is referred to as:
	     \begin{enumerate}
	        	\item a \textit{pre-order} if it is reflexive and transitive;
	        	\item a \textit{partial order} if it is an antisymmetric pre-order; and
	        	\item a \textit{total order} if it is a total partial order.
	     \end{enumerate}
	    \end{definition}

		\begin{definition}
			Let $X$ be set partially ordered by $\leq$ and $A$ be a subset of $X$. An element $a \in A$ is referred to as:
			\begin{enumerate}
				\item \textit{the greatest element of $A$} if $a$ is greater than or equal to every element of $A$; 
				\item \textit{the least element of $A$} if $a$ is less than or equal to every element of $A$;
				\item \textit{a minimal element of $A$} if there does not exist an element of $A$ that is greater than $a$; and
				\item \textit{a maximal element of $A$} if there does not exist another element of $A$ that is less than $a$.
			\end{enumerate}
			
			Furthermore an element $x \in X$ is referred to as:
			\begin{enumerate}
				\item \textit{an upper bound of $A$} if every element of $A$ is less than or equal to $x$;
				\item \textit{a lower bound of $A$} if every element of $A$ is greater than or equal to $x$;
				\item \textit{the least upper bound of $A$} if it is the least element of the upper bounds of $A$; and
				\item \textit{the greatest lower bound of $A$} if it is the greatest element of the lower bounds of $A$.
			\end{enumerate}
		\end{definition}

	\subsection{Lattices}
		\begin{definition}
			A partial order $\leq$ on a set $X$ is referred to as a:
			\begin{enumerate}
				\item \textit{join-semilattice} if each two-element subset $\{x, x'\} \subseteq X$ has a greatest lower bound, which is typically referred to as \textit{the join of $x$ and $x'$}, and denoted as $x \vee x'$;
				\item \textit{meet-semilattice} if each two-element subset $\{x, x'\} \subseteq X$ has a least upper bound, which is typically referred to as \textit{the meet of $x$ and $x'$}, and denoted as $x \wedge x'$; and
				\item \textit{lattice} if it is both a join-semilattice and a meet-semilattice.
			\end{enumerate}
	    \end{definition} 
	    
	    For example, for each $k \in \posints$, the subsets of $\{1, \ldots, k\}$ are partially ordered by $\subseteq$, set union gives us a lattice join operation and set intersection gives us a lattice meet operation. Hence the subsets of $\{1, \ldots, k\}$ form a lattice under $\subseteq$.
	    
	    Let $X$ be a set. The reader should note it follows inductively that:
	    \begin{enumerate}
		    \item given a join-semilattice $\leq$ on a set $X$, every finite subset $A \subseteq X$ has a least upper bound which is referred to as \textit{the join of the elements of $A$}; and
		    \item given a meet-semilattice $\leq$ on a set $X$, every finite subset $A \subseteq X$ has a greatest lower bound which is referred to as \textit{the meet of the elements of $A$}.
	    \end{enumerate}
	    
   		\begin{definition}
   			A property $P$ of binary relations is referred to as \textit{closable} when for each binary relation $\mathcal{R}$ on a set $X$, while ordered by $\subseteq$, there exists a least element in the set of binary relations that contains $\mathcal{R}$ and satisfies property $P$. The least element, when it exists, is typically referred to as either the \textit{$P$ relation generated by $\mathcal{R}$} or the \textit{$P$ closure of $\mathcal{R}$}.
   		\end{definition}
   		
   		For example transitivity is a closable property of binary relations. The transitive closure of a binary relation $\mathcal{R}$ is stated in Proposition \ref{prop:transitiveclosure}.
   		
		\begin{proposition} \label{prop:transitiveclosure}
			If $\mathcal{R}$ is a binary relation on a set $X$ then the transitive closure of $\mathcal{R}$ is the set of all $(x, x') \in X\times{X}$ such that there exist $a_1, \ldots, a_n \in X$ satisfying $x = a_1, x' = a_n$, and $(a_1, a_2), \ldots, (a_{n-1}, a_n) \in \mathcal{R}$.
			
			\begin{proof}
				See page 337 of \cite{book:Lidl:AppliedAbstractAlgebra}.
			\end{proof}
		\end{proposition}
   		
   		Given a closable property of binary relations $P$:
   		\begin{enumerate}
	   		\item a lattice-meet operation is defined by taking the $P$ relation generated by the intersection of two binary relations with property $P$;
	   		\item a lattice-join operation is defined by taking the $P$ relation generated by the union of two binary relations with property $P$.
   		\end{enumerate} 
	    
	\subsection{Equivalence relations}
	    
	    \begin{definition}
	    	A binary relation $\mathcal{R}$ on a set $X$ is referred to as an \textit{equivalence relation} if it is a symmetric pre-order.
	    \end{definition}
	    
	    \begin{definition}
	        Let $\sim$ be an equivalence relation on a set $X$. \textit{The equivalence class of $x \in X$} is the set $\{y \in X: y \sim x\}$, which we denote as	 $[x]$. The set of equivalence classes $\{[x]: x \in X\}$ is denoted by $X/{\sim}$.
	    \end{definition}
	    
	    Note that reflexivity, symmetry and transitivity are all preserved under intersections of relations, hence equivalence relations are also preserved under intersections, which forms a lattice meet operation on equivalence relations.
	    
	    However while reflexivity and symmetry are preserved under unions of binary relations, transitivity may not be, consequently the union of two equivalence relations may not be an equivalence relation. Nevertheless we already noted that transitivity is a closable property of relations, hence giving us Proposition \ref{prop:joinofequivrltns}.
	    
	    \begin{proposition} \label{prop:joinofequivrltns}
	    	Let $\sim_1$ and $\sim_2$ be two equivalence relations. The smallest equivalence relation under set inclusion containing the union $\sim_1\cup\sim_2$ is equal to the transitive closure of the union $\sim_1\cup\sim_2$. Furthermore, this operation forms a lattice join operation on equivalence relations.\qed
	    \end{proposition}
	    
	    We shall adopt the style of indicating that there will be no explicit proof (as in the case above) by terminating the statement itself with $\square$.
		
	\subsection{The finer and coarser than pre-orders}
		\begin{definition}
			Let $X$ be a set and $\mathcal{A}, \mathcal{B}$ be families of subsets of $X$. We say that \textit{$\mathcal{A}$ is finer than $\mathcal{B}$} and that \textit{$\mathcal{B}$ is coarser than $\mathcal{A}$}, which we denote as $\mathcal{A} \preceq \mathcal{B}$, if for each $A \in \mathcal{A}$ there exist $B \in \mathcal{B}$ such that $A \subseteq B$. We denote the finer than relation as $\preceq$.
		\end{definition}
		
		Note that the finer than relation $\preceq$ is trivially reflexive and transitive, and hence pre-orders families of subsets. Further note that $\preceq$ is not a partial order, for example the families $\{\{1,2\}\}$ and $\{\{1, 2\}, \{1\}\}$ are finer than and coarser than each other.	
		
	\subsection{Set partitions}
		\begin{definition} \label{def:partition}
			Let $X$ be a set and $\mathcal{A}$ be a collection of non-empty subsets of $X$. We say that $X$:
			\begin{enumerate}
				\item is \textit{pairwise disjoint} if for each $A, A' \in \mathcal{A}$, $A \cap A'$ is the empty set;
				\item \textit{covers $X$}, or \textit{is a cover of $X$}, if $\bigcup_{A \in \mathcal{A}} A$ is equal to $X$; and
				\item \textit{partitions $X$}, or \textit{is a partition of $X$}, if $\mathcal{A}$ is pairwise disjoint and covers $X$.
			\end{enumerate}
		\end{definition}
		
	 	\begin{definition} \label{def:graphicaldepictionofpartitions}
		 	Let $k \in \nonnegints$. Each partition $\alpha$ of $\{1, \ldots, k\}$ may be \textit{graphically depicted} as follows: 
		 	\begin{enumerate}
		 		\item for each $j \in \{1, \ldots, k\}$, the vertex $j$ is depicted as the point $(j, 0) \in \mathbb{R}^2$; and
		 		\item lines connecting points in $\{1, \ldots, k\}\times\{0\}$ are drawn non-linearly below the horizontal line $\{(x, 0): x \in \mathbb{R}\}$ and between the two vertical lines $\{(1, y), (k, y): y \in \mathbb{R}\}$ such that the connected components form the blocks of $\alpha$.
		 	\end{enumerate}
		\end{definition}
		
\begin{figure}[!ht]
			\caption[ ]{Given $k=11$, as outlined in Definition \ref{def:graphicaldepictionofpartitions}, the partition $\{\{1,5\}, \{2,3\}, \{4\}, \{6,7,8\}, \{9\},$ $\{10,11\}\}$ may be graphically depicted as:}
			\label{fig:partitioneg}
			\vspace{5pt}
			\centering
			\input{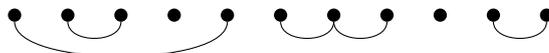}
		\end{figure}

		For example, Figure \ref{fig:partitioneg} illustrates a graphical depiction of the partition $\{\{1,5\}, \{2,3\}, \{4\}, \{6,7,8\}, \{9\},$ $\{10,11\}\}$.
		
		\begin{definition}
			The number of partitions of a set $X$ such that $|X| = k$ is referred to as \textit{the $k$-th Bell number}, which we denote as $\bellnos{k}$.
		\end{definition}
		
		Note the Bell numbers $\bellnos{k}$ are listed on the OEIS \cite{man:OEIS} as sequence $A000110$.

		When restricted to set partitions, the finer-than relation is additionally antisymmetric, hence set partitions are partially ordered by the finer-than relation.
		
		\begin{proposition} \label{prop:finestpartitioncoarserthanexistence}
			Let $X$ be a set. If $\mathcal{A}$ is a family of subsets of $X$ that covers $X$ then there exists a finest partition of $X$ that is coarser than $\mathcal{A}$.
			
			\begin{proof}
				Let $\sim$ be the relation on $\mathcal{A}$ where for each $A, A' \in \mathcal{A}$, $A \sim A'$ if there exist $k \in \intsge{2}$ and $A_1, \ldots, A_k \in \mathcal{A}$ such that $A_1 = A$, $A_k = A'$ and $A_i \cap A_{i+1} \neq \nullset$ for all $i \in \{1, \ldots, k-1\}$. It is trivially the case that $\sim$ is an equivalence relation and hence partitions the elements of $\mathcal{A}$. 
				
				Let $\mathcal{B} = \{\bigcup_{A' \in [A]}A': A \in \mathcal{A}\}$. It follows from how we defined $\sim$ that $\bigcup_{B \in \mathcal{B}}B = \bigcup_{A \in \mathcal{A}}A = X$ and for each $B, B' \in \mathcal{B}$, $B \cap B' = \nullset$, hence $\mathcal{B}$ partitions $X$. Furthermore for each $A \in \mathcal{A}$, $A \subseteq \bigcup_{A' \in [A]}A' \in \mathcal{B}$, hence $\mathcal{A}$ is finer than $\mathcal{B}$.
				
				Let $\mathcal{C}$ be a partition of $X$ that is coarser than $\mathcal{A}$.  Let $A, A' \in \mathcal{A}$ such that $A \sim A'$, therefore there exist $k \in \intsge{2}$, $A_1, \ldots, A_k \in \mathcal{A}$ and such that $A_1 = A$, $A_k = A'$ and $A_i \cap A_{i+1} \neq \nullset$ for all $i \in \{1, \ldots, k-1\}$. It follows from $\mathcal{C}$ being coarser than $\mathcal{A}$ that for each $i \in \{1, \ldots, k\}$, there exist $C_i \in \mathcal{C}$ such that $A_i \subseteq C_i$. Now for each $i \in \{1, \ldots, k-1\}$, $A_i \cap A_{i+1} \neq \nullset$ implies $C_i \cap C_{i+1} \neq \nullset$, which requires $C_i = C_{i+1}$ since $\mathcal{C}$ is a partition. Hence we must have $\bigcup_{A' \in [A]}A' \subseteq C_i \in \mathcal{C}$, giving that $\mathcal{B}$ is finer than $\mathcal{C}$. It follows from finer than partially ordering partitions that if $\mathcal{C}$ is also finer than $\mathcal{A}$ then $\mathcal{C} = \mathcal{A}$.
			\end{proof}
		\end{proposition}
		
		For example the finest partition coarser than $\{\{1, 2\}, \{2, 3\}\}$ is $\{\{1, 2, 3\}\}$.
		
		Given a set $X$, the intersection of two partitions of $X$ is also a partition of $X$, consequently the operation of taking the intersection of two partitions of $X$ forms a lattice meet operation. The union of two partitions of $X$ may not be a partition of $X$. However, the operation of taking the finest partition of $X$ coarser than the union of two partitions of $X$, which is well-defined as was established in Proposition \ref{prop:finestpartitioncoarserthanexistence}, forms a lattice join operation.
		
		\begin{proposition}
			Let $X$ be a set. With the above operations the set of partitions of $X$ forms a lattice. \qed
		\end{proposition}
		
		Furthermore, given an equivalence relation $\sim$ on $X$, the equivalence classes $\sfrac{X}{\sim}$ form a partition of $X$. Distinct equivalence relations form distinct partitions, each partition is formed by an equivalence relation, and the lattice operations match. Hence the lattice of partitions of $X$ and the lattice of equivalence relations on $X$ are isomorphic to each other.
		
		\begin{proposition}
			Let $X$ be a set. The lattice of partitions of $X$ and the lattice of equivalence relations on $X$ are isomorphic to each other. \qed
		\end{proposition}
		
		\begin{proposition}
			If $X$ is a set, $\mathcal{A}$ is a partition of $X$ and $Y \subseteq X$ then $\{A \cap Y: A \in \mathcal{A}\}$ is a partition of $Y$. \qed
		\end{proposition}
		
		\begin{definition}
			Let $X$ be a set, $\mathcal{A}$ be a partition of $X$ and $Y \subseteq X$. We refer to the partition $\{A \cap Y: A \in \mathcal{A}\}$ as \textit{the partition $\mathcal{A}$ restricted to $Y$} or \textit{the restriction to $Y$ of the partition $\mathcal{A}$}.
		\end{definition}
		
		\begin{definition} \label{def:partitiontypes}
			Let $k \in \nonnegints$ and $m \in \posints$. A partition $\alpha$ of $\{1, \ldots, k\}$ is referred to as:
			\begin{enumerate}
				\item \textit{non-crossing} or \textit{planar} if it may be graphically depicted, as in Definition \ref{def:graphicaldepictionofpartitions}, without connected components crossing; and
				\item \textit{m-divisible} if for each block $b \in \alpha$, the number of vertices in $b$ is divisible by $m$, that is $\card{b} \in m\posints$.
			\end{enumerate}
		\end{definition}
		 
		The lattice of non-crossing partitions of $\{1, \ldots, k\}$:
		\begin{enumerate}
			\item was first identified by \citet{art:Kreweras:PartitionsNonCrossing};
			\item is self-dual and admits a symmetric chain decomposition \citep{art:Simion:StructureLatticeNoncrossingPartitions}, with further structural and enumerative properties of chains in \cite{art:Edelman:ChainEnumerationAndNonCrossingPartitions}; and
			\item arise in the context of algebraic and geometric combinatorics, topological problems, questions in probability theory and mathematical biology \cite{art:Simion:NoncrossingPartitions}.
		\end{enumerate} 
		
		The poset of $m$-divisible non-crossing partitions of $\{1, \ldots, k\}$:
		\begin{enumerate}
			\item was first considered by \citet{art:Edelman:ChainEnumerationAndNonCrossingPartitions}; and
			\item has cardinality $\frac{\binom{(m+1)k}{k}}{mk+1}$ \citep{art:Arizmendi:StatsBlocksKDivisNCPs, art:Arizmendi:ProductsofRVsandKDivisNCPs}, \citep{art:Edelman:ChainEnumerationAndNonCrossingPartitions}, which is the Pfaff-Fuss-Catalan sequence \citep[page 347]{book:Graham:ConcreteMathematics}, \citep{art:Fuss:Solutio}.
		\end{enumerate}

		Note that a non-crossing partition may always be graphically depicted without individual connected components crossing, but some graphical depictions of planar partitions may require individual connected components to cross.
		
	\subsection{Modular arithmetic}
		\begin{definition}
			Given $m \in \nonnegints$ and $a, b \in \ints$, it is said that \textit{$a$ is congruent to $b$ modulo $m$}, which is typically denoted as $\congmod{a}{b}{m}$, if the difference $b-a$ between $a$ and $b$ is equal to a product of an integer and $m$, that is $b-a \in m\ints$.
		\end{definition}
		
		\begin{proposition}
			For each $a_1, a_2, b_1, b_2 \in \ints$ such that $\congmod{a_1}{b_1}{m}$ and $\congmod{a_2}{b_2}{m}$:
			\begin{enumerate}
				\item $\congmod{a_1 + a_2}{b_1 + b_2}{m}$; and
				\item $\congmod{a_1 - a_2}{b_1 - b_2}{m}$.
			\end{enumerate}
			
			\begin{proof}
				Since $\congmod{a_1}{b_1}{m}$ and $\congmod{a_2}{b_2}{m}$, there exists $x, y \in \ints$ such that $b_1 - a_1 = mx$ and $b_2 - a_2 = my$. It trivially follows that $(b_1 + b_2) - (a_1 + a_2) = m(x+y) \in m\ints$ and $(b_1 - b_2) - (a_1 - a_2) = m(x-y) \in m\ints$.
			\end{proof}
		\end{proposition}
		
		\begin{proposition}
			For each $m \in \posints$, congruence modulo $m$ is an equivalence relation on integers with $m$ equivalence classes $\{m\ints + 1, \ldots, m\ints + m\}$.
			
			\begin{proof}
				Let $a, b, c \in \ints$ such that $\congmod{a}{b}{m}$ and $\congmod{b}{c}{m}$, hence there exists $x, y \in \ints$ such that $b-a = mx$ and $c-b = my$. Now: 
				\begin{enumerate}
					\item $a-a = 0 \in m\ints$ giving reflexivity;
					\item $a-b = m(-x) \in m\ints$ giving symmetry; and
					\item $c-a = c-b+b-a = m(y+x) \in m\ints$ giving transitivity.
				\end{enumerate}
				Hence congruence modulo $m$ is an equivalence relation. Furthermore:
				\begin{enumerate}
					\item the integers congruent to $a$ modulo $m$ are trivially equal to $m\ints + a$; and
					\item there trivially exists a distinct $d \in \{1, \ldots, m\}$ such that $a \in m\ints + d$.
				\end{enumerate}
			\end{proof}
		\end{proposition}
		
		\begin{proposition} \label{prop:sumsofcongruentintegersarecongruent}
			For each $m, n \in \nonnegints$ and $a_1, \ldots, a_n, b_1, \ldots, b_n \in \ints$, if $\congmod{a_i}{b_i}{m}$ for all $i \in \{1, \ldots, n\}$ then $\congmod{\Sigma^n_{i=1}a_i}{\Sigma^n_{i=1}b_i}{m}$.
			
			\begin{proof}
				For each $i \in \{1, \ldots, n\}$, since $\congmod{a_i}{b_i}{m}$ there exists $x_i \in \ints$ such that $b_i - a_i = mx_i$. It trivially follows that $\Sigma^n_{i=1}b_i - \Sigma^n_{i=1}a_i = \Sigma^n_{i=1}(b_i - a_i) = \Sigma^n_{i=1}mx_i = m\Sigma^n_{i=1}x_i \in m\ints$.
			\end{proof}
		\end{proposition}
		
	\subsection{Integer partitions}
		\begin{definition} 
			For each $n \in \posints$ and $p_1, \ldots, p_n \in \posints$:
			\begin{enumerate}
				\item $\set{p_1, \ldots, p_n}$ is an \textit{integer partition of $\Sigma^n_{i=1}p_i$}; and
				\item $\tuple{p_1, \ldots, p_n}$ is an \textit{ordered integer partition of $\Sigma^n_{i=1}p_i$}.
			\end{enumerate}
		\end{definition}
		
		For example the five integer partitions of $4$ are $\set{1, 1, 1, 1}$, $\set{2, 1, 1}$, $\set{2, 2}$, $\set{3, 1}$ and $\{4\}$, while the eight ordered integer partitions of $4$ are $\tuple{1, 1, 1, 1}$, $\tuple{2, 1, 1}$, $\tuple{1, 2, 1}$, $\tuple{1, 1, 2}$, $\tuple{2, 2}$, $\tuple{3, 1}$, $\tuple{1, 3}$, $\tuple{4}$.
		
		\begin{definition} \label{def:numberofintpartsnotation}
			For each $k \in \posints$, we shall denote:
			\begin{enumerate}
				\item the number of integer partitions of $k$ as $\intparts{k}$;
				\item the number of integer partitions of $k$ into parts of size less than or equal to $m$ as $\intparts{m, k}$ (see Table \ref{table:nointparts} for example values);
				\item the number of ordered integer partitions of $k$ as $\ointparts{k}$; and
				\item the number of ordered integer partitions of $k$ into parts of size less than or equal to $m$ as $\ointparts{m, k}$ (see Table \ref{table:noorderedintparts} for example values).
			\end{enumerate}
		\end{definition}
		
        \begin{table}[!ht]
            \caption[ ]{Example values for $\intparts{m, k}$.}
            \label{table:nointparts}
            \centering
            \begin{tabular}{|c|rrrrrrrrrr|}
                \hline
                \diagbox{m}{k} & 1 & 2 & 3 & 4 & 5 & 6 & 7 & 8 & 9 & 10 \\
                \hline
				1 & 1 & 1 & 1 & 1 & 1 &  1 &  1 &  1 &  1 &  1 \\
				2 & 1 & 2 & 2 & 3 & 3 &  4 &  4 &  5 &  5 &  6 \\
				3 & 1 & 2 & 3 & 4 & 5 &  7 &  8 & 10 & 12 & 14 \\
				4 & 1 & 2 & 3 & 5 & 6 &  9 & 11 & 15 & 18 & 23 \\
				5 & 1 & 2 & 3 & 5 & 7 & 10 & 13 & 18 & 23 & 30 \\
                \hline
            \end{tabular}
        \end{table}
        
        \begin{table}[!ht]
            \caption[ ]{Example values for $\ointparts{m, k}$.}
            \label{table:noorderedintparts}
            \centering
            \begin{tabular}{|c|rrrrrrrrrr|}
                \hline
                \diagbox{m}{k} & 1 & 2 & 3 & 4 & 5 & 6 & 7 & 8 & 9 & 10 \\
                \hline
				1 & 1 & 1 & 1 & 1 &  1 &  1 &  1 &  1  &   1 &  1  \\
				2 & 1 & 2 & 3 & 5 &  8 & 13 & 21 & 34  &  55 & 89  \\
				3 & 1 & 2 & 4 & 7 & 13 & 24 & 44 & 81  & 149 & 274 \\
				4 & 1 & 2 & 4 & 8 & 15 & 29 & 56 & 108 & 208 & 401 \\
				5 & 1 & 2 & 4 & 8 & 16 & 31 & 61 & 120 & 236 & 464 \\
                \hline
            \end{tabular}
        \end{table}
		
		\begin{proposition}
			For each $m, k \in \posints$:
			\begin{enumerate}
				\item $\intparts{m, k} = \begin{cases}1 & m = 1, \\ \intparts{k, k} & m > k, \\ 1 + \intparts{m-1, k} & m = k, \\ \Sigma^{m-1}_{i=1}\intparts{i, k-i} & m < k; \end{cases}$
				\item $\intparts{k} = \intparts{k, k}$;
				\item $\ointparts{k} = 1 + \Sigma^{k-1}_{i=1}\ointparts{i} = 2^{\max\set{0, k-1}}$; and
				\item $\ointparts{m, k} = \begin{cases} \ointparts{k} & m \leq k, \\ \Sigma^{m}_{i=1}\ointparts{m, k-i} & m > k.\end{cases}$ \qed
			\end{enumerate}
		\end{proposition}
    \section{Semigroups and monoids}
	In this section we outline a number of elementary definitions and results from the theory of semigroups and monoids that will pop up at various times throughout this thesis. Vastly more comprehensive expositions of the fundamentals of semigroups may be found in \cite{book:Clifford:AlgThrySemigrpsVol1, book:Clifford:AlgThrySemigrpsVol2}, \cite{book:Higgins:TechniquesOfSemigrpThry} and \cite{ art:Howie:FundamentalsOfSemigrpThry}

	\begin{definition}
	    A set $S$ equipped with an associative binary operation, commonly denoted by $(s,s')\mapsto{ss'}$ for all $s, s' \in S$, is typically referred to as a \textit{semigroup}. 
	\end{definition}
	
	\begin{definition}
	      Given a semigroup $S$, for each $a \in S$ and $A, B \subseteq S$ we denote by:
	      \begin{enumerate}
		      \item $aB$ the set $\{ab: b \in B\}$;
		      \item $Ab$ the set $\{ab: a \in A\}$; and
		      \item $AB$ the set $\{ab: a \in A, b \in B\}$.
	      \end{enumerate}
	\end{definition}
		    
	\begin{definition}
	    Let $S$ be a semigroup. A subset $T \subseteq S$ is referred to as \textit{a subsemigroup of $S$} if $TT \subseteq T$. Note that, when convenient to do so, we denote such a relationship as $T \leq S$.
	\end{definition}
	
	\begin{definition}
		Given a semigroup $S$ and a subset $G \subseteq S$, \textit{the semigroup generated by $G$}, which is typically denoted by $\gens{G}$, is the smallest subsemigroup of $S$ for which $G$ is a subset, or equivalently the set of all finite combinations of elements from $G$ under the operation of $S$. 
	\end{definition}
		    
	\begin{definition}
	    	A semigroup $S$ is referred to as \textit{regular} if for each $a \in S$ there exists $b \in S$ such that $aba = a$ and $bab = b$.
	\end{definition}
	
	\begin{definition}
		Given a semigroup $S$, an element $e \in S$ is referred to as \textit{the identity of $S$} if for each $s \in S$, $se = s = es$.
	\end{definition}
	
	Given a semigroup $S$ with identity element $e \in S$, Proposition \ref{prop:identitiesareunique} justifies referring to $e$ as \textit{the identity of $S$} rather than as \textit{an identity of $S$}.
	
	\begin{proposition} \label{prop:identitiesareunique}
		Given a semigroup $S$, if $e, e' \in S$ are identity elements then $e = e'$.
		
		\begin{proof}
			$e = ee' = e'$.
		\end{proof}
	\end{proposition}
	
	\begin{definition}
		A semigroup $S$ containing an identity element $e \in S$, that is so that $es = s = se$ for all $s \in S$, is typically referred to as a \textit{monoid}.
	\end{definition}
	
	\begin{definition}
		Given a semigroup $S$, $S^1$ typically denotes:
		\begin{enumerate}
			\item $S$ when $S$ is already a monoid; and
			\item $S\cup\{e\}$ such that $se = s = es$ for all $s \in S$ when $S$ is not already a monoid.
		\end{enumerate}
	\end{definition}
	
	\begin{proposition}
		If $S$ is a semigroup then $S^1$ is a monoid. \qed
	\end{proposition}
	
	\begin{definition}
		Given a monoid $M$ with identity $e \in M$, a subset $N \subseteq M$ is referred to as \textit{a submonoid of $M$} if $N$ is a subsemigroup of $M$ and $e \in N$.
	\end{definition}
	
	Note it is possible for a monoid $M$ to contain a distinct monoid $N$ that is not a submonoid of $M$. For example given a monoid $M$ with identity $e$, form $M \cup \{f\}$ such that $fm=m=mf$ for all $m \in M$, then $M$ and $M \cup \{f\}$ are both monoids, and $M$ is a subsemigroup of $M \cup \{f\}$, while $M$ is not a submonoid of $M \cup \{f\}$.
	
	\subsection{Idempotents}
		\begin{definition}
			An \textit{idempotent} of a semigroup $S$ is an element $s \in S$ satisfying $s^2 = s$.
		\end{definition}
		
		We denote the idempotents of a semigroup $S$ as $E(S)$. 
		
	\subsection{Ideals}
	    \begin{definition}
	        A subset $I \subseteq S$ is referred to as:
	        \begin{enumerate}
	            \item a \textit{left ideal} if $I$ is closed when multiplying on the left by elements of $S$, that is if $SI \subseteq I$;
	            \item a \textit{right ideal} if $I$ is closed when multiplying on the right by elements of $S$, that is if $IS \subseteq I$; and
	            \item an \textit{ideal} if it is both a left and a right ideal, that is if $SI, IS \subseteq I$.
	        \end{enumerate}
	    \end{definition}
	    
	    
	    \begin{proposition}
	        For each $a \in S$:
	        \begin{enumerate}
	            \item $Sa \cup \{a\}$ is a left ideal;
	            \item $aS \cup \{a\}$ is a right ideal; and
	            \item $SaS \cup Sa \cup aS \cup \{a\}$ is an ideal.
	        \end{enumerate}
		        
	        \begin{proof}
	        	For each $a \in S$ we have:
	            \begin{enumerate}
	            	\item $S(Sa \cup \{a\}) = SSa \cup Sa \subseteq Sa \cup \{a\}$;
	            	\item $(aS \cup \{a\})S = a(SS) \cup aS \subseteq Sa \cup \{a\}$; and
	            	\item $S(SaS \cup Sa \cup aS \cup \{a\}) = SSaS \cup SSa \cup SaS \cup Sa \subseteq SaS \cup Sa \cup aS \cup \{a\}$ and \newline $(SaS \cup Sa \cup aS \cup \{a\})S = SaSS \cup SaS \cup aSS \cup aS \subseteq SaS \cup Sa \cup aS \cup \{a\}$.
	            \end{enumerate}
	        \end{proof}
	    \end{proposition}
	    
	    \begin{definition}
	        Let $S$ be a semigroup. For each $a \in S$:
	        \begin{enumerate}
		        \item $Sa \cup \{a\}$ is referred to as \textit{the principal left ideal generated by $a$};
		        \item $aS \cup \{a\}$ is referred to as \textit{the principal right ideal generated by $a$}; and
		        \item $SaS \cup Sa \cup aS \cup \{a\}$ is referred to as \textit{the principal ideal generated by $a$.}
	        \end{enumerate}
	    \end{definition}
	    
	    
	
	\subsection{Regular $^*$-semigroups}
	    \begin{definition}
	        Let $S$ be a semigroup. An \textit{involution} is a unary operation $^*:S\rightarrow{S}$ such that $^*$ is:
	        \begin{enumerate}
	            \item its own inverse, that is $(s^*)^* = s$ for all $s \in S$; and
	            \item an anti-automorphism, that is $(st)^* = t^*s^*$ for all $s, t \in S$.
	        \end{enumerate}
	    \end{definition}
	    
	    \begin{definition}
	    	A semigroup $S$ equipped with an involution $^*:S\rightarrow{S}$ is referred to as a \textit{regular $^*$-semigroup} if for each $s \in S$, $ss^*s = s$.
	    \end{definition}
	    
	    Regular $^*$-semigroups were introduced by \citet{proceeding:Nordahl:regularstarsemigroups}.
	    
		\begin{proposition}
			If a submonoid $T$ of a regular $^*$-semigroup $S$ is closed under the involution $^*$, that is $T^* = T$, then $T$ is also a regular $^*$-semigroup.
			
			\begin{proof}
				It trivially follows by definition that:
				\begin{enumerate}
					\item for each $a \in T$ we have $a^* \in T$, ${a^*}^* = a$ and $aa^*a = a$; and
					\item for each $a, b \in T$ we have $a^*, b^*, b^*a^* \in T$ and $(ab)^* = b^*a^*$.
				\end{enumerate}
			\end{proof}
		\end{proposition}
	    
	    \begin{definition}
	        Let $S$ be a regular $^*$-semigroup. An element $p \in S$ is referred to as a \textit{projection} if it satisfies $p^2 = p = p^*$.
	    \end{definition}
	    
	    Note that all projections are idempotents, which trivially follows from the condition for a projection containing the condition for idempotency. 
	    
	    \begin{proposition}
	        Let $S$ be a regular $^*$-semigroup. The projections, and more inclusively the idempotents, of $S$ are partially ordered by $p \leq q$ if and only if $pq = p = qp$ for all $p, q \in S$.
	        \begin{proof}
	        	Clearly $\leq$ is reflexive and antisymmetric leaving transitivity. Let $p, q, r \in S$ such that $p \leq q \leq r$. Then $pr = (pq)r = p(qr) = pq = p = qp = rqp = rp$.
	        \end{proof}
	    \end{proposition}
	    
	    \begin{proposition}
	    	Given a regular $^*$-semigroup $S$:
	        \begin{enumerate}
	            \item for each $s \in S$, $ss^*$ and $s^*s$ are projections; 
	            \item every projection can be written as $ss^*$ and $s^*s$ for some $s \in S$;
	            \item the set of idempotents $E(S)$ is closed under involution; and
	            \item every idempotent is the product of two projections.
	        \end{enumerate}
	        
	        \begin{proof}
	        	\begin{enumerate}
	        		\item For each $s \in S$, $(ss^*)(ss^*) = (ss^*s)s^* = ss^* = (s^*)^*s^* = (ss^*)^*$, similarly $(s^*s)(s^*s) = s^*s = (s^*s)^*$;
	        		\item For each $p \in S$ such that $p^2 = p = p^*$, $p = p^2 = pp^* = p^*p$;
	        		\item For each idempotent $e \in E(S)$, $e^*e^* = (ee)^* = e^*$, hence $e^* \in E(S)$; and
	        		\item For each idempotent $e \in E(S)$, $(ee^*)^*(e^*e) = ee^*e^*e = ee^*e = e$.
	        	\end{enumerate}
	        \end{proof}
	    \end{proposition}
	
	\subsection{Green's relations}
		Green's relations, introduced by James Alexander Green in \cite{art:Green:OnTheStructureOfSemigroups}, are five equivalence relations that partition the elements of a semigroup with respect to the principal ideals that may be generated. 
	
	    \begin{definition}
	        Given a semigroup $S$, for each $a, b \in S$:
	        \begin{enumerate}
		        \item $(a, b) \in \leq_{\mathcal{R}}$ if and only if $aS \cup \{a\} \subseteq bS \cup \{b\}$; and 
		        \item $(a, b) \in \leq_{\mathcal{L}}$ if and only if $Sa \cup \{a\} \subseteq Sb \cup \{b\}$.
	        \end{enumerate}
	    \end{definition}
		    
	    \begin{proposition}
	        $\leq_{\mathcal{R}}$ and $\leq_{\mathcal{L}}$ are pre-orders on $S$.
		        
	        \begin{proof}
	            Reflexivity is obvious and transitivity follows from $\subseteq$ being transitive on sets.
	        \end{proof}
	    \end{proposition}
		
	    \begin{definition}
	        Given a semigroup $S$, \textit{Green's relations on $S$} are defined as follows:
	        \begin{enumerate}
	            \item $\mathcal{R} = \leq_{\mathcal{R}} \cap \leq_{\mathcal{R}}^{-1} = \{(a, b) \in S \times S: aS \cup \{a\} = bS \cup \{b\}\}$;
	            \item $\mathcal{L} = \leq_{\mathcal{L}} \cap \leq_{\mathcal{L}}^{-1} = \{(a, b) \in S \times S: Sa \cup \{a\} = Sb \cup \{b\}\}$;
	            \item $\mathcal{H} = \mathcal{L} \cap \mathcal{R}$;
	            \item $\mathcal{D} = \{(a, b) \in S\times S: \text{there exists } c \in S \text{ such that } a \mathcal{L} c \text{ and } c \mathcal{R} b\}$; and
	            \item $\mathcal{J} = \{(a, b) \in S\times S: SaS \cup \{a\} = SbS \cup \{b\}\}$.
	        \end{enumerate}
	    \end{definition}
	    
		\begin{proposition}
		    Given a semigroup $S$, 
		    \begin{enumerate}
		    \item Green's $\mathcal{R}$, $\mathcal{L}$, $\mathcal{H}$, $\mathcal{D}$ and $\mathcal{J}$ relations are equivalence relations;
		    \item $\mathcal{D} = \mathcal{R} \circ \mathcal{L} = \mathcal{L} \circ \mathcal{R} = \mathcal{R} \vee \mathcal{L}$; and
		    \item if $S$ is finite then Green's $\mathcal{D}$ and $\mathcal{J}$ relations are equal. \qed
		    \end{enumerate}
	    \end{proposition}

		\begin{proposition} \label{prop:greensrltnsforregularsubmons} 
			Let $T$ be a subsemigroup of a semigroup $S$. If $T$ is regular then Green's $\mathcal{R}^T$, $\mathcal{L}^T$ and $\mathcal{H}^T$ relations on $T$ are equal to the restrictions to $T$ of Green's $\mathcal{R}^S$, $\mathcal{L}^S$ and $\mathcal{H}^S$ relations on $S$ respectively, that is:
			\begin{enumerate}
				\item $\mathcal{R}^T$ = $\mathcal{R}^S \cap (T\times{T})$;
				\item $\mathcal{L}^T$ = $\mathcal{L}^S \cap (T\times{T})$; and
				\item $\mathcal{H}^T$ = $\mathcal{H}^S \cap (T\times{T})$.
			\end{enumerate}
			
			\begin{proof}
				Proposition 2.4.2 in \cite{art:Howie:FundamentalsOfSemigrpThry} (also in \cite{art:Hall:CongsAndGreensRltnsOnRegSemigrps}).
			\end{proof}
		\end{proposition}		    
		    
	        
	        
	
		    
		    
	\subsection{Green's relations on regular $^*$-semigroups}
	    \begin{proposition}
	    	Given a regular $^*$-semigroup $S$, for each $a, b \in S$:
	        \begin{enumerate}
	        	\item $(a, b) \in \leq_{\mathcal{R}}$ if and only if $aa^* \leq bb^*$;
	            \item $(a, b) \in \leq_{\mathcal{L}}$ if and only if $a^*a \leq b^*b$;
	            \item $(a, b) \in \mathcal{R}$ if and only if $aa^* = bb^*$;
	            \item $(a, b) \in \mathcal{L}$ if and only if $a^*a = b^*b$; and
	            \item $(a, b) \in \mathcal{D}$ if and only if there exists $c \in S$ such that $a^*a = c^*c$ and $cc^* = bb^*$.
	        \end{enumerate}
	        
	        \begin{proof}
		        Suppose $(a, b) \in \leq_{\mathcal{R}}$, and hence that there exist $s \in S^1$ such that $as = b$. Then $aa^*(b)b^* = (aa^*a)sb^* = (as)b^* = bb^* = b(s^*a^*) = bs^*(a^*aa^*) = b(b^*)aa^*$, giving us $aa^* \leq bb^*$. Conversely suppose $aa^* \leq bb^*$. Then $a(a^*bb^*b) = (aa^*bb^*)b = bb^*b = b$, hence $(a, b) \in \leq_{\mathcal{R}}$, establishing Part (i). Part (ii) follows dually to Part (i). Parts (iii), (iv) and (v) follow from applying Parts (i) and (ii) to how Green's $\mathcal{R}$, $\mathcal{L}$ and $\mathcal{D}$ relations are defined.
	        \end{proof}
	    \end{proposition}
	    
	    \begin{proposition}
		    Given a regular $^*$-semigroup $S$:
		    \begin{enumerate}
			    \item each $\mathcal{R}$ class contains precisely one projection; and
			    \item each $\mathcal{L}$ class contains precisely one projection.
		    \end{enumerate}
		    
		    \begin{proof}
			    Let $a, b \in S$ and suppose $(aa^*, bb^*) \in \mathcal{R}$, and hence that $(aa^*)(aa^*)^* = (bb^*)(bb^*)^*$. Then $aa^* = aa^*aa^* = bb^*bb^* = bb^*$, establishing Part (i). Part (ii) follows dually to Part (i).	
		    \end{proof}
	    \end{proposition}
	    
	\subsection{Congruence relations}
		\begin{definition}
			Given a semigroup $S$, a \textit{congruence relation} is an equivalence relation $\sim\ \subseteq S\times{S}$ that is compatible with the operation for $S$. The equivalence classes of a congruence relation are often also referred to as \textit{congruence classes}.
		\end{definition}
		
		The congruence classes of a congruence relation form a monoid where for each $s, s' \in S$, the congruence class containing $s$ multiplied by the congruence class containing $s'$ is equal to the congruence class containing $ss'$, that is $[s][s'] = [ss']$. 
		
		Note that congruence relations are closed under intersections, hence given a relation $\mathcal{R}$ on a semigroup $S$, the intersection of all congruence relations containing $\mathcal{R}$ is itself a congruence relation, and hence is the smallest congruence relation containing $\mathcal{R}$.
		
		\begin{definition}
			Given $n \in \posints$ binary relations $\mathcal{R}_1$, \ldots, $\mathcal{R}_n$ on a semigroup $S$, \textit{the congruence relation generated by $\mathcal{R}_1$, \ldots, $\mathcal{R}_n$}, which we denote by $\congrel{\gens{\mathcal{R}_1, \ldots, \mathcal{R}_n}}$, is the intersection of all congruence relations that contain $\cup_{i=1}^n \mathcal{R}_i$.
		\end{definition}
	
    \subsection{Free semigroups and free monoids}
        \begin{definition}
        	Given a set $X$:
        	\begin{enumerate}
            	\item the \textit{free semigroup of $X$}, which is often denoted by $\freesemigrp{X}$, is the set of non-zero finite strings from $X$ together with string concatenation, that is $\freesemigrp{X} = \{\Pi^n_{i=1}x_i: n \in \posints \text{ and } x_1, \ldots, x_n \in S\}$; and
            	\item the \textit{free monoid of $X$}, which is often denoted by $\freemon{X}$, adjoins the empty string, which acts as the identity under string concatenation, to the free semigroup $\freesemigrp{X}$.
            \end{enumerate}
        \end{definition}

    \subsection{Presentations}
        \begin{definition}
            A \textit{presentation} of a semigroup $S$ consists of:
            \begin{enumerate}
                \item a set of generators $\Sigma$; and 
                \item a set of binary relations $\{\mathcal{R}_1, \ldots, \mathcal{R}_n\}$ on the free semigroup $\freesemigrp{\Sigma}$ generated by $\Sigma$,
            \end{enumerate}
            such that $S \cong \freesemigrp{\Sigma}/\congrel{\gens{\mathcal{R}_1, \ldots, \mathcal{R}_n}}$.
        \end{definition}
        
        Perhaps more descriptively, given a semigroup $S$, a generating set $\Sigma$ of $S$ along with $n \in \posints$ binary relations $\mathcal{R}_1$, \ldots, $\mathcal{R}_n$ form a presentation of $S$ when any two words in the free semigroup $\freesemigrp{\Sigma}$ whose products form the same element of $S$ are able to be shown as equivalent in an abstract manner using the binary relations $\mathcal{R}_1$, \ldots, $\mathcal{R}_n$.
    \section{Fundamentals of diagram semigroups}
	There are fewer times more appropriate than a doctoral thesis to give a complete coverage of the foundations for diagram semigroups, consequently the author has taken the opportunity to do just that.

	\subsection{Diagrams $\diagrams{k}$}
	
		\begin{definition}
			Let $k \in \nonnegints$. A \textit{$k$-diagram} is a reflexive and symmetric binary relation on $\{1, \ldots, k\}\cup\{1', \ldots, k'\}$. We denote by $\diagrams{k}$ the set of all $k$-diagrams, and will refer to $k$-diagrams more succinctly as \textit{diagrams} either when $k$ may be any positive integer or when the value of $k$ is contextually unambiguous.
		\end{definition}
		
		\begin{definition}
			Let $k \in \nonnegints$ and $\delta \in \diagrams{k}$. We refer to:
			\begin{enumerate}
				\item elements of $\{1, \ldots, k\}$ as upper vertices;
				\item elements of $\{1', \ldots, k'\}$ as lower vertices; and
				\item elements of $\delta$ as \textit{edges} or \textit{lines}.
			\end{enumerate}
		\end{definition}
		
		Note provided $k \in \nonnegints$ is contextually unambiguous, we may specify a diagram using an edge set that is neither reflexive nor symmetric with the implicit understanding that we are actually referring to the reflexive and symmetric closure. For example given $k = 2$, if a diagram is specified as $\{(1,2), (2, 2')\}$, we will implicitly mean the diagram $\{(1, 1), (2,2), (1', 1'), (2', 2'), (1, 2), (2, 1), (2, 2')$, $(2', 2)\} \in \diagrams{2}$.

\begin{figure}[!ht]
	  		\caption[ ]{$\{(2,1'), (3,4), (4,7), (7,7'), (7',8'), (5,6), (8,6'),(6',5'), (2',3'), (3',4')\} \in \diagrams{8}$ may be depicted, as outlined in Definition \ref{def:graphicaldepictionofdiagrams}, as follows:}
	  		\label{fig:P8example}
			\vspace{5pt}
			\centering
	  		\input{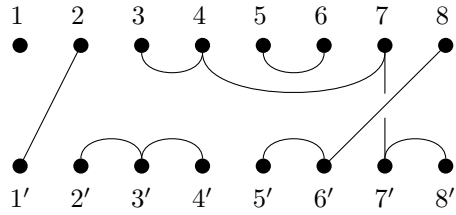}
		\end{figure}

	 	\begin{definition} \label{def:graphicaldepictionofdiagrams}
		 	Let $k \in \nonnegints$. Each diagram $\delta \in \diagrams{k}$ is typically \textit{depicted} as follows (see Figure \ref{fig:P8example} for an example): 
		 	\begin{enumerate}
		 		\item for each $j \in \{1, \ldots, k\}$, the upper vertex $j$ is depicted as the point $(j, 1)$ while the lower vertex $j'$ is depicted as the point $(j, 0)$; and
		 		\item $\delta$'s lines are drawn within the convex hull of the points $\{1, \ldots, k\}\times\{0, 1\}$.
		 	\end{enumerate}
		\end{definition}
	
		Note that vertex labels are trivially recoverable without their inclusion, consequently vertex labels are omitted from all figures proceeding Figure \ref{fig:P8example}.
		
		We will adopt the convention that lines between two upper vertices along with lines between two lower vertices will be non-linear, and that lines between an upper vertex and a lower vertex will be, inclusive of vertical sections, piecewise linear.
	
	\subsection{Diagram equivalence $\diagequiv$}
		Let $k \in \nonnegints$ and $\delta \in \diagrams{k}$. If we take the transitive closure of $\delta$ we get an equivalence relation on $\{1, \ldots, k\}\cup\{1', \ldots, k'\}$, the equivalence classes of which form a partition of $\{1, \ldots, k\}\cup\{1', \ldots, k'\}$.

		\begin{definition}
			Let $k \in \nonnegints$. The \textit{connected components} of a $k$-diagram $\delta \in \diagrams{k}$ are the equivalence classes of the transitive closure of $\delta$.
		\end{definition}
		
\begin{figure}[!ht]
	  		\caption[ ]{Given $k = 6$, the diagrams $\delta = \{(1,2), (2,3), (3,4), (4,5), (4,4'), (3',4), (4',5')\}$, $\delta' = \{(1,2),$ $(1,4), (3,4), (5,4'), (5,5'), (6,6'), (5',6')\} \in \diagrams{6}$ are depicted below.}
	  		\label{fig:equivdiagrams}
			\vspace{5pt}
			\begin{center}
				\input{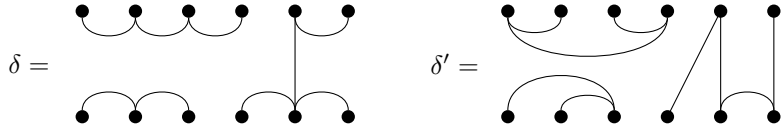}
			\end{center}
	  		While $\delta$ and $\delta'$ do not contain the same edges, the connected components of both $\delta$ and $\delta'$ are $\{\{1,2,3,4\}, \{5,6,4',5',6'\}, \{1',2',3'\}\}$ which is a partition of $\{1, \ldots, 6\}\cup\{1', \ldots, 6'\}$.
		\end{figure}

		For example the connected components of the diagram depicted in Figure \ref{fig:P8example} are: \[\{\{1\}, \{2, 1'\}, \{3,4,7,7',8'\}, \{5, 6\}, \{8,5',6'\}, \{2',3',4'\}\}.\]
		Given $k \in \nonnegints$, the connected components of each $k$-diagram form a partition of $\{1, \ldots, k\}\cup\{1', \ldots, k'\}$. Note however that there exist distinct diagrams whose connected components form precisely the same partition of $\{1, \ldots, k\}\cup\{1', \ldots, k'\}$. For example Figure \ref{fig:equivdiagrams} depicts two diagrams whose connected components form the same partition.
	
		\begin{definition}
			Let $k \in \nonnegints$. We say that two diagrams $\delta, \delta' \in \diagrams{k}$ are \textit{equivalent} if they have the same connected components. Furthermore we denote by $\diagequiv$ the binary relation that relates equivalent diagrams.
		\end{definition}
		
		For example the two diagrams depicted in Figure \ref{fig:equivdiagrams} are equivalent. Note that $\diagequiv$ is trivially an equivalence relation on $k$-diagrams. 
		
	\subsection{Bipartitions}
		\begin{definition} \label{def:bipartitions}
			Let $k \in \nonnegints$. A \textit{$k$-bipartition} is a set partition of $\{1, \ldots, k\}\cup\{1', \ldots, k'\}$. We denote by $\pttnmon{k}$ the set of all $k$-bipartitions, and will refer to $k$-bipartitions more succinctly as \textit{bipartitions} either when $k$ may be any positive integer or when the value of $k$ is contextually unambiguous. Furthermore we refer to elements of bipartitions as \textit{blocks}.
		\end{definition}
		
		Typically $\pttnmon{0}$ is treated as the singleton set containing the empty bipartition.
		
		\begin{definition} \label{def:diagramclassofbipartition}
			For each $k \in \nonnegints$ and $\alpha \in \pttnmon{k}$, we denote by $\diagclass{\alpha}$ the distinct equivalence class of diagrams in $\sfrac{\diagrams{k}}{\diagequiv}$ such that for each $\delta \in \diagclass{\alpha}$, the connected components of $\delta$ are the bipartition $\alpha$.
		\end{definition}
		
		To depict the bipartition $\alpha \in \pttnmon{k}$ we select any diagram $\delta \in \diagclass{\alpha}$, and treat blocks of bipartitions as being synonymous with connected components of diagrams.
		
		No consistent convention has been followed by the author for which diagram will be selected to depict a given bipartition, though the author has attempted to select diagrams that are aesthetically pleasing whilst ensuring a visible distance between distinct connected components.
 	
	\subsection{Block types}
		It is often convenient to establish and use various bits of notation and terminology for blocks of bipartitions. Terminology tends to differ between authors, though definitions are typically either identical or equivalent.
		
	 	\begin{definition}
	 		For each $k \in \nonnegints$, $\alpha \in \pttnmon{k}$ and $b \in \alpha$, we denote by:
		 	\begin{enumerate}
		 		\item $\upverts{b}$ the subset of upper vertices contained in $b$, that is $b \cap \set{1, \ldots, k}$; 
		 		\item $\noupverts{b}$ the number of upper vertices contained in $b$, that is $\card{\upverts{b}}$;
		 		\item $\lowverts{b}$ the subset of lower vertices contained in $b$, that is $b \cap \set{1', \ldots, k'}$; and
		 		\item $\nolowverts{b}$ the number of lower vertices contained in $b$, that is $\card{\lowverts{b}}$.
		 	\end{enumerate}
	 	\end{definition}

		\begin{definition}
			Given a bipartition $\alpha \in \pttnmon{k}$, the \textit{type} of a block $b \in \alpha$ is the pair $\big(u(b), l(b)\big)$, where $u(b), l(b) \in \set{0, \ldots, k}$ such that $u(b) + l(b) > 0$ are the number of upper vertices and number of lower vertices in block $b$ respectively. 
		\end{definition}
		
		It will be convenient for us to use some further terminology for blocks.
		
		\begin{definition} \label{def:blocktypes}
			Given $m, n \in \posints$ we refer to:
			\begin{enumerate}
				\item a block of type $(m, n)$ as a \textit{transversal};
				\item a block of type $(1, 1)$ as a \textit{transversal line};
				\item a block of type $(m, m)$ as \textit{uniform};
				\item a transversal line $\set{i, i'}$ where $i \in \set{1, \ldots, k}$ as a \textit{vertical line}; 
				\item a block of type $(m, 0)$ or $(0, m)$ as a \textit{non-transversal};
				\item a block of type $(2, 0)$ or $(0, 2)$ as a \textit{non-transversal line};
				\item a block of type $(1,1)$, $(2, 0)$ or $(0, 2)$ as a \textit{line};
				\item a non-transversal containing $m$ consecutive vertices as an \textit{$m$-apsis}; 
				\item $1$-apses as \textit{monapses}, $2$-apses as \textit{diapses} and $3$-apses as \textit{triapses}; and
				\item a transversal of type $(2, 2)$ containing the same two adjacent upper and lower vertices as a \textit{$(2,2)$-transapsis}.
			\end{enumerate}
		\end{definition}

	\subsection{Product graphs $\productgraph{\delta_1, \ldots, \delta_p}$}
		\begin{definition} \label{def:productgraph}
			For each $k \in \nonnegints$, $p \in \posints$ and $\delta_1, \ldots, \delta_p \in \diagrams{k}$:
			\begin{enumerate}
				\item by $\overbar{\delta}_1$ we denote the reflexive and symmetric relation on $\{1, \ldots, k\}\cup\{1^2, \ldots, k^2\}$ where for each $j \in \{1, \ldots, k\}$, every instance of the lower vertex $j'$ in $\delta_1$ has been relabelled to $j^2$ in $\overbar{\delta}_1$; 
				\item by $\overbar{\delta}_p$ we denote the reflexive and symmetric relation on $\{1^p, \ldots, k^p\}\cup\{1', \ldots, k'\}$ where for each $j \in \{1, \ldots, k\}$, every instance of the upper vertex $j$ in $\delta_p$ has been relabelled to $j^p$ in $\overbar{\delta}_p$; and
				\item for each $i \in \{2, \ldots, p-1\}$, by $\overbar{\delta}_i$ we denote the reflexive and symmetric relation on $\{1^i, \ldots, k^i\}\cup\{1^{i+1}, \ldots, k^{i+1}\}$ where for each $j \in \{1, \ldots, k\}$, every instance of the upper vertex $j$ in $\delta_i$ has been relabelled to $j^i$ in $\overbar{\delta}_i$ and every instance of the lower vertex $j'$ in $\delta_i$ has been relabelled to $j^{i+1}$ in $\overbar{\delta}_i$.
			\end{enumerate}
			\textit{The product graph of $\delta_1, \ldots, \delta_p$} is the union $\bigcup^p_{i=1}\overbar{\delta}_i$, which is a reflexive and symmetric relation on $\{1, \ldots, k\}\cup\{1', \ldots, k'\}\cup(\bigcup_{i=2}^{p}\{1^i, \ldots, k^i\})$ that is often denoted as $\productgraph{\delta_1, \ldots, \delta_p}$. 
		\end{definition}
		
		For each $k \in \nonnegints$ and $\delta \in \diagrams{k}$, it trivially follows by definition that $\productgraph{\delta} = \delta$, consequently we will implicitly interchange between the two whenever it is convenient to do so.
		
\begin{figure}[!ht]
			\caption[ ]{Let $k=5$, $\delta_{\alpha} = \{(1,1'), (2,2'), (3,4), (5,3'), (6,4'), (7,9'), (8,9), (5',8'), (6',7')\} \in \diagrams{5}$ and $\delta_{\beta} = \{(1,4), (1,1'), (1',2'), (2,3), (5,6'), (3',6'), (6,7), (8,9), (8,7'), (4',5'), (8',9')\} \in \diagrams{5}$. The product graph $\productgraph{\delta_{\alpha}, \delta_{\beta}}$, as outlined in Definition \ref{def:productgraphdepiction}, is depicted as:}
			\label{fig:productgraphexample}
			\vspace{-5pt}
			\begin{center}
				\input{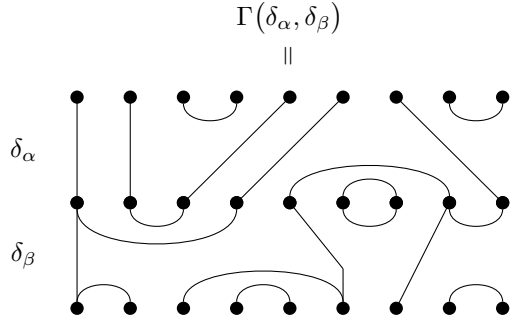}
			\end{center}
		\end{figure}

		\begin{definition} \label{def:productgraphdepiction}
			For each $k \in \nonnegints$, $p \in \posints$ and $\delta_1, \ldots, \delta_p \in \diagrams{k}$, the product graph $\productgraph{\delta_1, \ldots, \delta_p}$ is typically depicted by vertically stacking $\delta_1, \ldots, \delta_p$ such that for each $j \in \{1, \ldots, p\}$, $\delta_j$ is depicted as the $j$-th upper most diagram. Being more precise (see Figure \ref{fig:productgraphexample} for an example):
			\begin{enumerate}
				\item for each $j \in \{1, \ldots, k\}$, the upper vertex $j$ is depicted as the point $(j, p)$, the lower vertex $j'$ is depicted as the point $(j, 0)$, and for each $i \in \{2, \ldots, p\}$, the vertex $j^{i-1}$ is depicted as the point $(j, i-1)$; and
				\item for each $i \in \{1, \ldots, p\}$, each of $\overbar{\delta}_i$'s edges are drawn within the convex hull of the points $\{1, \ldots, k\}\times\{p-i, p-i+1\}$.
			\end{enumerate}
		\end{definition}
		
		We will often label the vertical layers of a product graph by whichever is more convenient out of the associated diagram or the underlying bipartition for each layer.
		
		\begin{definition}
			For each $k \in \nonnegints$, $p, q \in \posints$ and $\delta_1, \ldots, \delta_p, \eta_1, \ldots, \eta_q \in \diagrams{k}$, the \textit{product graph of $\productgraph{\delta_1, \ldots, \delta_p}$ and $\productgraph{\eta_1, \ldots, \eta_q}$}, which we denote by $\productgraph{\productgraph{\delta_1, \ldots, \delta_p}, \productgraph{\eta_1, \ldots, \eta_q}}$, is formed by vertically stacking $\productgraph{\delta_1, \ldots, \delta_p}$ above $\productgraph{\eta_1, \ldots, \eta_q}$ in an analogous way to how product graphs of diagrams are formed.
		\end{definition}
		
		\begin{proposition}
			For each $k \in \nonnegints$, $p, q \in \posints$ and $\delta_1, \ldots, \delta_p, \eta_1, \ldots, \eta_q \in \diagrams{k}$, \[\productgraph{\productgraph{\delta_1, \ldots, \delta_p}, \productgraph{\eta_1, \ldots, \eta_q}} = \productgraph{\delta_1, \ldots, \delta_p, \eta_1, \ldots, \eta_q}.\]
			
			\begin{proof}
				Vertically stacking $\delta_1, \ldots, \delta_p, \eta_1, \ldots, \eta_q$ is trivially equivalent to vertically stacking the vertical stack of $\delta_1, \ldots, \delta_p$ with the vertical stack of $\eta_1, \ldots, \eta_q$.
			\end{proof}
		\end{proposition}

		\begin{corollary}
			Forming product graphs of product graphs is associative.
			
			\begin{proof}
				For each $k \in \nonnegints$, $p, q, r \in \posints$ and $\delta_1, \ldots, \delta_p, \eta_1, \ldots, \eta_q, \xi_1, \ldots, \xi_r \in \diagrams{k}$,
				\begin{align*}
					\productgraph{\productgraph{\productgraph{\delta_1, \ldots, \delta_p}, \productgraph{\eta_1, \ldots, \eta_q}}, \productgraph{\xi_1, \ldots, \xi_r}} & = \productgraph{\productgraph{\delta_1, \ldots, \delta_p, \eta_1, \ldots, \eta_q}, \productgraph{\xi_1, \ldots, \xi_r}} \\
					& = \productgraph{\delta_1, \ldots, \delta_p, \eta_1, \ldots, \eta_q, \xi_1, \ldots, \xi_r} \\
					& = \productgraph{\productgraph{\delta_1, \ldots, \delta_p}, \productgraph{\eta_1, \ldots, \eta_q, \xi_1, \ldots, \xi_r}} \\
					& = \productgraph{\productgraph{\delta_1, \ldots, \delta_p}, \productgraph{\productgraph{\eta_1, \ldots, \eta_q}, \productgraph{\xi_1, \ldots, \xi_r}}}.
				\end{align*} 
			\end{proof}
		\end{corollary}
		
		\begin{definition}
			Let $k \in \nonnegints$, $p, q \in \posints$ and $\delta_1, \ldots, \delta_p, \eta_1, \ldots, \eta_q \in \diagrams{k}$. We say that the two product graphs $\productgraph{\delta_1, \ldots, \delta_p}$ and $\productgraph{\eta_1, \ldots, \eta_q}$ are \textit{equivalent}, which we denote by $\productgraph{\delta_1, \ldots, \delta_p} \diagequiv \productgraph{\eta_1, \ldots, \eta_q}$, if their connected components restricted to $\{1, \ldots, k\}\cup\{1', \ldots, k'\}$ are equal, with the implicit understanding that we forget about any then-empty components following our restriction.
		\end{definition}
		
		\begin{proposition} \label{prop:diagequivpreservesconnectedcomponents}
			For each $k \in \nonnegints$, $p, q, r, s \in \posints$ and $\delta_1, \ldots, \delta_p$, $\eta_1, \ldots, \eta_q$, $\xi_1, \ldots, \xi_r$, $\omega_1, \ldots, \omega_s$ $\in \diagrams{k}$, if:
			\begin{enumerate}
				\item $\productgraph{\delta_1, \ldots, \delta_p} \diagequiv \productgraph{\eta_1, \ldots, \eta_q}$; and
				\item $\productgraph{\xi_1, \ldots, \xi_r}$ $\diagequiv \productgraph{\omega_1, \ldots, \omega_s}$,
			\end{enumerate}
			then $\productgraph{\delta_1, \ldots, \delta_p, \xi_1, \ldots, \xi_r} \diagequiv \productgraph{\eta_1, \ldots, \eta_q, \omega_1, \ldots, \omega_s}$.
			
			\begin{proof}
				First note that 
				\[\begin{cases}\productgraph{\delta_1, \ldots, \delta_p, \xi_1, \ldots, \xi_r} = \productgraph{\productgraph{\delta_1, \ldots, \delta_p}, \productgraph{\xi_1, \ldots, \xi_r}}; \text{ and} \\ \productgraph{\eta_1, \ldots, \eta_q, \omega_1, \ldots, \omega_s} = \productgraph{\productgraph{\eta_1, \ldots, \eta_q}, \productgraph{\omega_1, \ldots, \omega_s}}.\end{cases}\]
				For each $j \in \{1, \ldots, k\}$, we relabel the vertex $j^{p+1}$ to $j''$ in the $(p+1)$-th vertex row of the product graph $\productgraph{\delta_1, \ldots, \delta_p, \xi_1, \ldots, \xi_r}$, and relabel the vertex $j^{r+1}$ to $j''$ in the $(r+1)$-th vertex row of the product graph $\productgraph{\eta_1, \ldots, \eta_q, \omega_1, \ldots, \omega_s}$.
				
				Let $i, j \in \{1, \ldots, k\}\cup\{1', \ldots, k'\}$ such that there exists a path $p_{i, j}$ from $i$ to $j$ in the product graph $\productgraph{\delta_1, \ldots, \delta_p, \xi_1, \ldots, \xi_r}$. Let $t \in \intsge{2}$ such that the minimum number of times $p_{i,j}$ crosses between $\productgraph{\delta_1, \ldots, \delta_p}$ and $\productgraph{\xi_1, \ldots, \xi_r}$ is $t-2$. Therefore there exists $t$ vertices $v_1, \ldots, v_t \in \{1, \ldots, k\}\cup\{1'', \ldots, k''\}\cup\{1', \ldots, k'\}$ and paths $p_{v_1, v_2}, \ldots, p_{v_{t-1}, v_t}$ such that $v_1 = i$, $v_t = j$ and for each $l \in \{1, \ldots, t-1\}$, $p_{v_{l}, v_{l+1}}$ is a path from $v_l$ to $v_{l+1}$ that is either a subset of $\productgraph{\delta_1, \ldots, \delta_p}$ or a subset of $\productgraph{\xi_1, \ldots, \xi_r}$. 

				For each $l \in \{1, \ldots, t-1\}$, it follows from $\productgraph{\delta_1, \ldots, \delta_p} \diagequiv \productgraph{\eta_1, \ldots, \eta_q}$ and $\productgraph{\xi_1, \ldots, \xi_r} \diagequiv \productgraph{\omega_1, \ldots, \omega_s}$ that there exists a path from $v_l$ to $v_{l+1}$, which we denote as $p'_{v_l, v_{l+1}}$, that is either a subset of $\productgraph{\eta_1, \ldots, \eta_q}$ or a subset of $\productgraph{\omega_1, \ldots, \omega_s}$ depending on whether $\productgraph{\delta_1, \ldots, \delta_p}$ or $\productgraph{\xi_1, \ldots, \xi_r}$ contains the path $p_{v_l, v_{l+1}}$.
				
				Finally $\bigcup_{l=1}^{t-1} p'_{v_l, v_{l+1}}$ forms a path from $i$ to $j$ and for each $l \in \{1, \ldots, t-1\}$, $p'_{v_l, v_{l+1}}$ is either a subset of $\productgraph{\eta_1, \ldots, \eta_q}$ or a subset of $\productgraph{\omega_1, \ldots, \omega_s}$. Hence $i$ and $j$, which share the same connected component of $\productgraph{\delta_1, \ldots, \delta_p, \xi_1, \ldots, \xi_r}$, also share the same connected component of $\productgraph{\eta_1, \ldots, \eta_q, \omega_1, \ldots, \omega_s}$. 
				
				The same argument in the opposite direction establishes that for any two vertices $i, j \in \{1, \ldots, k\}\cup\{1', \ldots, k'\}$ that share the same connected component of the product graph $\productgraph{\eta_1, \ldots, \eta_q, \omega_1, \ldots, \omega_s}$, $i$ and $j$ also share the same connected component of $\productgraph{\delta_1, \ldots, \delta_p, \xi_1, \ldots, \xi_r}$. Hence $\productgraph{\delta_1, \ldots, \delta_p, \xi_1, \ldots, \xi_r} \diagequiv \productgraph{\eta_1, \ldots, \eta_q, \omega_1, \ldots, \omega_s}$.
			\end{proof}
		\end{proposition}

	\subsection{The partition monoid $\pttnmon{k}$}
		The partition monoid consists of all bipartitions equipped with an associative binary operation, we review two equivalent approaches at defining this operation. The first approach forms the product of two bipartitions using the product graph of representative diagrams, admitting a rather intuitive way to depict products of bipartitions in the process. The second approach multiplies bipartitions directly. Since both products are equivalent, which we will later establish, we will only make a distinction between the two when it is relevant to do so.
	
		\begin{definition} \label{def:bipartitionproduct}
			Let $k \in \nonnegints$ and $\alpha, \beta \in \pttnmon{k}$. The \textit{product of $\alpha$ and $\beta$}, denoted as $\alpha\beta$, is formed as follows: 
			\begin{enumerate}
				\item Let $\delta_{\alpha} \in \diagrams{\alpha}$ and $\delta_{\beta} \in \diagrams{\beta}$;
				\item two vertices $i, j \in \{1, \ldots, k\}\cup\{1', \ldots, k'\}$ share the same block in the product $\alpha\beta$ if they share the same connected component of the product graph $\productgraph{\delta_{\alpha}, \delta_{\beta}}$, that is the blocks of the product $\alpha\beta$ are the connected components of $\productgraph{\delta_{\alpha}, \delta_{\beta}}$ restricted to $\{1, \ldots, k\}\cup\{1', \ldots, k'\}$ with any then-empty connected components implicitly removed.
			\end{enumerate}
		\end{definition}
		
		For example consider $\delta_{\alpha}, \delta_{\beta} \in \diagrams{8}$ from Figure \ref{fig:productgraphexample} where the product graph $\productgraph{\delta_{\alpha}, \delta_{\beta}}$ is depicted. It is easily verified that the underlying bipartition of $\delta_{\alpha}$ is $\alpha = \{\{1, 1'\}, \{2, 2'\}, \{3, 4\}, \{5, 3'\}, \{6, 4'\}, \{7, 9'\}, \{8, 9\},$ $\{5', 8'\}, \{7', 8'\}\} \in \pttnmon{8}$ and that the underlying bipartition of $\delta_{\beta}$ is $\beta = \{\{1, 4, 1', 2'\}, \{2, 3\}, \{5, 3', 6'\}$, $\{6, 7\}, \{8, 9, 7'\}, \{4', 5'\}, \{8', 9'\}\} \in \pttnmon{8}$, that is $\delta_{\alpha} \in \diagrams{\alpha}$ and $\delta_{\beta} \in \diagrams{\beta}$. Hence two vertices $i, j \in \{1, \ldots, k\}\cup\{1', \ldots, k'\}$ share the same block in $\alpha\beta$ if there is a path from $i$ to $j$ in the product graph $\productgraph{\delta_{\alpha}, \delta_{\beta}}$ that is depicted in Figure \ref{fig:productgraphexample}. Giving us $\alpha\beta = \{\{1,6,1',2'\}, \{2,5\}, \{3,4\}, \{7,3',6',7'\}, \{8,9\}, \{4',5'\}, \{8',9'\}\} \in \pttnmon{8}$.	
		
		Note that in order for the bipartition product to be well-defined, the choice of representative diagrams from $\diagrams{\alpha}$ and $\diagrams{\beta}$ must be arbitrary, which is conveniently a restricted case of Proposition \ref{prop:diagequivpreservesconnectedcomponents}.
		
		\begin{corollary} \label{prop:bipartitionproductwelldefined}
			For each $k \in \nonnegints$, $\alpha, \beta \in \pttnmon{k}$, $\delta_{\alpha}, \delta'_{\alpha} \in \diagrams{\alpha}$ and $\delta_{\beta}, \delta'_{\beta} \in \diagrams{\beta}$, $\productgraph{\delta_{\alpha}, \delta_{\beta}} \diagequiv \productgraph{\delta'_{\alpha}, \delta'_{\beta}}$.
			
			\begin{proof}
				A restricted case of Proposition \ref{prop:diagequivpreservesconnectedcomponents}.
			\end{proof}
		\end{corollary}	
		
		\begin{proposition} \label{prop:bipartitionproductassociative}
			For each $k \in \nonnegints$ and $\alpha, \beta, \gamma \in \pttnmon{k}$, $(\alpha\beta)\gamma = \alpha(\beta\gamma)$.
			
			\begin{proof}
				Let $\delta_{\alpha} \in \diagrams{\alpha}$, $\delta_{\beta} 
				\in \diagrams{\beta}$, $\delta_{\alpha\beta} \in \diagrams{\alpha\beta}$ and $\delta_{\gamma} \in \diagrams{\gamma}$. Note that the connected components of $\delta_{\alpha\beta}$ and the connected components of $\productgraph{\delta_{\alpha}, \delta_{\beta}}$ restricted to $\{1, \ldots, k\}\cup\{1', \ldots, k'\}$ are the blocks of $\alpha\beta$, and hence that the connected components of $\delta_{\alpha\beta}$ are equal to the connected components of $\productgraph{\delta_{\alpha}, \delta_{\beta}}$ restricted to $\{1, \ldots, k\}\cup\{1', \ldots, k'\}$, hence $\delta_{\alpha\beta} \diagequiv \productgraph{\delta_{\alpha}, \delta_{\beta}}$. 
				
				Note it follows from Proposition \ref{prop:diagequivpreservesconnectedcomponents} that the connected components of $\productgraph{\delta_{\alpha\beta}, \delta_{\gamma}}$ and the connected components of $\productgraph{\productgraph{\delta_{\alpha}, \delta_{\beta}}, \delta_{\gamma}} = \productgraph{\delta_{\alpha}, \delta_{\beta}, \delta_{\gamma}}$ are equal. Hence the blocks of $(\alpha\beta)\gamma$ are equal to the connected components of $\productgraph{\delta_{\alpha}, \delta_{\beta}, \delta_{\gamma}}$ restricted to $\{1, \ldots, k\}\cup\{1', \ldots, k'\}$.
				
				By an analogous argument the blocks of $\alpha(\beta\gamma)$ are also equal to the connected components of $\productgraph{\delta_{\alpha}, \delta_{\beta}, \delta_{\gamma}}$ restricted to $\{1, \ldots, k\}\cup\{1', \ldots, k'\}$, giving us $(\alpha\beta)\gamma = \alpha(\beta\gamma)$.
			\end{proof}
		\end{proposition}
		
		\begin{definition}
			For each $k \in \posints$, we denote by $\id{k}$ the bipartition $\Big\{\{j, j'\}: j \in \{1, \ldots, k\}\Big\} \in \pttnmon{k}$ (see Figure \ref{fig:id4} for a depiction of $\id{4}$). 
		\end{definition}
		
		\begin{figure}[!ht] 
			\caption[ ]{Given $k=4$,}
			\label{fig:id4}
			\vspace{5pt}
			\centering
			\input{chap_background/tikz/fig-id4.tex}
		\end{figure}
		
		Note $\id{k}$ trivially acts as both a left and right identity, and hence is the identity of the bipartition product. 
		
	 	\begin{definition}
	 		Given $k \in \nonnegints$, \textit{the partition monoid}, which is often also denoted as $\pttnmon{k}$ for convenience, consists of all $k$-bipartitions $\pttnmon{k}$ along with the bipartition product. 
	 	\end{definition}
		
		\begin{corollary}
			For each $k \in \nonnegints$, $p \in \intsge{2}$ $\alpha_1, \ldots, \alpha_p \in \pttnmon{k}$ and $\delta_{\alpha_j} \in \diagrams{\alpha_j}$ where $j \in \{1, \ldots, p\}$, the blocks of $\alpha_1\ldots\alpha_p$ are the connected components of $\productgraph{\delta_{\alpha_1}, \ldots, \delta_{\alpha_p}}$.
			
			\begin{proof}
				Follows inductively on $p$ in an analogous way to the proof of Proposition \ref{prop:bipartitionproductassociative}, where we established that the bipartition product is associative.
			\end{proof}
		\end{corollary}
		
		Let $p \in \intsge{2}$ and $\alpha_1, \ldots, \alpha_p \in \pttnmon{k}$. When depicting the formation of the product $\alpha_1\ldots\alpha_p$, we will equate the product graph $\Gamma(\alpha_1, \ldots, \alpha_p)$ to any representative diagram $\delta \in \diagrams{\alpha_1\ldots\alpha_p}$, with the implicit understanding that we really mean the connected components on the product graph $\Gamma(\alpha_1, \ldots, \alpha_p)$ restricted to $\{1, \ldots, k\}\cup\{1', \ldots, k'\}$ are equal to the connected components of the diagram depicting the product $\alpha_1\ldots\alpha_p$.
		
		Next we review how the product of two bipartitions may equivalently be formed directly without using representative diagrams. 
		
		\begin{proposition}
			For each $k \in \nonnegints$ and $\alpha, \beta \in \diagrams{k}$, $\alpha\beta$ is the finest partition coarser than $\alpha_{\vee}\cup\beta^{\wedge}$ where:
			\begin{enumerate}
				\item $\alpha_{\vee}$ denotes the partition of $\{1, \ldots, k\}\cup\{1'', \ldots, k''\}$ such that for each $j \in \set{1, \ldots, k}$, the lower vertex $j'$ in $\alpha$ has been relabelled to $j''$ in $\alpha_{\vee}$; and
				\item $\beta^{\wedge}$ denotes the partition of $\{1'', \ldots, k''\}\cup\{1', \ldots, k'\}$ where for each $j \in \set{1, \ldots, k}$, the upper vertex $j$ in $\beta$ has been relabelled to $j''$ in $\beta^{\wedge}$.
			\end{enumerate}
			
			\begin{proof}
				Let $\delta_{\alpha} \in \diagrams{\alpha}$ and $\delta_{\beta} \in \diagrams{\beta}$. The product graph $\productgraph{\delta_{\alpha}, \delta_{\beta}}$ is formed in an equivalent way to how $\alpha_{\vee}\cup\beta^{\wedge}$ is formed, and the transitive closure of $\productgraph{\delta_{\alpha}, \delta_{\beta}}$ is the finest equivalence relation containing $\productgraph{\delta_{\alpha}, \delta_{\beta}}$. Hence the associated partition must be the finest partition coarser than $\alpha_{\vee}\cup\beta^{\wedge}$.
			\end{proof}
		\end{proposition}
		
		For example, consider $\alpha = \{\{1,5,4',5'\}, \{2,3,4\}, \{1'\}, \{2', 3'\}\}$, $\beta = \{\{1,4,5,1',2',3'\}, \{2,3\}, \{4',5'\}\} \in \pttnmon{5}$. The formation of the product $\alpha\beta$ as outlined in Definition \ref{def:bipartitionproduct} is depicted in Figure \ref{fig:simplediagramproductexample}, alternatively it may be computed as follows:
		\begin{enumerate}
			\item $\alpha_{\vee} = \{\{1,5,4'',5''\}, \{2,3,4\}, \{1''\}, \{2'', 3''\}\}$;
			\item $\beta^{\wedge} = \{\{1'',4'',5'',1',2',3'\}, \{2'',3''\}, \{4',5'\}\}$;
			\item $\{\{1,5,1'',4'',5'',1',2',3'\}, \{2,3,4\}, \{2'',3''\}, \{4',5'\}\}$ is the finest partition coarser than $\alpha_{\vee}\cup\beta^{\wedge}$;
			\item $\{\{1,5,1',2',3'\}, \{2,3,4\}, \{\}, \{4',5'\}\}$ are the blocks of the finest partition coarser than $\alpha_{\vee}\cup\beta^{\wedge}$ restricted to $\{1, \ldots, k\}\cup\{1', \ldots, k'\}$; and
			\item $\{\{1,5,1',2',3'\}, \{2,3,4\}, \{4',5'\}\}$ are the non-empty blocks of the finest partition coarser than $\alpha_{\vee}\cup\beta^{\wedge}$ restricted to $\{1, \ldots, k\}\cup\{1', \ldots, k'\}$, which is the bipartition $\alpha\beta$ and depicted in Figure \ref{fig:simplediagramproductexample}.
		\end{enumerate}
		
		\begin{figure}[!ht]
			\caption[ ]{Given $k = 5$ and $\alpha = \{\{1,5,4',5'\}, \{2,3,4\}, \{1'\}, \{2', 3'\}\}, \beta = \{\{1,4,5,1',2',3'\}, \{2,3\}, \{4',5'\}\} \in \pttnmon{5}$,}
			\label{fig:simplediagramproductexample}
			\vspace{5pt}
			\centering
			\input{chap_background/tikz/fig-simplediagramproductexample.tex}
		\end{figure}
		
\begin{figure}[!ht]
			\caption[ ]{Given $\alpha = \Big\{\{1,2\}, \{3,4'\}, \{4,5'\}, \{5,6'\}, \{6\}, \{7,7'\}, \{8,13'\}, \{9,10,11,12,13\}, \{1',2',3'\}, \{8',12'\},$ $\{9',10'\}, \{11'\}\Big\}$, $\beta = \Big\{\{1,2,3\}, \{4,7,3'\}, \{5,6\}, \{8,4',7',8'\}, \{9\}, \{10,11\}, \{12,13,11'\}, \{1',2'\}, \{5',6'\}, \{9',10'\},$ $\{12',13'\}\Big\} \in \pttnmon{13}$,}
			\label{fig:bigdiagramproductexample}
			\vspace{5pt}
			\centering
			\input{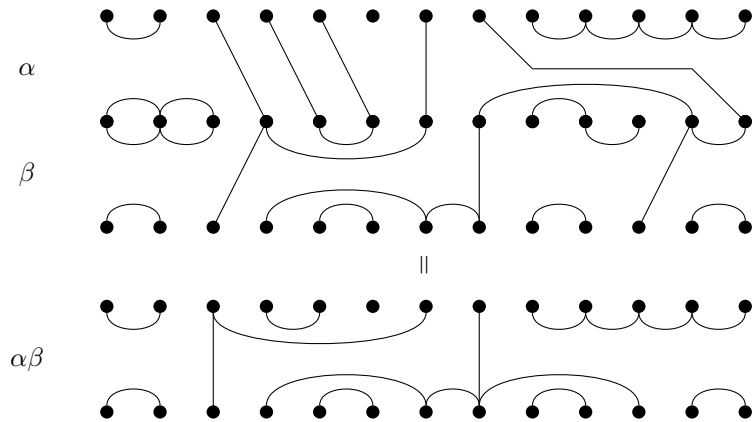}
		\end{figure}

		Figure \ref{fig:bigdiagramproductexample} depicts the formation of a more intricate product of bipartitions than we have considered up to now, which we proceed to discuss. Notice that when forming the product $\alpha\beta$:
		
		\begin{enumerate}
			\item the upper non-transversal blocks of $\alpha$ and lower non-transversal blocks of $\beta$ are preserved;
			\item the upper $2$-apsis block $\{5,6\} \in \beta$ joins to the transversal lines $\{4,5'\}$, $\{5,6'\} \in \alpha$, forming the upper non-transversal $\{4,5\} \in \alpha\beta$;
			\item the lower $3$-apsis block $\{1', 2', 3'\} \in \alpha$ and upper $3$-apsis block $\{1, 2, 3\} \in \beta$ cap each other off;
			\item the lower non-transversals $\{9', 10'\}, \{11'\} \in \alpha$ and upper non-transversals $\{9\}, \{10, 11\} \in \beta$ cap each other off;
			\item the transversal lines $\{3,4'\}, \{7, 7'\} \in \alpha$ join to the transversal $\{4,7,3'\} \in \beta$, forming the transversal $\{3,7,3'\} \in \alpha\beta$; and
			\item the blocks $\{8,13'\}, \{8', 12'\} \in \alpha$ join to the transversal blocks $\{8,4',7',8'\}$, $\{12,11'\} \in \beta$, forming the transversal block $\{8,4',7',8',11'\} \in \alpha\beta$.
		\end{enumerate}
		
	\subsection{Horizontal sum $\oplus:\pttnmon{k_1}\times\pttnmon{k_2}\rightarrow\pttnmon{k_1+k_2}$}
	
		\begin{definition} \label{def:horizontalsumofbipartitions}
			For each $k_1, k_2 \in \nonnegints$, $\alpha \in \pttnmon{k_1}$ and $\beta \in \pttnmon{k_2}$, \textit{the horizontal sum of $\alpha$ and $\beta$}, which we denote as $\alpha \oplus \beta$, is formed as follows:
			\begin{enumerate}
				\item let $\beta^{\rhd}$ denote the partition of $\{k_1+1, \ldots, k_1+k_2\}\cup\{(k_1+1)', \ldots, (k_1+k_2)'\}$ where for each $j \in \set{1, \ldots, k_2}$, the upper vertex $j$ in $\beta$ has been relabelled to $k_1+j$ and the lower vertex $j'$ in $\beta$ has been relabelled to $(k_1+j)'$ in $\beta^{\rhd}$; and
				\item the horizontal sum $\alpha \oplus \beta$ is the union $\alpha \cup \beta^{\rhd} \in \pttnmon{k_1 + k_2}$.
			\end{enumerate}
		\end{definition}
		
		The horizontal sum of diagrams is defined analogously as follows.
		
		\begin{definition} \label{def:horizontalsumofdiagrams}
			For each $k_1, k_2 \in \nonnegints$, $\delta \in \diagrams{k_1}$ and $\eta \in \diagrams{k_2}$, \textit{the horizontal sum of $\delta$ and $\eta$}, which we denote as $\delta \oplus \eta$, is formed as follows:
			\begin{enumerate}
				\item let $\eta^{\rhd}$ denote the translation of $\eta$ where for each $j \in \set{1, \ldots, k_2}$, the upper vertex $j$ in $\eta$ has been relabelled to $k_1+j$ and the lower vertex $j'$ in $\eta$ has been relabelled to $(k_1+j)'$ in $\eta^{\rhd}$; and
				\item the horizontal sum $\delta \oplus \eta$ is the union $\delta \cup \eta^{\rhd} \in \diagrams{k_1 + k_2}$.
			\end{enumerate}
		\end{definition}
		
		For each $k_1, k_2 \in \nonnegints$, $\alpha \in \pttnmon{k_1}$ and $\beta \in \pttnmon{k_2}$, it trivially follows by definition that $\diagrams{\alpha}\oplus\diagrams{\beta} = \diagrams{\alpha\oplus\beta}$. Hence if we select representative diagrams $\delta_{\alpha} \in \diagrams{\alpha}$ and $\delta_{\beta} \in \diagrams{\beta}$, we may depict $\alpha \oplus \beta$ by $\delta_{\alpha} \oplus \delta_{\beta}$ without having to explicitly determine $\alpha \oplus \beta$ directly to either verify that $\delta_{\alpha} \oplus \delta_{\beta} \in \diagrams{\alpha \oplus \beta}$ or to select a representative diagram from $\diagrams{\alpha\oplus\beta}$.
	
\begin{figure}[!ht]
			\caption[ ]{Given $k_1 = 6$, $k_2 = 8$,}
			\label{fig:horizontalsumeg}
			\vspace{5pt}
			\centering
			\input{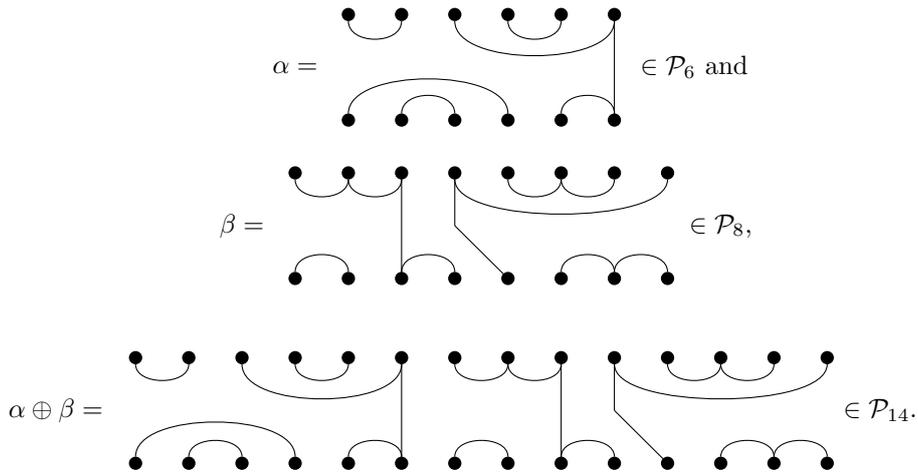}
		\end{figure}

		Figure \ref{fig:horizontalsumeg} depicts the horizontal sum of $\alpha = \{\{1,2\}, \{3,6,5',6'\}, \{4,5\}$, $\{1',4'\}, \{2',3'\}\} \in \pttnmon{6}$ and $\beta = \{\{1,2,3,3',4'\}, \{4,8,5'\}, \{1',2'\}, \{6',7',8'\}\}$.
		
		Note that the horizontal sum operation $\oplus$ is trivially associative.
	 	
	\subsection{Vertical flip involution $^*:\pttnmon{k}\rightarrow\pttnmon{k}$}
		\begin{definition} \label{def:verticalflipofbipartition}
			For each $k \in \nonnegints$ and $\alpha \in \pttnmon{k}$, \textit{the vertical flip of $\alpha$}, which is typically denoted as $\alpha^*$, is formed by for each $j \in \{1, \ldots, k\}$, switching the labels of $j$ and $j'$ in $\alpha$.
		\end{definition}
		
		The vertical flip of a diagram is defined analogously as follows.
	
		\begin{definition} \label{def:verticalflipofdiagram}
			For each $k \in \nonnegints$ and $\delta \in \diagrams{k}$, \textit{the vertical flip of $\delta$}, which is typically denoted as $\delta^*$, is formed by for each $j \in \{1, \ldots, k\}$, switching the labels of $j$ and $j'$ in $\delta$, or when depicted by vertically flipping $\delta$.
		\end{definition}
		
		For each $k \in \nonnegints$ and $\alpha \in \pttnmon{k}$, it trivially follows by definition that $\diagrams{\alpha}^* = \diagrams{\alpha^*}$. Hence if we select a representative diagram $\delta_{\alpha} \in \diagrams{\alpha}$, we may depict $\alpha^*$ by depicting $\delta_{\alpha}^*$, which is easy to do without having to explicitly determine $\alpha^*$ directly to either verify that $\delta_{\alpha}^* \in \diagrams{\alpha^*}$ or to select a representative diagram from $\diagrams{\alpha^*}$.
		
		Figure \ref{fig:verticalflipeg} depicts the vertical flip of the bipartition $\alpha = \{\{1,2\}, \{3,6,5',6'\}$, $\{4,5\}, \{1',4'\}, \{2',3'\}\} \in \pttnmon{6}$.
		
		\begin{figure}[!ht]
			\caption[ ]{Given $k = 6$ and}
			\label{fig:verticalflipeg}
			\vspace{5pt}
			\centering
			\input{chap_background/tikz/fig-verticalflipeg.tex}
		\end{figure}
		
\begin{figure}[!ht]
			\caption[ ]{Given $k = 6$ and $\alpha = \{\{1,2,3\}, \{4,7,3'\}, \{5,6\}, \{8,4',7',8'\}, \{9\}, \{10,11\}, \{12,13,11'\}, \{1',2'\},$ $\{5',6'\}, \{9',10'\}, \{12',13'\}\} \in \pttnmon{13}$, as established in Proposition \ref{prop:pttnmonisaregularstarsemigrp},}
			\label{fig:regularstarsemigroupeg}
			\vspace{5pt}
			\centering
			\input{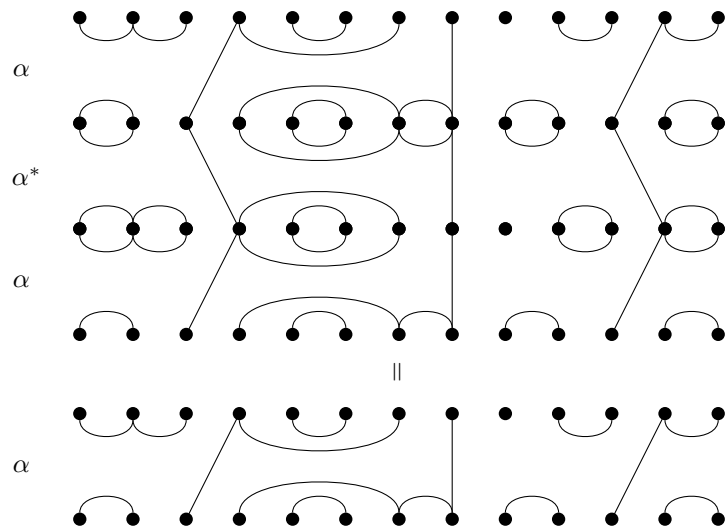}
		\end{figure}

		\begin{proposition} \label{prop:pttnmonisaregularstarsemigrp}
			For each $k \in \nonnegints$, the partition monoid $\pttnmon{k}$ is a regular $^*$-semigroup. That is for each $\alpha, \beta \in \pttnmon{k}$:
			\begin{enumerate}
				\item $\alpha^{**} = \alpha$;
				\item $(\alpha\beta)^* = \beta^*\alpha^*$; and
				\item $\alpha\alpha^*\alpha = \alpha$ (see Figure \ref{fig:regularstarsemigroupeg} for an example).
			\end{enumerate}
			Note that conditions (i) and (ii) give that $^*$ is an involution. 
			
			\begin{proof}
				Switching the upper and lower labels twice is trivially the same as never switching them to begin with, hence $\alpha^{**} = \alpha$. 
				
				Let $\delta_{\alpha} \in \diagrams{\alpha}$ and $\delta_{\beta} \in \diagrams{\beta}$. It is trivially the case that flipping the product graph of $\delta_{\alpha}$ and $\delta_{\beta}$ is equivalent to forming the product graph of $\delta_{\beta}$ flipped and $\delta_{\alpha}$ flipped, that is  $\productgraph{\delta_{\alpha}, \delta_{\beta}}^* = \productgraph{\delta_{\beta}^*, \delta_{\alpha}^*}$, hence $(\alpha\beta)^* = \beta^*\alpha^*$.
				
				Finally, when forming the product $\alpha\alpha^*\alpha$:
				\begin{enumerate}
					\item each upper non-transversal block in the left-most $\alpha$ is preserved;
					\item each lower non-transversal block in the right-most $\alpha$ is preserved;
					\item each lower non-transversal block in the left-most $\alpha$ and the corresponding upper non-transversal block in $\alpha^*$ join and are removed,
					similarly with each upper non-transversal block in the right-most $\alpha$ and the corresponding lower non-transversal block in $\alpha^*$; and
					\item each transversal block in the left-most $\alpha$ joins with the corresponding transversal block in $\alpha^*$, which then joins back up with the corresponding transversal block in $\alpha$.
				\end{enumerate}
				Hence $\alpha\alpha^*\alpha = \alpha$.
			\end{proof}
		\end{proposition}
		
		\begin{corollary}
			If a submonoid $S$ of the partition monoid $\pttnmon{k}$ is closed under the vertical flip involution $^*$, that is $S^* = S$, then $S$ is also a regular $^*$-semigroup.
		\end{corollary}
	
	\subsection{Patterns}	
		Patterns encode information about the non-transversal and transversal blocks for the upper or lower half of a bipartition. Note that while definitions are typically consistent between authors, terminology and notation frequently differs.
		
		\begin{definition}
			Let $k \in \nonnegints$. A \textit{$k$-pattern} is a pair $(\mathcal{A}, \mathcal{B})$ where $\mathcal{A}$ and $\mathcal{B}$ are families of subsets of $\{1, \ldots, k\}$, that is $\mathcal{A}, \mathcal{B} \subseteq \powerset{\{1, \ldots, k\}}$, such that:
			\begin{enumerate}
				\item $\mathcal{A}$ and $\mathcal{B}$ are pairwise disjoint; and
				\item the union of $\mathcal{A}$ and $\mathcal{B}$ is a partition of $\{1, \ldots, k\}$.
			\end{enumerate} 
			We denote by $\uppat{\pttnmon{k}}$ the set of all $k$-patterns, and will refer to $k$-patterns more succinctly as \textit{patterns} either when $k$ may be any positive integer or when the value of $k$ is contextually unambiguous. Given a pattern $(\mathcal{A}, \mathcal{B}) \in \uppat{\pttnmon{k}}$, we refer to elements of $\mathcal{A}$ as \textit{non-transversal blocks} and elements of $\mathcal{B}$ as \textit{transversal blocks}.
		\end{definition}
		
		Let $k \in \nonnegints$ and $(\mathcal{A}, \mathcal{B})$ be a pattern. The reader should note that requiring $\mathcal{A}$ and $\mathcal{B}$ be pairwise disjoint does not exclude the possibility of there existing $A \in \mathcal{A}$ and $B \in \mathcal{B}$ such that $A \cap B$ is non-empty, however the additional requirement that the union $\mathcal{A} \cup \mathcal{B}$ partition $\{1, \ldots, k\}$ does, that is for each $A \in \mathcal{A}$ and $B \in \mathcal{B}$, $A \cap B$ is the empty set. Similarly, requiring that the union $\mathcal{A} \cup \mathcal{B}$ be a partition of $\{1, \ldots, k\}$ does not exclude the possibility of there existing $A \in \mathcal{A}$ and $B \in \mathcal{B}$ such that $A = B$, however the additional requirement that $\mathcal{A}$ and $\mathcal{B}$ be pairwise disjoint does.
		
		We may depict patterns in a similar way to how we use diagrams to depict bipartitions.

	 	\begin{definition} \label{def:graphicaldepictionofpatterns}
		 	Let $k \in \nonnegints$. Each pattern $(\mathcal{A}, \mathcal{B}) \in \uppat{\pttnmon{k}}$ may be \textit{depicted} as follows: 
		 	\begin{enumerate}
		 		\item for each $j \in \{1, \ldots, k\}$, the vertex $j$ is depicted as the point $(j, 0)$;
		 		\item lines connecting points in $\{1, \ldots, k\}\times\{0\}$ are drawn non-linearly either all above or all below the horizontal line $\{(x, 0): x \in \mathbb{R}\}$ and between the two vertical lines $\{(1, y), (k, y): y \in \mathbb{R}\}$ such that the connected components form the blocks of $\mathcal{A} \cup \mathcal{B}$; and
		 		\item so that blocks of $\mathcal{A}$ may be distinguished from blocks of $\mathcal{B}$, which is needed for every pattern to be uniquely recoverable from some depictions, for each block $B \in \mathcal{B}$, a vertical line is drawn downwards from one of the vertices in $B$. Note some authors alternatively use a two-tonne colouring of the vertices to indicate whether each vertex is contained in a non-transversal or transversal block.
		 	\end{enumerate}
		 	
		 	For example, Figure \ref{fig:patterneg} contains a depiction of the pattern $(\{\{2,3\}, \{4\}, \{6,7,8\}\}, \{\{1,5\}, \{9\}, \{10,11\}\})$ $\in \uppat{\pttnmon{11}}$.
		\end{definition}
		
		\begin{figure}[!ht]
			\caption[ ]{Given $k=11$, as outlined in Definition \ref{def:graphicaldepictionofpatterns}, the pattern $(\{\{2,3\}, \{4\}, \{6,7,8\}\}, \{\{1,5\}, \{9\},$ $\{10,11\}\}) \in \uppat{\pttnmon{11}}$ may be depicted as:}
			\label{fig:patterneg}
			\vspace{5pt}
			\centering
			\input{chap_background/tikz/fig-patterneg.tex}
		\end{figure}
		
		Similar to how distinct diagrams may depict the same bipartition, it is possible for distinct depictions to depict the same pattern. However we need not go in to as much detail with depictions of patterns since:
		\begin{enumerate}
			\item the details for depicting patterns follow in an analogous fashion to depicting bipartitions using diagrams; and
			\item the same amount of detail for depicting patterns will not be needed in our later discussions.
		\end{enumerate} 
		
		Typically two patterns, a lower and an upper, are associated with each bipartition, which may be outlined as follows.
	
		\begin{definition}
			For each $k \in \nonnegints$ and $\alpha \in \pttnmon{k}$, we denote by:
			\begin{enumerate}
				\item $\upnontrans{\alpha}$ the set of all upper non-transversal blocks in $\alpha$, that is $\upnontrans{\alpha} = \{b \in \alpha: \noupverts{b} > 0 \text{ and } \nolowverts{b} = 0\}$;
				\item $\lownontrans{\alpha}$ the set of all upper non-transversal blocks in $\alpha^*$, that is $\lownontrans{\alpha} = \upnontrans{\alpha^*} = \{b \in \alpha^*: \noupverts{b} > 0 \text{ and } \nolowverts{b} = 0\}$; 
				\item $\trans{\alpha}$ the set of all transversal blocks in $\alpha$, that is $\trans{\alpha} = \{b \in \alpha: \noupverts{b}, \nolowverts{b} > 0\}$;
				\item $\uptrans{\alpha}$ the set $\{\upverts{b}: b \in \trans{\alpha}\}$;
				\item $\lowtrans{\alpha}$ the set $\{\upverts{b}: b \in \trans{\alpha^*}\}$;
				\item $\uppat{\alpha}$ the pair $(\upnontrans{\alpha}, \uptrans{\alpha})$, which we refer to as the \textit{upper pattern of $\alpha$}; and
				\item $\lowpat{\alpha}$ the pair $(\lownontrans{\alpha}, \lowtrans{\alpha})$, which we refer to as the \textit{lower pattern of $\alpha$}.
			\end{enumerate}
			Furthermore given $A \subseteq \pttnmon{k}$, we denote by:
			\begin{enumerate}
				\item $\uppat{A}$ the set $\{\uppat{\alpha}: \alpha \in A\}$, which we refer to as the \textit{$A$-admissible upper patterns}; and
				\item $\lowpat{A}$ the set $\{\lowpat{\alpha}: \alpha \in A\}$, which we refer to as the \textit{$A$-admissible lower patterns}.
			\end{enumerate}
		\end{definition}
		
		\begin{proposition}
			For each $k \in \posints$ and $\alpha \in \pttnmon{k}$, $\uppat{\alpha}, \lowpat{\alpha} \in \uppat{\pttnmon{k}}$.
			
			\begin{proof}
				It is trivially the case that:
				\begin{enumerate}
					\item $\upnontrans{\alpha}$ and $\uptrans{\alpha}$ are pairwise disjoint, as are $\lownontrans{\alpha}$ and $\lowtrans{\alpha}$; and
					\item the unions $\upnontrans{\alpha} \cup \uptrans{\alpha}$ and $\lownontrans{\alpha} \cup \lowtrans{\alpha}$ each form a partition of $\{1, \ldots, k\}$.
				\end{enumerate} 
				Hence $\uppat{\alpha} = (\upnontrans{\alpha}, \uptrans{\alpha}), \lowpat{\alpha} = (\lownontrans{\alpha}, \lowtrans{\alpha}) \in \uppat{\pttnmon{k}}$, that is the upper and lower patterns of $\alpha$ are indeed patterns, justifying the terminology.
			\end{proof}
		\end{proposition}
		
		\begin{proposition}
			If a subsemigroup $S$ of the partition monoid $\pttnmon{k}$ is closed under the vertical flip involution $^*$ then $S$-admissible upper patterns and $S$-admissible lower patterns coincide.
			
			\begin{proof}
				For each $p \in \lowpat{S}$ there exist $\alpha \in S$ such that $\lowpat{\alpha} = p$. Since $S$ is closed under $^*$, we have $\alpha^* \in S$ and $\uppat{\alpha^*} = \lowpat{\alpha} = p$, and hence that $p \in \uppat{S}$. The converse follows analogously.
			\end{proof}
		\end{proposition}
		
		\begin{definition}
			Given a subsemigroup $S$ of the partition monoid $\pttnmon{k}$ that is closed under the vertical flip involution $^*$, we refer to $\uppat{S}$, which is equal to $\uppat{S} \cup \lowpat{S}$, as \textit{the $S$-admissible patterns}.
		\end{definition}
		
		For example if $k=8$ and $\alpha = \big\{\{1,5,2',3',6',7'\}, \{2,3,4\}, \{6,7,8,8'\}, \{1'\}, \{4',5'\}\big\} \in \pttnmon{8}$ then:
		\begin{enumerate}
			\item $\upnontrans{\alpha} = \big\{\{2,3,4\}\big\}$;
			\item $\lownontrans{\alpha} = \big\{\{1\}, \{4,5\}\big\}$;
			\item $\trans{\alpha} = \big\{\{1,5,2',3',6',7'\}, \{6,7,8,8'\}\big\}$;
			\item $\uptrans{\alpha} = \big\{\{1,5\}, \{6,7,8\}\big\}$;
			\item $\lowtrans{\alpha} = \big\{\{2,3,6,7\}, \{8\}\big\}$;
			\item $\uppat{\alpha} = \big(\upnontrans{\alpha}, \uptrans{\alpha}\big) = \big(\big\{\{2,3,4\}\big\}, \big\{\{1,5\}, \{6,7,8\}\big\}\big)$; and
			\item $\lowpat{\alpha} = \big(\lownontrans{\alpha}, \lowtrans{\alpha}\big) = \big(\big\{\{1\}, \{4,5\}\big\}, \big\{\{2,3,6,7\}, \{8\}\big\}\big)$.
		\end{enumerate}
		
		Note that $\alpha$, $\uppat{\alpha}$ and $\lowpat{\alpha}$ from above are depicted in Figure \ref{fig:upperlowerpatterneg}.
		
		\begin{figure}[!ht]
			\caption[ ]{Given $k = 8$,}
			\label{fig:upperlowerpatterneg}
			\vspace{5pt}
			\centering
			\input{chap_background/tikz/fig-upperlowerpatterneg.tex}
		\end{figure}
		
		Let $k \in \nonnegints$ and $\alpha \in \pttnmon{k}$. Graphical depictions of the upper and lower patterns of $\alpha$ may essentially be thought of as being formed by taking any diagram from $\diagrams{\alpha}$ to depict $\alpha$ and cutting it in half horizontally.
		
		However the reader should note that there do exist $\alpha \in \pttnmon{k}$ and $\delta_{\alpha} \in \diagrams{\alpha}$ such that simply cutting each of the transversal lines in $\delta_{\alpha}$ does not leave us with depictions of the upper and lower patterns of $\alpha$. Consequently it is also implicitly meant that lines are added so that upper vertices in the same connected component of $\delta_{\alpha}$ are in the same connected component of the upper cut diagram, and analogously so that lower vertices in the same connected component of $\delta_{\alpha}$ are in the same connected component of the lower cut diagram.
		
\begin{figure}[!ht]
			\caption[ ]{Given $k = 2$, $\alpha = \{\{1,2,1',2'\}\} \in \pttnmon{2}$ and $\delta_{\alpha} = \{(1,1'), (1',2')$, $(2,2')\} \in \diagrams{\alpha}$. If we cut each of the transversal lines in the depiction of $\delta_{\alpha}$ without additionally joining both upper vertices, which share the same connected component of $\delta_{\alpha}$ or equivalently the same block of $\alpha$, then we do not get a depiction of $\uppat{\alpha}$.}
			\label{fig:cutdiagram}
			\vspace{5pt}
			\centering
			\input{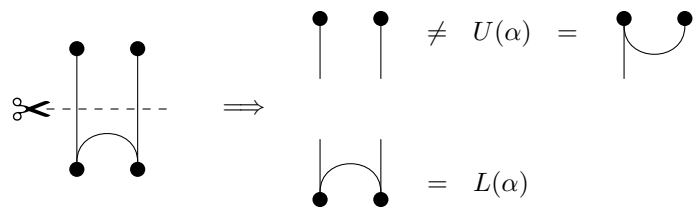}
		\end{figure}

		For example consider $k=2$, $\alpha = \{\{1,2,1',2'\}\} \in \pttnmon{2}$ and $\delta_{\alpha} = \{(1,1')$, $(1',2'), (2,2')\} \in \diagrams{\alpha}$. Figure $\ref{fig:cutdiagram}$ depicts that were we to simply cut each transversal line in the depiction of $\delta_{\alpha}$, the resulting depiction of an upper pattern is not a depiction of the upper pattern $\uppat{\alpha}$, though the resulting depiction of a lower pattern is a depiction of the lower pattern $\lowpat{\alpha}$. In order to obtain a depiction of the upper pattern $\uppat{\alpha}$ when cutting $\delta_{\alpha}$ in half horizontally, we additionally need to join both upper vertices.
		
		\begin{proposition}
			For each $k \in \nonnegints$ and $\alpha \in \pttnmon{k}$, $\uppat{\alpha^*} = \lowpat{\alpha}$ and dually $\lowpat{\alpha^*} = \uppat{\alpha}$.
			
			\begin{proof}
				Follows directly from how $\lowpat{\alpha}$ was defined.
			\end{proof}
		\end{proposition}
		
		Furthermore, since the partition monoid $\pttnmon{k}$ is closed under the vertical flip involution $^*$, that is $\pttnmon{k}^* = \pttnmon{k}$, it follows that $\lowpat{\pttnmon{k}} = \uppat{\pttnmon{k}}$, hence why we denote the set of all $k$-patterns simply as $\uppat{\pttnmon{k}}$ rather than $\uppat{\pttnmon{k}}\cup\lowpat{\pttnmon{k}}$.
		
		\begin{proposition}
			For each $k \in \nonnegints$ and $\alpha, \beta \in \pttnmon{k}$, $\upnontrans{\alpha} \subseteq \upnontrans{\alpha\beta}$ and $\lownontrans{\beta} \subseteq \lownontrans{\alpha\beta}$.
			
			\begin{proof}
				Trivially obvious that none of $\beta$'s blocks can interact with $\alpha$'s upper non-transversal blocks and that none of $\alpha$'s blocks can interact with $\beta$'s lower non-transversal blocks when forming the product $\alpha\beta$ as the connected components of two arbitrary diagrams $\delta_{\alpha} \in \diagrams{\alpha}$ and $\delta_{\beta} \in \diagrams{\beta}$, hence $\upnontrans{\alpha} \subseteq \upnontrans{\alpha\beta}$ and $\lownontrans{\beta} \subseteq \lownontrans{\alpha\beta}$.
			\end{proof}
		\end{proposition}

		\begin{proposition} \label{prop:productuplowpatwithsamemiddlepat}
			For each $k \in \nonnegints$ and $\alpha, \beta \in \pttnmon{k}$, if $\lowpat{\alpha} = \uppat{\beta}$ then $\uppat{\alpha\beta} = \uppat{\alpha}$ and $\lowpat{\alpha\beta} = \lowpat{\beta}$. \qed
		\end{proposition}

	\subsection{Ranks}
	 	\begin{definition}
			Let $k \in \nonnegints$ and $\alpha \in \pttnmon{k}$. The \textit{rank of $\alpha$}, which is often denoted as $\rank{\alpha}$, is the number of transversal blocks that $\alpha$ contains, that is $\rank{\alpha} = \card{\trans{\alpha}}$.
	 	\end{definition}
	 	
	 	\begin{proposition}
		 	For each $k \in \nonnegints$ and $\alpha, \beta \in \pttnmon{k}$,
		 	\begin{enumerate}
		 		\item $\rank{\alpha^*} = \rank{\alpha}$;
		 		\item $\rank{\alpha\oplus\beta} = \rank{\alpha} + \rank{\beta}$;
			 	\item $\rank{\alpha\beta} \leq \rank{\alpha}, \rank{\beta}$;
			 	\item $\rank{\alpha\beta} = \rank{\alpha}$ if and only if $\uppat{\alpha\beta} = \uppat{\alpha}$; and
			 	\item $\rank{\alpha\beta} = \rank{\beta}$ if and only if $\lowpat{\alpha\beta} = \lowpat{\beta}$.
		 	\end{enumerate}
		 	
		 	\begin{proof}
		 		Vertically flipping a bipartition trivially does not alter the rank, that is $\rank{\alpha^*} = \rank{\alpha}$, and the rank of the vertical sum of bipartitions is trivially equal to the sum of their ranks, that is $\rank{\alpha \oplus \beta} = \rank{\alpha} + \rank{\beta}$, establishing Conditions (i) and (ii). There are two ways in which $\rank{\alpha\beta}$ may differ from $\rank{\alpha}$:
		 		\begin{enumerate}[label=\Roman*]
					\item some of $\beta$'s transversal blocks, possibly along with some of $\beta$'s upper non-transversal blocks and $\alpha$'s lower non-transversal blocks, may join some of $\alpha$'s transversal blocks; and
                   	\item some of $\beta$'s upper non-transversal blocks, along with possibly some of $\alpha$'s lower non-transversal blocks, may cap some of $\alpha$'s transversal blocks.
               	\end{enumerate}
               	Neither I or II allows us to split a transversal block of $\alpha$ so that $\rank{\alpha} < \rank{\alpha\beta}$, giving Condition (iii). If $\rank{\alpha\beta} = \rank{\alpha}$ then neither I or II could have occurred, hence we must have $\uppat{\alpha\beta} = \uppat{\alpha}$. Conversely if $\uppat{\alpha\beta} = \uppat{\alpha}$ then we trivially must have $\rank{\alpha\beta} = \rank{\alpha}$, giving Condition (iv). Condition (v) is trivially the dual of Condition (iv).
		 	\end{proof} 
	 	\end{proposition}
	 	
	\subsection{Green's relations}

		\begin{theorem} \label{thm:greensrltnsforpttnmon}  (see \citep{art:Fitzgerald:OnThePartitionMonoid} \citep{art:Wilcox:CullularityDiagramAlgebras})
			For each $k \in \posints$ and $\alpha, \beta \in \pttnmon{k}$:
			\begin{enumerate}
				\item $(\alpha, \beta) \in \mathcal{R}$ if and only if $\uppat{\alpha} = \uppat{\beta}$;
				\item $(\alpha, \beta) \in \mathcal{L}$ if and only if $\lowpat{\alpha} = \lowpat{\beta}$;
				\item $(\alpha, \beta) \in \mathcal{H}$ if and only if $\uppat{\alpha} = \uppat{\beta}$ and $\lowpat{\alpha} = \lowpat{\beta}$; and
				\item $(\alpha, \beta) \in \mathcal{J}$ if and only if $\uppat{\alpha} = \uppat{\beta}$.
			\end{enumerate}
		\end{theorem}

		\begin{corollary} \label{cor:greensrltnsforregulardiagramsemigrps}
			If $S$ is a subsemigroup of the partition monoid that is closed under $^*$ then for each $\alpha, \beta \in S$:
			\begin{enumerate}
				\item $(\alpha, \beta) \in \mathcal{R}$ if and only if $\uppat{\alpha} = \uppat{\beta}$;
				\item $(\alpha, \beta) \in \mathcal{L}$ if and only if $\lowpat{\alpha} = \lowpat{\beta}$;
				\item $(\alpha, \beta) \in \mathcal{H}$ if and only if $\uppat{\alpha} = \uppat{\beta}$ and $\lowpat{\alpha} = \uppat{\beta}$;
				\item $\uppat{\alpha^*\alpha} = \uppat{\beta^*\beta}$ if and only if $\lowpat{\alpha^*\alpha} = \lowpat{\beta^*\beta}$ if and only if $\alpha^*\alpha = \beta^*\beta$; and
				\item $\uppat{\alpha\alpha^*} = \uppat{\beta\beta^*}$ if and only if $\lowpat{\alpha\alpha^*} = \lowpat{\beta\beta^*}$ if and only if $\alpha\alpha^* = \beta\beta^*$.
			\end{enumerate}
		
			\begin{proof}
				Coniditions (i), (ii) and (iii) follow by applying Proposition \ref{prop:greensrltnsforregularsubmons} to Theorem \ref{thm:greensrltnsforpttnmon}. It follows from Condition (i) that $\uppat{\alpha^*\alpha} = \uppat{\beta^*\beta}$ implies $(\alpha^*\alpha, \beta^*\beta) \in \mathcal{R}$. Since $\alpha^*\alpha$ and $\beta^*\beta$ are projections, and each $\mathcal{R}$ class contains precisely one projection, we must have $\alpha^*\alpha = \beta^*\beta$. The converse is trivial and equivalence with $\lowpat{\alpha^*\alpha} = \lowpat{\beta^*\beta}$ follows dually, from which Condition (iv) follows. Condition (v) follows analogously to Condition (iv).
			\end{proof}
		\end{corollary}
	 	
\section{Contextually relevant diagram semigroups}
	\begin{definition}
		By a \textit{diagram semigroup} we shall mean any subsemigroup of the partition monoid.
	\end{definition}
	
	A number of diagram semigroups will either be directly relevant at various stages or simply useful as having similar properties to the diagram semigroups that will be introduced and investigated later in our discussion. 	
	
	\subsection{Generating sets}
		So that readers may conveniently refer back to one place, we next define all generating sets that will appear at some point through the remainder of the thesis. 
		
\begin{figure}[!ht]
			\caption[ ]{Given $k = 4$, the three transposition generators may be depicted as:}
			\label{fig:k=4transpositiongens}
			\vspace{5pt}
			\centering
			\input{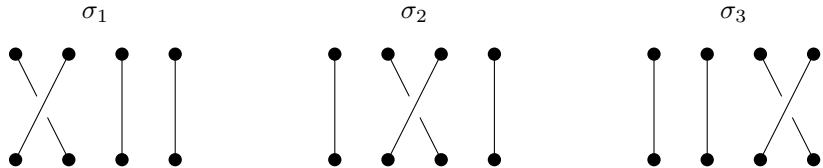}
		\end{figure}

		\begin{definition} \label{def:transpositiongens}
			For each $k \in \posints$ and $i \in \{1, \ldots, k-1\}$, we denote by $\symgen{i}$ the bipartition containing:
			\begin{enumerate}
				\item the two transversal lines $\{i, (i+1)'\}$ and $\{i+1, i'\}$; and
				\item for each $j \in \{1, \ldots, i-1, i+2, \ldots, k\}$, the vertical line $\{j, j'\}$.
			\end{enumerate}
			We refer to the $k-1$ bipartitions $\Big\{\symgen{i}: i \in \{1, \ldots, k-1\}\Big\}$ as \textit{the transposition generators} (see Figure \ref{fig:k=4transpositiongens} for an example).
		\end{definition}
		
		\begin{definition} \label{def:transapgens}
			For each $k \in \posints$ and $i \in \{1, \ldots, k-1\}$, we denote by $\transapgen{j}$ the bipartition containing:
			\begin{enumerate}
				\item the $(2,2)$-transapsis $\{i, i+1, i', (i+1)'\}$; and
				\item for each $j \in \{1, \ldots, i-1, i+2, \ldots, k\}$, the vertical line $\{j, j'\}$.
			\end{enumerate}
			We refer to the $k-1$ bipartitions $\Big\{\transapgen{i}: i \in \{1, \ldots, k-1\}\Big\}$ as \textit{the $(2,2)$-transapsis generators} (see Figure \ref{fig:k=4(2,2)-transapgens} for an example).
		\end{definition}
		
		\begin{figure}[!ht]
			\caption[ ]{Given $k = 4$, the three $(2,2)$-transapsis generators may be depicted as:}
			\label{fig:k=4(2,2)-transapgens}
			\vspace{5pt}
			\centering
			\input{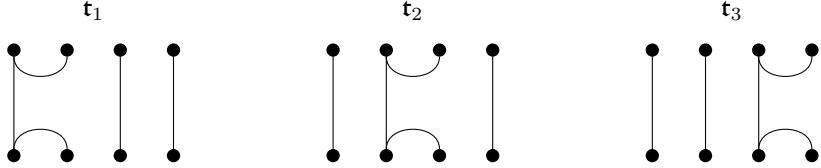}
		\end{figure}
		
\begin{figure}[!ht]
			\caption[ ]{Given $k = 4$, the four monapsis generators may be depicted as:}
			\label{fig:k=4monapsisgens}
			\vspace{5pt}
			\centering
			\input{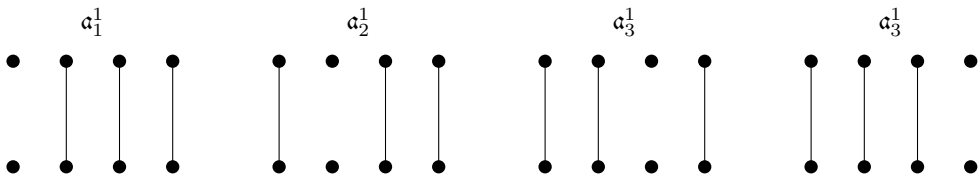}
		\end{figure}

\begin{figure}[!ht]
			\caption[ ]{Given $k = 5$, the three triapsis generators may be depicted as:}
			\label{fig:k=5triapsisgens}
			\vspace{5pt}
			\centering
			\input{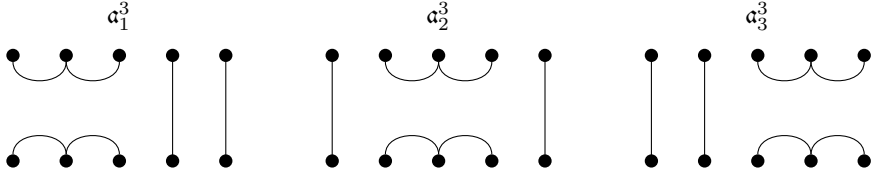}
		\end{figure}

		\begin{definition} \label{def:mapsisgens}
			For each $m \in \posints$, $k \in \intsge{m}$ and $i \in \{1, \ldots, k-m+1\}$, we denote by $\apgen{m}{i}$ the bipartition containing:
			\begin{enumerate}
				\item the two $m$-apses $\{i, \ldots, i+m-1\}$ and $\{i', \ldots, (i+m-1)'\}$; and
				\item for each $j \in \{1, \ldots, i-1, i+m, \ldots, k\}$, the vertical line $\{j, j'\}$.
			\end{enumerate}
			We refer to the $k-m+1$ bipartitions $\Big\{\apgen{m}{i}: i \in \{1, \ldots, k-m+1\}\Big\}$ as \textit{the $m$-apsis generators}. Recall from Definition \ref{def:blocktypes} that we are referring to $1$-apses as monapses, $2$-apses as diapses and $3$-apses as triapses, consequently we shall refer to $1$-apsis generators as \textit{monapsis generators}, $2$-apsis generators as \textit{diapsis generators} and $3$-apsis generators as \textit{triapsis generators} (see Figures \ref{fig:k=4monapsisgens}, \ref{fig:k=4diapsisgens} and \ref{fig:k=5triapsisgens} for examples). Furthermore, when convenient we denote the $i$th diapsis generator $\apgen{2}{i}$ as $\hookgen{i}$.
		\end{definition}
		
		\begin{figure}[!ht]
			\caption[ ]{Given $k = 4$, the three diapsis generators may be depicted as:}
			\label{fig:k=4diapsisgens}
			\vspace{5pt}
			\centering
			\input{chap_background/tikz/fig-k=4diapsisgens.tex}
		\end{figure}
		
		\begin{definition} \label{def:partialtranspositiongens}
			For each $k \in \posints$ and $i \in \{1, \ldots, k-1\}$, we denote by:
			\begin{enumerate}
				\item $\ftransgen{i}$ the bipartition containing:
				\begin{enumerate}
					\item for each $j \in \{1, \ldots, i-1\}$, the vertical line $\{j, j'\}$;
					\item the monapses $\{k\}$ and $\{i'\}$; and
					\item for each $j \in \{i, \ldots, k-1\}$, the line $\{j, (j+1)'\}$, and
				\end{enumerate}
				\item $\btransgen{i}$ the bipartition containing:
				\begin{enumerate}
					\item for each $j \in \{1, \ldots, i-1\}$, the vertical line $\{j, j'\}$;
					\item the monapses $\{i\}$ and $\{k'\}$; and
					\item for each $j \in \{i, \ldots, k-1\}$, the line $\{j+1, j'\}$.
				\end{enumerate}
			\end{enumerate}
			We refer to the $2k-2$ bipartitions $\Big\{\ftransgen{i}, \btransgen{i}: i \in \{1, \ldots, k-1\}\Big\}$ as \textit{the $\psyminvmon{k}$ generators} (see Figure \ref{fig:k=2partialtranspositiongens} for an example), where $\psyminvmon{k}$ is the planar symmetric inverse monoid, defined later in Definition \ref{def:psyminvmon}.
		\end{definition}
		
		\begin{figure}[!ht]
			\caption[ ]{Given $k = 3$, the four $\psyminvmon{3}$ generators may be depicted as:}
			\label{fig:k=2partialtranspositiongens}
			\vspace{5pt}
			\centering
			\input{chap_background/tikz/fig-k=3partialtranspositiongens.tex}
		\end{figure}

	\subsection{The planar partition monoid $\ppttnmon{k}$} \label{sec:ppttnmon}
		\begin{definition}
			Let $k \in \posints$. A diagram $\delta \in \diagrams{k}$ is referred to as \textit{planar} if it may be depicted without two distinct connected components crossing and \textit{non-planar} otherwise. We denote by $\pdiagrams{k}$ the set of all planar $k$-diagrams.
		\end{definition}
		
		\begin{definition}
			Let $k \in \posints$. A bipartition $\alpha \in \pttnmon{k}$ is referred to as \textit{planar} when the diagrams $\diagrams{\alpha}$ that depict $\alpha$ are planar, and \textit{non-planar} when the diagrams $\diagrams{\alpha}$ that depict $\alpha$ are non-planar. We denote by $\ppttnmon{k}$ the set of all planar bipartitions.
		\end{definition}
		
		\begin{proposition} \label{prop:diagramsalphaeitherplanarornonplanar}
			Let $k \in \posints$ and $\alpha \in \pttnmon{k}$. Either $\alpha$ is planar or $\alpha$ is non-planar, that is either every diagram in $\diagrams{\alpha}$ is planar or every diagram in $\diagrams{\alpha}$ is non-planar. \qed
		\end{proposition}
		
		
		Without giving all the gory details, Proposition \ref{prop:diagramsalphaeitherplanarornonplanar} may be established by arguing that given $k \in \posints$ and $b \subseteq \{1, \ldots, k\}\cup\{1', \ldots, k'\}$, every connected graph on $b$ induces the same partition of $\setdiff{\{1, \ldots, k\}\cup\{1', \ldots, k'\}}{b}$, and the partition induced by each depiction of the block $b$ within the convex hull of $\{1, \ldots, k\}\times\{0, 1\}$ is the coarsest partition finer than the partitions induced by distinct $b_1, b_2 \in b$. 
		
		\begin{proposition}
			For each $k \in \posints$, the set of planar bipartitions $\ppttnmon{k}$ is a submonoid of the partition monoid $\pttnmon{k}$.
			
			\begin{proof}
				Let $\alpha, \beta \in \ppttnmon{k}$, $\delta_{\alpha} \in \diagrams{\alpha}$ and $\delta_{\beta} \in \diagrams{\beta}$. In the product graph $\productgraph{\delta_{\alpha}, \delta_{\beta}}$, none of the connected components cross in the upper or lower half since $\delta_{\alpha}$ and $\delta_{\beta}$ are planar, and any connected components that join at one of the vertices in the middle row become part of the same connected component in the product $\alpha\beta$. 
			\end{proof}
		\end{proposition}
		
		\begin{definition} \label{def:ppttnmon}
			For each $k \in \posints$, the monoid of planar bipartitions $\ppttnmon{k}$ is referred to as \textit{the planar partition monoid}.
		\end{definition}
		
		\begin{proposition} \label{prop:ppttnmonpresentation} (see \cite{art:Halverson:PartitionAlgebras})
			For each $k \in \nonnegints$, the planar partition monoid $\ppttnmon{k}$ is characterised by the generators $\set{\monapgen{k}, \monapgen{i}, \transapgen{i}: i = 1, \ldots, k-1}$ along with the relations:
			\begin{enumerate}
				\item $\monapgen{i}\monapgen{i} = \monapgen{i}$;
				\item $\monapgen{j}\monapgen{i} = \monapgen{i}\monapgen{j}$ for all $|j - i| \geq 1$;
				\item $\transapgen{i}^2 = \transapgen{i}$;
				\item $\transapgen{j}\transapgen{i} = \transapgen{i}\transapgen{j}$ for all $|j - i| \geq 1$; and
				\item $\monapgen{i}\transapgen{j}\monapgen{i} = \monapgen{i}$ for all $|j - i| = 1$. \qed
			\end{enumerate}
		\end{proposition}

	\subsection{The symmetric group $\symgrp{k}$}
		The symmetric group, which really ought to need no introduction at all, is most commonly known as either the permutations of a set or equivalently the bijections from a set to itself. In the context of diagram semigroups, the symmetric group may be defined as follows.
	
		\begin{definition}
			For each $k \in \nonnegints$, \textit{the symmetric group}, which we denote as $\symgrp{k}$, is the set of bipartitions such that every block is a transversal line.
		\end{definition}
		
		\begin{proposition} (see \cite{art:Moore:SymGrpPresentation})
			For each $k \in \posints$, the symmetric group $\symgrp{k}$ is characterised by the transposition generators $\Big\{\symgen{i}: i \in \{1, \ldots, k-1\}\Big\}$ along with the relations:
			\begin{enumerate}
				\item $\symgen{i}^2 = 1$; 
				\item $\symgen{i}\symgen{i+1}\symgen{i} = \symgen{i+1}\symgen{i}\symgen{i+1}$; and
				\item $\symgen{i}\symgen{k} = \symgen{k}\symgen{i}$ for all $\card{k-i} > 1$. \qed
			\end{enumerate}
		\end{proposition}
		
		\begin{definition}
			For each $k, x \in \posints$ and $\mu_1, \ldots, \mu_x, \gamma_1, \ldots, \gamma_x \in \{0, \ldots, k\}$ such that $\Sigma^x_{i=1}\mu_i = \Sigma^x_{i=1}\gamma_i = k$ and $\mu_i + \gamma_i > 0$ for all $i \in \{1, \ldots, x\}$, we denote by $\Psi^{(\mu_1, \gamma_1), \ldots, (\mu_x, \gamma_x)}$ the set of all bipartitions $\alpha = \{b_1, \ldots, b_x\} \in \pttnmon{k}$ such that for each $i \in \{1, \ldots, x\}$, $b_i$ is a block of type $(\mu_i, \gamma_i)$.
		\end{definition}
		
		\begin{proposition} \label{prop:symgrpedgetypes}
			For each $k, x \in \posints$ and $\mu_1, \ldots, \mu_x, \gamma_1, \ldots, \gamma_x \in \{0, \ldots, k\}$ such that $\Sigma^x_{i=1}\mu_i = \Sigma^x_{i=1}\gamma_i = k$ and $\mu_i + \gamma_i > 0$ for all $i \in \{1, \ldots, x\}$:
			\begin{enumerate}
			\item $\Psi^{(\mu_1, \gamma_1), \ldots, (\mu_x, \gamma_x)} \cap \ppttnmon{k}$ is non-empty; and
			\item $\Psi^{(\mu_1, \gamma_1), \ldots, (\mu_x, \gamma_x)} = \symgrp{k}\psi\symgrp{k}$ for all $\psi \in \Psi^{(\mu_1, \gamma_1), \ldots, (\mu_x, \gamma_x)}$. \qed
			\end{enumerate}
			
			\begin{proof}
				\begin{enumerate}
					\item $\{\{1, \ldots, \mu_1, 1', \ldots, \gamma_1'\}, \ldots, \{\Sigma^{x-1}_{i=1}\mu_i + 1, \ldots, k, (\Sigma^{x-1}_{i=1}\gamma_i + 1)', k'\}\} \in \Psi^{(\mu_1, \gamma_1), \ldots, (\mu_x, \gamma_x)} \cap \ppttnmon{k}$.
					\item It is trivially the case that multiplication on the left by permutations is equivalent to permuting the upper vertices, and similarly that multiplication on the right by permutations is equivalent to permuting the lower vertices. 
				\end{enumerate}
			\end{proof}
		\end{proposition}

	\subsection{The Jones monoid $\jonesmon{k}$}
		\begin{definition} \label{def:jonesmon}
			For each $k \in \nonnegints$, \textit{the Jones monoid}, which is often denoted as $\jonesmon{k}$ and also commonly known as the Temperley-Lieb monoid, is the set of planar bipartitions such that every block contains precisely two vertices.
		\end{definition}
		
		\begin{proposition}
			The cardinality of the Jones monoid $\card{\jonesmon{k}}$ is equal to the $k$th Catalan number (sequence $A000108$ on the OEIS \cite{man:OEIS}). \qed
		\end{proposition}
		
		Recall from Definition \ref{def:mapsisgens} that we will be denoting the $i$th diapsis generator $\apgen{2}{i}$ as $\hookgen{i}$ whenever it is convenient to do so.
	
		\begin{proposition} \label{prop:jonesmonpresentation}
			For each $k \in \nonnegints$, the Jones monoid $\jonesmon{k}$ is characterised by the generators $\set{\hookgen{i}: i = 1, \ldots, k-1}$ along with the relations:
			\begin{enumerate}
				\item $\hookgen{i}^2 = \hookgen{i}$;
				\item $\hookgen{i}\hookgen{j}\hookgen{i} = \hookgen{j}\hookgen{i}\hookgen{j}$ for all $|j-i| = 1$; and
				\item $\hookgen{j}\hookgen{i} = \hookgen{i}\hookgen{j}$ for all $|j-i| \geq 2$. \qed
			\end{enumerate}
		\end{proposition}
		
		First we review how Ridout and Saint-Aubin \cite{art:Ridout:StandardModules} establish Proposition \ref{prop:jonesmonpresentation}.
		
		\begin{definition}
			We refer to elements of the free semigroup of the Jones monoid $\freesemigrp{{\jonesmon{k}}}$ as \textit{$\jonesmon{k}$-words}, and say that a $\jonesmon{k}$-word is \textit{reduced} if it may not be written with fewer generators using the relations from Proposition \ref{prop:jonesmonpresentation}.
		\end{definition}
		
		\begin{proposition}
			In any reduced $\jonesmon{k}$-word $\hookgen{i_1}\ldots\hookgen{i_n}$, the maximal index $m = \max\big\{i_j: j \in \{1, \ldots, n\}\big\}$ occurs precisely once.
			\begin{proof}
				See Lemma $2.2$ of \cite{art:Ridout:StandardModules}.
			\end{proof}
		\end{proposition}
		
		\begin{proposition} \label{prop:firststepnormalisingreducedjonesword}
			If $W$ is a reduced $\jonesmon{k}$-word with maximal index $m$ then $W$ may be rewritten as $W = W'\hookgen{m}\ldots\hookgen{l}$ where $W'$ is a reduced $\pmodmon{2}{k}$-word with maximal index less than $m$ and $l \in \set{1, \ldots, m}$.
			
			\begin{proof}
				See \cite{art:Ridout:StandardModules}.
			\end{proof}
		\end{proposition}
		
		\begin{proposition} \label{prop:jonesnormalform}
			For each $k \in \intsge{2}$, any reduced $\jonesmon{k}$-word $W$ may be rewritten as $W = r_{j_1, i_1}\ldots r_{j_n, i_n}$ where:
			\begin{enumerate}
				\item $n \in \posints$ and $i_1, \ldots, i_n, j_1, \ldots, j_n \in \set{1, \ldots, k-1}$ such that: \begin{enumerate}
					\item for each $l \in \set{1, \ldots, n}$, $j_l \geq i_l$; and
					\item for each $l \in \set{1, \ldots, n-1}$, $i_l < i_{l+1}$ and $j_l < j_{l+1}$, and
				\end{enumerate}
				\item for each $l \in \set{1, \ldots, n}$, $r_{j_l, i_l} = \hookgen{j_l}\ldots\hookgen{i_l}$.
			\end{enumerate}
			
			\begin{proof}
				See Proposition 2.3 (Jones' Normal Form) in \cite{art:Ridout:StandardModules}.
			\end{proof}
		\end{proposition}
		
		\begin{definition} \label{def:Jonesnormalform}
			$\jonesmon{k}$-words written in the form for Proposition \ref{prop:jonesnormalform} are referred to as being in \textit{normal form}.
		\end{definition}
		
		To establish Proposition \ref{prop:jonesmonpresentation}, it is now sufficient to establish that the number of $\jonesmon{k}$-words in normal form is equal to the cardinality of the Jones monoid $\card{\jonesmon{k}}$. Ridout and Saint-Aubin did so by establishing a bijection between $\jonesmon{k}$-words in normal form and lattice paths from the origin to $(k, k)$ that do not step above $y=x$ using the step-set $\set{(1, 0), (0, 1)}$, while also noting that the number of such lattice paths also form the Catalan numbers. We provide an alternative counting argument here, which is similar to how we shall bound reduced $\pmodmon{2}{k}$-words in Subsection \ref{subsec:pmod2boundednormals} later on.
		
		\begin{definition}
			For each $k \in \intsge{2}$ and $i, j \in \set{1, ..., k-1}$ such that $i < j$, we denote by $r_{j, i}$ the product $\hookgen{j}\ldots\hookgen{i}$ of diapsis genators, and refer to $r_{j, i}$ as \textit{the run of diapsis generators from $j$ to $i$}.
		\end{definition}
		
\begin{figure}[!ht]
			\caption[ ]{Runs of diapsis generators ordered by $\prec$:}
			\label{fig:runHassediagram}
			\vspace{5pt}
			\centering
			\input{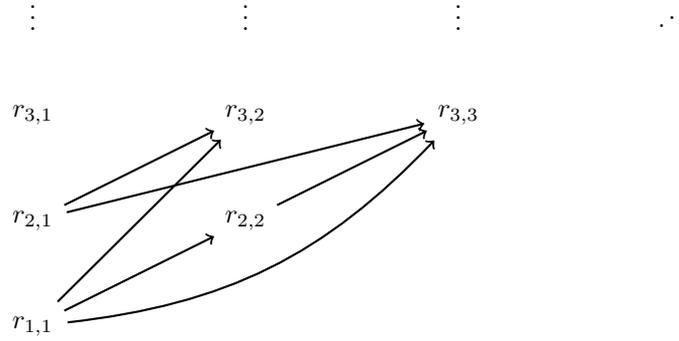}
		\end{figure}

		\begin{definition}
			For each $k \in \posints$, we define a binary relation $\prec$ on runs of diapsis generators by $r_{j, i} \prec r_{j', i'}$ if and only if $j < j'$ and $i < i'$ (see Figure \ref{fig:runHassediagram} for a depiction).
		\end{definition}
		
		Note that $\prec$ is trivially transitive, that is if $r_{j, i} \prec r_{j', i'}$ and $r_{j', i'} \prec r_{j'', i''}$ then $r_{j, i} \prec r_{j'', i''}$. Further note that a product of runs of diapsis generators $r_{j_1, i_1}\ldots r_{j_n, i_n}$ is a $\jonesmon{k}$-word in normal form if and only if for each $l \in \set{1, ..., n-1}$, $r_{j_l, i_l} \prec r_{j_{l+1}, i_{l+1}}$.
		
		\begin{definition} \label{def:noJknormalformsendingwithrun}
			For each $i, j \in \posints$ such that $i < j$, we denote by $R_{j, i}$ the number of $\jonesmon{k}$-words in normal form that end with the run $r_{j, i}$.
		\end{definition}
		
		For example the $\jonesmon{k}$-words in normal form that end with the run $r_{3,3}$ are $\{r_{3, 3}, r_{2,1}r_{3,3}, r_{2, 2}r_{3,3}, r_{1,1}r_{2,2}r_{3,3},$ $r_{1,1}r_{3,3}\}$, and hence $R_{3,3} = 5$.
		
		\begin{proposition}
			For each $i, j \in \posints$ such that $i \leq j$, 
			\[R_{j, i} = \begin{cases}
				1 & i = 1; \text{ and} \\
				1 + \Sigma^{j-1}_{j'=1}\Sigma^{j'}_{i'=1}R_{j', i'} & i > 1. \\
			\end{cases}\]
			
			\begin{proof}
				Since no $\jonesmon{k}$-words in normal form may have a run before $r_{j, 1}$, $r_{j, 1}$ is the only $\jonesmon{k}$-word in normal form that ends with run $r_{j, 1}$, hence $R_{j, 1} = 1$. When $i > 1$, $r_{j, i}$ is trivially a $\jonesmon{k}$-word in normal form, then for each $j' \in \set{1, \ldots, j-1}$ and $i' \in \set{1, \ldots, j'}$, each $\jonesmon{k}$-word in normal form ending with the run $r_{j', i'}$ is still a $\jonesmon{k}$-word in normal form when multiplied on the right by $r_{j, i}$, and each word formed when doing so is unique.
			\end{proof}
		\end{proposition}
		
		Table \ref{table:Rjivalues} contains computed values of $R_{j, i}$ for $i, j \in \set{1, \ldots, 10}$ such that $i \leq j$, which together form Catalan's triangle (sequence $A009766$ on the OEIS \cite{man:OEIS}) with the right border removed (sequence $A030237$ on the OEIS \cite{man:OEIS}). 
		
		\begin{table}[!ht]
			\caption[ ]{Number of planar mod-$2$ normal form words ending with a run from $\Run{j}{i}$.}
			\label{table:Rjivalues}
			\centering
			\begin{tabular}{| c | r r r r r r r r r r |}
				\hline
				\diagbox{$j$}{$i$} & 1 & 2 & 3 & 4 & 5 & 6 & 7 & 8 & 9 & 10 \\
				\hline
				1  & 1 &    &    &     &     &      &      &      &       &       \\
				2  & 1 & 2  &    &     &     &      &      &      &       &       \\
				3  & 1 & 3  & 5  &     &     &      &      &      &       &       \\
				4  & 1 & 4  & 9  & 14  &     &      &      &      &       &       \\
				5  & 1 & 5  & 14 & 28  & 42  &      &      &      &       &       \\
				6  & 1 & 6  & 20 & 48  & 90  & 132  &      &      &       &       \\
				7  & 1 & 7  & 27 & 75  & 165 & 297  & 429  &      &       &       \\
				8  & 1 & 8  & 35 & 110 & 275 & 572  & 1001 & 1430 &       &       \\
				9  & 1 & 9  & 44 & 154 & 429 & 1001 & 2002 & 3432 & 4862  &       \\
				10 & 1 & 10 & 54 & 208 & 637 & 1638 & 3640 & 7072 & 11934 & 16796 \\
				\hline
			\end{tabular}
		\end{table}
		
		Note that the diagonal of Catalan's triangle is the Catalan numbers. Further note that the number of $\jonesmon{k}$-words in normal form is trivially equal to $1 + \Sigma^{k-1}_{j=1}\Sigma^{j}_{i=1}R_{j,i} = R_{k, k}$, which provides an alternative argument than that given by Ridout and Saint-Aubin \cite{art:Ridout:StandardModules} for the number of $\jonesmon{k}$-words in normal form being equal to the cardinality of the Jones monoid $\card{\jonesmon{k}}$, establishing Proposition \ref{prop:jonesmonpresentation}. This also provides us with a constructive way to recursively generate the $\jonesmon{k}$-words in normal form in an abstract manner.

	\subsection{The Brauer monoid $\brauermon{k}$} \label{subsec:brauermon}
		\begin{definition} \label{def:brauermon}
			For each $k \in \nonnegints$, \textit{the Brauer monoid}, which is often denoted as $\brauermon{k}$, is the set of bipartitions such that every block contains precisely two vertices.
		\end{definition}
		
		\begin{proposition} \label{prop:brauermonpresentation1} (see \cite{art:Kudryavtseva:PresentationsBrauerTypeMonoids})
			For each $k \in \nonnegints$, the Brauer monoid $\brauermon{k}$ is characterised by the generators $\set{\symgen{i}, \hookgen{i}: i = 1, \ldots, k-1}$ along with the relations:
			\begin{enumerate}
				\item $\symgen{i}^2 = \id{k}$;
				\item $\symgen{i+1}\symgen{i}\symgen{i+1} = \symgen{i}\symgen{i+1}\symgen{i}$;
				\item $\symgen{j}\symgen{i} = \symgen{i}\symgen{j}$ for all $j-i \geq 2$;
				\item $\hookgen{i}^2 = \hookgen{i}$;
				\item $\hookgen{i}\hookgen{j}\hookgen{i} = \hookgen{i}$ for all $|j-i| = 1$;
				\item $\hookgen{j}\hookgen{i} = \hookgen{i}\hookgen{j}$ for all $j-i \geq 2$;
				\item $\symgen{i}\hookgen{i} = \hookgen{i} = \hookgen{i}\symgen{i}$;
				\item $\hookgen{i}\hookgen{j}\symgen{i} = \hookgen{i}\symgen{j}$ for all $|j-i| = 1$;
				\item $\symgen{i}\hookgen{j}\hookgen{i} = \symgen{j}\hookgen{i}$ for all $|j-i| = 1$; and
				\item $\hookgen{j}\symgen{i} = \symgen{i}\hookgen{j}$ for all $|j-i| \geq 2$. \qed
			\end{enumerate}
		\end{proposition}
	
		It turns out that Relations (vii), (viii), (ix) and (x) from Proposition \ref{prop:brauermonpresentation1} are stricter than required, which we establish in the following proposition.
		
		\begin{proposition} \label{prop:brauermonpresentation2}
			For each $k \in \nonnegints$, the Brauer monoid $\brauermon{k}$ is characterised by the generators $\set{\symgen{i}, \hookgen{i}: i = 1, \ldots, k-1}$ along with the relations:
			\begin{enumerate}
				\item $\symgen{i}^2 = \id{k}$;
				\item $\symgen{i+1}\symgen{i}\symgen{i+1} = \symgen{i}\symgen{i+1}\symgen{i}$;
				\item $\symgen{j}\symgen{i} = \symgen{i}\symgen{j}$ for all $j-i \geq 2$;
				\item $\hookgen{i}^2 = \hookgen{i}$;
				\item $\hookgen{i}\hookgen{j}\hookgen{i} = \hookgen{i}$ for all $|j-i| = 1$;
				\item $\hookgen{j}\hookgen{i} = \hookgen{i}\hookgen{j}$ for all $j-i \geq 2$;
				\item $\symgen{1}\hookgen{1} = \hookgen{1} = \hookgen{1}\symgen{1}$;
				\item $\hookgen{i}\hookgen{i+1}\symgen{i} = \hookgen{i}\symgen{i+1}$;
				\item $\symgen{i}\hookgen{i+1}\hookgen{i} = \symgen{i+1}\hookgen{i}$; and
				\item $\hookgen{j}\symgen{i} = \symgen{i}\hookgen{j}$ either for all $j-i \geq 2$ or for all $i-j \geq 2$.
			\end{enumerate}
			
			\begin{proof}
				It is sufficient for us to show that Relations (vii), (viii), (ix) and (x) from Proposition \ref{prop:brauermonpresentation1} may be deduced from Relations (i)-(x) above.
				
				Let $i \in \set{1, \ldots, k-2}$. First:
				\begin{enumerate}
					\item applying Relations (ix) and (v), $(\symgen{i+1}\hookgen{i})\hookgen{i+1} = \symgen{i}(\hookgen{i+1}\hookgen{i}\hookgen{i+1}) = \symgen{i}\hookgen{i+1}$; 
					\item applying Relations (viii) and (v), $\hookgen{i+1}(\hookgen{i}\symgen{i+1}) = (\hookgen{i+1}\hookgen{i}\hookgen{i+1})\symgen{i} = \hookgen{i+1}\symgen{i}$; and
					\item applying Relations (ix), (i), (viii), (v) and (i) again, $\symgen{i}(\symgen{i+1}\hookgen{i})\symgen{i+1}\symgen{i} = (\symgen{i}\symgen{i})\hookgen{i+1}(\hookgen{i}\symgen{i+1})\symgen{i} = \id{k}(\hookgen{i+1}\hookgen{i}\hookgen{i+1})(\symgen{i}\symgen{i}) = \hookgen{i+1}$.
				\end{enumerate} 
				
				To establish that Relation (vii) from Proposition \ref{prop:brauermonpresentation1} may be deduced we proceed by induction. Suppose that $\symgen{i}\hookgen{i} = \hookgen{i} = \hookgen{i}\symgen{i}$, noting that for induction the base case of $i=1$ is still given by Relation (vii), then we have $\hookgen{i+1}\symgen{i+1} = \symgen{i}\symgen{i+1}\hookgen{i}(\symgen{i+1}\symgen{i}\symgen{i+1}) = \symgen{i}\symgen{i+1}(\hookgen{i}\symgen{i})\symgen{i+1}\symgen{i} = \symgen{i}\symgen{i+1}\hookgen{i}\symgen{i+1}\symgen{i} = \hookgen{i+1}$.
				
				It remains for us to show that Relation (x) from Proposition \ref{prop:brauermonpresentation1} may be deduced. First note that $\symgen{i+1}\hookgen{i}\symgen{i+1} = \symgen{i+1}\hookgen{i}\hookgen{i+1}\symgen{i} = \symgen{i}\hookgen{i+1}\hookgen{i}\hookgen{i+1}\symgen{i} = \symgen{i}\hookgen{i+1}\symgen{i}$. Now for each $j \in \set{i+1, \ldots, k-1}$,
				\begin{align*}
				 	\symgen{j-1}\symgen{j}\hookgen{j-1}\symgen{j}\symgen{j-1} &= \symgen{j-1}\symgen{j-1}\hookgen{j}\symgen{j-1}\symgen{j-1} \\
				 	&= \hookgen{j}; \text{ and} \\ 
				 	\symgen{j-1}\ldots\symgen{i}\symgen{j}\ldots\symgen{i+1}\hookgen{i}\symgen{i+1}&\ldots\symgen{j}\symgen{i}\ldots\symgen{j-1} \\
				 	&= \symgen{j-1}\ldots\symgen{i}\symgen{j}\ldots\symgen{i+2}\symgen{i}\hookgen{i+1}\symgen{i}\symgen{i+2}\ldots\symgen{j}\symgen{i}\ldots\symgen{j-1} \\
				 	&= \symgen{j-1}\ldots\symgen{i+1}\symgen{j}\ldots\symgen{i+2}\hookgen{i+1}\symgen{i+2}\ldots\symgen{j}\symgen{i+1}\ldots\symgen{j-1} \\
				 	&= \hookgen{j}.
				\end{align*}
				
				Finally for each $j \in \set{i+2, \ldots, k-1}$,
				\begin{align*}
					\symgen{i}\hookgen{j} &= \symgen{i}\symgen{j-1}\ldots\symgen{i}\symgen{j}\ldots\symgen{i+1}\hookgen{i}\symgen{i+1}\ldots\symgen{j}\symgen{i}\ldots\symgen{j-1} \\
					&= \symgen{j-1}\ldots\symgen{i}\symgen{i+1}\symgen{j}\ldots\symgen{i+1}\hookgen{i}\symgen{i+1}\ldots\symgen{j}\symgen{i}\ldots\symgen{j-1} \\
					&=
					\symgen{j-1}\ldots\symgen{i}\symgen{j}\ldots\symgen{i+1}\symgen{i+2}\hookgen{i}\symgen{i+1}\ldots\symgen{j}\symgen{i}\ldots\symgen{j-1} \\
					&=
					\symgen{j-1}\ldots\symgen{i}\symgen{j}\ldots\symgen{i+1}\hookgen{i}\symgen{i+2}\symgen{i+1}\ldots\symgen{j}\symgen{i}\ldots\symgen{j-1} \\
					&=
					\symgen{j-1}\ldots\symgen{i}\symgen{j}\ldots\symgen{i+1}\hookgen{i}\symgen{i+1}\ldots\symgen{j}\symgen{i+1}\symgen{i}\ldots\symgen{j-1} \\
					&=
					\symgen{j-1}\ldots\symgen{i}\symgen{j}\ldots\symgen{i+1}\hookgen{i}\symgen{i+1}\ldots\symgen{j}\symgen{i}\ldots\symgen{j-1}\symgen{i} \\
					&= \hookgen{j}\symgen{i}.
				\end{align*}
			\end{proof}
		\end{proposition}
		
		\newpage
		
		Next we establish that the relations used by Kosuda \cite{art:Kosuda:CharacterizationModularPartyAlgebra} on diapsis and $(2,2)$-transapsis generators when giving a presentation of the mod-$2$ monoid $\modmon{2}{k}$, outlined later in Proposition \ref{prop:modmonpresentation}, also form a presentation of the Brauer monoid $\brauermon{k}$.

		\begin{proposition} \label{prop:brauermonpresentation3}
			For each $k \in \nonnegints$, the Brauer monoid $\brauermon{k}$ is characterised by the generators $\set{\symgen{i}, \hookgen{i}: i = 1, \ldots, k-1}$ along with the relations:
			\begin{enumerate}
				\item $\symgen{i}^2 = \id{k}$;
				\item $\symgen{i+1}\symgen{i}\symgen{i+1} = \symgen{i}\symgen{i+1}\symgen{i}$;
				\item $\symgen{j}\symgen{i} = \symgen{i}\symgen{j}$ for all $j-i \geq 2$;
				\item $\hookgen{i}^2 = \hookgen{i}$;
				\item $\hookgen{j}\hookgen{i} = \hookgen{i}\hookgen{j}$ for all $j-i \geq 2$;
				\item $\symgen{1}\hookgen{1} = \hookgen{1} = \hookgen{1}\symgen{1}$;
				\item $\hookgen{i+1} = \symgen{i}\symgen{i+1}\hookgen{i}\symgen{i+1}\symgen{i}$;
				\item $\hookgen{1}\symgen{2}\hookgen{1} = \hookgen{1}$; and
				\item $\hookgen{j}\symgen{i} = \symgen{i}\hookgen{j}$ either for all $j-i \geq 2$ or for all $i-j \geq 2$.
			\end{enumerate}
			
			\begin{proof}
				It is sufficient for us to show that Relations (v), (viii) and (ix) in Proposition \ref{prop:brauermonpresentation2} may be deduced from Relations (i)-(ix) above.
				
				First for each $i \in \set{1, ..., k-1}$, 
				\begin{align*}
					\symgen{i+1}(\hookgen{i+1}) &= (\symgen{i+1}\symgen{i}\symgen{i+1})\hookgen{i}\symgen{i+1}\symgen{i} \\
					&= \symgen{i}\symgen{i+1}(\symgen{i}\hookgen{i})\symgen{i+1}\symgen{i} \\
					&= (\symgen{i}\symgen{i+1}\hookgen{i}\symgen{i+1}\symgen{i}) \\
					&= \hookgen{i+1} \\
					&= (\symgen{i}\symgen{i+1}\hookgen{i}\symgen{i+1}\symgen{i}) \\
					&= \symgen{i}\symgen{i+1}(\hookgen{i}\symgen{i})\symgen{i+1}\symgen{i} \\
					&= \symgen{i}\symgen{i+1}\hookgen{i}(\symgen{i+1}\symgen{i}\symgen{i+1}) = (\hookgen{i+1})\symgen{i+1}.
				\end{align*}
				
				Now for each $i \in \set{2, ..., k-2}$ (note that our proof uses the first choice for Relation (ix) above), 
				\begin{align*}
					(\hookgen{i})\symgen{i+1}(\hookgen{i}) &= \symgen{i-1}\symgen{i}\hookgen{i-1}\symgen{i}(\symgen{i-1}\symgen{i+1})\symgen{i-1}\symgen{i}\hookgen{i-1}\symgen{i}\symgen{i-1} \\
					&= \symgen{i-1}\symgen{i}\hookgen{i-1}\symgen{i}\symgen{i+1}(\symgen{i-1}\symgen{i-1})\symgen{i}\hookgen{i-1}\symgen{i}\symgen{i-1} \\
					&= \symgen{i-1}\symgen{i}\hookgen{i-1}(\symgen{i}\symgen{i+1}\symgen{i})\hookgen{i-1}\symgen{i}\symgen{i-1} \\
					&= \symgen{i-1}\symgen{i}(\hookgen{i-1}\symgen{i+1})\symgen{i}(\symgen{i+1}\hookgen{i-1})\symgen{i}\symgen{i-1} \\
					&= \symgen{i-1}\symgen{i}\symgen{i+1}(\hookgen{i-1}\symgen{i}\hookgen{i-1})\symgen{i+1}\symgen{i}\symgen{i-1} \\
					&= \symgen{i-1}\symgen{i}(\symgen{i+1}\hookgen{i-1})\symgen{i+1}\symgen{i}\symgen{i-1} \\
					&= \symgen{i-1}\symgen{i}\hookgen{i-1}(\symgen{i+1}\symgen{i+1})\symgen{i}\symgen{i-1} \\
					&= (\symgen{i-1}\symgen{i}\hookgen{i-1}\symgen{i}\symgen{i-1}) = \hookgen{i}. \\
				\end{align*}
				
				Finally for each $i \in \set{1, ..., k-2}$:
				\begin{align*}
					\hookgen{i}(\hookgen{i+1})\hookgen{i} &= (\hookgen{i}\symgen{i})\symgen{i+1}\hookgen{i}\symgen{i+1}(\symgen{i}\hookgen{i}) \\
					&= (\hookgen{i}\symgen{i+1}\hookgen{i})\symgen{i+1}\hookgen{i} \\
					&= (\hookgen{i}\symgen{i+1}\hookgen{i}) = \hookgen{i}; \\
					(\hookgen{i+1})\hookgen{i}(\hookgen{i+1}) &= \symgen{i}\symgen{i+1}\hookgen{i}\symgen{i+1}(\symgen{i}\hookgen{i})\symgen{i}\symgen{i+1}\hookgen{i}\symgen{i+1}\symgen{i} \\
					&= \symgen{i}\symgen{i+1}\hookgen{i}\symgen{i+1}(\hookgen{i}\symgen{i})\symgen{i+1}\hookgen{i}\symgen{i+1}\symgen{i} \\
					&= \symgen{i}\symgen{i+1}(\hookgen{i}\symgen{i+1}\hookgen{i})\symgen{i+1}\hookgen{i}\symgen{i+1}\symgen{i} \\
					&= \symgen{i}\symgen{i+1}(\hookgen{i}\symgen{i+1}\hookgen{i})\symgen{i+1}\symgen{i} \\
					&= (\symgen{i}\symgen{i+1}\hookgen{i}\symgen{i+1}\symgen{i}) = \hookgen{i+1};\\
					\hookgen{i}(\hookgen{i+1})\symgen{i} &= (\hookgen{i}\symgen{i})\symgen{i+1}\hookgen{i}\symgen{i+1}(\symgen{i}\symgen{i}) \\
					&= (\hookgen{i}\symgen{i+1}\hookgen{i})\symgen{i+1} = \hookgen{i}\symgen{i+1}; \text{ and} \\
					\symgen{i}(\hookgen{i+1})\hookgen{i} &= (\symgen{i}\symgen{i})\symgen{i+1}\hookgen{i}\symgen{i+1}(\symgen{i}\hookgen{i}) \\
					&= \symgen{i+1}(\hookgen{i}\symgen{i+1}\hookgen{i}) =  \symgen{i+1}\hookgen{i}.
				\end{align*}
			\end{proof}
		\end{proposition}
		
		It was pointed out to the author by James East that by rewriting Relation (vii) from Proposition \ref{prop:brauermonpresentation3} as $\hookgen{i+1}\symgen{i}\symgen{i+1} = \symgen{i}\symgen{i+1}\hookgen{i}$ we get one of the relations for the singular braid monoid \cite{art:Baez:LinkInvariantsFiniteType, art:Birman:NewPointsOfViewInKnowTheory}. The author notes that in fact all of the relations for the singular braid monoid \cite{art:Baez:LinkInvariantsFiniteType, art:Birman:NewPointsOfViewInKnowTheory} may be deduced from the relations in Proposition \ref{prop:brauermonpresentation3}.
		
	\subsection{The symmetric inverse semigroup $\syminvmon{k}$}
		\begin{definition} \label{def:syminvmon}
			For each $k \in \nonnegints$, \textit{the symmetric inverse semigroup}, which is often denoted as $\syminvmon{k}$, is the set of bipartitions such that every block is either a transversal line or a monapsis.
		\end{definition}
		
		\begin{proposition} \label{prop:syminvmonpresentation} (see \cite{art:Easdown:InverseBraidMonoid, art:East:FactorizableInvMonsEtc, art:Popova:DefiningRltnsEtc})
			For each $k \in \nonnegints$, the symmetric inverse semigroup $\syminvmon{k}$ is characterised by the generators $\set{\symgen{i}, \apgen{1}{i}, \apgen{1}{k}: i = 1, \ldots, k-1}$ along with the relations:
			\begin{enumerate}
				\item $\symgen{i}^2 = \symgen{i}$;
				\item $\symgen{i}\symgen{i+1}\symgen{i} = \symgen{i+1}\symgen{i}\symgen{i+1}$;
				\item $\symgen{j}\symgen{i} = \symgen{i}\symgen{j}$ for all $|j-i| \geq 2$;
				\item ${(\apgen{1}{i})}^2 = \apgen{1}{i}$;
				\item $\apgen{1}{j}\apgen{1}{i} = \apgen{1}{i}\apgen{1}{j}$ for all $|j-i| \geq 1$;
				\item $\symgen{i}\apgen{1}{i} = \apgen{1}{i+1}\symgen{i}$;
				\item $\apgen{1}{i+1}\symgen{i} = \symgen{i}\apgen{1}{i}$; and
				\item $\symgen{j}\apgen{1}{i} = \apgen{1}{i}\symgen{j}$ for all $j \neq i, i+1$. \qed
			\end{enumerate}
		\end{proposition}
	
	\subsection{The planar symmetric inverse semigroup $\psyminvmon{k}$}
	
		\begin{definition} \label{def:psyminvmon}
			For each $k \in \nonnegints$, \textit{the planar symmetric inverse semigroup}, which we shall denote as $\psyminvmon{k}$, is the set of planar bipartitions such that every block is either a transversal line or a monapsis (see Figure \ref{fig:psyminvmon3} for examples), that is the meet of the planar partition monoid $\ppttnmon{k}$ and the symmetric inverse semigroup $\syminvmon{k}$.
		\end{definition}
	
		\begin{figure}[!ht]
			\caption[ ]{Given $k = 3$, the twenty elements of $\psyminvmon{3}$ may be depicted as:}
			\label{fig:psyminvmon3}
			\vspace{5pt}
			\centering
			\input{chap_background/tikz/fig-psyminvmon3.tex}
		\end{figure}
		
		\begin{proposition} \label{prop:psyminvmonmonpresentation} (see \cite{art:East:PresentationOfSingularPartSymInvMon})
			For each $k \in \nonnegints$, the planar symmetric inverse semigroup $\psyminvmon{k}$ is characterised by the generators $\set{\ftransgen{i}, \btransgen{i}: i = 1, \ldots, k}$ along with the relations:
			\begin{enumerate}
				\item $\ftransgen{j}\ftransgen{i} = \ftransgen{i}\ftransgen{j+1}$ for all $1 \leq i \leq j \leq k-1$; 
				\item $\ftransgen{k}\ftransgen{i} = \ftransgen{i}$ for all $1 \leq i \leq k$;
				\item $\btransgen{i}\btransgen{j} = \btransgen{j+1}\btransgen{i}$ for all $1 \leq i \leq j \leq k-1$; 
				\item $\btransgen{i}\btransgen{k} = \btransgen{i}$ for all $1 \leq i \leq k$; and
				\item $\btransgen{i}\ftransgen{j} = \begin{cases} \ftransgen{k}\ftransgen{j-1}\btransgen{i} & \text{ for all } 1 \leq i < j \leq k, \\ \ftransgen{k} = \btransgen{k} & \text{ for all } 1 \leq i=j \leq k, \text{ and} \\ \ftransgen{k}\ftransgen{j}\btransgen{i-1} & \text{ for all } 1 \leq j < i \leq k. \end{cases}$ \qed
			\end{enumerate}
		\end{proposition}
		
	\subsection{The monoid of planar uniform block bijections $\puniblockbijmon{k}$}
		\begin{definition} \label{def:puniblockbijmon}
			For each $k \in \nonnegints$, \textit{the planar monoid of uniform block bijections}, which is often denoted as $\puniblockbijmon{k}$, is the set of planar bipartitions such that every block is uniform.
		\end{definition}
		
		\begin{proposition} \label{prop:puniblockbijmonpresentation}
			For each $k \in \nonnegints$, the monoid of planar uniform block bijections $\puniblockbijmon{k}$ is characterised by the generators $\set{\transapgen{i}: i = 1, \ldots, k-1}$ along with the relations:
			\begin{enumerate}
				\item $\transapgen{i}^2 = \transapgen{i}$; and
				\item $\transapgen{j}\transapgen{i} = \transapgen{i}\transapgen{j}$ for all $|j-i| \geq 1$. \qed
			\end{enumerate}
		\end{proposition}

	\subsection{The monoid of uniform block bijections $\uniblockbijmon{k}$}
		\begin{definition} \label{def:uniblockbijmon}
			For each $k \in \nonnegints$, \textit{the monoid of uniform block bijections}, which is often denoted as $\uniblockbijmon{k}$, is the set of bipartitions such that every block is uniform.
		\end{definition}
		
		\begin{proposition} \label{prop:uniblockbijmonpresentation1} (see \cite{art:Fitzgerald:PresentationMonUniBlockBij, art:Kosuda:CharacterizationPartAlgebras})
			For each $k \in \nonnegints$, the monoid of uniform block bijections $\uniblockbijmon{k}$ is characterised by the generators $\set{\symgen{i}, \transapgen{1}: i = 1, \ldots, k-1}$ along with the relations:
			\begin{enumerate}
				\item $\symgen{i}^2 = \id{k}$;
				\item $\symgen{i+1}\symgen{i}\symgen{i+1} = \symgen{i}\symgen{i+1}\symgen{i}$;
				\item $\symgen{j}\symgen{i} = \symgen{i}\symgen{j}$ for all $|j-i| \geq 2$;
				\item $\transapgen{1}^2 = \transapgen{1}$;
				\item $\symgen{1}\transapgen{1} = \transapgen{1} = \transapgen{1}\symgen{1}$;
				\item $\symgen{2}\transapgen{1}\symgen{2}\transapgen{1} = \transapgen{1}\symgen{2}\transapgen{1}\symgen{2}$;
				\item $\symgen{2}\symgen{1}\symgen{3}\symgen{2}\transapgen{1}\symgen{2}\symgen{1}\symgen{3}\symgen{2}\transapgen{1} = \transapgen{1}\symgen{2}\symgen{1}\symgen{3}\symgen{2}\transapgen{1}\symgen{2}\symgen{1}\symgen{3}\symgen{2}$; and
				\item $\symgen{i}\transapgen{1} = \transapgen{1}\symgen{i}$ for all $i = 3, \ldots, k-1$. \qed
			\end{enumerate}
		\end{proposition}
		
		\begin{proposition} \label{prop:uniblockbijmonpresentation2} (see \cite{art:East:FactorizableBraidMon, art:Kosuda:CharacterizationPartAlgebras, art:Kudryavtseva:PresentationsBrauerTypeMonoids})
			For each $k \in \nonnegints$, the monoid of uniform block bijections $\uniblockbijmon{k}$ is characterised by the generators $\set{\symgen{i}, \transapgen{i}: i = 1, \ldots, k-1}$ along with the relations:
			\begin{enumerate}
				\item $\symgen{i}^2 = \id{k}$;
				\item $\symgen{i+1}\symgen{i}\symgen{i+1} = \symgen{i}\symgen{i+1}\symgen{i}$;
				\item $\symgen{j}\symgen{i} = \symgen{i}\symgen{j}$ for all $|j-i| \geq 2$;
				\item $\transapgen{i}^2 = \transapgen{i}$;
				\item $\transapgen{j}\transapgen{i} = \transapgen{i}\transapgen{j}$;
				\item $\symgen{i}\transapgen{i} = \transapgen{i} = \transapgen{i}\symgen{i}$;
				\item either $\symgen{i+1}\transapgen{i}\symgen{i+1} = \symgen{i}\transapgen{i+1}\symgen{i}$ \cite{art:Kosuda:CharacterizationPartAlgebras} or $\symgen{i+1}\symgen{i}\transapgen{i+1} = \transapgen{i}\symgen{i+1}\symgen{i}$ \cite{art:East:FactorizableBraidMon}; and
				\item $\symgen{j}\transapgen{i} = \transapgen{i}\symgen{j}$ for all $|j-i| \geq 2$. \qed
			\end{enumerate}
		\end{proposition}
		
	\subsection{The modular partition monoid $\modmon{m}{k}$}
		\begin{definition} \label{def:modmon}
			For each $m \in \posints$ and $k \in \intsge{m}$, \textit{the modular partition monoid} or more simply \textit{the mod-$m$ monoid} when ambiguity is of no concern in the context the term is used, which we denote as $\modmon{m}{k}$, is the set of all bipartitions $\alpha \in \pttnmon{k}$ such that for each block $b \in \alpha$, $\congmod{\noupverts{b}}{\nolowverts{b}}{m}$ (see Figure \ref{fig:modmonelmnts} for some examples).
		\end{definition}
		
		Note that the mod-$m$ monoid $\modmon{m}{k}$ was referred to by Kosuda as \textit{the party algebra of type B} in \cite{art:Kosuda:PartyAlgTypeBandIrreps, art:Kosuda:StructurePartyAlgebraTypeB} and as \textit{the $m$-modular party algebra} in \cite{art:Kosuda:StdExpForPartyAlgebra, art:Kosuda:CharacterizationModularPartyAlgebra}, and referred to by Ahmed, Martin and Mazorchuk as \textit{the $d$-tonal partition monoid} in \cite{art:Ahmed:OnTheNoOfPrincipalIdealsInDTonalPartMons}.
		
		\begin{figure}[!ht]
			\caption[ ]{Examples of bipartitions in the mod-$m$ monoid $\modmon{m}{k}$.}
			\label{fig:modmonelmnts}
			\vspace{5pt}
			\centering
			\input{chap_characterisations/tikz/fig-modmonelmnts.tex}
		\end{figure} 
		
		\begin{proposition} \label{prop:modmonisamonoid}
			$\modmon{m}{k}$ is a monoid.
			
			\begin{proof}
				Let $\alpha, \beta \in \modmon{m}{k}$. Suppose that when forming the product $\alpha\beta$, $x \in \posints$ blocks $a_1, \ldots, a_x \in \alpha$ join with $y \in \posints$ blocks $b_1, \ldots, b_y \in \beta$, forming the block $\bigcup_{i=1}^x\upverts{a_i}\cup\bigcup_{j=1}^y\lowverts{b_j} \in \alpha\beta$. It follows from $\congmod{\Sigma^x_{i=1}u(a_i)}{\Sigma^x_{i=1}l(a_i)}{m}$, $\Sigma^x_{i=1}l(a_i) = \Sigma^y_{j=1}u(b_j)$ and $\congmod{\Sigma^y_{j=1}u(b_j)}{\Sigma^y_{j=1}l(b_j)}{m}$ that $\congmod{\Sigma^x_{i=1}u(a_i)}{\Sigma^y_{j=1}l(b_j)}{m}$.
			\end{proof}
		\end{proposition}
		
		Note it follows by definition that 
		\begin{enumerate}
			\item the mod-$1$ monoid $\modmon{1}{k}$ is equal to the partition monoid $\pttnmon{k}$; and
			\item if $k < m$ then the mod-$m$ monoid $\modmon{m}{k}$ is equal to the monoid of uniform block bijections $\uniblockbijmon{k}$.
		\end{enumerate}
		
		\begin{proposition} \label{prop:modmonpresentation} (see \cite{art:Kosuda:CharacterizationModularPartyAlgebra})
			For each $k \in \nonnegints$, the mod-$m$ monoid $\modmon{m}{k}$ is characterised by the generators $\set{\symgen{i}, \transapgen{i}, \apgen{m}{j}: i = 1, \ldots, k-1, j = 1, \ldots, k-m+1}$ along with the relations:
			\begin{enumerate}
				\item $\symgen{i}^2 = \id{k}$;
				\item $\symgen{i}\symgen{i+1}\symgen{i} = \symgen{i+1}\symgen{i}\symgen{i+1}$;
				\item $\symgen{j}\symgen{i} = \symgen{i}\symgen{j}$ for all $j-i \geq 2$;
				\item $\transapgen{i}^2 = \transapgen{i}$;
				\item $\transapgen{j}\transapgen{i} = \transapgen{i}\transapgen{j}$ for all $|j-i| \geq 1$;
				\item $\apgen{m}{i}\apgen{m}{i} = \apgen{m}{i}$;
				\item $\apgen{m}{j}\apgen{m}{i} = \apgen{m}{i}\apgen{m}{j}$ for all $|j-i| \geq m$;
				\item $\transapgen{i}\symgen{i} = \transapgen{i} = \symgen{i}\transapgen{i}$;
				\item $\symgen{i}\symgen{i+1}\transapgen{i}\symgen{i+1}\symgen{i} = \transapgen{i+1}$;
				\item $\transapgen{j}\symgen{i} = \symgen{i}\transapgen{j}$ for all $|j-i| \geq 2$;
				\item $\apgen{m}{i}\symgen{j} = \apgen{m}{i} = \symgen{j}\apgen{m}{i}$ for all $j \in \set{i, ..., i+m-2}$;
				\item $\apgen{m}{i}\symgen{i+m-1}\apgen{m}{i} = \apgen{m}{i} = \apgen{m}{i}\symgen{i-1}\apgen{m}{i}$;
				\item $\symgen{i}\ldots\symgen{i+m-1}\apgen{m}{i}\symgen{i+m-1}\ldots\symgen{i} = \apgen{m}{i+1}$;
				\item $\apgen{m}{i}\symgen{j} = \symgen{j}\apgen{m}{i}$ for all $j \in \set{1, ..., i-2, i+m, ..., k-1}$;
				\item $\apgen{m}{i}\transapgen{i} = \apgen{m}{i} = \transapgen{i}\apgen{m}{i}$;
				\item $\transapgen{i}\apgen{m}{i+1}\transapgen{i} = \transapgen{i+m-1}\apgen{m}{i}\transapgen{i+m-1} = \transapgen{i}\ldots\transapgen{i+m-1}$; and
				\item $\apgen{m}{i}\transapgen{j} = \transapgen{j}\apgen{m}{i}$ for all $j \in \set{1, ..., i-2, i+m, ..., k-1}$. \qed
			\end{enumerate}
		\end{proposition}
		
		Recall from Subsection \ref{subsec:brauermon} that Kudryavtseva and Mazorchuk used a slightly different set of relations to give a presentation of the Brauer monoid $\brauermon{k}$ than the relations used on diapsis and transposition generators in Proposition \ref{prop:modmonpresentation}. Switching these relations again gives us a second equivalent set of relations that present the mod-$2$ monoid $\modmon{2}{k}$.
		
		\begin{proposition} \label{prop:modmonpresentation2}
			For each $k \in \nonnegints$, the mod-$m$ monoid $\modmon{m}{k}$ is characterised by the generators $\set{\symgen{i}, \transapgen{i}, \apgen{m}{j}: i = 1, \ldots, k-1, j = 1, \ldots, k-m+1}$ along with the relations:
			\begin{enumerate}
				\item $\symgen{i}^2 = \id{k}$;
				\item $\symgen{i}\symgen{i+1}\symgen{i} = \symgen{i+1}\symgen{i}\symgen{i+1}$;
				\item $\symgen{j}\symgen{i} = \symgen{i}\symgen{j}$ for all $j-i \geq 2$;
				\item $\transapgen{i}^2 = \transapgen{i}$;
				\item $\transapgen{j}\transapgen{i} = \transapgen{i}\transapgen{j}$ for all $|j-i| \geq 1$;
				\item $\hookgen{i}\hookgen{i} = \hookgen{i}$;
				\item $\hookgen{i}\hookgen{j}\hookgen{i} = \hookgen{i}$ for all $|j-i| = 1$;
				\item $\hookgen{j}\hookgen{i} = \hookgen{i}\hookgen{j}$ for all $|j-i| \geq m$;
				\item $\transapgen{i}\symgen{i} = \transapgen{i} = \symgen{i}\transapgen{i}$;
				\item $\symgen{i}\symgen{i+1}\transapgen{i}\symgen{i+1}\symgen{i} = \transapgen{i+1}$;
				\item $\transapgen{j}\symgen{i} = \symgen{i}\transapgen{j}$ for all $|j-i| \geq 2$;
				\item $\hookgen{i}\symgen{i} = \hookgen{i} = \symgen{i}\hookgen{i}$;
				\item $\hookgen{i}\hookgen{i+1}\symgen{i} = \hookgen{i}\symgen{i+1}$;
				\item $\symgen{i}\hookgen{i+1}\hookgen{i} = \symgen{i+1}\hookgen{i}$;
				\item $\hookgen{i}\symgen{j} = \symgen{j}\hookgen{i}$ for all $|j-i| \geq 2$;
				\item $\hookgen{i}\transapgen{i} = \hookgen{i} = \transapgen{i}\hookgen{i}$;
				\item $\transapgen{i}\hookgen{i+1}\transapgen{i} = \transapgen{i+1}\hookgen{i}\transapgen{i+1} = \transapgen{i}\transapgen{i+1}$; and
				\item $\hookgen{i}\transapgen{j} = \transapgen{j}\hookgen{i}$ for all $|j-i| \geq 2$.
			\end{enumerate}
			
			\begin{proof}
				It is sufficient for us to show that Relations (vii), (xiii) and (xiv) from above follow from the relations in Proposition \ref{prop:modmonpresentation}. Note that all but Relations (vii), (xiii) and (xiv) from above already appear in Proposition \ref{prop:modmonpresentation}. For each $i, j \in \set{1, \ldots, k-1}$ such that $j-i = 1$:
				\begin{enumerate}
					\item $\hookgen{i}(\hookgen{j})\hookgen{i} = (\hookgen{i}\symgen{i})\symgen{j}\hookgen{i}\symgen{j}(\symgen{i}\hookgen{i}) = (\hookgen{i}\symgen{j}\hookgen{i})\symgen{j}\hookgen{i} = (\hookgen{i}\symgen{j}\hookgen{i}) = \hookgen{i}$;
					\item $\hookgen{i}(\hookgen{j})\symgen{i} = (\hookgen{i}\symgen{i})\symgen{j}\hookgen{i}\symgen{j}(\symgen{i}\symgen{i}) = (\hookgen{i}\symgen{j}\hookgen{i})\symgen{j} = \hookgen{i}\symgen{j}$;
					\item $\symgen{i}(\hookgen{j})\hookgen{i} = (\symgen{i}\symgen{i})\symgen{j}\hookgen{i}\symgen{j}(\symgen{i}\hookgen{i}) = \symgen{j}(\hookgen{i}\symgen{j}\hookgen{i}) = \symgen{j}\hookgen{i}$; and
					\item $(\hookgen{j})\hookgen{i}(\hookgen{j}) = \symgen{i}\symgen{j}\hookgen{i}\symgen{j}(\symgen{i}\hookgen{i}\symgen{i})\symgen{j}\hookgen{i}\symgen{j}\symgen{i} = \symgen{i}\symgen{j}(\hookgen{i}\symgen{j}\hookgen{i})\symgen{j}\hookgen{i}\symgen{j}\symgen{i} = \symgen{i}\symgen{j}(\hookgen{i}\symgen{j}\hookgen{i})\symgen{j}\symgen{i} = (\symgen{i}\symgen{j}\hookgen{i}\symgen{j}\symgen{i}) = \hookgen{j}$.
				\end{enumerate}
			\end{proof}
		\end{proposition}
		
	\subsection{The planar modular partition monoid $\pmodmon{m}{k}$}
\begin{figure}[!ht]
			\caption[ ]{Given $m=3$ and $k=5$,}
			\label{fig:pmodmonelmnts}
			\vspace{5pt}
			\centering
			\input{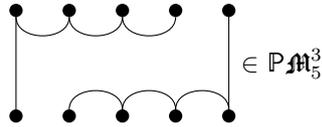}
		\end{figure}

		\begin{definition} \label{def:pmodmon}
			For each $m \in \posints$ and $k \in \intsge{m}$, \textit{the planar modular partition monoid} or more simply \textit{the planar mod-$m$ monoid} when ambiguity is of no concern in the context the term is used, which we denote as $\pmodmon{m}{k}$, is the planar analogue of the mod-$m$ monoid $\pmodmon{m}{k}$, that is the set of all planar bipartitions $\alpha \in \ppttnmon{k}$ such that for each block $b \in \alpha$, $\congmod{\noupverts{b}}{\nolowverts{b}}{m}$ (see Figure \ref{fig:pmodmonelmnts} for an example).
		\end{definition}
		
		Note it follows by definition that 
		\begin{enumerate}
			\item the planar mod-$1$ monoid $\pmodmon{1}{k}$ is equal to the planar partition monoid $\ppttnmon{k}$; and
			\item if $k < m$ then the planar mod-$m$ monoid $\pmodmon{m}{k}$ is equal to the monoid of planar uniform block bijections $\puniblockbijmon{k}$.
		\end{enumerate}
		
		Surprisingly, to the best of the author's knowledge, the planar mod-$m$ monoid $\pmodmon{m}{k}$ has only been examined elsewhere in \cite{art:Ahmed:OnTheNoOfPrincipalIdealsInDTonalPartMons}, the work for which was undertaken simultaneously and independently to the author.

    \clearpage{\pagestyle{empty}\cleardoublepage}
\chapter{Characterisations}	\label{chap:characterisations}	
	\section{The $m$-apsis generated diagram monoid $\apsismon{m}{k}$} 
	\subsection{$m$-apsis generators}
		For each $m \in \posints$ and $k \in \intsge{m}$:
		\begin{enumerate}
			\item recall from Definition \ref{def:blocktypes} that an $m$-apsis is a non-transversal block containing $m$ consecutive vertices, a monapsis is a $1$-apsis, a diapsis is a $2$-apsis, and a triapsis is a $3$-apsis; and
			\item recall from Definition \ref{def:mapsisgens} that for each $m \in \posints$ and $k \in \intsge{m}$, the $m$-apsis generators consist of the $k-m+1$ bipartitions $\apgen{m}{1}, \ldots, \apgen{m}{k-m+1}$ such that for each $i \in \{1, \ldots, k-m+1\}$, $\apgen{m}{i}$ contains:
				\begin{enumerate}
					\item the two $m$-apses $\{i, \ldots, i+m-1\}$ and $\{i', \ldots, (i+m-1)'\}$; and
					\item for each $j \in \{1, \ldots, i-1, i+m, \ldots, k\}$, the vertical line $\{j, j'\}$.
				\end{enumerate}
		\end{enumerate}
		
\begin{figure}[!ht]
			\caption[ ]{Given $m=1$ and $k = 4$, there are four monapsis generators depicted as follows:}
			\label{fig:A14generators}
			\vspace{5pt}
			\centering
			\input{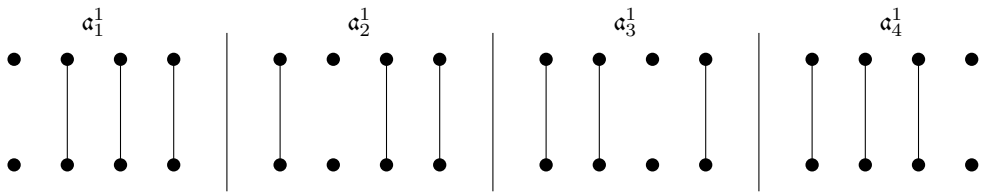}
		\end{figure}

\begin{figure}[!ht]
			\caption[ ]{Given $m=3$ and $k = 5$, there are three triapsis generators depicted as follows:}
			\label{fig:A35Generators}
			\vspace{5pt}
			\centering
			\input{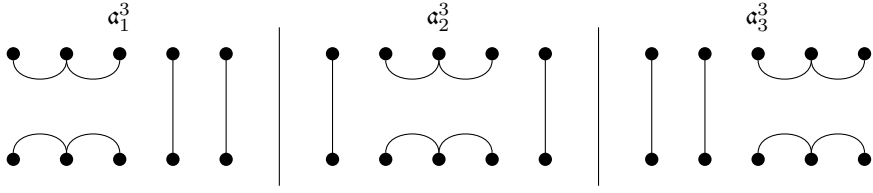}
		\end{figure}

		For example there are four monapsis generators when $k=4$ as depicted in Figure \ref{fig:A14generators}, and three triapsis generators when $k=5$ as depicted in Figure \ref{fig:A35Generators}.
		
		\begin{definition} \label{def:apsismon}
			Given $m \in \posints$ and $k \in \ints_{\geq m}$, we denote by $\apsismon{m}{k}$ the monoid generated by the $m$-apsis generators $\set{\apgen{m}{j}: j \in \set{1, \ldots, k-m+1}}$ along with the identity bipartition $\id{k}$, that is $\apsismon{m}{k}$ denotes the monoid $\langle\apgen{m}{j}, \id{k}: j \in \{1, \ldots, k-m+1\}\rangle$. We shall refer to the monoid $\apsismon{m}{k}$ as \textit{the $m$-apsis generated diagram monoid} or, for the sake of being succinct when ambiguity is of no concern in the context the term is used, simply as \textit{the $m$-apsis monoid}.
		\end{definition}
	
		Monapsis generators are idempotent and trivially commute with each other, hence the monapsis monoid $\uniapmon{k}$ is trivially isomorphic to the join-semilattice of subsets of $\set{1, \ldots, k}$ under $\uniapmon{i} \mapsto \set{i}$ for all $i \in \set{1, \ldots, k}$.
		
		The diapsis monoid $\diapmon{k}$ is by definition the Jones monoid $\jonesmon{k}$, that is $\apsismon{2}{k} = \jonesmon{k}$. We proceed to characterise the bipartitions contained in the $m$-apsis monoid $\apsismon{m}{k}$ for $m \in \intsge{3}$.   

	\subsection{Runs of $m$-apsis generators}
		\begin{definition}
    		For each $m \in \posints$, $k \in \intsge{m}$ and $i, j \in \set{1, \ldots, k-m+1}$:
	    	\begin{enumerate}
	    		\item the \textit{run of $m$-apsis generators from $i$ to $j$} is the product $\apgen{m}{i}\ldots\apgen{m}{j}$ of $m$-apsis generators where indices increase by one throughout the product when $i < j$, and decrease by one when $i > j$; and
		    	\item we denote by $\run{i}{j}$ the planar bipartition containing:
		    	\begin{enumerate}
			    	\item the upper $m$-apsis $\set{i, \ldots, i+m-1}$;
			    	\item the lower $m$-apsis $\set{j', \ldots, (j+m-1)'}$; and
			    	\item transversal lines connecting the remaining vertices in a planar fashion.
		    	\end{enumerate}
	    		
	    	\end{enumerate}
    	\end{definition}
    	
    	For example, given $m=3$ and $k=6$, $\run{4}{2} \in \ppttnmon{6}$ is depicted in Figure \ref{fig:run42}.
 
		\begin{figure}[!ht]
			\caption[ ]{Given $m=3$ and $k = 6$,}
			\label{fig:run42}
			\vspace{5pt}
			\centering
			\input{chap_characterisations/tikz/fig-run42.tex}
  		\end{figure}

		When reading the proof of Proposition \ref{prop:runcharacterisation}, where we establish that $\run{i}{j}$ is precisely the run $\apgen{m}{i}\ldots\apgen{m}{j}$ of $m$-apsis generators, the reader may find it useful to refer to Figure \ref{fig:run42=g4g3g2} which illustrates that, given $m=3$ and $k=6$, $\run{4}{2} = \apgen{3}{4}\apgen{3}{3}\apgen{3}{2} \in \triapmon{6}$, along with Figure \ref{fig:rim-1gm} which illustrates that given $m=3$ and $i, x \in \set{1, \ldots, k-2}$ such that $i < x$, $\run{i}{x} = \run{i}{x-1}\apgen{3}{x}$.

		\begin{proposition} \label{prop:runcharacterisation}
			For each $m \in \intsge{2}$, $k \in \intsge{m}$ and $i, j \in \{1, \ldots, k-m+1\}$, $\run{i}{j}$ is equal to the run $\apgen{m}{i}\ldots\apgen{m}{j} \in \apsismon{m}{k}$ of $m$-apsis generators. 
			
			\begin{proof}
				When $i=j$, $\run{i}{j} = \apgen{m}{i}$ by definition. When $i < j$, let $x \in \set{i+1, \ldots, j}$ and suppose $\run{i}{x-1} = \apgen{m}{i}\ldots\apgen{m}{x-1}$. When forming the product $\run{i}{x-1}\apgen{m}{x}$ (see Figure \ref{fig:rim-1gm} for a depiction): 
				\begin{enumerate}
					\item The upper $m$-apsis $\set{i, \ldots, i+m-1}$ in $\run{i}{x-1}$ along with the lower $m$-apsis $\set{x', \ldots, (x+m-1)'}$ in $\apgen{m}{x}$ are preserved;
					\item The lower $m$-apsis $\set{(x-1)', \ldots, (x+m-2)'}$ in $\run{i}{x-1}$ and the upper $m$-apsis $\set{x, \ldots, x+m-1}$ in $\apgen{m}{x}$ connect at the middle row vertices $x, \ldots, x+m-2$, also joining to both the block $\{x+m-1,$ $(x+m-1)'\} $ in $\run{i}{x-1}$ along with the block $\set{x-1, (x-1)'}$ in $\apgen{m}{x}$, all together forming the block $\set{x+m-1, (x-1)'}$; and
					\item Finally there is no option other than for the remaining transversal lines of $\run{i}{x-1}$ to connect to the remaining vertical lines of $\apgen{m}{x}$.
				\end{enumerate}  
				Hence $\run{i}{x} = \run{i}{x-1}\apgen{m}{x} = \apgen{m}{i}\ldots\apgen{m}{x}$, by induction it follows that $\run{i}{j}$ is the run of $m$-apsis generators $\apgen{m}{i}\ldots\apgen{m}{j}$. Finally when $i > j$, $\run{i}{j} = (\run{j}{i})^* = (\apgen{m}{j}\ldots\apgen{m}{i})^* = \apgen{m}{i}\ldots\apgen{m}{j}$.
			\end{proof}
		\end{proposition} 
		
		\begin{figure}[!ht]
        \caption[ ]{Given $m=3$ and $k = 6$, $\run{4}{2} = \apgen{3}{4}\apgen{3}{3}\apgen{3}{2} \in \triapmon{6}$,}
        \label{fig:run42=g4g3g2}
        \vspace{5pt}
        \centering
        \input{chap_characterisations/tikz/fig-run42=g4g3g2.tex}
    \end{figure}

		\begin{figure}[!ht]
      \caption[ ]{Given $m=3$ and $i, x \in \set{1, \ldots, k-2}$ such that $i < x$, $\run{i}{x} = \run{i}{x-1}\apgen{m}{x}$.}
      \label{fig:rim-1gm}
      \vspace{5pt}
      \centering
      \input{chap_characterisations/tikz/proof-rim-1gm.tex}
  \end{figure}

		Next we establish in Proposition \ref{prop:runsaretrans} that for each $i, j, x \in \{1, \ldots, k-m+1\}$, $\run{i}{x}\run{x}{j} = \run{i}{j}$. The reader may find it useful to refer to Figure \ref{fig:run35run51=run31} which illustrates, given $m=3$ and $k = 8$, that $\run{3}{5}\run{5}{1} = \run{3}{1}$.

		\begin{proposition} \label{prop:runsaretrans}
			For each $m \in \intsge{2}$, $k \in \intsge{m}$ and $i, j, x \in \set{1, \ldots, k-m+1}$, $\run{i}{x}\run{x}{j} = \run{i}{j}$.
			
			\begin{proof}
				In the product $\run{i}{x}\run{x}{j}$:
				\begin{enumerate}
					\item The upper $m$-apsis $\set{i, \ldots, i+m-1}$ in $\run{i}{x}$ is preserved;
					\item The lower $m$-apsis $\{j', \ldots, (j+m-1)'\}$ in $\run{x}{j}$ is preserved;
					\item The lower $m$-apsis of $\run{i}{x}$ and the upper $m$-apsis of $\run{x}{j}$ cap each other off; and
					\item Finally there is no option other than for the remaining transversal lines of $\run{i}{x}$ and $\run{x}{j}$ to connect to each other.
				\end{enumerate}
				Hence $\run{i}{x}\run{x}{j} = \run{i}{j}$.
			\end{proof}
		\end{proposition}
		
		\begin{figure}[!ht]
			\caption[ ]{Given $m=3$ and $k=8$, $\run{3}{5}\run{5}{1} = \run{3}{1}$.}
			\label{fig:run35run51=run31}
			\vspace{5pt}
			\centering
			\input{chap_characterisations/tikz/fig-run35run51=run31.tex}
		\end{figure}
 
 	\subsection{$m$-apmorphisms} \label{subsec:m-apmorphisms}
    	\begin{definition} \label{def:m-apmorphism}
  			Given $m \in \intsge{3}$ and $k \in \intsge{m}$ by an \textit{$m$-apmorphism} we shall mean a planar bipartition $\theta \in \ppttnmon{k}$ such that each block $b \in \theta$ is either an $m$-apsis or a transversal line. Furthermore we denote by $\apmorphs{m}{k}$ the set of all $m$-apmorphisms.
  		\end{definition}
    	
    	\begin{definition}
			For each $t \in \set{0, \ldots, \floor{\frac{k}{m}}}$ and $u_1, \ldots, u_t, l_1, \ldots, l_t \in \{1, \ldots$, $k-m+1\}$ such that $u_j+m \leq u_{j+1}$ and $l_j+m \leq l_{j+1}$ for all $j \in \set{1, \ldots, t-1}$, we shall denote by $\apmorph{u_1, .., u_t}{l_1, \ldots, l_t}$ the $m$-apmorphism containing:
			\begin{enumerate}
				\item $t$ upper $m$-apses $\set{u_i, \ldots, u_i+m-1}$ where $i \in \set{1, \ldots, t}$;
				\item $t$ lower $m$-apses $\set{l_i', \ldots, (l_i+m-1)'}$ where $i \in \set{1, \ldots, t}$; and
				\item transversal lines connecting the remaining vertices in a planar fashion.
			\end{enumerate}  
		\end{definition}
		For example $\apmorph{\emptyset}{\emptyset} = \id{k}$ and given $m=3$ and $k=9$, $\apmorph{3,6}{2,7} \in \apmorphs{3}{9}$ is depicted in Figure \ref{fig:sigma3627}.
		
		\begin{figure}[!ht]
			\caption[ ]{Given $m=3$ and $k=9$,}
			\label{fig:sigma3627}
			\vspace{5pt}
			\centering
			\input{chap_characterisations/tikz/fig-sigma3627.tex}
		\end{figure}

		We proceed by establishing that $\apmorph{u_1, .., u_t}{l_1, \ldots, l_t}$ is equal to the product of runs $\run{u_1}{1}\ldots\run{u_t}{m(t-1)+1}\run{m(t-1)+1}{l_t}\ldots\run{1}{l_1}$. The reader may find it useful to refer to Figure \ref{fig:sigma3627factorisation} which illustrates, given $m=3$ and $k = 9$, that $\apmorph{3,6}{2,7} = \run{3}{1}\apmorph{6}{7}\run{1}{2} = \run{3}{1}\run{6}{4}\run{4}{7}\run{1}{2} \in \triapmon{9}$.

		\begin{proposition} \label{prop:apmorphismsinapsismon}
			For each $t \in \set{0, \ldots, \floor{\frac{k}{m}}}$ and $u_1, \ldots, u_t, l_1, \ldots, l_t \in \{1,$ $\ldots, k-m+1\}$ such that $u_j+m \leq u_{j+1}$ and $l_j+m \leq l_{j+1}$ for all $j \in \set{1, \ldots, t-1}$, $\apmorph{u_1, .., u_t}{l_1, \ldots, l_t} = \run{u_1}{1}\ldots\run{u_t}{m(t-1)+1}\run{m(t-1)+1}{l_t}\ldots\run{1}{l_1} \in \apsismon{m}{k}$.
			
			\begin{proof}
				It follows by definition that $\apmorph{u_t}{l_t} = \run{u_t}{l_t}$ and from Proposition \ref{prop:runsaretrans} that $\run{u_t}{l_t} = \run{u_t}{m(t-1)+1}\run{m(t-1)+1}{l_t}$.
				
				Let $j \in \set{1, \ldots, t-1}$ and suppose $\apmorph{u_{j+1}, \ldots, u_t}{l_{j+1}, \ldots, l_t} = \run{u_{j+1}}{mj+1}\ldots\run{u_t}{m(t-1)+1}\run{m(t-1)+1}{l_t}\ldots\run{mj+1}{l_{j+1}}$. In order to show that $\apmorph{u_j, \ldots, u_t}{l_j, \ldots, l_t}$ is equal to $\run{u_j}{m(j-1)+1}\apmorph{u_{j+1}, \ldots, u_t}{l_{j+1}, \ldots, l_t}\run{m(j-1)+1}{l_j} = \run{u_j}{m(j-1)+1}\ldots\run{u_t}{m(t-1)+1}\run{m(t-1)+1}{l_t}\ldots\run{m(j-1)+1}{l_j}$, it is sufficient for us to show that the product $\run{u_j}{m(j-1)+1}\apmorph{u_{j+1}, \ldots, u_t}{l_{j+1}, \ldots, l_t}\run{m(j-1)+1}{l_j}$ joins and removes the lower $m$-apsis of $\run{u_j}{m(j-1)+1}$ and upper $m$-apsis of $\run{m(j-1)+1}{l_j}$, while also preserving: all $m$-apses in $\apmorph{u_{j+1}, \ldots, u_t}{l_{j+1}, \ldots, l_t}$ along with the upper $m$-apsis of $\run{u_j}{m(j-1)+1}$ and lower $m$-apsis of $\run{m(j-1)+1}{l_j}$.
				
				When forming the product $\run{u_j}{m(j-1)+1}\apmorph{u_{j+1}, \ldots, u_t}{l_{j+1}, \ldots, l_t}\run{m(j-1)+1}{l_j}$:
				\begin{enumerate}
					\item the upper $m$-apsis $\set{u_j, \ldots, u_j+m-1}$ in $\run{u_j}{m(j-1)+1}$ and the lower $m$-apsis $\set{l_j', \ldots, (l_j+m-1)'}$ in $\run{m(j-1)+1}{l_j}$ are preserved;
					
					\item For each $i \in \{1, \ldots, mj\}$, it follows from the left-most $m$-apsis of $\apmorph{u_{j+1}, \ldots, u_t}{l_{j+1}, \ldots, l_t}$ starting at $\min\{u_{j+1},$ $l_{j+1}\}$ and $u_{j+1}, l_{j+1} > mj$ that  $\apmorph{u_{j+1}, \ldots, u_t}{l_{j+1}, \ldots, l_t}$ contains the vertical line $\set{i, i'}$. Hence the lower $m$-apsis $\set{(m(j-1)+1)', \ldots, (mj)'} \in \run{u_j}{m(j-1)+1}$ and the upper $m$-apsis $\set{m(j-1)+1, \ldots, mj} \in \run{m(j-1)+1}{l_j}$ join and are removed; and
					\item Since $u_j+m-1$ is the right-most vertex in the right-most $m$-apsis of $\run{u_j}{m(j-1)+1}$, for each $i \in \set{u_j+m, \ldots, k-m+1}$, $\run{u_j}{m(j-1)+1}$ has vertical lines $\set{i, i'}$. Now $u_{j+1}$ is the left-most vertex in the left-most upper $m$-apsis in $\apmorph{u_{j+1}, \ldots, u_t}{l_{j+1}, \ldots, l_t}$. It follows from $u_{j+1} > u_j+m-1$ that the upper $m$-apses in $\apmorph{u_{j+1}, \ldots, u_t}{l_{j+1}, \ldots, l_t}$ are preserved. By an analogous argument the lower $m$-apses in $\apmorph{u_{j+1}, \ldots, u_t}{l_{j+1}, \ldots, l_t}$ are also preserved.
				\end{enumerate}
				
				Hence it follows that \[\apmorph{u_j, \ldots, u_t}{l_j, \ldots, l_t}= \run{u_j}{m(j-1)+1}\apmorph{u_{j+1}, \ldots, u_t}{l_{j+1}, \ldots, l_t}\run{m(j-1)+1}{l_j} = \run{u_j}{m(j-1)+1}\ldots\run{u_t}{m(t-1)+1}\run{m(t-1)+1}{l_t}\ldots\run{m(j-1)+1}{l_j}.\] Finally, by induction $\apmorph{u_1, \ldots, u_t}{l_1, \ldots, l_t} = \run{u_1}{1}\ldots\run{u_t}{m(t-1)+1}\run{m(t-1)+1}{l_t}\ldots\run{1}{l_1}$.
			\end{proof}
		\end{proposition}

		\begin{figure}[!ht]
			\caption[ ]{$\apmorph{3,6}{2,7} = \run{3}{1}\apmorph{6}{7}\run{1}{2} = \run{3}{1}\run{6}{4}\run{4}{7}\run{1}{2} \in \triapmon{9}$.}
			\label{fig:sigma3627factorisation}
			\vspace{5pt}
			\centering
			\input{chap_characterisations/tikz/fig-sigma3627factorisation.tex}
		\end{figure}
		
		\begin{corollary}
			$\apmorphs{m}{k} \subseteq \apsismon{m}{k}$.
		\end{corollary}
		
		\begin{definition} \label{def:apmorphismofbipartitions}
			Given bipartitions $\alpha, \beta \in \pttnmon{k}$ such that the number of upper $m$-apses in $\alpha$ is equal to the number of lower $m$-apses in $\beta$, we shall denote by $\apmorph{\alpha}{\beta}$ the $m$-apmorphism with precisely the same upper $m$-apses as $\alpha$ and precisely the same lower $m$-apses as $\beta$.
		\end{definition}
		
	\subsection{Blocks that must appear}
		Our first step towards characterising the bipartitions that are elements of the $m$-apsis monoid $\apsismon{m}{k}$ is to establish which blocks must appear in all products of $m$-apsis generators.
	
		\begin{proposition}
			For each $m \in \intsge{3}$, $k \in \intsge{m}$, with the exception of the identity, every element of the $m$-apsis monoid $\apsismon{m}{k}$ must contain at least one upper $m$-apsis and at least one lower $m$-apsis. 
			
			\begin{proof}
				Trivially follows from upper non-transversals being preserved under right-multiplication along with lower non-transversals being preserved under left-multiplication.
			\end{proof}
		\end{proposition}

		Additionally, note that when $k=m$ there is a single $m$-apsis generator containing exactly one upper $m$-apsis and one lower $m$-apsis. Hence when considering all $k \in \posints$, no further conditions may be placed on blocks that must appear in products of $m$-apsis generators.	
		
	\subsection{Feasible block types} \label{sec:feasibleblocktypes}
 \subsubsection{Block types that do appear}
    	In this subsection we seek to establish given $m \in \ints_{\geq 3}$ and $k \in \ints_{\geq m}$, for which combinations of $\mu, \gamma \in \nonnegints$ such that $\mu + \gamma > 0$ do there exist products of $m$-apsis generators containing a block of type $(\mu, \gamma)$?
    	
    	Note a product of $m$-apsis generators must contain at least one upper $m$-apsis and at least one lower $m$-apsis. Hence after excluding one upper $m$-apsis and one lower $m$-apsis to cover cases where $k < 2m$, which is done implicitly from this point, blocks in a product of $m$-apsis generators may contain at most $k-m$ upper vertices and at most $k-m$ lower vertices. Consequently we must have $\mu, \gamma \leq k-m$.
    	
    	For $k = m$ the $m$-apsis monoid trivially consists of the identity and the single $m$-apsis generator $\apgen{m}{m}$, one way in which this may be seen is observing it is trivially the case that each $m$-apsis generator is idempotent, ie. ${(\apgen{m}{i})}^2 = \apgen{m}{i}$ for all $m \in \posints$, $k \in \intsge{m}$ and $i \in \set{1, \ldots, k-m+1}$. 
    	
    	For $k > m$ we have seen that $m$-apses, which satisfy $\set{\mu, \gamma} = \set{0, m}$, and transversal lines, which have $\mu = \gamma = 1$, appear in products of $m$-apsis generators. To form an idea of which other combinations of $\mu$ and $\gamma$ do appear we proceed by examining block types that appear in products of triapsis generators for smaller values of $k$ greater than three:
    	\begin{enumerate}
	    	\item When $k=4$ there are two triapsis generators, which generate $\triapmon{4} = \set{\id{4}, \run{1}{1}, \run{1}{2}, \run{2}{1}, \run{2}{2}}$. Since runs contain only $m$-apses and transversal lines, no further block types appear;
	    	\item When $k=5$ transversals appear of type $(2,2)$, for example in the product $\triapgen{1}\triapgen{3}$ as illustrated in Figure \ref{fig:trans22}. Hence for combinations of $\mu, \gamma \in \set{0, 1, 2}$ such that $\mu + \gamma > 0$, blocks of type $(\mu,\gamma)$ appear in elements of $\triapmon{5}$ when $\mu = \gamma$; and
	    	\item When $k=6$ transversals also appear of type $(3,3)$, for example in the product $\triapgen{1}\triapgen{3}\triapgen{2}\triapgen{4}$ as illustrated in Figure $\ref{fig:trans33}$. Hence for combinations of $\mu, \gamma \in \set{0, 1, 2, 3}$ such that $\mu + \gamma > 0$, blocks of type $(\mu,\gamma)$ appear in elements of $\triapmon{6}$ when $\mu = \gamma$ or $\set{\mu, \gamma} = \set{0, 3}$.
    	\end{enumerate}
		
		\begin{figure}[!ht]
			\caption[ ]{Given $m=3$ and $k=5$, $\triapgen{1}\triapgen{3}$ has a type $(2,2)$ transversal.}
			\label{fig:trans22}
			\vspace{5pt}
			\centering
			\input{chap_characterisations/tikz/fig-trans22.tex}
		\end{figure} 
		
		\begin{figure}[!ht]
			\caption[ ]{Given $m=3$ and $k=6$, $\triapgen{1}\triapgen{3}\triapgen{2}\triapgen{4}$ has a type $(3,3)$ transversal.}
			\label{fig:trans33}
			\vspace{5pt}
			\centering
			\input{chap_characterisations/tikz/fig-trans33.tex}
		\end{figure} 
		
		Before examining when $k \geq 7$, at this point it is reasonable to conjecture that, given $m \in \ints_{\geq 3}$, $k \in \ints_{\geq m}$ and $\mu \in \set{1, \ldots, k-m}$, there exist products of $m$-apsis generators containing a block of type $(\mu, \mu)$. 

		\begin{definition}
			Let $\tomegabar{(\mu, \mu)}{k}$ denote the planar bipartition consisting of: 
			\begin{enumerate}
				\item the upper $m$-apsis $\set{1, \ldots, m}$ containing the $m$ left-most upper vertices;
				\item the lower $m$-apsis $\set{(\mu+1)', \ldots, (\mu+m)'}$;
				\item the $(\mu,\mu)$-transversal $\set{m+1, \ldots, m+\mu, 1', \ldots, \mu'}$; and
				\item the vertical lines $\set{j, j'}$ where $j \in \set{\mu+m+1, \ldots, k}$.
			\end{enumerate} 
		\end{definition}
		
\begin{figure}[!ht]
			\caption[ ]{Given $m=2$ and $k=5$,}
			\label{fig:omegabar22}
			\vspace{5pt}
			\centering
			\input{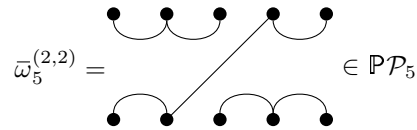}
		\end{figure}

		For example given $m=3$, $\tomegabar{(2,2)}{5} \in \ppttnmon{5}$ is depicted in Figure \ref{fig:omegabar22}.
		
		We proceed by establishing in Proposition \ref{prop:tomegabarmumuinapsismon} that $\tomegabar{(\mu,\mu)}{k} \in \apsismon{m}{k}$. The reader may find it useful to refer to Figure \ref{fig:omegabarmmfactorisation} which illustrates that given $m=3$, $k \in \ints_{\geq m}$, and $\lambda \in \set{2, \ldots, k-m}$, $\tomegabar{(\lambda,\lambda)}{k} = \tomegabar{(\lambda-1,\lambda-1)}{k}\triapgen{\lambda-1}\triapgen{\lambda+1}$.
		
		\begin{proposition} \label{prop:tomegabarmumuinapsismon}
			For each $m \in \ints_{\geq 3}$, $k \in \ints_{\geq m}$ and $\mu \in \set{1, \ldots, k-m}$,
			\[\tomegabar{(\mu,\mu)}{k} = \run{1}{2}\prod^\mu_{i=2}\apgen{m}{i-1}\apgen{m}{i+1} \in \apsismon{m}{k}.\]
			
			\begin{proof}
				That $\tomegabar{(1,1)}{k}$ is equal to the run $\run{1}{2}$ follows by definition. Let $\lambda \in \set{2, \ldots, \mu}$ and suppose that $\tomegabar{(\lambda-1,\lambda-1)}{k} = \run{1}{2}\prod^{\lambda-1}_{i=2}\apgen{m}{i-1}\apgen{m}{i+1}$. Note that $\apgen{m}{\lambda-1}\apgen{m}{\lambda+1} = \id{\lambda-2}\oplus\tomegabar{(2,2)}{k-\lambda+2}$ consists of: 
				\begin{enumerate}
					\item two $m$-apses $\set{\lambda-1, \ldots, \lambda+m-2}$ and $\set{(\lambda+1)', \ldots, (\lambda+m)'}$;
					\item a $(2,2)$-transversal $\set{\lambda+m-1, \lambda+m, (\lambda-1)', \lambda'}$; and
					\item vertical lines $\set{j, j'}$ where $j \in \set{1, \ldots, \lambda-2, \lambda+m+1, \ldots, k}$.
				\end{enumerate}
				When forming the product $\tomegabar{\lambda-1}{\lambda-1}\apgen{m}{\lambda-1}\apgen{m}{\lambda+1}$:
				\begin{enumerate}
					\item The upper $m$-apsis $\set{1, \ldots, m}$ in $\tomegabar{\lambda-1}{\lambda-1}$ is preserved;
					\item The lower $m$-apsis $\set{(\lambda+1)', \ldots, (\lambda+m)'}$ in $\apgen{m}{\lambda-1}\apgen{m}{\lambda+1}$ is preserved; and
					\item The vertical line $\set{\lambda+m, (\lambda+m)'}$ in $\tomegabar{\lambda-1}{\lambda-1}$ joins to the type $(2,2)$ block $\{\lambda+m-1, \lambda+m, (\lambda-1)',$ $\lambda'\}$ in $\apgen{m}{\lambda-1}\apgen{m}{\lambda+1}$, which joins to the lower $m$-apsis $\set{\lambda', \ldots, (\lambda+m-1)'}$ in $\tomegabar{\lambda-1}{\lambda-1}$, which joins to the upper $m$-apsis $\set{\lambda-1, \ldots, \lambda+m-2}$ in $\apgen{m}{\lambda-1}\apgen{m}{\lambda+1}$, which joins to the type $(\lambda-1,\lambda-1)$ block $\set{m+1, \ldots, m+\lambda-1, 1', \ldots, (\lambda-1)'}$ in $\tomegabar{\lambda-1}{\lambda-1}$, which joins to the vertical lines $\set{j, j'}$ in $\apgen{m}{\lambda-1}\apgen{m}{\lambda+1}$ where $j \in \set{1, \ldots, \lambda-2}$. This forms the block $\set{m+1, \ldots, m+\lambda, 1', \ldots, \lambda'}$ when collecting the upper vertices in the aforementioned blocks of $\tomegabar{\lambda-1}{\lambda-1}$ along with the lower vertices in the aforementioned blocks of $\apgen{m}{\lambda-1}\apgen{m}{\lambda+1}$. 
				\end{enumerate}
				Therefore we have $\tomegabar{(\lambda, \lambda)}{k} = \tomegabar{(\lambda-1, \lambda-1)}{k}\apgen{m}{\lambda-1}\apgen{m}{\lambda+1} = \run{1}{2}\prod^{\lambda}_{i=2}\apgen{m}{i-1}\apgen{m}{i+1}$, and hence it follows by induction that $\tomegabar{(\mu, \mu)}{k} = \run{1}{2}\prod^{\mu}_{i=2}\apgen{m}{i-1}\apgen{m}{i+1} \in \apsismon{m}{k}$.
			\end{proof}
		\end{proposition}
		
		\begin{figure}[!ht]
			\caption[ ]{Given $m=3$, $k \in \ints_{\geq m}$, and $\lambda \in \set{2, \ldots, k-m}$, $\tomegabar{(\lambda, \lambda)}{k} = \tomegabar{(\lambda-1, \lambda-1)}{k}\apgen{3}{\lambda-1}\apgen{3}{\lambda+1}$.}
			\label{fig:omegabarmmfactorisation}
			\vspace{5pt}
			\centering
			\input{chap_characterisations/tikz/fig-omegabarmmfactorisation.tex}
		\end{figure} 

		Next we seek to identify whether for any further combinations of $\mu, \gamma \in \set{0, \ldots, k-m}$ such that $\mu + \gamma > 0$, there exist products of $m$-apsis generators containing a block of type $(\mu, \gamma)$. To do so we examine products of triapsis generators when $k \geq 7$:
		\begin{enumerate}
			\item When $k=7$ the transversal types $(4,1)$ and $(1,4)$ appear, for example in the products $\tomegabar{4}{4}\triapgen{4}$ and $\triapgen{4}\tomegabar{4}{4}$;
			\item When $k=8$ the transversal types $(5,2)$ and $(2,5)$ appear, for example in the products $\tomegabar{5}{5}\triapgen{4}$ and $\triapgen{4}\tomegabar{5}{5}$; and
			\item When $k=9$ the transversal types $(6,3)$ and $(3,6)$ along with the non-transversal types $(6,0)$ and $(0,6)$ appear, for example in the products $\tomegabar{6}{6}\triapgen{4}$, $\triapgen{4}\tomegabar{6}{6}$, $\tomegabar{6}{6}\triapgen{4}\triapgen{7}$ and $\triapgen{4}\triapgen{7}\tomegabar{6}{6}$.
		\end{enumerate}
		
		At this point it would be reasonable to conjecture that for each $m \in \intsge{3}$, $k \in \intsge{m}$, and $\mu, \gamma \in \set{0, \ldots, k-m}$ such that $\mu + \gamma > 0$ and $\congmod{\mu}{\gamma}{m}$, there exist products of $m$-apsis generators that contain a block of type $(\mu, \gamma)$.

		\begin{definition}
			For each $m \in \ints_{\geq 3}$, $k \in \ints_{\geq m}$ and $\mu, \gamma \in \set{0, \ldots, k-m}$ such that $\mu + \gamma > 0$ and $\congmod{\mu}{\gamma}{m}$, let $\overbar{k} = \max\set{\mu, \gamma} + m$ and let $\tomegabar{(\mu, \gamma)}{k}$ denote the planar bipartition that contains: 
			\begin{enumerate}
				\item for each $j \in \set{1, \ldots, \frac{\overbar{k} - \mu}{m}}$, the upper $m$-apsis $\{m(j-1) + 1, \ldots, mj\}$;
				\item for each $j \in \set{1, \ldots, \frac{\overbar{k} - \gamma}{m}}$, the lower $m$-apsis $\{(\gamma+m(j-1)+1)', \ldots, (\gamma+mj)'\}$;
				\item a type $(\mu, \gamma)$ transversal $\set{\overbar{k}-\mu + 1, \ldots, \overbar{k}, 1', \ldots, \gamma'}$; and
				\item for each $j \in \set{\overbar{k}+1, \ldots, k}$, the vertical line $\set{j, j'}$.
			\end{enumerate} 
		\end{definition}
		
		For example, given $m=3$, $\tomegabar{(5,2)}{9} \in \ppttnmon{9}$ is depicted in Figure \ref{fig:omega52}.
		
		\begin{figure}[!ht]
			\caption[ ]{Given $m=3$ and $k=9$,}
			\label{fig:omega52}
			\vspace{5pt}
			\centering
			\input{chap_characterisations/tikz/fig-omega52.tex}
		\end{figure} 
		
\begin{figure}[!ht]
			\caption[ ]{Given $m=3$ and $k=12$,}
			\label{fig:omega82=omega88sigma4747}
			\vspace{5pt}
			\centering
			\input{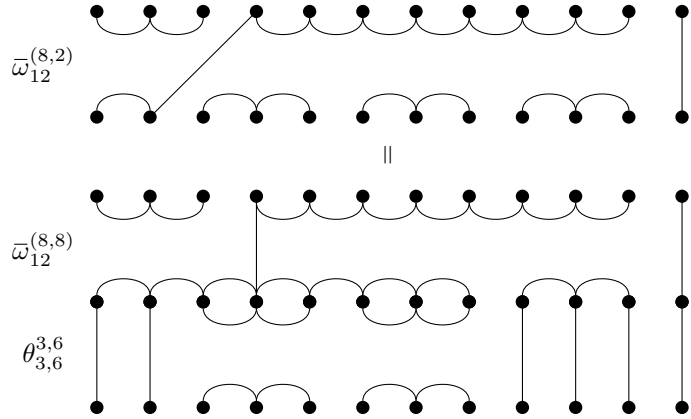}
		\end{figure}

		We proceed by establishing in Proposition \ref{prop:tomegabarmugammainapsismon} that $\tomegabar{(\mu, \gamma)}{k} \in \apsismon{m}{k}$. The reader may find it useful to refer to Figure \ref{fig:omega82=omega88sigma4747} which illustrates that, given $m=3$ and $k=12$, $\tomegabar{(8,2)}{12} = \tomegabar{(8,8)}{12}\apmorph{3,6}{3,6}$.
		
		\begin{proposition} \label{prop:tomegabarmugammainapsismon}
			For each $m \in \ints_{\geq 3}$, $k \in \ints_{\geq m}$ and $\mu, \gamma \in \set{0, \ldots, k-m}$ such that $\mu + \gamma > 0$ and $\congmod{\mu}{\gamma}{m}$,
			
			$\tomegabar{(\mu, \gamma)}{k} = \begin{cases}\tomegabar{(\mu, \mu)}{k}\apmorph{\gamma+1, \gamma+m+1, \ldots, \mu-m+1}{\gamma+1, \gamma+m+1, \ldots, \mu-m+1} & \text{ if } \mu > \gamma; \\ \run{1}{2}\prod^\mu_{i=2}\apgen{m}{i-1}\apgen{m}{i+1} & \text{ if } \mu = \gamma; \text{ and} \\ \apmorph{m+1, 2m+1, \ldots, \gamma-\mu+1}{m+1, 2m+1, \ldots, \gamma-\mu+1}\tomegabar{(\gamma, \gamma)}{k} & \text{ if } \mu < \gamma, \\ \end{cases}$
			
			and hence $\tomegabar{(\mu, \gamma)}{k} \in \apsismon{m}{k}$.
			
			\begin{proof}
				We already established in Proposition \ref{prop:tomegabarmumuinapsismon} that if $\mu = \gamma$ then $\tomegabar{(\mu, \gamma)}{k} = \run{1}{2} \prod^\mu_{i=2} \apgen{m}{i-1} \apgen{m}{i+1} \in \apsismon{m}{k}$. 
				
				If $\mu > \gamma$ then when forming the product $\tomegabar{(\mu, \mu)}{k}\apmorph{\gamma+1, \gamma+m+1, \ldots, \mu-m+1}{\gamma+1, \gamma+m+1, \ldots, \mu-m+1}$:
				\begin{enumerate}
					\item the upper $m$-apsis $\set{1, \ldots, m}$ in $\tomegabar{(\mu, \mu)}{k}$ is preserved;
					\item the lower $m$-apsis $\set{(\mu+1)', \ldots, (\mu+m)'}$ in $\tomegabar{(\mu,\mu)}{k}$ is preserved since it connects to the vertical lines $\set{j, j'}$ in $\apmorph{\gamma+1, \gamma+m+1, \ldots, \mu-m+1}{\gamma+1, \gamma+m+1, \ldots, \mu-m+1}$ where $j \in \set{\mu+1, \ldots, \mu+m}$;
					\item the $(\mu, \gamma)$-transversal $\set{m+1, \ldots, m+\mu, 1', \ldots, \mu'}$ in $\tomegabar{(\mu,\mu)}{k}$ connects to 
					\begin{enumerate}
						\item the $m$-apses $\set{\gamma+(j-1)m+1, \ldots, \gamma+jm}$ in $\apmorph{\gamma+1, \gamma+m+1, \ldots, \mu-m+1}{\gamma+1, \gamma+m+1, \ldots, \mu-m+1}$ where $j \in \set{1, \ldots, \frac{\mu-\gamma}{m}}$; and
						\item the lines $\set{j, j'}$ in $\apmorph{\gamma+1, \gamma+m+1, \ldots, \mu-m+1}{\gamma+1, \gamma+m+1, \ldots, \mu-m+1}$ where $j \in \set{1, \ldots, \gamma}$, 
					\end{enumerate}
					collectively forming the block $\set{m+1, \ldots, m+\mu, 1', \ldots, \gamma'}$; and
					\item for each $j \in \set{\mu+m+1, \ldots, k}$, the vertical line $\set{j, j'}$ in $\tomegabar{(\mu, \mu)}{k}$ connects to the vertical line $\set{j, j'}$ in $\apmorph{\gamma+1, \gamma+m+1, \ldots, \mu-m+1}{\gamma+1, \gamma+m+1, \ldots, \mu-m+1}$.
				\end{enumerate}
				Hence $\tomegabar{(\mu, \gamma)}{k} = \tomegabar{(\mu, \mu)}{k}\apmorph{\gamma+1, \gamma+m+1, \ldots, \mu-m+1}{\gamma+1, \gamma+m+1, \ldots, \mu-m+1} \in \apsismon{m}{k}$. It follows analogously that if $\mu < \gamma$ then $\tomegabar{(\mu, \gamma)}{k} = \apmorph{m+1, 2m+1, \ldots, \gamma-\mu+1}{m+1, 2m+1, \ldots, \gamma-\mu+1}\tomegabar{(\gamma, \gamma)}{k} \in \apsismon{m}{k}$.
			\end{proof}
		\end{proposition}
		
		\begin{corollary}
			For each $m \in \intsge{3}$, $k \in \intsge{m}$ and $\mu, \gamma \in \set{0, \ldots, k-m}$ such that $\mu + \gamma > 0$ and $\congmod{\mu}{\gamma}{m}$, blocks of type $(\mu, \gamma)$ do appear in some elements of the $m$-apsis monoid $\apsismon{m}{k}$.
			
			\begin{proof}
				Trivially follows from Proposition \ref{prop:tomegabarmugammainapsismon} where we established that the well-defined bipartition $\tomegabar{(\mu,\gamma)}{k}$, which contains a block of type $(\mu, \gamma)$, may be factorised into a product of $m$-apsis generators.
			\end{proof}
		\end{corollary}

	\subsubsection{Block types that must not appear}
		Given $m \in \ints_{\geq 3}$, $k \in \ints_{\geq m}$ and $\mu, \gamma \in \set{0, \ldots, k-m}$ such that $\ncongmod{\mu}{\gamma}{m}$, it remains to establish whether there exist products of $m$-apsis generators containing a block of type $(\mu, \gamma)$. 
		
		\begin{proposition} \label{prop:apsismonsubsetneqmodmon}
			For each $m \in \intsge{2}$ and $k \in \intsge{m}$, the $m$-apsis monoid $\apsismon{m}{k}$ is a proper submonoid of the planar mod-$m$ monoid $\pmodmon{m}{k}$.
			
			\begin{proof}
				Containment follows from $m$-apsis generators trivially sitting inside the planar mod-$m$ monoid $\pmodmon{m}{k}$, that is $\apgen{m}{1}, \ldots, \apgen{m}{k-m+1} \in \pmodmon{m}{k}$. Inequality trivially follows from $(2,2)$-transapsis generators not being elements of the $m$-apsis generated monoid $\apsismon{m}{k}$.
			\end{proof}
		\end{proposition}
		
		\begin{corollary}
			For each $m \in \intsge{3}$, $k \in \intsge{m}$ and $\mu, \gamma \in \set{0, \ldots, k-m}$ such that $\ncongmod{\mu}{\gamma}{m}$, blocks of type $(\mu, \gamma)$ do not appear in any elements of the $m$-apsis monoid $\apsismon{m}{k}$.
		\end{corollary}

	\subsection{Bounding $\apsismon{m}{k}$ above by $\apsisbound{m}{k}$}
		\begin{definition}
			For each $m \in \intsge{3}$ and $k \in \intsge{m}$, we shall denote by $\apsisbound{m}{k}$ the set of all $\alpha \in \pmodmon{m}{k}$ such that either:
			\begin{enumerate}
				\item $\alpha$ is the identity $\id{k}$; or 
				\item $\alpha$ contains at least one upper $m$-apsis and at least one lower $m$-apsis.
			\end{enumerate} 
		\end{definition}
		
		\begin{proposition}
			For each $m \in \intsge{3}$ and $k \in \intsge{m}$, $\apsisbound{m}{k}$ is a monoid.
			
			\begin{proof}
				Closure under multiplication follows from:
				\begin{enumerate}
					\item the planar mod-$m$ monoid $\pmodmon{m}{k}$ being closed under multiplication; and
					\item upper $m$-apses and dually lower $m$-apses being preserved under right-multiplication and left-multiplication respectively.
				\end{enumerate}
			\end{proof}
		\end{proposition}
		
		\begin{proposition} \label{prop:apsismonpropersubmonapsisbound}
			For each $m \in \intsge{3}$ and $k \in \intsge{m}$, the $m$-apsis monoid $\apsismon{m}{k}$ is a submonoid of $\apsisbound{m}{k}$.
			
			\begin{proof}
				We have already established that products of $m$-apsis generators must:
				\begin{enumerate}
					\item sit inside the planar mod-$m$ monoid; and
					\item contain at least one upper $m$-apsis and at least one lower $m$-apsis.
				\end{enumerate}
				Hence containment follows.
			\end{proof}
		\end{proposition}

	\subsection{Building blocks}
		It remains to establish whether $\apsisbound{m}{k}$ is equal to the $m$-apsis monoid $\apsismon{m}{k}$, that is whether each element of $\apsisbound{m}{k}$ may be factorised into a product of $m$-apsis generators. Before doing so, it will be convenient for us to first identify some further well-defined elements and subsets of $\apsisbound{m}{k}$, then show that they are contained within the $m$-apsis monoid $\apsismon{m}{k}$.
		
		Both for the sake of succinctness and due to there being no notational ambiguity in doing so, when we need to require $\congmod{\mu}{\gamma}{m}$ and $\congmod{\mu}{k}{m}$ we simply state that we require $\tcongmod{\mu}{\gamma}{k}{m}$.
		
		Note that for each $m \in \intsge{3}$, $k \in \intsge{m}$, and $\mu, \gamma \in \set{0, \ldots, k-m}$ such that $\mu + \gamma > 0$ and $\tcongmod{\mu}{\gamma}{k}{m}$, in order to be able to form a bipartition containing precisely one block of type $(\mu, \gamma)$ then only $m$-apses, we additionally require that $\congmod{\mu}{k}{m}$, or equivalently that $\congmod{\gamma}{k}{m}$. 
		
		\begin{definition} \label{def:tomegamugamma}
			For each $m \in \intsge{3}$, $k \in \intsge{m}$ and $\mu, \gamma \in \set{0, \ldots, k-m}$ such that $\mu + \gamma > 0$ and $\tcongmod{\mu}{\gamma}{k}{m}$, let $\tomega{(\mu, \gamma)}{k}$ denote the bipartition in $\apsisbound{m}{k}$ containing:
			\begin{enumerate}
				\item the type $(\mu, \gamma)$ block $\set{1, \ldots, \mu, 1', \ldots, \gamma'}$;
				\item for each $j \in \set{1, \ldots, \frac{k-\mu}{m}}$, the upper $m$-apsis $\{\mu + m(j-1) + 1, \ldots, \mu + mj\}$; and
				\item for each $j \in \set{1, \ldots, \frac{k-\gamma}{m}}$, the lower $m$-apsis $\{(\gamma + m(j-1) + 1)', \ldots, (\gamma + mj)'\}$.
			\end{enumerate} 
		\end{definition}
		
\begin{figure}[!ht]
			\caption[ ]{Given $m=4$ and $k=9$,}
			\label{fig:omega51}
			\vspace{5pt}
			\centering
			\input{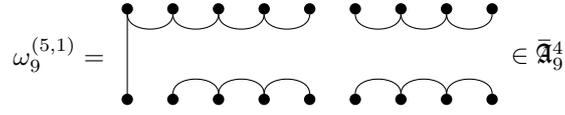}
		\end{figure}

		For example given $m=4$ and $k=9$, Figure \ref{fig:omega51} depicts $\tomega{(5, 1)}{9} \in \apsisbound{4}{9}$.
		
\begin{figure}[!ht]
			\caption[ ]{Given $m=3$ and $k=11$,}
			\label{fig:omega25_11=omegabar25_11sigma369169}
			\vspace{5pt}
			\centering
			\input{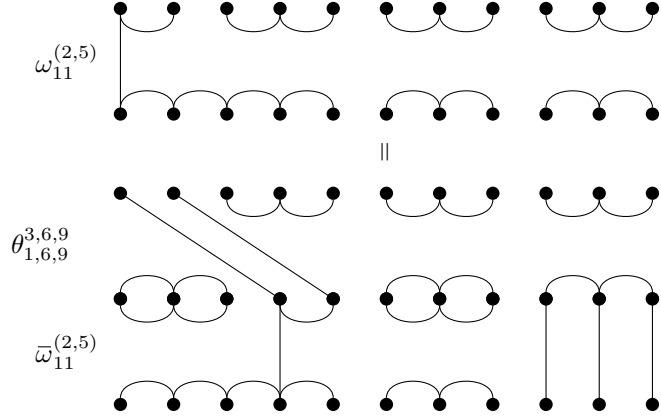}
		\end{figure}

		We proceed by establishing in Proposition \ref{prop:tomegamugammainapsismon} that $\tomega{(\mu, \gamma)}{k} \in \apsismon{m}{k}$. The reader may find it useful to refer to Figure \ref{fig:omega25_11=omegabar25_11sigma369169} which illustrates that, given $m=3$ and $k=11$, $\tomega{(2,5)}{11} = \apmorph{3,6,9}{1,6,9}\tomegabar{(2,5)}{11} \in \apsismon{3}{11}$.
		
    	\begin{proposition} \label{prop:tomegamugammainapsismon}
    		For each $m \in \intsge{3}$, $k \in \intsge{m}$ and $\mu, \gamma \in \set{0, \ldots, k-m}$ such that $\mu + \gamma > 0$ and $\tcongmod{\mu}{\gamma}{k}{m}$, \[\tomega{(\mu, \gamma)}{k} = \apmorph{\mu+1, \mu+m+1, \mu+2m+1, \ldots, k-m+1}{\ \ \ 1, \mu+m+1, \mu+2m+1, \ldots, k-m+1}\tomegabar{(\mu, \gamma)}{k} \in \apsismon{m}{k}.\]
    		
    		\begin{proof}
    			Let $\overbar{k} = \max\set{\mu, \gamma} +m$ and recall that for each $j \in \set{\overbar{k}+1, \ldots, k}$, $\tomegabar{(\mu, \gamma)}{k}$ contains the vertical line $\set{j, j'}$. 
    			
    			When forming the product $\apmorph{\mu+1, \mu+m+1, \mu+2m+1, \ldots, k-m+1}{\ \ \ 1, \mu+m+1, \mu+2m+1, \ldots, k-m+1}\tomegabar{(\mu, \gamma)}{k}$:
    			\begin{enumerate}
    				\item the type $(\mu, \gamma)$ block $\set{1, \ldots, \mu, 1', \ldots, \gamma'}$ in $\tomegabar{(\mu, \gamma)}{k}$ connects to the transversal lines $\set{j, (m+j)'}$ in $\apmorph{\mu+1, \mu+m+1, \mu+2m+1, \ldots, k-m+1}{\ \ \ 1, \mu+m+1, \mu+2m+1, \ldots, k-m+1}$ where $j \in \set{1, \ldots, \mu}$, forming the type $(\mu, \gamma)$ block $\{1, \ldots, \mu, 1', \ldots,$ $\gamma'\}$;
    				\item the upper $m$-apses in $\apmorph{\mu+1, \mu+m+1, \mu+2m+1, \ldots, k-m+1}{\ \ \ 1, \mu+m+1, \mu+2m+1, \ldots, k-m+1}$, which are upper $m$-apses in $\tomega{(\mu, \gamma)}{k}$, are preserved;
    				\item the lower $m$-apsis $\set{1', \ldots, m'}$ in $\apmorph{\mu+1, \mu+m+1, \mu+2m+1, \ldots, k-m+1}{\ \ \ 1, \mu+m+1, \mu+2m+1, \ldots, k-m+1}$ and upper $m$-apsis $\set{1, \ldots, m}$ in $\tomegabar{(\mu, \gamma)}{k}$ join and are removed;
    				\item for each $j \in \set{1, \ldots, \frac{k-\overbar{k}}{m}}$, the lower $m$-apsis $\{(\overbar{k} + m(j-1) + 1)', \ldots, (\overbar{k} + mj)'\}$ from the $m$-apmorphism $\apmorph{\mu+1, \mu+m+1, \mu+2m+1, \ldots, k-m+1}{\ \ \ 1, \mu+m+1, \mu+2m+1, \ldots, k-m+1}$, which joins to vertical lines in $\tomegabar{(\mu, \gamma)}{k}$, is preserved;
					\item for each $j \in \set{2, \ldots, \frac{\overbar{k}-\mu}{m}}$, the upper $m$-apsis $\{\mu+m(j-1) + 1, \ldots, \mu+mj\}$ in $\tomegabar{(\mu, \gamma)}{k}$ and the lower $m$-apsis $\{(\mu+m(j-1)+1)', \ldots, (\mu+mj)'\}$ in $\apmorph{\mu+1, \mu+m+1, \mu+2m+1, \ldots, k-m+1}{\ \ \ 1, \mu+m+1, \mu+2m+1, \ldots, k-m+1}$ join and are removed; and
					\item for each $j \in \set{1, \ldots, \frac{\overbar{k}-\gamma}{m}}$, the lower $m$-apsis $\{(\gamma+m(j-1)+1)', \ldots, (\gamma+mj)'\}$ in $\tomegabar{(\mu, \gamma)}{k}$ is preserved.
    			\end{enumerate}
    			
    			Hence $\tomega{(\mu, \gamma)}{k} = \apmorph{\mu+1, \mu+m+1, \mu+2m+1, \ldots, k-m+1}{\ \ \ 1, \mu+m+1, \mu+2m+1, \ldots, k-m+1}\tomegabar{(\mu, \gamma)}{k} \in \apsismon{m}{k}$.
    		\end{proof}
    	\end{proposition}
    	
    \subsection{Transversal building blocks}        	
       	We now turn our attention to bipartitions in $\apsisbound{m}{k}$ where non-transversals must be $m$-apses. 
       	
       	\begin{definition}
       		For each $m \in \intsge{3}$ and $k \in \intsge{m}$, let $\aptrans{m}{k}$ denote the set of all $\alpha \in \apsisbound{m}{k}$ such that every non-transversal $b \in \alpha$ is an $m$-apsis. We shall refer to elements of $\aptrans{m}{k}$ as \textit{transversal building blocks}.
       	\end{definition}
       	
       	We proceed in this subsection to establish that $\aptrans{m}{k}$ is contained within the $m$-apsis monoid $\apsismon{m}{k}$. In order to do so we partition $\aptrans{m}{k}$ based on transversal types from left to right then establish our desired result inductively.
       	
       	Note that for each $m \in \intsge{3}$, $k \in \intsge{m}$, $r \in \set{0, \ldots, k-m}$, and $\mu_1, \ldots, \mu_r$, $\gamma_1, \ldots, \gamma_r \in \set{1, \ldots, k-m}$ such that $\congmod{\mu_j}{\gamma_j}{k}$ for all $j \in \set{1, \ldots, r}$ and $\Sigma^r_{j=1}\mu_j, \Sigma^r_{j=1}\gamma_j \leq k-m$, in order to be able to form a bipartition of rank $r$ with a distinct block of type $(\mu_j, \gamma_j)$ designated for each $j \in \set{1, \ldots, r}$, we additionally require that $\congmod{\Sigma^r_{j=1}\mu_j}{k}{m}$, or equivalently that $\congmod{\Sigma^r_{j=1}\gamma_j}{k}{m}$. 
       	
       	\begin{definition}
       		For each $m \in \intsge{3}$, $k \in \intsge{m}$, $r \in \set{0, \ldots, k-m}$, and $\mu_1, \ldots, \mu_r, \gamma_1, \ldots, \gamma_r \in \set{1, \ldots, k-m}$ such that $\congmod{\mu_j}{\gamma_j}{k}$ for all $j \in \set{1, \ldots, r}$, $\Sigma^r_{j=1}\mu_j, \Sigma^r_{j=1}\gamma_j \leq k-m$ and $\congmod{\Sigma^r_{j=1}\mu_j}{k}{m}$, let $\aptrans{(\mu_1, \gamma_1),\ldots,(\mu_r, \gamma_r)}{k}$ denote the set of all $\omega \in \aptrans{m}{k}$ such that $\omega$ contains precisely:
			\begin{enumerate}
				\item $\frac{k - \Sigma^r_{j=1}\mu_j}{m}$ upper $m$-apses;
				\item $\frac{k - \Sigma^r_{j=1}\gamma_j}{m}$ lower $m$-apses; and
				\item $r$ transversals of type $(\mu_1, \gamma_1)$, \ldots, $(\mu_r, \gamma_r)$ from left to right.
			\end{enumerate} 
       	\end{definition}
       	
       	Note that by definition $\aptrans{m}{k} = \bigcup_{\substack{r \in \set{0, \ldots, k-m} \\ \mu_1, \ldots, \mu_r, \gamma_1, \ldots, \gamma_r \in \set{1, \ldots, k-m} \\ \congmod{\mu_j}{\gamma_j}{m}\ \forall\ j \in \set{1, \ldots, r} \\ \Sigma^r_{j=1}\mu_j, \Sigma^r_{j=1}\gamma_j \leq k-m \\ \congmod{\Sigma^r_{j=1}\mu_j}{k}{m}}} \aptrans{(\mu_1, \gamma_1), \ldots, (\mu_r, \gamma_r)}{k}$.
		
		Next we seek to establish whether $\aptrans{(\mu_1, \gamma_1),\ldots,(\mu_r, \gamma_r)}{k} \subseteq \apsismon{m}{k}$. We begin by establishing in Proposition \ref{prop:aptransmujgammajgeneratesitself} that for each $\omega \in \aptrans{(\mu_1, \gamma_1),\ldots,(\mu_r, \gamma_r)}{k}$, each element of $\aptrans{(\mu_1, \gamma_1),\ldots,(\mu_r, \gamma_r)}{k}$ may be factorised into a product containing $\omega$ and $m$-apmorphisms. An example of the method outlined for factorisation is depicted in Figure \ref{fig:aptransprop}.
		
		\begin{proposition} \label{prop:aptransmujgammajgeneratesitself}
			For each $\omega \in \aptrans{(\mu_1, \gamma_1),\ldots,(\mu_r, \gamma_r)}{k}$, $\aptrans{(\mu_1, \gamma_1),\ldots,(\mu_r, \gamma_r)}{k} \subseteq \apmorphs{m}{k}\omega\apmorphs{m}{k}$.
			
			\begin{proof}
				Let $\psi \in \aptrans{(\mu_1, \gamma_1),\ldots,(\mu_r, \gamma_r)}{k}$. Since $\gamma$ and $\psi$ have the same number of upper $m$-apses as well as the same number of lower $m$-apses, both of the $m$-apmorphisms $\apmorph{\psi}{\omega^*}$ and $\apmorph{\omega^*}{\psi}$ are well-defined. When forming the product $\apmorph{\psi}{\omega^*}\omega\apmorph{\omega^*}{\psi}$:
				\begin{enumerate}
					\item the upper $m$-apses in $\apmorph{\psi}{\omega^*}$, which are identical to the upper $m$-apses in $\psi$, are preserved;
					\item the lower $m$-apses in $\apmorph{\omega^*}{\psi}$, which are identical to the lower $m$-apses in $\psi$, are preserved;
					\item the upper $m$-apses of $\omega$ and lower $m$-apses of $\apmorph{\psi}{\omega^*}$ join and are removed;
					\item the lower $m$-apses of $\omega$ and upper $m$-apses of $\apmorph{\omega^*}{\psi}$ join and are removed; and
					\item each of the upper vertices in the transversals of $\omega$ connects to a transversal line of $\apmorph{\psi}{\omega^*}$, and each of the lower vertices in the transversals of $\omega$ connects to a transversal line of $\apmorph{\omega^*}{\psi}$, consequently block types of transversals are preserved.
				\end{enumerate}
				Hence $\psi = \apmorph{\psi}{\omega^*}\omega\apmorph{\omega^*}{\psi}$, consequently $\aptrans{(\mu_1, \gamma_1),\ldots,(\mu_r, \gamma_r)}{k} \subseteq \apmorphs{m}{k}\omega\apmorphs{m}{k}$.
			\end{proof}
		\end{proposition}
		
		\begin{figure}[!ht]
			\caption[ ]{Given $m=3$,}
			\label{fig:aptransprop}
			\vspace{10pt}
			\centering
			\input{chap_characterisations/tikz/fig-aptransprop.tex}
		\end{figure}
		
		\begin{corollary} \label{cor:whenaptransinapsismon}
			If the intersection of $\aptrans{(\mu_1, \gamma_1),\ldots,(\mu_r, \gamma_r)}{k}$ with the $m$-apsis monoid $\apsismon{m}{k}$ is non-empty then $\aptrans{(\mu_1, \gamma_1),\ldots,(\mu_r, \gamma_r)}{k} \subseteq \apsismon{m}{k}$.
			
			\begin{proof}
				Suppose there exists $\omega \in \aptrans{(\mu_1, \gamma_1),\ldots,(\mu_r, \gamma_r)}{k} \cap \apsismon{m}{k}$, then it trivially follows from Proposition \ref{prop:aptransmujgammajgeneratesitself} that $\aptrans{(\mu_1, \gamma_1),\ldots,(\mu_r, \gamma_r)}{k} \subseteq \apmorphs{m}{k}\omega\apmorphs{m}{k} \subseteq \apsismon{m}{k}$.
			\end{proof}
		\end{corollary}
    	
    	\begin{proposition} \label{prop:aptransmugammainapsismon}
    		For each $m \in \intsge{3}$, $k \in \intsge{m}$, $\mu, \gamma \in \set{1, \ldots, k-m}$ such that $\tcongmod{\mu}{\gamma}{k}{m}$, $\aptrans{(\mu, \gamma)}{k} \subseteq \apsismon{m}{k}$.
    		
    		\begin{proof}
    			It follows by definition that $\tomega{(\mu, \gamma)}{k} \in \aptrans{(\mu, \gamma)}{k}$ and we established in Proposition \ref{prop:tomegamugammainapsismon} that $\tomega{(\mu, \gamma)}{k} \in \apsismon{m}{k}$, hence $\tomega{(\mu, \gamma)}{k} \in \aptrans{(\mu, \gamma)}{k} \cap \apsismon{m}{k}$. Finally, employing Corollary \ref{cor:whenaptransinapsismon} we have $\aptrans{(\mu, \gamma)}{k} \subseteq \apsismon{m}{k}$.
    		\end{proof}
    	\end{proposition}

		\begin{definition}
	    	For each $m \in \intsge{3}$, $k \in \intsge{m}$, $r \in \set{0, \ldots, k-m}$, and $\mu_1, \ldots, \mu_r, \gamma_1, \ldots, \gamma_r \in \set{1, \ldots, k-m}$ such that $\congmod{\mu_j}{\gamma_j}{k}$ for all $j \in \set{1, \ldots, r}$, $\Sigma^r_{j=1}\mu_j, \Sigma^r_{j=1}\gamma_j \leq k-m$ and $\congmod{\Sigma^r_{j=1}\mu_j}{k}{m}$, let $\tomega{(\mu_1, \gamma_1), \ldots, (\mu_r, \gamma_r)}{k}$ denote the element of $\aptrans{(\mu_1, \gamma_1), \ldots, (\mu_r, \gamma_r)}{k}$ containing: 
		    	\begin{enumerate}
		    		\item the type $(\mu_1, \gamma_1)$ transversal $\set{1, \ldots, \mu_1, 1', \ldots, \gamma_1'}$;
		    		\item for each $j \in \set{1, \ldots, \frac{k - \Sigma^r_{j=1}\mu_j}{m}}$, the upper $m$-apsis $\{\mu_1 + m(j-1) + 1, \ldots, \mu_1 + mj\}$;
		    		\item for each $j \in \set{1, \ldots, \frac{k - \Sigma^r_{j=1}\gamma_j}{m}}$, the lower $m$-apsis $\{(\gamma_1+m(j-1)+1)', \ldots, (\gamma_1 + mj)'\}$;
		    		\item for each $j \in \set{2, \ldots, r}$, the type $(\mu_j, \gamma_j)$ transversal $\big\{k - \Sigma^{r}_{i=j}\mu_i + 1, \ldots, k - \Sigma^{r}_{i=j+1}\mu_i, (k - \Sigma^{r}_{i=j}\gamma_i + 1)', \ldots, (k - \Sigma^{r}_{i=j+1}\gamma_i)'\big\}$,
		    	\end{enumerate} 
			and let $\tupsilon{(\mu_1, \gamma_1), \ldots, (\mu_r, \gamma_r)}{k}$ denote the element of $\aptrans{(\mu_1, \gamma_1), \ldots, (\mu_r, \gamma_r)}{k}$ containing: 
		    	\begin{enumerate}
		    		\item for each $j \in \set{1, \ldots, \frac{k - \Sigma^r_{j=1}\mu_j}{m}}$, the upper $m$-apsis $\{m(j-1) + 1, \ldots, mj\}$;
		    		\item for each $j \in \set{1, \ldots, \frac{k - \Sigma^r_{j=1}\gamma_j}{m}}$, the lower $m$-apsis $\{(m(j-1) + 1)', \ldots, (mj)'\}$;
		    		\item for each $j \in \set{1, \ldots, r}$, the type $(\mu_j, \gamma_j)$ transversal $\big\{k - \Sigma^{r}_{i=j}\mu_i + 1, \ldots, k - \Sigma^{r}_{i=j+1}\mu_i, (k - \Sigma^{r}_{i=j}\gamma_i + 1)', \ldots, (k - \Sigma^{r}_{i=j+1}\gamma_i)'\big\}$.
		    	\end{enumerate} 
		\end{definition}
    	
    	For example given $m=3$, $\tomega{(1,4),(1,1),(4,1)}{12}, \tupsilon{(1,4),(1,1),(4,1)}{12} \in \aptrans{(1,4),(1,1),(4,1)}{12}$ are depicted in Figure \ref{fig:upsilon411114}.
    	
    	\begin{figure}[!ht]
			\caption[ ]{Given $m=3$,}
			\label{fig:upsilon411114}
			\vspace{10pt}
			\centering
			\input{chap_characterisations/tikz/fig-upsilon411114.tex}
		\end{figure}
		
		We proceed by establishing that for each $m \in \intsge{3}$, $k \in \intsge{m}$, $r \in \set{0, \ldots, k-m}$ and $\mu_1, \ldots, \mu_r, \gamma_1, \ldots, \gamma_r$ $\in \set{1, \ldots, k-m}$ such that $\congmod{\mu_j}{\gamma_j}{k}$ for all $j \in \set{1, \ldots, r}$, $\Sigma^r_{j=1}\mu_j, \Sigma^r_{j=1}\gamma_j \leq k-m$ and $\congmod{\Sigma^r_{j=1}\mu_j}{k}{m}$, \[\tomega{(\mu_1, \gamma_1), \ldots, (\mu_r, \gamma_r)}{k} \in \begin{cases}(\aptrans{(\mu_1, \gamma_1)}{k-\Sigma^r_{j=2}\mu_j} \oplus \id{\Sigma^r_{j=2}\mu_j})(\id{\gamma_1} \oplus \aptrans{(\mu_2,\gamma_2), \ldots, (\mu_r, \gamma_r)}{k-\gamma_1}) & \mu_1 \geq \gamma_1; \\ (\id{\mu_1} \oplus \aptrans{(\mu_2,\gamma_2), \ldots, (\mu_r, \gamma_r)}{k-\mu_1})(\aptrans{(\mu_1, \gamma_1)}{k-\Sigma^r_{j=2}\gamma_j} \oplus \id{\Sigma^r_{j=2}\gamma_j}) & \mu_1 \leq \gamma_1. \end{cases}\] 
		
		For example, given $m=3$, Figure \ref{fig:aptrans411114element} illustrates that $\tomega{(4,1),(1,1),(1,4)}{12} = (\tomega{(4,1)}{10}\oplus\id{2})(\id{1}\oplus\tupsilon{(1,1), (1, 4)}{11})$ $\in \aptrans{(4,1),(1,1),(1,4)}{12}$ and Figure \ref{fig:aptrans1771element} illustrates that $\tomega{(1,7),(7,1)}{11} = (\id{1}\oplus\tupsilon{(7,1)}{10})(\tomega{(1,7)}{11}\oplus\id{1}) \in \aptrans{(1,7),(7,1)}{11}$.

		\begin{proposition} \label{prop:aptransfactorisation}
			For each $m \in \intsge{3}$, $k \in \intsge{m}$, $r \in \set{0, \ldots, k-m}$, and $\mu_1, \ldots, \mu_r, \gamma_1, \ldots, \gamma_r \in \set{1, \ldots, k-m}$ such that $\congmod{\mu_j}{\gamma_j}{k}$ for all $j \in \set{1, \ldots, r}$, $\Sigma^r_{j=1}\mu_j, \Sigma^r_{j=1}\gamma_j \leq k-m$ and $\congmod{\Sigma^r_{j=1}\mu_j}{k}{m}$, 
			\[\tomega{(\mu_1, \gamma_1),\ldots,(\mu_r,\gamma_r)}{k} = \begin{cases} (\tomega{(\mu_1, \gamma_1)}{k - \Sigma^r_{j=2}\mu_j} \oplus \id{\Sigma^r_{j=2}\mu_j})(\id{\gamma_1} \oplus \tupsilon{(\mu_2, \gamma_2), \ldots, (\mu_r, \gamma_r)}{k - \gamma_1}) & \mu_1 \geq \gamma_1; \\ (\id{\mu_1} \oplus \tupsilon{(\mu_2, \gamma_2), \ldots, (\mu_r, \gamma_r)}{k-\mu_1})(\tomega{(\mu_1, \gamma_1)}{k-\Sigma^r_{j=2}\gamma_j} \oplus \id{\Sigma^r_{j=2}\gamma_j}) & \mu_1 \leq \gamma_1. \end{cases}\]
			
			\begin{proof}
				Suppose $\mu_1 \geq \gamma_1$. Note that since $k \geq \Sigma^r_{j=1}\mu_j + m \geq \gamma_1 + \Sigma^r_{j=2}\mu_j + m$ we have:
				\begin{enumerate}
					\item $k - \Sigma^r_{j=2}\mu_j \geq \mu_1 + m = \max\set{\mu_1, \gamma_1} + m$, ensuring $\tomega{(\mu_1, \gamma_1)}{k - \Sigma^r_{j=2}\mu_j}$ is well-defined;
					\item $k-\gamma_1-m \geq \Sigma^r_{j=2}\mu_j, \Sigma^r_{j=2}\gamma_j$, ensuring $\tupsilon{(\mu_2, \gamma_2), \ldots, (\mu_r, \gamma_r)}{k - \gamma_1}$ is well-defined; and
					\item $\frac{k - \gamma_1 - \Sigma^r_{j=2}\mu_j}{m}$, which is both the number of lower $m$-apses in $\tomega{(\mu_1, \gamma_1)}{k - \Sigma^r_{j=2}\mu_j}$ and the number of upper $m$-apses in $\tupsilon{(\mu_2, \gamma_2), \ldots, (\mu_r, \gamma_r)}{k - \gamma_1}$, is a positive integer.
				\end{enumerate}
					
				Let $t^* = \frac{k - \Sigma^r_{j=1}\mu_j}{m}, t = \frac{k - \gamma_1 - \Sigma^r_{j=2}\mu_j}{m}, t_* = \frac{k - \Sigma^r_{j=1}\gamma_j}{m} \in \posints$. When forming the product $(\tomega{(\mu_1, \gamma_1)}{k - \Sigma^r_{j=2}\mu_j} \oplus \id{\Sigma^r_{j=2}\mu_j})(\id{\gamma_1} \oplus \tupsilon{(\mu_2, \gamma_2), \ldots, (\mu_r, \gamma_r)}{k - \gamma_1})$:
				\begin{enumerate}
					\item for each $j \in \set{1, \ldots, t^*}$, the upper $m$-apsis $\{\mu_1 + m(j-1) + 1, \ldots, \mu_1 + mj\}$ in $\tomega{(\mu_1, \gamma_1)}{k - \Sigma^r_{j=2}\mu_j}$ is preserved; 
					\item for each $j \in \set{1, \ldots, t_*}$, the lower $m$-apsis $\{(\gamma_1 + m(j-1) + 1)', \ldots, (\gamma_1 + mj)'\}$ in $\id{\gamma_1}\oplus\tupsilon{(\mu_2, \gamma_2), \ldots, (\mu_r, \gamma_r)}{k - \gamma_1}$ is preserved; 
					\item for each $j \in \set{1, \ldots, t}$, the lower $m$-apsis $\{(\gamma_1 + m(j-1) + 1)', \ldots, (\gamma_1 + mj)'\}$ in $\tomega{(\mu_1, \gamma_1)}{k - \Sigma^r_{j=2}\mu_j}$ and upper $m$-apsis $\{\gamma_1 + m(j-1) + 1, \ldots, \gamma_1 + mj\}$ in $\id{\gamma_1}\oplus\tupsilon{(\mu_2, \gamma_2), \ldots, (\mu_r, \gamma_r)}{k - \gamma_1}$ join and are removed; 
					\item the transversal $\set{1, \ldots, \mu_1, 1', \ldots, \gamma_1'}$ in $\tomega{(\mu_1, \gamma_1)}{k - \Sigma^r_{j=2}\mu_j}$ is preserved since it joins to the vertical lines in $\id{\gamma_1}$; and
					\item for each $j \in \set{2, \ldots, r}$, the transversal $\{\gamma_1 + mt + \Sigma^{j-1}_{i=2}\mu_i + 1, \ldots, \gamma_1 + mt + \Sigma^{j}_{i=2}\mu_j, (mt_* + \Sigma^{j-1}_{i=1}\gamma_j + 1)', \ldots, (mt_* + \Sigma^j_{i=1}\gamma_j)'\}$ in $\id{\gamma_1}\oplus\tupsilon{(\mu_2, \gamma_2), \ldots, (\mu_r, \gamma_r)}{k - \gamma_1}$ is preserved since it joins to the vertical lines $\set{l, l'}$ in $\tomega{(\mu_1, \gamma_1)}{k - \Sigma^r_{j=2}\mu_j} \oplus \id{\Sigma^r_{j=2}\mu_j}$ where $l \in \{\gamma_1 + mt + \Sigma^{j-1}_{i=2}\mu_i + 1, \ldots, \gamma_1 + mt + \Sigma^{j}_{i=2}\mu_j\}$.
				\end{enumerate}
				Hence $\tomega{(\mu_1, \gamma_1),\ldots,(\mu_r,\gamma_r)}{k} = (\tomega{(\mu_1, \gamma_1)}{k - \Sigma^r_{j=2}\mu_j} \oplus \id{\Sigma^r_{j=2}\mu_j})(\id{\gamma_1} \oplus \tupsilon{(\mu_2, \gamma_2), \ldots, (\mu_r, \gamma_r)}{k - \gamma_1})$. The case when $\mu_1 \leq \gamma_1$ follows analogously. 
			\end{proof}
		\end{proposition}
		
		\begin{figure}[!ht]
			\caption[ ]{Given $m=3$,}
			\label{fig:aptrans411114element}
			\vspace{10pt}
			\centering
			\input{chap_characterisations/tikz/fig-aptrans411114element.tex}
		\end{figure}
		
		\begin{figure}[!ht]
			\caption[ ]{Given $m=3$,}
			\label{fig:aptrans1771element}
			\vspace{10pt}
			\centering
			\input{chap_characterisations/tikz/fig-aptrans1771element.tex}
		\end{figure}

		\begin{proposition} \label{prop:aptransmujgammajinapsismon}
			For each $m \in \intsge{3}$, $k \in \intsge{m}$, $r \in \set{0, \ldots, k-m}$, and $\mu_1, \ldots, \mu_r, \gamma_1, \ldots, \gamma_r \in \set{1, \ldots, k-m}$ such that $\congmod{\mu_j}{\gamma_j}{k}$ for all $j \in \set{1, \ldots, r}$, $\Sigma^r_{j=1}\mu_j, \Sigma^r_{j=1}\gamma_j \leq k-m$ and $\congmod{\Sigma^r_{j=1}\mu_j}{k}{m}$, we have $\aptrans{(\mu_1, \gamma_1),\ldots,(\mu_r, \gamma_r)}{k} \subseteq \apsismon{m}{k}$.
			
			\begin{proof}
				If $r=0$ then we trivially have $k \in m\posints$ and $\aptrans{\nullset}{k} = \set{\tomega{(m, 0)}{k}} \subseteq \apsismon{m}{k}$, where $\tomega{(m, 0)}{k}$ was described in Definition \ref{def:tomegamugamma} and established as an element of the $m$-apsis monoid in Proposition \ref{prop:tomegamugammainapsismon}. If $r=1$ then we already established in Proposition \ref{prop:aptransmugammainapsismon} that $\aptrans{(\mu_1, \gamma_1)}{k} \subseteq \apsismon{m}{k}$. 
				
				If $r > 1$ then let $\overbar{r} \in \set{2, \ldots, r}$. Suppose $\mu_1 \geq \gamma_1$ and let $\overbar{k} = k-\Sigma^r_{j=\overbar{r}+1}\mu_j$. It  follows from $\aptrans{(\mu_1, \gamma_1)}{k-\Sigma^{r}_{j=2}\mu_j} \subseteq \apsismon{m}{k - \Sigma^{r}_{j=2}\mu_j}$, as was established in Proposition \ref{prop:aptransmugammainapsismon}, that $\tomega{(\mu_1, \gamma_1)}{k - \Sigma^r_{j=2}\mu_j} \oplus \id{\Sigma^{\overbar{r}}_{j=2}\mu_j} \in \apsismon{m}{\overbar{k}}$.
				
				Suppose $\aptrans{(\mu_2, \gamma_2),\ldots,(\mu_{\overbar{r}}, \gamma_{\overbar{r}})}{\overbar{k}-\gamma_1} \subseteq \apsismon{m}{\overbar{k}-\gamma_1}$ giving us $\id{\gamma_1} \oplus \tupsilon{(\mu_2, \gamma_2), \ldots, (\mu_{\overbar{r}}, \gamma_{\overbar{r}})}{\overbar{k} - \gamma_1} \in \apsismon{m}{\overbar{k}}$, then $\tomega{(\mu_1, \gamma_1), \ldots, (\mu_{\overbar{r}}, \gamma_{\overbar{r}})}{\overbar{k}} = (\tomega{(\mu_1, \gamma_1)}{\overbar{k} - \Sigma^{\overbar{r}}_{j=2}\mu_j} \oplus \id{\Sigma^{\overbar{r}}_{j=2}\mu_j})(\id{\gamma_1} \oplus \tupsilon{(\mu_2, \gamma_2), \ldots, (\mu_{\overbar{r}}, \gamma_{\overbar{r}})}{\overbar{k} - \gamma_1}) \in \apsismon{m}{\overbar{k}} \cap \aptrans{(\mu_1, \gamma_1), \ldots, (\mu_{\overbar{r}}, \gamma_{\overbar{r}})}{\overbar{k}}$. Therefore $\apsismon{m}{\overbar{k}} \cap \aptrans{(\mu_1, \gamma_1), \ldots, (\mu_{\overbar{r}}, \gamma_{\overbar{r}})}{\overbar{k}}$ is non-empty, employing Corollary \ref{cor:whenaptransinapsismon} we have $\aptrans{(\mu_1, \gamma_1), \ldots, (\mu_{\overbar{r}}, \gamma_{\overbar{r}})}{\overbar{k}} \subseteq \apsismon{m}{\overbar{k}}$. It follows by induction that $\aptrans{(\mu_1, \gamma_1), \ldots, (\mu_r, \gamma_r)}{k} \subseteq \apsismon{m}{k}$. 
				
				The case when $\mu_1 \leq \gamma_1$ follows analogously using the second half of Proposition \ref{prop:aptransfactorisation}.
			\end{proof}
		\end{proposition}

		\begin{proposition}
			For each $m \in \intsge{3}$ and $k \in \intsge{m}$, $\aptrans{m}{k} \subseteq \apsismon{m}{k}$.
			
			\begin{proof}
				Recall that by definition $\aptrans{m}{k} = \bigcup_{\substack{r \in \set{0, \ldots, k-m} \\ \mu_1, \ldots, \mu_r, \gamma_1, \ldots, \gamma_r \in \set{1, \ldots, k-m} \\ \congmod{\mu_j}{\gamma_j}{m}\ \forall\ j \in \set{1, \ldots, r} \\ \Sigma^r_{j=1}\mu_j, \Sigma^r_{j=1}\gamma_j \leq k-m \\ \congmod{\Sigma^r_{j=1}\mu_j}{k}{m}}} \aptrans{(\mu_1, \gamma_1), \ldots, (\mu_r, \gamma_r)}{k}$, hence it follows from Proposition \ref{prop:aptransmujgammajinapsismon} that $\aptrans{m}{k} \subseteq \apsismon{m}{k}$.
			\end{proof}
     \end{proposition}

	\subsection{Non-transversal building blocks}
    	We now turn our attention to bipartitions in $\apsisbound{m}{k}$ where transversals must be lines, and either lower non-transversals or upper non-transversals must be $m$-apses. 
    
		\begin{definition}
			For each $m \in \intsge{3}$ and $k \in \intsge{m}$, let $\apupnontrans{m}{k}$ denote the set of all $\alpha \in \apsisbound{m}{k}$ such that every:
			\begin{enumerate}
				\item transversal in $\alpha$ is a line; and
				\item lower non-transversal in $\alpha$ is an $m$-apsis, 
			\end{enumerate}
			furthermore let $\aplownontrans{m}{k}$ denote the set of all $\alpha \in \apsisbound{m}{k}$ such that every:
			\begin{enumerate}
				\item transversal in $\alpha$ is a line; and
				\item upper non-transversal in $\alpha$ is an $m$-apsis.
			\end{enumerate}
			We shall refer to elements of $\apupnontrans{m}{k}$ as \textit{upper non-transversal building blocks}, and elements of $\aplownontrans{m}{k}$ as \textit{lower non-transversal building blocks}. By a \textit{non-transversal building block} we shall mean an upper or lower non-transversal building block.
		\end{definition}
     
     	We proceed in this subsection to establish that both $\apupnontrans{m}{k}$ and $\aplownontrans{m}{k}$ are contained within the $m$-apsis monoid $\apsismon{m}{k}$. Containment of $\aplownontrans{m}{k}$ will trivially follow from containment of $\apupnontrans{m}{k}$. In order to establish the latter containment, first we partition $\apupnontrans{m}{k}$ based on the number of upper non-transversals that are not $m$-apses, then we establish our desired result inductively.
       	
       	Note that $\floor{\frac{k}{m}}$ is the maximum number of upper non-transversals an element of $\apupnontrans{m}{k}$ may contain, however in such a case at least one upper non-transversal must be an $m$-apsis, hence $\apupnontrans{m}{k}$ may contain at most $\floor{\frac{k-m}{m}}$ upper non-transversals that are not $m$-apses.
       	
       	\begin{definition}
       		For each $x \in \set{0, \ldots, \floor{\frac{k-m}{m}}}$, let $\apupnontrans{}{x}$ denote the subset of all bipartitions $\eta \in \apupnontrans{m}{k}$ such that $\eta$ contains precisely $x$ upper non-transversals that are not $m$-apses. 
       	\end{definition}
       		
       	Note it follows by definition that $\apupnontrans{}{0} = \apmorphs{m}{k}$ and that $\apupnontrans{m}{k} = \bigcup^{\floor{\frac{k-m}{m}}}_{x=0}\apupnontrans{}{x}$. 
       	
       	\begin{definition} \label{def:eta-andomegaeta}
					For each $x \in \set{1, \ldots, \floor{\frac{k-m}{m}}}$ and $\eta \in \apupnontrans{}{x}$, by definition there exists:
       			\begin{enumerate}
       				\item $\mu \in m\posints$ and $b = \set{b_1, \ldots, b_{\mu}} \in \eta$, where $b_1, \ldots, b_{\mu} \in \set{1, \ldots, k}$, such that $\mu \leq k-m$, $b$ is a type $(\mu, 0)$ non-transversal that is not an $m$-apsis and no upper non-transversals of $\eta$ pass underneath $b$;
       				\item $r_1 \in \set{0, \ldots, b_1-1}$ and $r_2 \in \set{0, \ldots, k-b_{\mu}}$ such that the number of transversal lines in $\eta$ whose upper vertex sits to the left or right of $b \in \eta$ is $r_1$ and $r_2$ respectively; and
       				\item upper vertices $u_1, \ldots, u_{r_1} \in \set{1, \ldots, b_1-1}$, $u_{r_1+1}, \ldots, u_{r_1+r_2} \in \{b_{\mu}+1,$ $\ldots, k\}$ and lower vertices $l_1, \ldots, l_{r_1+r_2} \in \set{1, \ldots, k}$ such that for each $j \in \{1, \ldots, r_1+r_2\}$, $\eta$ contains the transversal line $\set{u_j, l_j}$.
       			\end{enumerate}
       			
       			\begin{itemize}
       				\item Let $\omega_{\eta}$ denote $ \id{r_1}\oplus\tomega{(\mu,0)}{k-r_1-r_2}\oplus\id{r_2}$ where $\tomega{(\mu,0)}{k-r_1-r_2}$ was outlined in Definition \ref{def:tomegamugamma} and shown to be an element of the $m$-apsis monoid $\apsismon{m}{k-r_1-r_2}$ in Proposition \ref{prop:tomegamugammainapsismon}; and
       			
	         		\item let $\eta^-$  denote the upper non-transversal building block in $\apupnontrans{}{x-1}$ containing:
		         	\begin{enumerate}
		         		\item excluding $b$, the upper non-transversals of $\eta$;
		         		\item for each $j \in \set{1, \ldots, \mu}$, the transversal line $\set{b_j, r_1+j}$;
		         		\item for each $j \in \set{1, \ldots, r_1}$, the transversal line $\set{u_j, j'}$;
		         		\item for each $j \in \set{1, \ldots, r_2}$, the transversal line $\set{u_j, (k-r_2+j)'}$; and
		         		\item for each $j \in \set{1, \ldots, \frac{k-\mu-r_1-r_2}{m}}$, the lower $m$-apsis $\{(r_1+\mu+m(j-1)+1)', \ldots, (r_1+\mu+mj)'\}$.
		         	\end{enumerate}
        \end{itemize}
       	\end{definition}

		Figure \ref{fig:apupnontranseg} contains an example of $\omega_{\eta}$ and $\eta^-$ for $m=3$, $k=13$ and a given $\eta \in \apupnontrans{}{1}$.

		\begin{figure}[!ht]
			\caption[ ]{Let $m=3$, $k=13$.}
			\label{fig:apupnontranseg}
			\vspace{10pt}
			\centering
			\input{chap_characterisations/tikz/fig-apupnontranseg.tex}
		\end{figure}
       	
       	We proceed to establish that each $\apupnontrans{}{x}$ is contained within the $m$-apsis monoid $\apsismon{m}{k}$. To do so we will show that each $\eta \in \apupnontrans{}{x}$, where $x \in \set{1, \ldots, \floor{\frac{k-m}{m}}}$, may be factorised into the product $\eta^-\omega_{\eta}\apmorph{\omega_{\eta}^*}{\eta}$, an example of which is depicted in Figure \ref{fig:apupnontransfactorisation}.
       	
       	\begin{proposition} \label{prop:apupnontransinapsismon}
       		For each $m \in \intsge{3}$, $k \in \intsge{m}$ and $x \in \set{0, \ldots, \floor{\frac{k-m}{m}}}$, $\apupnontrans{}{x} \subseteq \apsismon{m}{k}$.
       		
       		\begin{proof}
					By definition $\apupnontrans{}{0} = \apmorphs{m}{k}$ and we established in Proposition \ref{prop:apmorphismsinapsismon} that $\apmorphs{m}{k} \subseteq \apsismon{m}{k}$, hence $\apupnontrans{}{0} \subseteq \apsismon{m}{k}$.
       			 
        			Let $x \in \set{1, \ldots, \floor{\frac{k-m}{m}}}$ and suppose $\apupnontrans{}{x-1} \subseteq \apsismon{m}{k}$. For each $\eta \in \apupnontrans{}{x}$, as previously noted in Definition \ref{def:eta-andomegaeta}, by definition there exists:
       			\begin{enumerate}
       				\item $\mu \in m\posints$ and $b = \set{b_1, \ldots, b_{\mu}} \in \eta$, where $b_1, \ldots, b_{\mu} \in \set{1, \ldots, k}$, such that $\mu \leq k-m$, $b$ is a type $(\mu, 0)$ non-transversal that is not an $m$-apsis and no upper non-transversals of $\eta$ pass underneath $b$;
       				\item $r_1 \in \set{0, \ldots, b_1-1}$ and $r_2 \in \set{0, \ldots, k-b_{\mu}}$ such that the number of transversal lines in $\eta$ whose upper vertex sits to the left or right of $b \in \eta$ is $r_1$ and $r_2$ respectively; and
       				\item upper vertices $u_1, \ldots, u_{r_1} \in \set{1, \ldots, b_1-1}$, $u_{r_1+1}, \ldots, u_{r_1+r_2} \in \{b_{\mu}+1,$ $\ldots, k\}$ and lower vertices $l_1, \ldots, l_{r_1+r_2} \in \set{1, \ldots, k}$ such that for each $j \in \{1, \ldots, r_1+r_2\}$, $\eta$ contains the transversal line $\set{u_j, l_j}$.
       			\end{enumerate}
        			        			
				Note $\omega_{\eta}$ and $\eta$ both have $\frac{k-\mu-r_1-r_2}{m}$ lower $m$-apses, hence the $m$-apmorphism $\apmorph{\omega_{\eta}^*}{\eta}$ is well-defined. When forming the product $\eta^-\omega_{\eta}\apmorph{\omega_{\eta}^*}{\eta}$:
       			\begin{enumerate}
        			\item the upper $m$-apses in $\eta^-$ are preserved, which are identical to the $m$-apses in $\eta$ after excluding $b$;
        			\item $b$ is formed by the $\mu$-apsis $\set{r_1 + 1, \ldots, r_1+\mu} \in \omega_{\eta}$ joining to each of $\eta^-$'s transversal lines $\set{b_j, r_1+j}$ where $j \in \set{1, \ldots, \mu}$;
        			\item the lower $m$-apses in $\apmorph{\omega_{\eta}^*}{\eta}$, which are identical to the lower $m$-apses in $\eta$, are preserved;
        			\item for each $j \in \set{1, \ldots, \frac{k-\mu-r_1-r_2}{m}}$, the lower $m$-apsis $\{r_1 + \mu + m(j-1)+1, \ldots, r_1+\mu+mj\} \in \eta^-$ and the upper $m$-apsis $\{(r_1 + \mu + m(j-1)+1)', \ldots, (r_1+\mu+mj)\} \in \omega_{\eta}$ join and are removed;
        			\item for each $j \in \set{1, \ldots, \frac{k-r_1-r_2}{m}}$, the lower $m$-apsis $\{r_1 + m(j-1)+1, \ldots, r_1+mj\} \in \omega_{\eta}$ and upper $m$-apsis $\{(r_1 + m(j-1)+1)', \ldots, (r_1+mj)\} \in \apmorph{\omega_{\eta}^*}{\eta}$ join and are removed;
        			\item for each $j \in \set{1, \ldots, r_1}$, the transversal line $\set{u_j, l_j'}$ is formed by the lines $\set{u_j, j'} \in \eta^-$, $\set{j, j'} \in \omega_{\eta}$ and $\set{j, l_j'} \in \apmorph{\omega_{\eta}^*}{\eta}$ joining; and
        			\item for each $j \in \set{k-r_2+1, \ldots, k}$, the line $\set{u_j, l_j'}$ is formed by the lines $\set{u_j, (k-r_2+j)'} \in \eta^-$, $\set{k-r_2+j, (k-r_2+j)'} \in \omega_{\eta}$ and $\set{k-r_2+j, l_j'} \in \apmorph{\omega_{\eta}^*}{\eta}$ joining.
       			\end{enumerate}
       			
       			Hence $\eta = \eta^-\omega_{\eta}\apmorph{\omega_{\eta}^*}{\eta} \in \apsismon{m}{k}$. It follows by induction that $\apupnontrans{}{x} \subseteq \apsismon{m}{k}$ for all $x \in \set{0, \ldots, \floor{\frac{k-m}{m}}}$.
       		\end{proof}
       	\end{proposition}
		
		\begin{figure}[!ht]
			\caption[ ]{Given $m=3$, $k=13$ and $\eta \in \apupnontrans{}{1}$ from Figure \ref{fig:apupnontranseg},}
			\label{fig:apupnontransfactorisation}
			\vspace{10pt}
			\centering
			\input{chap_characterisations/tikz/fig-apupnontransfactorisation.tex}
		\end{figure}
		
		\begin{proposition}
			For each $m \in \intsge{3}$ and $k \in \intsge{m}$, $\apupnontrans{m}{k}, \aplownontrans{m}{k} \subseteq \apsismon{m}{k}$.
			
			\begin{proof}
				$\apupnontrans{m}{k} = \bigcup_{x=0}^{\floor{\frac{k-m}{m}}}\apupnontrans{}{x} \subseteq \apsismon{m}{k}$. Containment of $\aplownontrans{m}{k}$ may be established in a dual fashion to $\apupnontrans{m}{k}$, or more succinctly by noting that $\aplownontrans{m}{k} = \tuple{\apupnontrans{m}{k}}^* \subseteq \apsismon{m}{k}$.
			\end{proof}
		\end{proposition}

   	\subsection{Factorising elements of the $m$-apsis monoid $\apsismon{m}{k}$}
    	To establish that the $m$-apsis monoid $\apsismon{m}{k}$ is in fact the monoid $\apsisbound{m}{k}$, which will bring our characterisation of the $m$-apsis monoid $\apsismon{m}{k}$ to a conclusion, it remains for us to establish that each element of $\apsisbound{m}{k}$ is a product of $m$-apsis generators. 
		
		We will reach our desired conclusion by establishing that each element of $\apsisbound{m}{k}$ may be factorised into a product containing an upper non-transversal building block, a transversal building block and a lower non-transversal building block.
		
		\begin{definition} \label{def:transandnontransblockelmnts}
			For each $m \in \intsge{3}$, $k \in \intsge{m}$ and $\alpha \in \apsisbound{m}{k}$:
			\begin{itemize}
				\item let $u_{\alpha}$ denote the upper non-transversal building block containing:
				\begin{enumerate}
					\item the upper non-transversals in $\alpha$;
					\item the vertical line $\set{j, j'}$ for all $j \in \set{1, \ldots, k}$ such that $j$ is an element of a transversal in $\alpha$;
					\item lower $m$-apses containing the remaining vertices,
				\end{enumerate}
				\item let $t_{\alpha}$ denote the transversal building block containing:
				\begin{enumerate}
					\item the transversals contained in $\alpha$;
					\item $m$-apses replacing the non-transversals in $\alpha$, and
				\end{enumerate}
				\item let $l_{\alpha}$ denote the lower non-transversal building block containing:
				\begin{enumerate}
					\item the lower non-transversals contained in $\alpha$;
					\item for each $j \in \set{1, \ldots, k}$ such that $j'$ is an element of a transversal in $\alpha$, the vertical line $\set{j, j'}$;
					\item upper $m$-apses containing the remaining vertices.
				\end{enumerate}
			\end{itemize}
		\end{definition}
		
\begin{figure}[!ht]
			\caption[ ]{Let $m=3$ and $k =16$.}
			\label{fig:apsismonintoblocks}
			\vspace{5pt}
			\centering
			\input{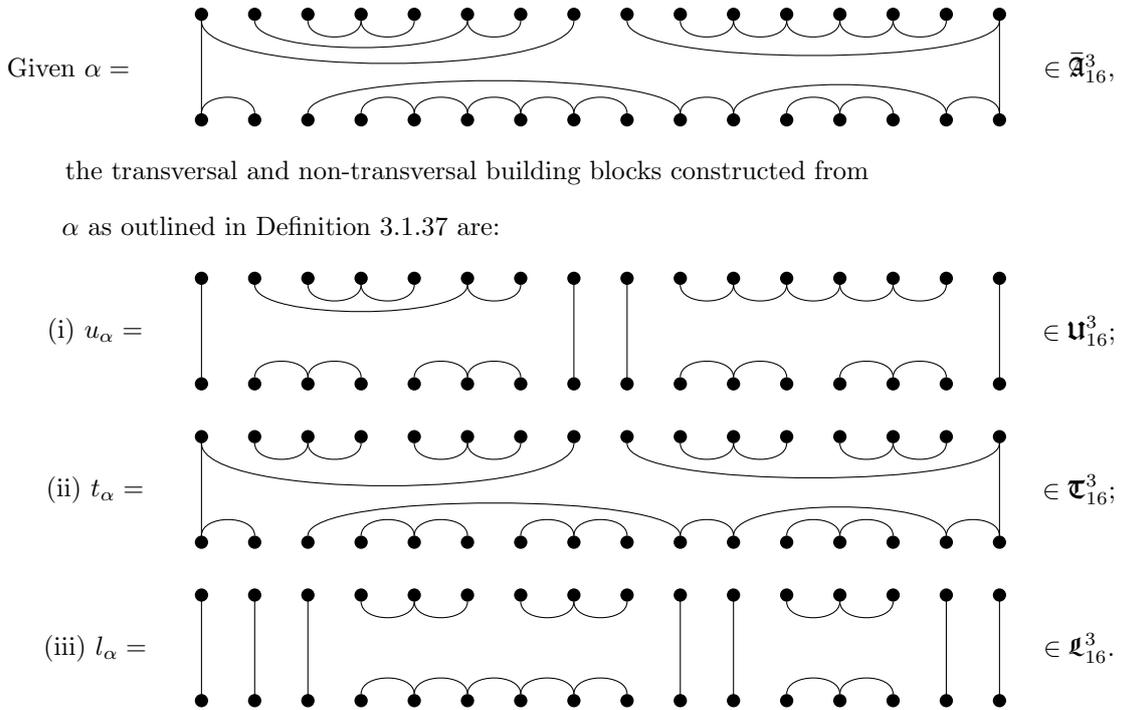}
		\end{figure}

		Figure \ref{fig:apsismonintoblocks} contains an example of $u_{\alpha}$, $t_{\alpha}$ and $l_{\alpha}$ for a given $\alpha \in \apsisbound{3}{16}$.
		
		\begin{theorem} \label{thm:apsismon=apsisbound}
			For each $m \in \intsge{3}$ and $k \in \intsge{m}$, the $m$-apsis monoid $\apsismon{m}{k}$ is equal to the monoid $\apsisbound{m}{k}$.
			
			\begin{proof}
				We already established in Proposition \ref{prop:apsismonpropersubmonapsisbound} that the $m$-apsis monoid $\apsismon{m}{k}$ is contained within $\apsisbound{m}{k}$. Conversely for each $\alpha \in \apsisbound{m}{k}$, when forming the product $u_{\alpha}t_{\alpha}l_{\alpha}$ (see Figure \ref{fig:apsismonfactorisation} for an example):
				\begin{enumerate}
					\item the upper non-transversals in $u_{\alpha}$, which are identical to the upper non-transversals in $\alpha$, are preserved;
					\item the transversals in $t_{\alpha}$, which are identical to the transversals in $\alpha$, are preserved as they join to vertical lines in both $u_{\alpha}$ and $l_{\alpha}$ by construction;
					\item the lower non-transversals in $l_{\alpha}$, which are identical to the lower non-transversals in $\alpha$, are preserved; and
					\item the lower $m$-apses in $u_{\alpha}$ and upper $m$-apses in $t_{\alpha}$ join and are removed, similarly with the lower $m$-apses in $t_{\alpha}$ and upper $m$-apses in $l_{\alpha}$.
				\end{enumerate}
				
				Hence $\alpha = u_{\alpha}t_{\alpha}l_{\alpha}$, giving us $\apsisbound{m}{k} \subseteq \apupnontrans{m}{k}\aptrans{m}{k}\aplownontrans{m}{k} \subseteq \apsismon{m}{k}$.
			\end{proof}
		\end{theorem}
		
		\begin{figure}[!ht]
			\caption[ ]{Given $m=3$, $k=16$ and $\alpha \in \apsisbound{3}{16}$ from Figure \ref{fig:apsismonintoblocks},}
			\label{fig:apsismonfactorisation}
			\vspace{5pt}
			\centering
			\input{chap_characterisations/tikz/fig-apsismonfactorisation.tex}
		\end{figure} 

	\section{The crossed $m$-apsis generated diagram monoid $\crossedmon{\apsismon{m}{k}}$} 
	\begin{definition} \label{def:capsismon}
		We shall refer to the join of the $m$-apsis monoid $\apsismon{m}{k}$ and the symmetric group $\symgrp{k}$ as \textit{the crossed $m$-apsis generated diagram monoid}, or more succinctly as \textit{the crossed $m$-apsis monoid}, and shall denote it as $\capsismon{m}{k}$.
	\end{definition}

	\begin{proposition}
		The crossed $m$-apsis monoid $\capsismon{m}{k}$ is equal to the set of all bipartitions $\alpha \in \modmon{m}{k}$ such that either $\alpha \in \symgrp{k}$ or $\alpha$ contains at least one upper and at least one lower non-transversal that each contain precisely $m$ vertices.
		
		\begin{proof}
			Let $\capsisbound{m}{k}$ denote the set of all bipartitions $\alpha \in \modmon{m}{k}$ such that either $\alpha \in \symgrp{k}$ or $\alpha$ contains at least one upper and at least one lower non-transversal that each contain precisely $m$ vertices.
			
			It follows from Proposition \ref{prop:symgrpedgetypes} that $\capsisbound{m}{k} = \symgrp{k}\apsismon{m}{k}\symgrp{k}$. That $\capsisbound{m}{k}$ is a monoid follows trivially from upper non-transversals being preserved when multiplying on the right, and from lower non-transversals dually being preserved when multiplying on the left.
		\end{proof}
	\end{proposition}
	
	Note that when $m \geq 3$ and $k \geq 3m$, taking the meet of the crossed $m$-apsis monoid and the planar partition monoid does not get us back to the $m$-apsis-monoid (see Figure \ref{fig:planarbutnotinapsismon} for an example of a planar element of $\capsismon{3}{9}$ that is not an element of $\apsismon{3}{9}$). Within this thesis, the $m$-apsis monoid $\apsismon{m}{k}$ is the only submonoid of the planar partition monoid that is not equal to the meet of the planar partition monoid with the join of the symmetric group and itself.
	
	\begin{figure}[!ht]
		\caption[ ]{A planar element of $\capsismon{3}{9}$ that is not an element of $\apsismon{3}{9}$.}
		\label{fig:planarbutnotinapsismon}
		\vspace{5pt}
		\centering
		\input{chap_characterisations/tikz/fig-planarbutnotinapsismon.tex}
	\end{figure}

	\section{The planar modular partition monoid $\pmodmon{m}{k}$} \label{sec:pmodmon}
	Recall from Definition \ref{def:pmodmon} that the planar mod-$m$ monoid $\pmodmon{m}{k}$ is the set of all planar bipartitions $\alpha \in \ppttnmon{k}$ such that for each block $b \in \alpha$, $\congmod{\noupverts{b}}{\nolowverts{b}}{m}$. In this section we establish generators for the planar mod-$m$ monoid $\pmodmon{m}{k}$, our objective will be obtained with the following approach:
	\begin{enumerate}
		\item First we consider the bipartitions that sit inside the planar mod-$m$ monoid $\pmodmon{m}{k}$ that are not products of $m$-apsis generators, that is consider the bipartitions $\pmodmon{m}{k}-\apsismon{m}{k}$, along with known generators of the planar partition monoid $\ppttnmon{k}$ which is equal to the planar mod-$1$ monoid $\pmodmon{1}{k}$ by definition;
		\item Second we bound the planar mod-$m$ monoid $\pmodbound{m}{k}$ below by conjecturing a sufficient generating set based on our characterisation of the $m$-apsis monoid $\apsismon{m}{k}$; and
		\item Finally we establish equivalence between the planar mod-$m$ monoid $\pmodmon{m}{k}$ and our lower bound $\pmodbound{m}{k}$ in an analogous fashion to how we established equivalence between the $m$-apsis monoid $\apsismon{m}{k}$ and our upper bound $\apsisbound{m}{k}$.
	\end{enumerate} 
	
	Note that this was the approach originally taken by the author to identify a generating set for the planar mod-$m$ monoid $\pmodmon{m}{k}$. If he had been familiar with the work of Kosuda \cite{art:Kosuda:StructurePartyAlgebraTypeB, art:Kosuda:StdExpForPartyAlgebra, art:Kosuda:CharacterizationModularPartyAlgebra} one would have likely guessed the same generating set, though the bulk of our characterisation would have remained original. 

    \subsection{Bounding $\pmodmon{m}{k}$ below by $\pmodbound{m}{k}$}
    	Recall that when $m \in \intsge{3}$, the $m$-apsis monoid $\apsismon{m}{k}$ consists of all bipartitions $\alpha \in \pmodmon{m}{k}$ such that $\alpha$ contains at least one upper $m$-apsis and at least one lower $m$-apsis. When generating the planar mod-$m$ monoid $\pmodmon{m}{k}$, at each step we need to lose the restriction that the bipartitions we may genearate must contain at least one upper $m$-apsis and at least one lower $m$-apsis.
    	
    	Further recall that our first step when characterising the $m$-apsis monoid $\apsismon{m}{k}$ was to, for each $m \in \intsge{3}$, $k \in \intsge{m}$ and $\mu \in \{1, \ldots, k-m\}$, establish a well-defined product of $m$-apsis generators containing a block of type $(\mu, \mu)$. The analogous step for characterising the planar mod-$m$ monoid $\pmodmon{m}{k}$ will require for each $m \in \posints$, $k \in \intsge{m}$ and $\mu \in \{1, \ldots, k\}$, to establish a well-defined product of generators containing a block of type $(\mu, \mu)$. 
    	
    	The reader may be able to recall that:
    	\begin{enumerate}
	    	\item the monoid of planar uniform block bijections $\uniblockbijmon{k}$, which is generated by $(2,2)$-transapsis generators, consists of all planar diagrams $\alpha \in \ppttnmon{k}$ such that each of $\alpha$'s blocks is uniform, that is $\noupverts{b} = \nolowverts{b}$ for all $b \in \alpha$; and
	    	\item Halverson and Ram \cite{art:Halverson:PartitionAlgebras} gave a presentation of planar partition monoid $\ppttnmon{k}$, which is by definition the planar mod-$1$ monoid $\pmodmon{1}{k}$, showing that it was generated by $(2,2)$-transapsis generators and monapsis generators.
    	\end{enumerate}
    	
    	Furthermore, it may trivially be checked that $\pmodmon{2}{2} - \apsismon{2}{2} = \{\transapgen{1}\}$. Based on the information we have now established, it is reasonable to conjecture that the planar mod-$m$ monoid may be generated by the $(2,2)$-transapsis generators, the $m$-apsis generators and the identity bipartition.
    	
		\begin{definition} \label{def:pmodbound}
			Given $m \in \intsge{0}$ and $k \in \ints_{\geq m}$, we denote by $\pmodbound{m}{k}$ the monoid generated by:
			\begin{enumerate}
				\item the $(2,2)$-transapsis generators $\set{\transapgen{j}: j \in \set{1, \ldots, k-1}}$;
				\item the $m$-apsis generators $\set{\apgen{m}{j}: j \in \set{1, \ldots, k-m+1}}$; and
				\item the identity bipartition $\id{k}$.
			\end{enumerate}
		\end{definition}
		
		The remainder of Section \ref{sec:pmodmon} will be spent establishing that the planar mod-$m$ monoid $\pmodmon{m}{k}$ is generated by the $(2,2)$-transapsis generators, the $m$-apsis generators and the identity bipartition, that is $\pmodmon{m}{k} = \pmodbound{m}{k}$.
    
    \subsection{Monapmorphisms}
    	Recall from Subsection \ref{def:m-apmorphism} that the $m$-apmorphisms $\apmorphs{m}{k}$ are the set of all planar bipartitions $\theta \in \ppttnmon{k}$ such that each block $b \in \theta$ is either an $m$-apsis or a transversal line.

		We established in Proposition \ref{prop:apmorphismsinapsismon} that for each $m \in \intsge{2}$ and $k \in \intsge{m}$, the $m$-apmorphisms $\apmorphs{m}{k}$ are contained in the $m$-apsis monoid $\apsismon{m}{k}$, from which it trivially follows that the $m$-apmorphisms $\apmorphs{m}{k}$ are also contained in the planar mod-$m$ monoid $\pmodmon{m}{k}$.
		
		When $m=1$, recall that the monapsis monoid $\uniapmon{k}$ is trivially isomorphic to the join-semilattice of subsets of $\set{1, \ldots, k}$ under $\monapgen{i} \mapsto \set{i}$ for all $i \in \set{1, \ldots, k}$. It follows by definition that $\apmorphs{1}{k} = \psyminvmon{k}$, which for each $m \in \intsge{2}$ does not contain the monapmorphisms $\apmorphs{1}{k} = \psyminvmon{k}$, that is $\apmorphs{1}{k} = \psyminvmon{k} \nsubseteq \apsismon{m}{k}$.
		
		To establish that $\apmorphs{1}{k} = \psyminvmon{k} \subseteq \pmodbound{m}{k}$ we will establish that each generator of $\psyminvmon{k}$ may be formed as a product of generators for $\pmodbound{m}{k}$. The reader may find it useful to refer to Figure \ref{fig:lrtransapgenfrommoanpgensandtransapgens} which illustrates that for each $k \in \intsge{2}$ and $i \in \{1, \ldots, k-1\}$, $\ftransgen{i} = \ftransgen{i+1}\transapgen{i}\monapgen{i}$ and $\btransgen{i} = \monapgen{i}\transapgen{i}\btransgen{i+1}$.
		
		\begin{proposition} \label{prop:lrtransapgenfrommoanpgensandtransapgens}
			For each $k \in \intsge{2}$ and $i \in \{1, \ldots, k-1\}$:
			\begin{enumerate}
				\item $\ftransgen{i} = \begin{cases}\monapgen{i+1} & i = k; \\ \ftransgen{i+1}\transapgen{i}\monapgen{i} = \monapgen{k}\left(\prod_{j=k-1}^{i}\monapgen{i}\transapgen{i}\right) & i < k, \text{ and}\end{cases}$
				\item $\btransgen{i} = \begin{cases}\monapgen{i} & i = k; \\ \monapgen{i}\transapgen{i}\btransgen{i+1} = \left(\prod_{j=i}^{k-1}\monapgen{i}\transapgen{i}\right)\monapgen{k} & i < k.\end{cases}$
			\end{enumerate}
			
			\begin{proof}
				When forming the product $\ftransgen{i+1}\transapgen{i}\monapgen{i}$:
				\begin{enumerate}
					\item for each $j \in \{1, \ldots, i-1\}$, $\monapgen{i+1}$, $\transapgen{i}$ and $\monapgen{i}$ all contain the vertical line $\{j, j'\}$, which trivially join and form the vertical line $\{j, j'\} \in \ftransgen{i+1}\transapgen{i}\monapgen{i}$;
					\item the upper monapsis $\{k\} \in \ftransgen{i+1}$ and lower monapsis $\{i\} \in \monapgen{i}$ are preserved; 
					\item the vertical line $\{i, i'\} \in \ftransgen{i+1}$ joins to the $(2,2)$-transapsis $\{i, i+1, i', (i+1)'\} \in \transapgen{i}$ which joins to the transversal line $\{i+1, (i+1)'\} \in \monapgen{i}$, forming the transversal line $\{i, (i+1)'\} \in \ftransgen{i+1}\transapgen{i}\monapgen{i}$; and
					\item for each $j \in \{i+2, \ldots, k\}$, the line $\set{j-1, j'} \in \ftransgen{i+1}$ joins to the line $\set{j, j'} \in \transapgen{i}$, which joins to the vertical line $\set{j, j'} \in \monapgen{i}$, forming the line $\{j-1, j'\} \in \ftransgen{i+1}\transapgen{i}\monapgen{i}$.
				\end{enumerate}
				Hence $\ftransgen{i+1}\transapgen{i}\monapgen{i} = \ftransgen{i}$. It follows by induction that $\ftransgen{i+1}\transapgen{i}\monapgen{i} = \monapgen{k}\left(\prod_{j=k-1}^{i}\monapgen{i}\transapgen{i}\right)$. That $\btransgen{i} = \monapgen{i}\transapgen{i}\monapgen{i+1} = \left(\prod_{j=i}^{k-1}\monapgen{i}\transapgen{i}\right)\monapgen{k}$ follows dually.
			\end{proof}
		\end{proposition}
		
		\begin{figure}[!ht]
			\caption[ ]{For each $k \in \intsge{2}$ and $i \in \{1, \ldots, k-1\}$,}
			\label{fig:lrtransapgenfrommoanpgensandtransapgens}
			\vspace{5pt}
			\centering
			\input{chap_characterisations/tikz/fig-lrtransapgenfrommoanpgensandtransapgens.tex}
		\end{figure} 
		
		\begin{corollary} \label{cor:monapmorphsinpmodbound}
			For each $k \in \posints$, $\apmorphs{1}{k} = \psyminvmon{k} \subseteq \pmodbound{1}{k}$.
			
			\begin{proof}
				When $k=1$, we trivially have $\apmorphs{1}{1} = \psyminvmon{1} =
				\{\id{1}, \monapgen{1}\} = \pmodbound{1}{1}$. When $k >1$, $\apmorphs{1}{k} = \psyminvmon{k} \subseteq \pmodbound{1}{k}$ trivially follows from Proposition \ref{prop:lrtransapgenfrommoanpgensandtransapgens} when recalling that:
				\begin{enumerate}
					\item $\psyminvmon{k}$ is generated by $\Big\{\id{k}, \ftransgen{i}, \btransgen{i} : i \in \{1, \ldots, k-1\}\Big\}$; and
					\item $\pmodbound{1}{k}$ is generated by $\Big\{\id{k}, \transapgen{i}, \monapgen{i}, \monapgen{k} : i \in \{1, \ldots, k-1\}\Big\}$.
				\end{enumerate}
			\end{proof}
		\end{corollary}
    
    \subsection{Generating feasible block types}
		\begin{definition}
			For each $m \in \posints$, $k \in \ints_{\geq m}$ and $\mu, \gamma \in \set{0, \ldots, k}$ such that $\mu + \gamma > 0$ and $\congmod{\mu}{\gamma}{m}$, let $\overbar{k} = \max\set{\mu, \gamma}$ and let $\tomegabar{(\mu, \gamma)}{k}$ denote the bipartition in the planar mod-$m$ monoid $\pmodmon{m}{k}$ containing: 
			\begin{enumerate}
				\item the type $(\mu, \gamma)$ transversal $\set{1, \ldots, \mu, 1', \ldots, \gamma'}$;
				\item for each $j \in \set{1, \ldots, \frac{\overbar{k} - \mu}{m}}$, the upper $m$-apsis $\{\mu+m(j-1) + 1, \ldots, \mu+mj\}$;
				\item for each $j \in \set{1, \ldots, \frac{\overbar{k} - \gamma}{m}}$, the lower $m$-apsis $\{(\gamma+m(j-1)+1)', \ldots, (\gamma+mj)'\}$; and
				\item for each $j \in \set{\overbar{k}+1, \ldots, k}$, the vertical line $\set{j, j'}$.
			\end{enumerate} 
		\end{definition}
		
		For example, given $m=2$, $\tomegabar{(2,6)}{8} \in \pmodmon{2}{8}$ is depicted in Figure \ref{fig:tomegabar26pmodmon}.
		
		\begin{figure}[!ht]
			\caption[ ]{Given $m=2$ and $k=8$,}
			\label{fig:tomegabar26pmodmon}
			\vspace{5pt}
			\centering
			\input{chap_characterisations/tikz/fig-tomegabar26pmodmon.tex}
		\end{figure} 
		
\begin{figure}[!ht]
			\caption[ ]{Given $m=3$ and $k=10$,}
			\label{fig:tomegabar82=tomegabar88sigma357357}
			\vspace{5pt}
			\centering
			\input{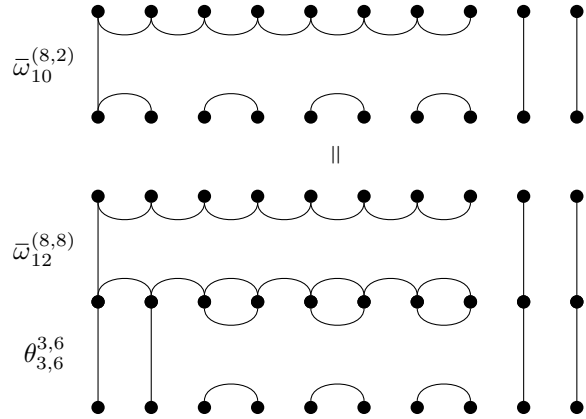}
		\end{figure}

		We proceed by establishing in Proposition \ref{prop:tomegabarmugammainpmodbound} that $\tomegabar{(\mu, \gamma)}{k} \in \pmodbound{m}{k}$. The reader may find it useful to refer to Figure \ref{fig:tomegabar82=tomegabar88sigma357357} which illustrates that, given $m=2$ and $k=10$, $\tomegabar{(8,2)}{10} = \tomegabar{(8,8)}{10}\apmorph{3,5,7}{3,5,7}$.
		
		\begin{proposition} \label{prop:tomegabarmugammainpmodbound}
			For each $m \in \posints$, $k \in \intsge{m}$ and $\mu, \gamma \in \set{0, \ldots, k}$ such that $\mu + \gamma > 0$ and $\congmod{\mu}{\gamma}{m}$,
			
			$\tomegabar{(\mu, \gamma)}{k} = \begin{cases}\tomegabar{(\mu, \mu)}{k}\apmorph{\gamma+1, \gamma+m+1, \ldots, \mu-m+1}{\gamma+1, \gamma+m+1, \ldots, \mu-m+1} & \text{ if } \mu > \gamma; \\ \prod^{\mu-1}_{i=1}\transapgen{i} & \text{ if } \mu = \gamma; \text{ and} \\\apmorph{\mu+1, \mu+m+1, \ldots, \gamma-m+1}{\mu+1, \mu+m+1, \ldots, \gamma-m+1}\tomegabar{(\gamma, \gamma)}{k} & \text{ if } \mu < \gamma, \\ \end{cases}$
			
			and hence $\tomegabar{(\mu, \gamma)}{k} \in \pmodbound{m}{k}$.
			
			\begin{proof}
				If $\mu = \gamma$ then we trivially have $\tomegabar{(\mu, \gamma)}{k} = \prod^{\mu-1}_{i=1} \transapgen{i} \in \pmodbound{m}{k}$. 
				
				If $\mu > \gamma$ then when forming the product $\tomegabar{(\mu, \mu)}{k}\apmorph{\gamma+1, \gamma+m+1, \ldots, \mu-m+1}{\gamma+1, \gamma+m+1, \ldots, \mu-m+1}$:
				\begin{enumerate}
					\item the $(\mu, \gamma)$-transversal $\set{1, \ldots, \mu, 1', \ldots, \mu'}$ in $\tomegabar{(\mu,\mu)}{k}$ connects to the $m$-apses $\{\gamma+(j-1)m+1, \ldots,$ $\gamma+jm\}$ in $\apmorph{\gamma+1, \gamma+m+1, \ldots, \mu-m+1}{\gamma+1, \gamma+m+1, \ldots, \mu-m+1}$ where $j \in \set{1, \ldots, \frac{\mu-\gamma}{m}}$ along with the vertical lines $\set{j, j'}$ in $\apmorph{\gamma+1, \gamma+m+1, \ldots, \mu-m+1}{\gamma+1, \gamma+m+1, \ldots, \mu-m+1}$ where $j \in \set{1, \ldots, \gamma}$, collectively forming the block $\set{1, \ldots, \mu, 1', \ldots, \gamma'}$; and
					\item for each $j \in \set{\mu+1, \ldots, k}$, the vertical line $\set{j, j'}$ in $\tomegabar{(\mu, \mu)}{k}$ connects to the vertical line $\set{j, j'}$ in $\apmorph{\gamma+1, \gamma+m+1, \ldots, \mu-m+1}{\gamma+1, \gamma+m+1, \ldots, \mu-m+1}$.
				\end{enumerate}
				Hence $\tomegabar{(\mu, \gamma)}{k} = \tomegabar{(\mu, \mu)}{k}\apmorph{\gamma+1, \gamma+m+1, \ldots, \mu-m+1}{\gamma+1, \gamma+m+1, \ldots, \mu-m+1} \in \pmodbound{m}{k}$. It follows analogously that if $\mu < \gamma$ then $\tomegabar{(\mu, \gamma)}{k} = \apmorph{\mu+1, \mu+m+1, \ldots, \gamma-m+1}{\mu+1, \mu+m+1, \ldots, \gamma-m+1}\tomegabar{(\gamma, \gamma)}{k} \in \pmodbound{m}{k}$.
			\end{proof}
		\end{proposition}
		
		\begin{corollary}
			For each $m \in \posints$, $k \in \intsge{m}$ and $\mu, \gamma \in \set{0, \ldots, k}$ such that $\mu + \gamma > 0$ and $\congmod{\mu}{\gamma}{m}$, blocks of type $(\mu, \gamma)$ do appear in some elements of $\pmodbound{m}{k}$.
			
			\begin{proof}
				Trivially follows from Proposition \ref{prop:tomegabarmugammainpmodbound} where we established that the well-defined bipartition $\tomegabar{(\mu,\gamma)}{k}$, which contains a block of type $(\mu, \gamma)$, may be factorised into a product containing only $(2,2)$-transapsis generators and $m$-apsis generators.
			\end{proof}
		\end{corollary}	

    \subsection{Building blocks}
		We established our characterisation of the $m$-apsis monoid $\apsismon{m}{k}$ by first establishing transversal and non-transversal building blocks, then factorised each element of the $m$-apsis monoid $\apsismon{m}{k}$ into a product containing transversal and non-transversal building blocks. Establishing that the planar mod-$m$ monoid $\pmodmon{m}{k}$ is equal to the submonoid generated by $(2,2)$-transapsis and $m$-apsis generators $\pmodbound{m}{k}$ may be done in a similar fashion. 
		
		Note that for each $m \in \posints$, $k \in \intsge{m}$, and $\mu, \gamma \in \set{0, \ldots, k}$ such that $\mu + \gamma > 0$ and $\tcongmod{\mu}{\gamma}{k}{m}$, in order to be able to form a bipartition containing precisely one block of type $(\mu, \gamma)$ then only $m$-apses, we additionally require that $\congmod{\mu}{k}{m}$, or equivalently that $\congmod{\gamma}{k}{m}$. 
		
		Both for the sake of succinctness and due to there being no notational ambiguity in doing so, when we need to require $\congmod{\mu}{\gamma}{m}$ and $\congmod{\mu}{k}{m}$ we simply state that we require $\tcongmod{\mu}{\gamma}{k}{m}$.
		
		\begin{definition} \label{def:tomegamugammapmod}
			For each $m \in \posints$, $k \in \intsge{m}$ and $\mu, \gamma \in \set{0, \ldots, k}$ such that $\mu + \gamma > 0$ and $\tcongmod{\mu}{\gamma}{k}{m}$, let $\tomega{(\mu, \gamma)}{k}$ denote the bipartition in $\pmodmon{m}{k}$ containing:
			\begin{enumerate}
				\item the type $(\mu, \gamma)$ block $\set{1, \ldots, \mu, 1', \ldots, \gamma'}$;
				\item for each $j \in \set{1, \ldots, \frac{k-\mu}{m}}$, the upper $m$-apsis $\{\mu + m(j-1) + 1, \ldots, \mu + mj\}$; and
				\item for each $j \in \set{1, \ldots, \frac{k-\gamma}{m}}$, the lower $m$-apsis $\{(\gamma + m(j-1) + 1)', \ldots, (\gamma + mj)'\}$.
			\end{enumerate} 
		\end{definition}
		
		For example given $m=2$ and $k=9$, Figure \ref{fig:omega15pmodmon} depicts $\tomega{(1, 5)}{9} \in \pmodmon{2}{9}$.
		
		\begin{figure}[!ht]
			\caption[ ]{Given $m=2$ and $k=9$,}
			\label{fig:omega15pmodmon}
			\vspace{5pt}
			\centering
			\input{chap_characterisations/tikz/fig-omega15pmodmon.tex}
		\end{figure} 
		
		We proceed by establishing in Proposition \ref{prop:tomegamugammainpmodbound} that $\tomega{(\mu, \gamma)}{k} \in \pmodbound{m}{k}$. The reader may find it useful to refer to Figure \ref{fig:omega51_9=omegabar51_9sigma6868pmodmon} which illustrates that, given $m=2$ and $k=9$, $\tomega{(5,1)}{9} = \tomegabar{(5,1)}{9}\apmorph{2,4,6,8}{2,4,6,8} \in \pmodbound{2}{9}$.
		
    	\begin{proposition} \label{prop:tomegamugammainpmodbound}
    		For each $m \in \posints$, $k \in \intsge{m}$ and $\mu, \gamma \in \set{0, \ldots, k}$ such that $\mu + \gamma > 0$ and $\tcongmod{\mu}{\gamma}{k}{m}$, \[\tomega{(\mu, \gamma)}{k} = \tomegabar{(\mu, \gamma)}{k}\apmorph{\gamma+1, \gamma+m+1, \gamma+2m+1, \ldots, k-m+1}{\gamma+1, \gamma+m+1, \gamma+2m+1, \ldots, k-m+1} \in \pmodbound{m}{k}.\]
    		
    		\begin{proof}
    			Let $\overbar{k} = \max\set{\mu, \gamma}$ and recall that for each $j \in \set{\overbar{k}+1, \ldots, k}$, $\tomegabar{(\mu, \gamma)}{k}$ contains the vertical line $\set{j, j'}$. 
    			
    			When forming the product $\tomegabar{(\mu, \gamma)}{k}\apmorph{\gamma+1, \gamma+m+1, \gamma+2m+1, \ldots, k-m+1}{\gamma+1, \gamma+m+1, \gamma+2m+1, \ldots, k-m+1}$:
    			\begin{enumerate}
 					\item the type $(\mu, \gamma)$ block $\set{1, \ldots, \mu, 1', \ldots, \gamma'}$ in $\tomegabar{(\mu, \gamma)}{k}$, which connects to the vertical lines $\set{j, j'}$ in $\apmorph{\gamma+1, \gamma+m+1, \gamma+2m+1, \ldots, k-m+1}{\gamma+1, \gamma+m+1, \gamma+2m+1, \ldots, k-m+1}$ where $j \in \set{1, \ldots, \gamma}$, is preserved;
    				\item the lower $m$-apses in $\apmorph{\gamma+1, \gamma+m+1, \gamma+2m+1, \ldots, k-m+1}{\gamma+1, \gamma+m+1, \gamma+2m+1, \ldots, k-m+1}$, which are lower $m$-apses in $\tomega{(\mu, \gamma)}{k}$, are preserved;
    				\item for each $j \in \set{1, \ldots, \frac{k-\overbar{k}}{m}}$, the $m$-apsis $\{\overbar{k} + m(j-1) + 1, \ldots, \overbar{k} + mj\}$ in $\apmorph{\gamma+1, \gamma+m+1, \gamma+2m+1, \ldots, k-m+1}{\gamma+1, \gamma+m+1, \gamma+2m+1, \ldots, k-m+1}$, which joins to vertical lines in $\tomegabar{(\mu, \gamma)}{k}$, is preserved;
					\item for each $j \in \set{2, \ldots, \frac{\overbar{k}-\gamma}{m}}$, the lower $m$-apsis $\{(\gamma+m(j-1)+1)', \ldots, (\gamma+mj)'\}$ in $\tomegabar{(\mu, \gamma)}{k}$ and the upper $m$-apsis $\{\gamma+m(j-1) + 1, \ldots, \gamma+mj\}$ in $\apmorph{\gamma+1, \gamma+m+1, \gamma+2m+1, \ldots, k-m+1}{\gamma+1, \gamma+m+1, \gamma+2m+1, \ldots, k-m+1}$ join and are removed; and
					\item for each $j \in \set{1, \ldots, \frac{\overbar{k}-\mu}{m}}$, the upper $m$-apsis $\{\mu+m(j-1)+1, \ldots, \mu+mj\}$ in $\tomegabar{(\mu, \gamma)}{k}$ is preserved.
    			\end{enumerate}
    			
    			Hence $\tomega{(\mu, \gamma)}{k} = \tomegabar{(\mu, \gamma)}{k}\apmorph{\gamma+1, \gamma+m+1, \gamma+2m+1, \ldots, k-m+1}{\gamma+1, \gamma+m+1, \gamma+2m+1, \ldots, k-m+1} \in \pmodbound{m}{k}$.
    		\end{proof}
    	\end{proposition}
    	
		\begin{figure}[!ht]
			\caption[ ]{Given $m=2$ and $k=9$,}
			\label{fig:omega51_9=omegabar51_9sigma6868pmodmon}
			\vspace{5pt}
			\centering
			\input{chap_characterisations/tikz/fig-omega51_9=omegabar51_9sigma6868pmodmon.tex}
		\end{figure} 
    
    \subsection{Transversal building blocks}
       	We now turn our attention to bipartitions in $\pmodmon{m}{k}$ where non-transversals must be $m$-apses. 
       	
       	\begin{definition}
       		For each $m \in \intsge{0}$ and $k \in \intsge{m}$, let $\pmodtrans{m}{k}$ denote the set of all $\alpha \in \pmodmon{m}{k}$ such that every non-transversal block $b \in \alpha$ is an $m$-apsis. We shall refer to elements of $\pmodtrans{m}{k}$ as \textit{transversal building blocks}.
       	\end{definition}
       	
       	We proceed in this subsection to establish that $\pmodtrans{m}{k}$ is contained within the monoid generated by $(2,2)$-transapsis generators and $m$-apsis generators $\pmodbound{m}{k}$. In order to do so we partition $\pmodtrans{m}{k}$ based on transversal types from left to right then establish our desired result inductively.
       	
       	Note that for each $m \in \posints$, $k \in \intsge{m}$, $r \in \set{0, \ldots, k}$, and $\mu_1, \ldots, \mu_r$, $\gamma_1, \ldots, \gamma_r \in \set{1, \ldots, k}$ such that $\congmod{\mu_j}{\gamma_j}{k}$ for all $j \in \set{1, \ldots, r}$ and $\Sigma^r_{j=1}\mu_j, \Sigma^r_{j=1}\gamma_j \leq k$, in order to be able to form a bipartition of rank $r$ with a distinct block of type $(\mu_j, \gamma_j)$ designated for each $j \in \set{1, \ldots, r}$, we additionally require that $\congmod{\Sigma^r_{j=1}\mu_j}{k}{m}$, or equivalently that $\congmod{\Sigma^r_{j=1}\gamma_j}{k}{m}$. 
       	
       	\begin{definition}
       		For each $m \in \posints$, $k \in \intsge{m}$, $r \in \set{0, \ldots, k}$, and $\mu_1, \ldots, \mu_r$, $\gamma_1, \ldots, \gamma_r \in \set{1, \ldots, k}$ such that $\congmod{\mu_j}{\gamma_j}{k}$ for all $j \in \set{1, \ldots, r}$, $\Sigma^r_{j=1}\mu_j$, $\Sigma^r_{j=1}\gamma_j \leq k$ and $\congmod{\Sigma^r_{j=1}\mu_j}{k}{m}$, let $\pmodtrans{(\mu_1, \gamma_1),\ldots,(\mu_r, \gamma_r)}{k}$ denote the set of all $\omega \in \pmodtrans{m}{k}$ such that $\omega$ contains precisely:
			\begin{enumerate}
				\item $\frac{k - \Sigma^r_{j=1}\mu_j}{m}$ upper $m$-apses;
				\item $\frac{k - \Sigma^r_{j=1}\gamma_j}{m}$ lower $m$-apses; and
				\item $r$ transversals of type $(\mu_1, \gamma_1)$, \ldots, $(\mu_r, \gamma_r)$ from left to right.
			\end{enumerate} 
       	\end{definition}
       	
       	Note that by definition $\pmodtrans{m}{k} = \bigcup_{\substack{r \in \set{0, \ldots, k} \\ \mu_1, \ldots, \mu_r, \gamma_1, \ldots, \gamma_r \in \set{1, \ldots, k} \\ \congmod{\mu_j}{\gamma_j}{m}\ \forall\ j \in \set{1, \ldots, r} \\ \Sigma^r_{j=1}\mu_j, \Sigma^r_{j=1}\gamma_j \leq k \\ \congmod{\Sigma^r_{j=1}\mu_j}{k}{m}}} \pmodtrans{(\mu_1, \gamma_1), \ldots, (\mu_r, \gamma_r)}{k}$.
		
\begin{figure}[!ht]
			\caption[ ]{Given $m=2$ and $k=13$,}
			\label{fig:pmodtransprop}
			\vspace{10pt}
			\centering
			\input{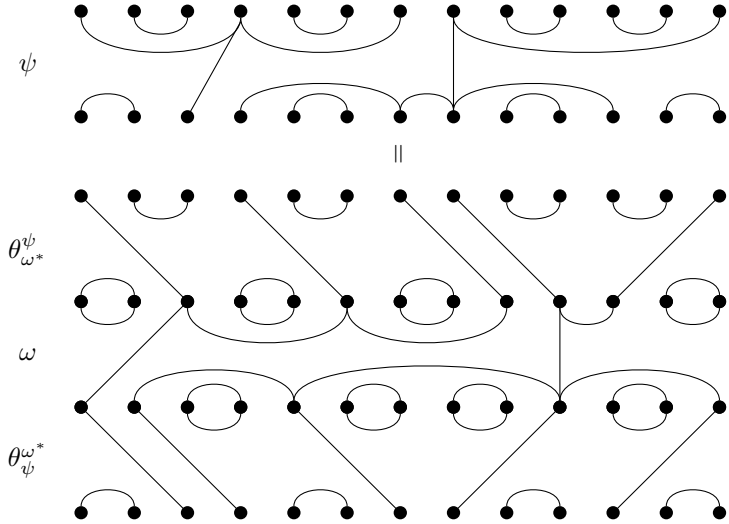}
		\end{figure}

		Next we seek to establish whether $\pmodtrans{(\mu_1, \gamma_1),\ldots,(\mu_r, \gamma_r)}{k} \subseteq \pmodbound{m}{k}$, which may be done in an almost identical fashion to how we established our characterisation of the $m$-apsis monoid $\apsismon{m}{k}$. We again begin by establishing that for each $\omega \in \pmodtrans{(\mu_1, \gamma_1),\ldots,(\mu_r, \gamma_r)}{k}$, each element of $\pmodtrans{(\mu_1, \gamma_1),\ldots,(\mu_r, \gamma_r)}{k}$ may be factorised into a product containing $\omega$ and $m$-apmorphisms. An example of such a factorisation when $m=2$ is depicted in Figure \ref{fig:pmodtransprop}.
		
		\begin{proposition} \label{prop:pmodtransmujgammajinpmodbound}
			For each $\omega \in \pmodtrans{(\mu_1, \gamma_1),\ldots,(\mu_r, \gamma_r)}{k}$,\[\pmodtrans{(\mu_1, \gamma_1),\ldots,(\mu_r, \gamma_r)}{k} \subseteq \apmorphs{m}{k}\omega\apmorphs{m}{k}.\]
			
			\begin{proof}
				Follows identically to the proof of Proposition \ref{prop:aptransmujgammajinapsismon}, which was our analogous result when characterising the $m$-apsis monoid $\apsismon{m}{k}$, to show that for each $\psi \in \pmodtrans{(\mu_1, \gamma_1),\ldots,(\mu_r, \gamma_r)}{k}$, $\psi = \apmorph{\psi}{\omega^*}\omega\apmorph{\omega^*}{\psi}$ and hence $\pmodtrans{(\mu_1, \gamma_1),\ldots,(\mu_r, \gamma_r)}{k} \subseteq \apmorphs{m}{k}\omega\apmorphs{m}{k}$.
			\end{proof}
		\end{proposition}
		
		\begin{corollary} \label{cor:whenpmodtransinpmodbound}
			If the intersection of $\pmodtrans{(\mu_1, \gamma_1),\ldots,(\mu_r, \gamma_r)}{k}$ with the monoid generated by $(2,2)$-transapsis generators and $m$-apsis generators $\pmodbound{m}{k}$ is non-empty then $\pmodtrans{(\mu_1, \gamma_1),\ldots,(\mu_r, \gamma_r)}{k} \subseteq \pmodbound{m}{k}$.
			
			\begin{proof}
				Suppose there exists $\omega \in \pmodtrans{(\mu_1, \gamma_1),\ldots,(\mu_r, \gamma_r)}{k} \cap \pmodbound{m}{k}$, then it trivially follows from Proposition \ref{prop:pmodtransmujgammajinpmodbound} that $\pmodtrans{(\mu_1, \gamma_1),\ldots,(\mu_r, \gamma_r)}{k} \subseteq \apmorphs{m}{k}\omega\apmorphs{m}{k} \subseteq \pmodbound{m}{k}$.
			\end{proof}
		\end{corollary}
    	
    	\begin{proposition} \label{prop:pmodtransmugammainpmodmon}
    		For each $m \in \posints$, $k \in \intsge{m}$, $\mu, \gamma \in \set{1, \ldots, k}$ such that $\tcongmod{\mu}{\gamma}{k}{m}$, $\pmodtrans{(\mu, \gamma)}{k} \subseteq \pmodbound{m}{k}$.
    		
    		\begin{proof}
    			It follows by definition that $\tomega{(\mu, \gamma)}{k} \in \pmodtrans{(\mu, \gamma)}{k}$, and we established in Proposition \ref{prop:tomegamugammainpmodbound} that $\tomega{(\mu, \gamma)}{k} \in \pmodbound{m}{k}$. Hence $\tomega{(\mu, \gamma)}{k} \in \pmodtrans{(\mu, \gamma)}{k} \cap \pmodbound{m}{k}$, employing Corollary \ref{cor:whenpmodtransinpmodbound} we have $\pmodtrans{(\mu, \gamma)}{k} \subseteq \pmodbound{m}{k}$.
    		\end{proof}
    	\end{proposition}
    	
		\begin{definition}
	    	For each $m \in \posints$, $k \in \intsge{m}$, $r \in \set{0, \ldots, k}$, and $\mu_1, \ldots, \mu_r$, $\gamma_1, \ldots, \gamma_r \in \set{1, \ldots, k}$ such that $\congmod{\mu_j}{\gamma_j}{k}$ for all $j \in \set{1, \ldots, r}$, $\Sigma^r_{j=1}\mu_j$, $\Sigma^r_{j=1}\gamma_j \leq k$ and $\congmod{\Sigma^r_{j=1}\mu_j}{k}{m}$, let $\tomega{(\mu_1, \gamma_1), \ldots, (\mu_r, \gamma_r)}{k}$ denote the element of $\pmodtrans{(\mu_1, \gamma_1), \ldots, (\mu_r, \gamma_r)}{k}$ containing: 
		    	\begin{enumerate}
		    		\item the type $(\mu_1, \gamma_1)$ transversal $\set{1, \ldots, \mu_1, 1', \ldots, \gamma_1'}$;
		    		\item for each $j \in \set{1, \ldots, \frac{k - \Sigma^r_{j=1}\mu_j}{m}}$, the upper $m$-apsis $\{\mu_1 + m(j-1) + 1, \ldots, \mu_1 + mj\}$;
		    		\item for each $j \in \set{1, \ldots, \frac{k - \Sigma^r_{j=1}\gamma_j}{m}}$, the lower $m$-apsis $\{(\gamma_1+m(j-1)+1)', \ldots, (\gamma_1 + mj)'\}$;
		    		\item for each $j \in \set{2, \ldots, r}$, the type $(\mu_j, \gamma_j)$ transversal $\big\{k - \Sigma^{r}_{i=j}\mu_i + 1, \ldots, k - \Sigma^{r}_{i=j+1}\mu_i, (k - \Sigma^{r}_{i=j}\gamma_i + 1)', \ldots, (k - \Sigma^{r}_{i=j+1}\gamma_i)'\big\}$,
		    	\end{enumerate} 
			and let $\tupsilon{(\mu_1, \gamma_1), \ldots, (\mu_r, \gamma_r)}{k}$ denote the element of $\pmodtrans{(\mu_1, \gamma_1), \ldots, (\mu_r, \gamma_r)}{k}$ containing: 
		    	\begin{enumerate}
		    		\item for each $j \in \set{1, \ldots, \frac{k - \Sigma^r_{j=1}\mu_j}{m}}$, the upper $m$-apsis $\{m(j-1) + 1, \ldots, mj\}$;
		    		\item for each $j \in \set{1, \ldots, \frac{k - \Sigma^r_{j=1}\gamma_j}{m}}$, the lower $m$-apsis $\{(m(j-1) + 1)', \ldots,$ $(mj)'\}$;
		    		\item for each $j \in \set{1, \ldots, r}$, the type $(\mu_j, \gamma_j)$ transversal $\big\{k - \Sigma^{r}_{i=j}\mu_i + 1, \ldots, k - \Sigma^{r}_{i=j+1}\mu_i, (k - \Sigma^{r}_{i=j}\gamma_i + 1)', \ldots, (k - \Sigma^{r}_{i=j+1}\gamma_i)'\big\}$.
		    	\end{enumerate} 
		\end{definition}
    	
    	For example given $m=2$, $\tomega{(3,1),(1,3)}{4} = \tupsilon{(3,1),(1,3)}{4} \in \pmodtrans{(3,1),(1,3)}{4}$ along with $\tomega{(4,2),(3,1)}{7}, \tupsilon{(4,2),(3,1)}{7} \in \pmodtrans{(4,2),(3,1)}{7}$ are depicted in Figure \ref{fig:tomegatupsilonpmod}.
    	
    	\begin{figure}[!ht]
			\caption[ ]{Given $m=2$,}
			\label{fig:tomegatupsilonpmod}
			\vspace{10pt}
			\centering
			\input{chap_characterisations/tikz/fig-tomegatupsilonpmod.tex}
		\end{figure}
		
		Analogously to when characterising the $m$-apsis monoid $\apsismon{m}{k}$, we proceed by establishing in Proposition \ref{prop:pmodtransfactorisation} that for each $m \in \posints$, $k \in \intsge{m}$, $r \in \set{0, \ldots, k}$, and $\mu_1, \ldots, \mu_r, \gamma_1, \ldots, \gamma_r \in \set{1, \ldots, k}$ such that $\congmod{\mu_j}{\gamma_j}{k}$ for all $j \in \set{1, \ldots, r}$, $\Sigma^r_{j=1}\mu_j, \Sigma^r_{j=1}\gamma_j \leq k$ and $\congmod{\Sigma^r_{j=1}\mu_j}{k}{m}$, \[\tomega{(\mu_1, \gamma_1), \ldots, (\mu_r, \gamma_r)}{k} \in \begin{cases}(\pmodtrans{(\mu_1, \gamma_1)}{k-\Sigma^r_{j=2}\mu_j} \oplus \id{\Sigma^r_{j=2}\mu_j})(\id{\gamma_1} \oplus \pmodtrans{(\mu_2,\gamma_2), \ldots, (\mu_r, \gamma_r)}{k-\gamma_1}) & \mu_1 \geq \gamma_1; \\ (\id{\mu_1} \oplus \pmodtrans{(\mu_2,\gamma_2), \ldots, (\mu_r, \gamma_r)}{k-\mu_1})(\pmodtrans{(\mu_1, \gamma_1)}{k-\Sigma^r_{j=2}\gamma_j} \oplus \id{\Sigma^r_{j=2}\gamma_j}) & \mu_1 \leq \gamma_1. \end{cases}\] 
		
		For example, given $m=2$, Figure \ref{fig:pmodtrans3113element} illustrates that $\tomega{(3,1),(1,3)}{4} = (\tomega{(3,1)}{3}\oplus\id{1})(\id{1}\oplus\tupsilon{(1,3)}{4}) \in \pmodtrans{(3,1),(1,3)}{4}$ and Figure \ref{fig:pmodtrans156213element} illustrates that $\tomega{(1,5),(6,2),(1,3)}{10} = (\id{1}\oplus\tupsilon{(6,2),(1,3)}{9})(\tomega{(1,5)}{5}\oplus\id{5}) \in \pmodtrans{(1,5),(6,2),(1,3)}{10}$.

		\begin{proposition} \label{prop:pmodtransfactorisation}
			For each $m \in \posints$, $k \in \intsge{m}$, $r \in \set{0, \ldots, k}$, and $\mu_1, \ldots, \mu_r$, $\gamma_1, \ldots, \gamma_r \in \set{1, \ldots, k}$ such that $\congmod{\mu_j}{\gamma_j}{k}$ for all $j \in \set{1, \ldots, r}$, $\Sigma^r_{j=1}\mu_j$, $\Sigma^r_{j=1}\gamma_j \leq k$ and $\congmod{\Sigma^r_{j=1}\mu_j}{k}{m}$, 
			\[\tomega{(\mu_1, \gamma_1),\ldots,(\mu_r,\gamma_r)}{k} = \begin{cases} (\tomega{(\mu_1, \gamma_1)}{k - \Sigma^r_{j=2}\mu_j} \oplus \id{\Sigma^r_{j=2}\mu_j})(\id{\gamma_1} \oplus \tupsilon{(\mu_2, \gamma_2), \ldots, (\mu_r, \gamma_r)}{k - \gamma_1}) & \mu_1 \geq \gamma_1; \\ (\id{\mu_1} \oplus \tupsilon{(\mu_2, \gamma_2), \ldots, (\mu_r, \gamma_r)}{k-\mu_1})(\tomega{(\mu_1, \gamma_1)}{k-\Sigma^r_{j=2}\gamma_j} \oplus \id{\Sigma^r_{j=2}\gamma_j}) & \mu_1 \leq \gamma_1. \end{cases}\]
			
			\begin{proof}
				Suppose $\mu_1 \geq \gamma_1$. Note that since $k \geq \Sigma^r_{j=1}\mu_j \geq \gamma_1 + \Sigma^r_{j=2}\mu_j$ we have:
				\begin{enumerate}
					\item $k - \Sigma^r_{j=2}\mu_j \geq \mu_1 = \max\set{\mu_1, \gamma_1}$, ensuring $\tomega{(\mu_1, \gamma_1)}{k - \Sigma^r_{j=2}\mu_j}$ is well-defined;
					\item $k-\gamma_1 \geq \Sigma^r_{j=2}\mu_j, \Sigma^r_{j=2}\gamma_j$, ensuring $\tupsilon{(\mu_2, \gamma_2), \ldots, (\mu_r, \gamma_r)}{k - \gamma_1}$ is well-defined; and
					\item $\frac{k - \gamma_1 - \Sigma^r_{j=2}\mu_j}{m}$, which is both the number of lower $m$-apses in $\tomega{(\mu_1, \gamma_1)}{k - \Sigma^r_{j=2}\mu_j}$ and the number of upper $m$-apses in $\tupsilon{(\mu_2, \gamma_2), \ldots, (\mu_r, \gamma_r)}{k - \gamma_1}$, is a non-negative integer.
				\end{enumerate}
					
				The remainder of the proof follows identically to the proof in Proposition \ref{prop:aptransfactorisation}, which was our analogous result when characterising the $m$-apsis monoid $\apsismon{m}{k}$.
			\end{proof}
		\end{proposition}
		
		\begin{figure}[!ht]
			\caption[ ]{Given $m=2$,}
			\label{fig:pmodtrans3113element}
			\vspace{10pt}
			\centering
			\input{chap_characterisations/tikz/fig-pmodtrans3113element.tex}
		\end{figure}
		
		\begin{figure}[!ht]
			\caption[ ]{Given $m=3$,}
			\label{fig:pmodtrans156213element}
			\vspace{10pt}
			\centering
			\input{chap_characterisations/tikz/fig-pmodtrans156213element.tex}
		\end{figure}

		\begin{proposition} \label{prop:pmodtransmujgammajinpmodmon}
			For each $m \in \posints$, $k \in \intsge{m}$, $r \in \set{0, \ldots, k}$, and $\mu_1, \ldots, \mu_r$, $\gamma_1, \ldots, \gamma_r \in \set{1, \ldots, k}$ such that $\congmod{\mu_j}{\gamma_j}{k}$ for all $j \in \set{1, \ldots, r}$, $\Sigma^r_{j=1}\mu_j$, $\Sigma^r_{j=1}\gamma_j \leq k$ and $\congmod{\Sigma^r_{j=1}\mu_j}{k}{m}$, we have $\pmodtrans{(\mu_1, \gamma_1),\ldots,(\mu_r, \gamma_r)}{k} \subseteq \pmodbound{m}{k}$.
			\begin{proof}
				Follows identically to the proof of Proposition \ref{prop:aptransmujgammajinapsismon}, which was our analogous result when characterising the $m$-apsis monoid $\apsismon{m}{k}$.
			\end{proof}
		\end{proposition}

		\begin{proposition}
			For each $m \in \posints$ and $k \in \intsge{m}$, $\pmodtrans{m}{k} \subseteq \pmodbound{m}{k}$.
			
			\begin{proof}
				Recall that by definition $\pmodtrans{m}{k} = \bigcup_{\substack{r \in \set{0, \ldots, k} \\ \mu_1, \ldots, \mu_r, \gamma_1, \ldots, \gamma_r \in \set{1, \ldots, k} \\ \congmod{\mu_j}{\gamma_j}{m}\ \forall\ j \in \set{1, \ldots, r} \\ \Sigma^r_{j=1}\mu_j, \Sigma^r_{j=1}\gamma_j \leq k \\ \congmod{\Sigma^r_{j=1}\mu_j}{k}{m}}} \pmodtrans{(\mu_1, \gamma_1), \ldots, (\mu_r, \gamma_r)}{k}$, hence it follows from Proposition \ref{prop:pmodtransmujgammajinpmodmon} that $\pmodtrans{m}{k} \subseteq \pmodbound{m}{k}$.
			\end{proof}
     \end{proposition}

    \subsection{Non-transversal building blocks}
       	We now turn our attention to bipartitions in the planar mod-$m$ monoid $\pmodmon{m}{k}$ where transversals must be lines, and either lower non-transversals or upper non-transversals must be $m$-apses. 
    
		\begin{definition}
			For each $m \in \posints$ and $k \in \intsge{m}$, let $\pmodupnontrans{m}{k}$ denote the set of all $\alpha \in \pmodmon{m}{k}$ such that every:
			\begin{enumerate}
				\item transversal in $\alpha$ is a line; and
				\item lower non-transversal in $\alpha$ is an $m$-apsis, 
			\end{enumerate}
			furthermore let $\pmodlownontrans{m}{k}$ denote the set of all $\alpha \in \pmodmon{m}{k}$ such that every:
			\begin{enumerate}
				\item transversal in $\alpha$ is a line; and
				\item upper non-transversal in $\alpha$ is an $m$-apsis.
			\end{enumerate}
			We shall refer to elements of $\pmodupnontrans{m}{k}$ as \textit{upper non-transversal building blocks}, and elements of $\pmodlownontrans{m}{k}$ as \textit{lower non-transversal building blocks}. By a \textit{non-transversal building block} we shall mean an upper or lower non-transversal building block.
		\end{definition}
     
		We proceed in this subsection to establish that both $\pmodupnontrans{m}{k}$ and $\pmodlownontrans{m}{k}$ are contained within the monoid generated by $(2,2)$-transapsis genrators and $m$-apsis generators $\pmodbound{m}{k}$. Containment of $\pmodlownontrans{m}{k}$ will trivially follow from containment of $\pmodupnontrans{m}{k}$. In order to establish the latter containment, first we partition $\pmodupnontrans{m}{k}$ based on the number of upper non-transversals that are not $m$-apses, then we establish our desired result inductively.
       	
       	Note that $\floor{\frac{k}{m}}$ is the maximum number of upper non-transversals an element of $\pmodupnontrans{m}{k}$ may contain, however in such a case at least one upper non-transversal must be an $m$-apsis, hence $\pmodupnontrans{m}{k}$ may contain at most $\floor{\frac{k-m}{m}}$ upper non-transversals that are not $m$-apses.
       	
       	\begin{definition}
       		For each $x \in \set{0, \ldots, \floor{\frac{k-m}{m}}}$, let $\pmodupnontrans{}{x}$ denote the subset of all bipartitions $\eta \in \pmodupnontrans{m}{k}$ such that $\eta$ contains precisely $x$ upper non-transversals that are not $m$-apses. 
       	\end{definition}
       		
       	Note it follows by definition that $\pmodupnontrans{}{0} = \apmorphs{m}{k}$ and that $\pmodupnontrans{m}{k} = \bigcup^{\floor{\frac{k-m}{m}}}_{x=0}\pmodupnontrans{}{x}$. 
       	
       	\begin{definition} \label{def:eta-andomegaetapmod}
			For each $x \in \set{1, \ldots, \floor{\frac{k-m}{m}}}$ and $\eta \in \pmodupnontrans{}{x}$, by definition there exists:
   			\begin{enumerate}
   				\item $\mu \in m\posints$ and $b = \set{b_1, \ldots, b_{\mu}} \in \eta$, where $b_1, \ldots, b_{\mu} \in \set{1, \ldots, k}$, such that $\mu \leq k$, $b$ is a type $(\mu, 0)$ non-transversal that is not an $m$-apsis and no upper non-transversals of $\eta$ pass underneath $b$;
   				\item $r_1 \in \set{0, \ldots, b_1-1}$ and $r_2 \in \set{0, \ldots, k-b_{\mu}}$ such that the number of transversal lines in $\eta$ whose upper vertex sits to the left or right of $b \in \eta$ is $r_1$ and $r_2$ respectively; and
   				\item upper vertices $u_1, \ldots, u_{r_1} \in \set{1, \ldots, b_1-1}$, $u_{r_1+1}, \ldots, u_{r_1+r_2} \in \{b_{\mu}+1,$ $\ldots, k\}$ and lower vertices $l_1, \ldots, l_{r_1+r_2} \in \set{1, \ldots, k}$ such that for each $j \in \{1, \ldots, r_1+r_2\}$, $\eta$ contains the transversal line $\set{u_j, l_j}$.
   			\end{enumerate}
      			
   			\begin{itemize}
   				\item Let $\omega_{\eta}$ denote $ \id{r_1}\oplus\tomega{(\mu,0)}{k-r_1-r_2}\oplus\id{r_2}$ where $\tomega{(\mu,0)}{k-r_1-r_2}$ was outlined in Definition \ref{def:tomegamugammapmod} and shown to be an element of $\pmodbound{m}{k-r_1-r_2}$ in Proposition \ref{prop:tomegamugammainpmodbound}; and
   			
      			\item let $\eta^-$  denote the upper non-transversal building block in $\pmodupnontrans{}{x-1}$ containing:
	         	\begin{enumerate}
	         		\item excluding $b$, the upper non-transversals of $\eta$;
	         		\item for each $j \in \set{1, \ldots, \mu}$, the transversal line $\set{b_j, r_1+j}$;
	         		\item for each $j \in \set{1, \ldots, r_1}$, the transversal line $\set{u_j, j'}$;
	         		\item for each $j \in \set{1, \ldots, r_2}$, the transversal line $\set{u_j, (k-r_2+j)'}$; and
	         		\item for each $j \in \set{1, \ldots, \frac{k-\mu-r_1-r_2}{m}}$, the lower $m$-apsis $\{(r_1+\mu+m(j-1)+1)', \ldots, (r_1+\mu+mj)'\}$.
	         	\end{enumerate}
        	\end{itemize}
       	\end{definition}
       	
       	Figure \ref{fig:pmodupnontranseg} contains an example of $\omega_{\eta}$ and $\eta^-$ for $m=2$, $k=12$ and a given $\eta \in \apupnontrans{}{3}$.

		\begin{figure}[!ht]
			\caption[ ]{Let $m=2$, $k=12$.}
			\label{fig:pmodupnontranseg}
			\vspace{10pt}
			\centering
			\input{chap_characterisations/tikz/fig-pmodupnontranseg.tex}
		\end{figure}
       	
       	We proceed to establish that each $\pmodupnontrans{}{x}$ is contained within the monoid generated by $(2,2)$-transapsis generators and $m$-apsis generators $\pmodbound{m}{k}$. To do so we will show that each $\eta \in \pmodupnontrans{}{x}$, where $x \in \set{1, \ldots \floor{\frac{k-m}{m}}}$, may be factorised into the product $\eta^-\omega_{\eta}\apmorph{\omega_{\eta}^*}{\eta}$, an example of which is depicted in Figure \ref{fig:pmodupnontransfactorisation}.
       	
       	\begin{proposition} \label{prop:pmodupnontransinpmodmon}
       		For each $m \in \posints$, $k \in \intsge{m}$ and $x \in \set{0, \ldots, \floor{\frac{k-m}{m}}}$, $\pmodupnontrans{}{x} \subseteq \pmodbound{m}{k}$.
       		
       		\begin{proof}
					By definition $\pmodupnontrans{}{0} = \apmorphs{m}{k}$. We established in Proposition \ref{prop:apmorphismsinapsismon} that for each $m \in \intsge{2}$ and $k \in \intsge{m}$, $\apmorphs{m}{k} \subseteq \apsismon{m}{k}$. Since $\apsismon{m}{k} \subseteq \pmodbound{m}{k}$, it follows that we also have $\apmorphs{m}{k} \subseteq \pmodbound{m}{k}$. We also noted in Corollary \ref{cor:monapmorphsinpmodbound} that for each $k \in \posints$, $\apmorphs{1}{k} \subseteq \pmodbound{1}{k}$. Hence $\pmodupnontrans{}{0} \subseteq \pmodbound{m}{k}$ for all $m \in \posints$ and $k \in \intsge{m}$.
       			 
        			Let $x \in \set{1, \ldots, \floor{\frac{k-m}{m}}}$ and suppose $\pmodupnontrans{}{x-1} \subseteq \pmodbound{m}{k}$. For each $\eta \in \pmodupnontrans{}{x}$, as previously noted in Definition \ref{def:eta-andomegaetapmod}, by definition there exists:
       			\begin{enumerate}
       				\item $\mu \in m\posints$ and $b = \set{b_1, \ldots, b_{\mu}} \in \eta$, where $b_1, \ldots, b_{\mu} \in \set{1, \ldots, k}$, such that $\mu \leq k$, $b$ is a type $(\mu, 0)$ non-transversal that is not an $m$-apsis and no upper non-transversals of $\eta$ pass underneath $b$;
       				\item $r_1 \in \set{0, \ldots, b_1-1}$ and $r_2 \in \set{0, \ldots, k-b_{\mu}}$ such that the number of transversal lines in $\eta$ whose upper vertex sits to the left or right of $b \in \eta$ is $r_1$ and $r_2$ respectively; and
       				\item upper vertices $u_1, \ldots, u_{r_1} \in \set{1, \ldots, b_1-1}$, $u_{r_1+1}, \ldots, u_{r_1+r_2} \in \{b_{\mu}+1,$ $\ldots, k\}$ and lower vertices $l_1, \ldots, l_{r_1+r_2} \in \set{1, \ldots, k}$ such that for each $j \in \{1, \ldots, r_1+r_2\}$, $\eta$ contains the transversal line $\set{u_j, l_j}$.
       			\end{enumerate}
        			        			
				Note $\omega_{\eta}$ and $\eta$ both have $\frac{k-\mu-r_1-r_2}{m}$ lower $m$-apses, hence the $m$-apmorphism $\apmorph{\omega_{\eta}^*}{\eta}$ is well-defined. When forming the product $\eta^-\omega_{\eta}\apmorph{\omega_{\eta}^*}{\eta}$:
       			\begin{enumerate}
        			\item the upper $m$-apses in $\eta^-$ are preserved, which are identical to the $m$-apses in $\eta$ after excluding $b$;
        			\item $b$ is formed by the $\mu$-apsis $\set{r_1 + 1, \ldots, r_1+\mu} \in \omega_{\eta}$ joining to each of $\eta^-$'s transversal lines $\set{b_j, r_1+j}$ where $j \in \set{1, \ldots, \mu}$;
        			\item the lower $m$-apses in $\apmorph{\omega_{\eta}^*}{\eta}$, which are identical to the lower $m$-apses in $\eta$, are preserved;
        			\item for each $j \in \set{1, \ldots, \frac{k-\mu-r_1-r_2}{m}}$, the lower $m$-apsis $\{r_1 + \mu + m(j-1)+1, \ldots, r_1+\mu+mj\} \in \eta^-$ and the upper $m$-apsis $\{(r_1 + \mu + m(j-1)+1)', \ldots, (r_1+\mu+mj)\} \in \omega_{\eta}$ join and are removed;
        			\item for each $j \in \set{1, \ldots, \frac{k-r_1-r_2}{m}}$, the lower $m$-apsis $\{r_1 + m(j-1)+1, \ldots, r_1+mj\} \in \omega_{\eta}$ and upper $m$-apsis $\{(r_1 + m(j-1)+1)', \ldots, (r_1+mj)\} \in \apmorph{\omega_{\eta}^*}{\eta}$ join and are removed;
        			\item for each $j \in \set{1, \ldots, r_1}$, the transversal line $\set{u_j, l_j'}$ is formed by the lines $\set{u_j, j'} \in \eta^-$, $\set{j, j'} \in \omega_{\eta}$ and $\set{j, l_j'} \in \apmorph{\omega_{\eta}^*}{\eta}$ joining; and
        			\item for each $j \in \set{k-r_2+1, \ldots, k}$, the line $\set{u_j, l_j'}$ is formed by the lines $\set{u_j, (k-r_2+j)'} \in \eta^-$, $\set{k-r_2+j, (k-r_2+j)'} \in \omega_{\eta}$ and $\set{k-r_2+j, l_j'} \in \apmorph{\omega_{\eta}^*}{\eta}$ joining.
       			\end{enumerate}
       			
       			Hence $\eta = \eta^-\omega_{\eta}\apmorph{\omega_{\eta}^*}{\eta} \in \pmodbound{m}{k}$. It follows by induction that $\pmodupnontrans{}{x} \subseteq \pmodbound{m}{k}$ for all $x \in \set{0, \ldots, \floor{\frac{k-m}{m}}}$.
       		\end{proof}
       	\end{proposition}
		
		\begin{figure}[!ht]
			\caption[ ]{Given $m=2$, $k=12$ and $\eta \in \pmodupnontrans{}{3}$ from Figure \ref{fig:pmodupnontranseg},}
			\label{fig:pmodupnontransfactorisation}
			\vspace{10pt}
			\centering
			\input{chap_characterisations/tikz/fig-pmodupnontransfactorisation.tex}
		\end{figure}
		
		\begin{proposition}
			For each $m \in \posints$ and $k \in \intsge{m}$, $\pmodupnontrans{m}{k}, \pmodlownontrans{m}{k} \subseteq \pmodbound{m}{k}$.
			
			\begin{proof}
				$\pmodupnontrans{m}{k} = \bigcup_{x=0}^{\floor{\frac{k-m}{m}}}\pmodupnontrans{}{x} \subseteq \pmodbound{m}{k}$. Containment of $\pmodlownontrans{m}{k}$ may be established in a dual fashion to $\pmodupnontrans{m}{k}$, or more succinctly by noting that $\pmodlownontrans{m}{k} = \tuple{\pmodupnontrans{m}{k}}^* \subseteq \pmodbound{m}{k}$.
			\end{proof}
		\end{proposition}

    \subsection{Factorising $\pmodmon{m}{k}$}
    	To establish that the planar mod-$m$ monoid $\pmodmon{m}{k}$ is in fact the monoid generated by $(2,2)$-transapsis generators and $m$-apsis generators $\pmodbound{m}{k}$, it remains for us to establish that each element of $\pmodmon{m}{k}$ may be factorised into a product of $(2,2)$-transapsis generators and $m$-apsis generators. 
		
		We will reach our desired conclusion by establishing that each element of $\pmodmon{m}{k}$ may be factorised into a product containing an upper non-transversal building block, a transversal building block and a lower non-transversal building block.
		
		\begin{definition} \label{def:transandnontransblockelmntspmod}
			For each $m \in \posints$, $k \in \intsge{m}$ and $\alpha \in \pmodmon{m}{k}$:
			\begin{itemize}
				\item let $u_{\alpha}$ denote the upper non-transversal building block containing:
				\begin{enumerate}
					\item the upper non-transversals in $\alpha$;
					\item for each $j \in \set{1, \ldots, k}$ such that $j$ is an element of a transversal in $\alpha$, the vertical line $\set{j, j'}$;
					\item lower $m$-apses containing the remaining vertices,
				\end{enumerate}
				\item let $t_{\alpha}$ denote the transversal building block containing:
				\begin{enumerate}
					\item the transversals contained in $\alpha$;
					\item $m$-apses replacing the non-transversals in $\alpha$, and
				\end{enumerate}
				\item let $l_{\alpha}$ denote the lower non-transversal building block containing:
				\begin{enumerate}
					\item the lower non-transversals contained in $\alpha$;
					\item for each $j \in \set{1, \ldots, k}$ such that $j'$ is an element of a transversal in $\alpha$, the vertical line $\set{j, j'}$;
					\item upper $m$-apses containing the remaining vertices.
				\end{enumerate}
			\end{itemize}
		\end{definition}
		
\begin{figure}[!ht]
			\caption[ ]{Let $m=2$ and $k =16$.}
			\label{fig:pmodmonintoblocks}
			\vspace{5pt}
			\centering
			\input{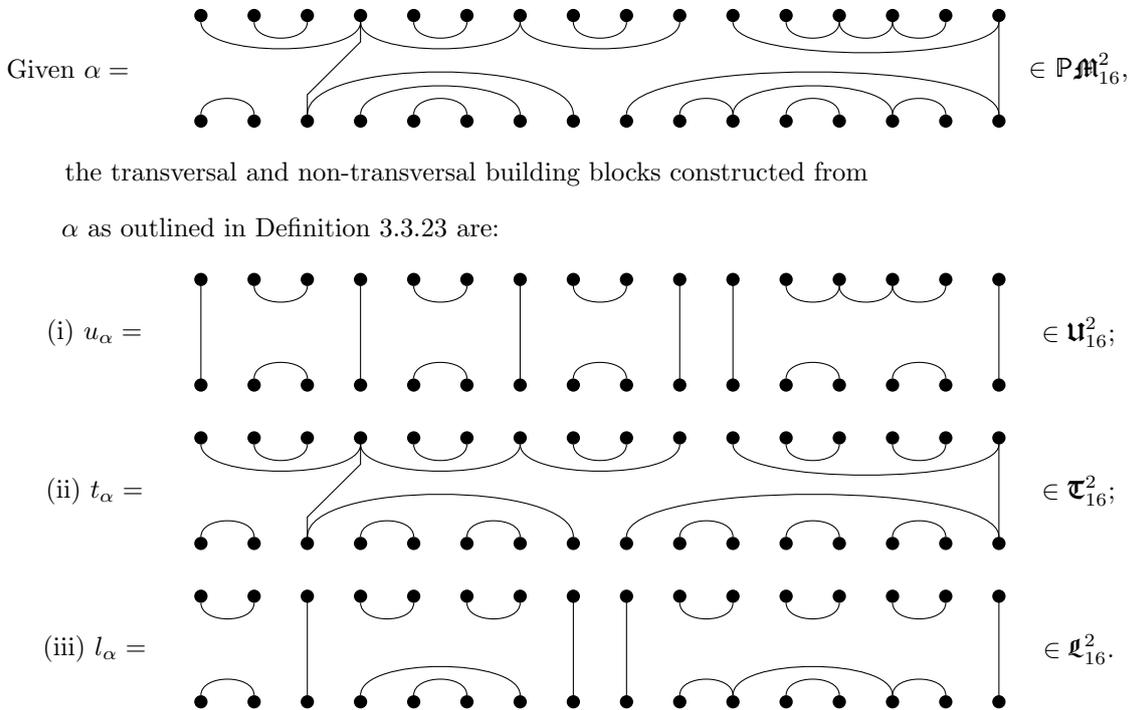}
		\end{figure}

		Figure \ref{fig:pmodmonintoblocks} contains an example of $u_{\alpha}$, $t_{\alpha}$ and $l_{\alpha}$ for a given $\alpha \in \apsisbound{2}{16}$.
		
		\begin{theorem} \label{thm:pmodmoncharacterisation}
			For each $m \in \posints$ and $k \in \intsge{m}$, the planar mod-$m$ monoid $\pmodmon{m}{k}$ is equal to the monoid generated by $(2,2)$-transapsis generators and $m$-apsis generators $\pmodbound{m}{k}$.
			
			\begin{proof}
				It follows trivially from $(2,2)$-transapsis generators and $m$-apsis generators sitting inside the planar mod-$m$ monoid $\pmodmon{m}{k}$ that $\pmodbound{m}{k}$ is contained within the planar mod-$m$ monoid $\pmodmon{m}{k}$. Conversely for each $\alpha \in \pmodmon{m}{k}$, when forming the product $u_{\alpha}t_{\alpha}l_{\alpha}$ (see Figure \ref{fig:pmodmonfactorisation} for an example):
				\begin{enumerate}
					\item the upper non-transversals in $u_{\alpha}$, which are identical to the upper non-transversals in $\alpha$, are preserved;
					\item the transversals in $t_{\alpha}$, which are identical to the transversals in $\alpha$, are preserved as they join to vertical lines in both $u_{\alpha}$ and $l_{\alpha}$ by construction;
					\item the lower non-transversals in $l_{\alpha}$, which are identical to the lower non-transversals in $\alpha$, are preserved; and
					\item the lower $m$-apses in $u_{\alpha}$ and upper $m$-apses in $t_{\alpha}$ join and are removed, similarly with the lower $m$-apses in $t_{\alpha}$ and upper $m$-apses in $l_{\alpha}$.
				\end{enumerate}
				
				Hence $\alpha = u_{\alpha}t_{\alpha}l_{\alpha}$, giving us $\pmodmon{m}{k} \subseteq \pmodupnontrans{m}{k}\pmodtrans{m}{k}\pmodlownontrans{m}{k} \subseteq \pmodbound{m}{k}$.
			\end{proof}
		\end{theorem}
		
		Note that no diapsis or $(2,2)$-transapsis generator may be generated by a product of other diapsis and $(2,2)$-transapsis generators. Hence the generating set of the planar mod-$m$ monoid established in Theorem \ref{thm:pmodmoncharacterisation} is minimal under set inclusion.
		
		\begin{figure}[!ht]
			\caption[ ]{Given $m=2$, $k=16$ and $\alpha \in \pmodmon{2}{16}$ from Figure \ref{fig:pmodmonintoblocks},}
			\label{fig:pmodmonfactorisation}
			\vspace{5pt}
			\centering
			\input{chap_characterisations/tikz/fig-pmodmonfactorisation.tex}
		\end{figure}

    \clearpage{\pagestyle{empty}\cleardoublepage}
\chapter{Cardinalities} \label{chap:cardinalities}
	\section{The planar modular partition monoid, $\card{\pmodmon{m}{k}}$}
		Recall from Definition \ref{def:pmodmon} that the planar mod-$m$ monoid $\pmodmon{m}{k}$ is the planar analogue of the mod-$m$ monoid $\modmon{m}{k}$, and hence consists of all planar bipartitions $\alpha \in \ppttnmon{k}$ such that each block $b \in \alpha$ satisfies $\congmod{\noupverts{b}}{\nolowverts{b}}{m}$.
		
		For each $m \in \posints$ and $k \in \intsge{m}$, we proceed to establish the cardinality of the planar mod-$m$ monoid $\pmodmon{m}{k}$ by:
		\begin{enumerate}
			\item for each $t \in m\nonnegints$ such that $t \leq k$, establishing recurrence relations for the number of $\pmodmon{m}{k}$-feasible upper non-transversal patterns such that the number of vertices contained in non-transversal blocks is equal to $t$;
			\item for each $k_1, k_2 \in \nonnegints$, establishing recurrence relations for the number of planar ways to connect $k_1$ upper vertices to $k_2$ lower vertices using only transversal blocks $b$ such that $\congmod{\noupverts{b}}{\nolowverts{b}}{m}$; and
			\item establish that the cardinality of the planar mod-$m$ monoid $\pmodmon{m}{k}$ can be calculated directly using the numbers generated from the recurrence relations above in (i) and (ii).
		\end{enumerate}
	
		\begin{definition} \label{def:pfeastrans}
			For each $m \in \posints$ and $k_1, k_2 \in \nonnegints$, by $\pfeastrans{m}{k_1}{k_2}$ we denote the number of planar ways to connect $k_1$ upper vertices to $k_2$ lower vertices using only transversal blocks $b$ such that $\congmod{\noupverts{b}}{\nolowverts{b}}{m}$ (see Figure \ref{fig:pfeastrans(3,6,3)} for an example).
		\end{definition}
		
		\begin{figure}[!ht]
			\caption[ ]{Given $m = 3$, $k_1 = 6$ and $k_2 = 3$, it may be verified with an exhaustive search that $\pfeastrans{3}{6}{3} = 8$, that is there are eight ways to connect six upper vertices to three lower vertices using only transversal blocks $b$ such that $\congmod{\noupverts{b}}{\nolowverts{b}}{3}$, which may be depicted as follows:}
			\label{fig:pfeastrans(3,6,3)}
			\vspace{5pt}
			\centering
			\input{chap_cardinalities/tikz/fig-pfeastransLB3-6-3RB.tex}
		\end{figure}
		
		\begin{proposition} \label{prop:pfeastransrecurrences}
			For each $m \in \posints$ and $k_1, k_2 \in \nonnegints$, $\pfeastrans{m}{k_1}{k_2}$ satisfies the recurrence (see Tables \ref{table:PT(2,k1,k2)}, \ref{table:PT(3,k1,k2)} and \ref{table:PT(4,k1,k2)} for example values): \[\pfeastrans{m}{k_1}{k_2} = \begin{cases}
				1 & \text{ if } k_1 = k_2 = 0 \\
				0 & \text{ if } \substack{k_1 = 0, k_2 > 0 \text{ or } k_1 > 0, k_2 = 0 \\ \text{ or } \ncongmod{k_1}{k_2}{m}} \\
				\sum_{\substack{(k_1', k_2') \leq (k_1, k_2)\\ \congmod{k_1'}{k_2'}{m}}} \pfeastrans{m}{k_1-k_1'}{k_2-k_2'} & \text{ if } \congmod{k_1}{k_2}{m}. 
			\end{cases}\]
			
			\begin{proof}
			 	Clearly we may not connect a positive number of vertices to zero vertices using transversal blocks, nor given $\ncongmod{k_1}{k_2}{m}$ may we connect $k_1$ upper vertices with $k_2$ lower vertices using only transversal blocks $b$ such that $\congmod{\noupverts{b}}{\nolowverts{b}}{m}$. Suppose $\congmod{k_1}{k_2}{m}$. If the transversal block containing the upper-right vertex $k_1$ connects the $k_1'$ upper right-most vertices and $k_2'$ lower right-most vertices such that $k_1' \in \{1, \ldots, k_1\}$ and $k_2' \in \{1, \ldots, k_2\}$, then there are $\pfeastrans{m}{k_1-k_1'}{k_2-k_2'}$ ways in which the remaining $k_1-k_1'$ upper vertices and $k_2-k_2'$ lower vertices may be connected in a planar fashion using only transversal blocks $b$ such that $\congmod{\noupverts{b}}{\nolowverts{b}}{m}$. Since there is one way for a single transversal to connect the $k_1$ upper vertices and $k_2$ lower vertices together it is convenient to have $\pfeastrans{m}{0}{0} = 1$.
			\end{proof}
		\end{proposition}
		
		\begin{table}[!ht]
			\caption[ ]{Example values for $\pfeastrans{2}{k_1}{k_2}$ calculated using the recurrences in Proposition \ref{prop:pfeastransrecurrences}.}
			\label{table:PT(2,k1,k2)}
			\centering
			\small
			\begin{tabular}{| c | r r r r r r r r r r r r r |}
				\hline
				\diagbox{$k_1$}{$k_2$} & 0 & 1 & 2 & 3 & 4 & 5 & 6 & 7 & 8 & 9 & 10 & 11 & 12 \\
				\hline
				0  & 1 &   &   &    &    &     &     &      &      &      &       &       &       \\
				1  &   & 1 &   & 1  &    & 1   &     & 1    &      & 1    &       & 1     &       \\
				2  &   &   & 2 &    & 3  &     & 4   &      & 5    &      & 6     &       & 7     \\
				3  &   & 1 &   & 4  &    & 8   &     & 13   &      & 19   &       & 26    &       \\
				4  &   &   & 3 &    & 10 &     & 22  &      & 40   &      & 65    &       & 98    \\
				5  &   & 1 &   & 8  &    & 26  &     & 61   &      & 120  &       & 211   &       \\
				6  &   &   & 4 &    & 22 &     & 70  &      & 171  &      & 356   &       & 665   \\
				7  &   & 1 &   & 13 &    & 61  &     & 192  &      & 483  &       & 1050  &       \\
				8  &   &   & 5 &    & 40 &     & 171 &      & 534  &      & 1373  &       & 3088  \\
				9  &   & 1 &   & 19 &    & 120 &     & 483  &      & 1500 &       & 3923  &       \\
				10 &   &   & 6 &    & 65 &     & 356 &      & 1373 &      & 4246  &       & 11257 \\
				11 &   & 1 &   & 26 &    & 211 &     & 1050 &      & 3923 &       & 12092 &       \\
				12 &   &   & 7 &    & 98 &     & 665 &      & 3088 &      & 11257 &       & 34606 \\
				\hline
			\end{tabular}
		\end{table}

		\begin{table}[!ht]
			\caption[ ]{Example values for $\pfeastrans{3}{k_1}{k_2}$ calculated using the recurrences in Proposition \ref{prop:pfeastransrecurrences}.}
			\label{table:PT(3,k1,k2)}
			\centering
			\small
			\begin{tabular}{| c | r r r r r r r r r r r r r |}
				\hline
				\diagbox{$k_1$}{$k_2$} & 0 & 1 & 2 & 3 & 4 & 5 & 6 & 7 & 8 & 9 & 10 & 11 & 12 \\
				\hline
				0  & 1 &   &   &    &    &     &     &     &     &      &      &      &      \\
				1  &   & 1 &   &    & 1  &     &     & 1   &     &      & 1    &      &      \\
				2  &   &   & 2 &    &    & 3   &     &     & 4   &      &      & 5    &      \\
				3  &   &   &   & 4  &    &     & 8   &     &     & 13   &      &      & 19   \\
				4  &   & 1 &   &    & 8  &     &     & 20  &     &      & 38   &      &      \\
				5  &   &   & 3 &    &    & 18  &     &     & 50  &      &      & 106  &      \\
				6  &   &   &   & 8  &    &     & 42  &     &     & 125  &      &      & 288  \\
				7  &   & 1 &   &    & 20 &     &     & 100 &     &      & 313  &      &      \\
				8  &   &   & 4 &    &    & 50  &     &     & 242 &      &      & 786  &      \\
				9  &   &   &   & 13 &    &     & 125 &     &     & 592  &      &      & 1979 \\
				10 &   & 1 &   &    & 38 &     &     & 313 &     &      & 1460 &      &      \\
				11 &   &   & 5 &    &    & 106 &     &     & 786 &      &      & 3624 &      \\
				12 &   &   &   & 19 &    &     & 288 &     &     & 1979 &      &      & 9042 \\
				\hline
			\end{tabular}
		\end{table}
		
		\begin{table}[!ht]
			\caption[ ]{Example values for $\pfeastrans{4}{k_1}{k_2}$ calculated using the recurrences in Proposition \ref{prop:pfeastransrecurrences}.}
			\label{table:PT(4,k1,k2)}
			\centering
			
			\begin{tabular}{| c | r r r r r r r r r r r r r |}
				\hline
				\diagbox{$k_1$}{$k_2$} & 0 & 1 & 2 & 3 & 4 & 5 & 6 & 7 & 8 & 9 & 10 & 11 & 12 \\
				\hline
				0 & 1 &  &  &  &  &  &  &  &  &  &  &  &  \\
				1 &  & 1 &  &  &  & 1 &  &  &  & 1 &  &  &  \\
				2 &  &  & 2 &  &  &  & 3 &  &  &  & 4 &  &  \\
				3 &  &  &  & 4 &  &  &  & 8 &  &  &  & 13 &  \\
				4 &  &  &  &  & 8 &  &  &  & 20 &  &  &  & 38 \\
				5 &  & 1 &  &  &  & 16 &  &  &  & 48 &  &  &  \\
				6 &  &  & 3 &  &  &  & 34 &  &  &  & 114 &  &  \\
				7 &  &  &  & 8 &  &  &  & 74 &  &  &  & 269 &  \\
				8 &  &  &  &  & 20 &  &  &  & 164 &  &  &  & 633 \\
				9 &  & 1 &  &  &  & 48 &  &  &  & 368 &  &  &  \\
				10 &  &  & 4 &  &  &  & 114 &  &  &  & 834 &  &  \\
				11 &  &  &  & 13 &  &  &  & 269 &  &  &  & 1904 &  \\
				12 &  &  &  &  & 38 &  &  &  & 633 &  &  &  & 4372 \\
				\hline
			\end{tabular}
		\end{table}
		
		\begin{definition}
			For each $m \in \posints$, $k \in \intsge{m}$, $x \in \floor{\frac{k}{m}}$ and $t_1, \ldots, t_x \in \nonnegints$ such that $\Sigma_{i=1}^ximt_i \leq k$, we denote by $\pfeasnontrans{m, k, t_1, \ldots, t_x}$ the number of $\pmodmon{m}{k}$-feasible upper non-transversal patterns containing precisely $t_1$ type $(m, 0)$, \ldots, and $t_x$ type $(xm, 0)$ non-transversal blocks (see Figure \ref{fig:pfeasnontrans(2,6,1,1)} for an example).
		\end{definition}
		
		\begin{figure}[!ht]
			\caption[ ]{Given $m = 2$, $k = 6$, and $t_1 = t_2 = 1$, it may be verified with an exhaustive search that $\pfeasnontrans{2, 6, 1, 1} = 6$, that is there are six $\pmodmon{2}{6}$-feasible upper non-transversal patterns containing precisely $1$ type $(2, 0)$ block and $1$ type $(4,0)$ block, which may be depicted as follows:}
			\label{fig:pfeasnontrans(2,6,1,1)}
			\vspace{5pt}
			\centering
			\input{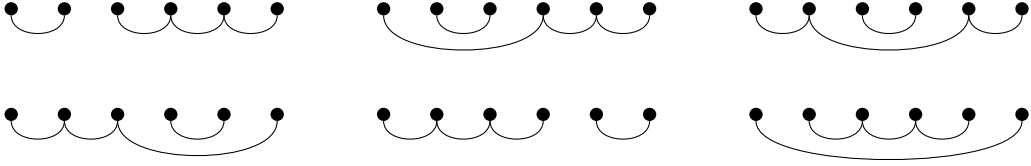}
		\end{figure}
        
        Note it trivially follows from the planar mod-$m$ monoid $\pmodmon{m}{k}$ being closed under the vertical flip involution $^*$, that is ${\pmodmon{m}{k}}^* = \pmodmon{m}{k}$, that the number of $\pmodmon{m}{k}$-feasible lower non-transversal patterns containing precisely $t_1$ type $(0, m)$, $\ldots$, and $t_x$ type $(0, xm)$ non-transversal blocks is equal to the number of $\pmodmon{m}{k}$-feasible upper non-transversal patterns containing precisely $t_1$ type $(m, 0)$, \ldots, and $t_x$ type $(xm, 0)$ non-transversal blocks.

        \begin{proposition} \label{prop:pfeasnontransrecurrences}
        	For each $m \in \posints$, $k \in \intsge{m}$, $x \in \floor{\frac{k}{m}}$ and $t_1, \ldots, t_x \in \nonnegints$ such that $\Sigma_{i=1}^ximt_i \leq k$, $\pfeasnontrans{m, k, t_1, \ldots, t_x}$ satisfies the recurrence (see Tables \ref{table:PN(2,k,t)}, \ref{table:PN(3,k,t)} and \ref{table:PN(4,k,t)} for example values):
            \label{prop:noNonTransRecurrencesPlanarPttn}
            \begin{align*}
                \pfeasnontrans{m, k, 0, \ldots, 0} =&\ 1; \\
                \pfeasnontrans{m, k, t_1, \ldots, t_x} =&\ \pfeasnontrans{m, k-1, t_1, \ldots, t_x} + \sum_{\substack{i \in \{1, \ldots, x\} \\ t_i > 0}}\pfeasnontrans{m, k-1, t_1, \ldots, t_i-1, \ldots, t_x}.
            \end{align*}

            \begin{proof}
                There is only one way to place no non-transversal blocks, giving us $\pfeasnontrans{m, k, 0, \ldots, 0} = 1$. Suppose there is at least one non-transversal block to be placed. There are $\pfeasnontrans{m, k-1, t_1, \ldots, t_x}$ $\pmodmon{m}{k}$-feasible upper non-transversal patterns that contain precisely $t_1$ type $(m, 0)$, \ldots, and $t_x$ type $(xm, 0)$ non-transversal blocks along the $k-1$ left-most vertices, that is so that no non-transversal blocks contain the right-most vertex $k$. Let $i \in \{1, \ldots, x\}$. If a type $(im, 0)$ non-transversal block $a$ contains the right-most vertex $k$, then there are $\pfeasnontrans{m, k-1, t_1, \ldots, t_i-1, \ldots, t_x}$ $\pmodmon{m}{k}$-feasible upper non-transversal patterns containing precisely $t_1$ type $(m, 0)$, \ldots, $t_i-1$ type $(im, 0)$, \ldots, and $t_x$ type $(xm, 0)$ non-transversal blocks along the $k-1$ left-most vertices, after which the type $(im, 0)$ block $a$ is placed along the $im$ right-most unused vertices. 
            \end{proof}
        \end{proposition}
        
        \begin{table}[!ht]
            \caption[ ]{Example values for $\pfeasnontrans{2, k, \vec{t}}$ where $\vec{t} = (t_1, t_2, t_3, t_4)$, calculated using the recurrences in Proposition \ref{prop:pfeasnontransrecurrences}.}
            \label{table:PN(2,k,t)}
            \centering 
            \begin{tabular}{| c | r r r r r r r r r r r r |}
                \hline
                \diagbox{$k$}{$\vec{t}$} & \rotatebox{90}{(0,0,0,0) } & \rotatebox{90}{(1,0,0,0)} & \rotatebox{90}{(2,0,0,0)} & \rotatebox{90}{(3,0,0,0)} & \rotatebox{90}{(4,0,0,0)} & \rotatebox{90}{(0,1,0,0)} & \rotatebox{90}{(1,1,0,0)} & \rotatebox{90}{(2,1,0,0)} & \rotatebox{90}{(0,2,0,0)} & \rotatebox{90}{(0,0,1,0)} & \rotatebox{90}{(1,0,1,0)} & \rotatebox{90}{(0,0,0,1)} \\
                \hline
				0  & 1 &    &     &     &      &    &     &      &     &    &     &    \\
				1  & 1 &    &     &     &      &    &     &      &     &    &     &    \\
				2  & 1 & 1  &     &     &      &    &     &      &     &    &     &    \\
				3  & 1 & 2  &     &     &      &    &     &      &     &    &     &    \\
				4  & 1 & 3  & 2   &     &      & 1  &     &      &     &    &     &    \\
				5  & 1 & 4  & 5   &     &      & 2  &     &      &     &    &     &    \\
				6  & 1 & 5  & 9   & 5   &      & 3  & 6   &      &     & 1  &     &    \\
				7  & 1 & 6  & 14  & 14  &      & 4  & 14  &      &     & 2  &     &    \\
				8  & 1 & 7  & 20  & 28  & 14   & 5  & 24  & 28   & 4   & 3  & 8   & 1  \\
				9  & 1 & 8  & 27  & 48  & 42   & 6  & 36  & 72   & 9   & 4  & 18  & 2  \\
				10 & 1 & 9  & 35  & 75  & 90   & 7  & 50  & 135  & 15  & 5  & 30  & 3  \\
				11 & 1 & 10 & 44  & 110 & 165  & 8  & 66  & 220  & 22  & 6  & 44  & 4  \\
				12 & 1 & 11 & 54  & 154 & 275  & 9  & 84  & 330  & 30  & 7  & 60  & 5  \\
				13 & 1 & 12 & 65  & 208 & 429  & 10 & 104 & 468  & 39  & 8  & 78  & 6  \\
				14 & 1 & 13 & 77  & 273 & 637  & 11 & 126 & 637  & 49  & 9  & 98  & 7  \\
				15 & 1 & 14 & 90  & 350 & 910  & 12 & 150 & 840  & 60  & 10 & 120 & 8  \\
				16 & 1 & 15 & 104 & 440 & 1260 & 13 & 176 & 1080 & 72  & 11 & 144 & 9  \\
                \hline
            \end{tabular}
        \end{table}
        
        \begin{table}[!ht]
            \caption[ ]{Example values for $\pfeasnontrans{3, k, \vec{t}}$ where $\vec{t} = (t_1, t_2, t_3, t_4)$, calculated using the recurrences in Proposition \ref{prop:pfeasnontransrecurrences}.}
            \label{table:PN(3,k,t)}
            \centering
            \begin{tabular}{| c | r r r r r r r r r r r r |}
                \hline
                \diagbox{$k$}{$\vec{t}$} & \rotatebox{90}{(0,0,0,0) } & \rotatebox{90}{(1,0,0,0)} & \rotatebox{90}{(2,0,0,0)} & \rotatebox{90}{(3,0,0,0)} & \rotatebox{90}{(4,0,0,0)} & \rotatebox{90}{(0,1,0,0)} & \rotatebox{90}{(1,1,0,0)} & \rotatebox{90}{(2,1,0,0)} & \rotatebox{90}{(0,2,0,0)} & \rotatebox{90}{(0,0,1,0)} & \rotatebox{90}{(1,0,1,0)} & \rotatebox{90}{(0,0,0,1)} \\
                \hline
                0  & 1 &    &    &     &     &    &     &     &    &   &    &   \\
                1  & 1 &    &    &     &     &    &     &     &    &   &    &   \\
                2  & 1 &    &    &     &     &    &     &     &    &   &    &   \\
                3  & 1 & 1  &    &     &     &    &     &     &    &   &    &   \\
                4  & 1 & 2  &    &     &     &    &     &     &    &   &    &   \\
                5  & 1 & 3  &    &     &     &    &     &     &    &   &    &   \\
                6  & 1 & 4  & 3  &     &     & 1  &     &     &    &   &    &   \\
                7  & 1 & 5  & 7  &     &     & 2  &     &     &    &   &    &   \\
                8  & 1 & 6  & 12 &     &     & 3  &     &     &    &   &    &   \\
                9  & 1 & 7  & 18 & 12  &     & 4  & 9   &     &    & 1 &    &   \\
                10 & 1 & 8  & 25 & 30  &     & 5  & 20  &     &    & 2 &    &   \\
                11 & 1 & 9  & 33 & 55  &     & 6  & 33  &     &    & 3 &    &   \\
                12 & 1 & 10 & 42 & 88  & 55  & 7  & 48  & 66  & 6  & 4 & 12 & 1 \\
                13 & 1 & 11 & 52 & 130 & 143 & 8  & 65  & 156 & 13 & 5 & 26 & 2 \\
                14 & 1 & 12 & 63 & 182 & 273 & 9  & 84  & 273 & 21 & 6 & 42 & 3 \\
                15 & 1 & 13 & 75 & 245 & 455 & 10 & 105 & 420 & 30 & 7 & 60 & 4 \\
                \hline
            \end{tabular}
        \end{table}
        
        \begin{table}[!ht]
            \caption[ ]{Example values for $\pfeasnontrans{4, k, \vec{t}}$ where $\vec{t} = (t_1, t_2, t_3, t_4)$, calculated using the recurrences in Proposition \ref{prop:pfeasnontransrecurrences}.}
            \label{table:PN(4,k,t)}
            \centering
            \begin{tabular}{| c | r r r r r r r r r r r r |}
                \hline
                \diagbox{$k$}{$\vec{t}$} & \rotatebox{90}{(0,0,0,0) } & \rotatebox{90}{(1,0,0,0)} & \rotatebox{90}{(2,0,0,0)} & \rotatebox{90}{(3,0,0,0)} & \rotatebox{90}{(4,0,0,0)} & \rotatebox{90}{(0,1,0,0)} & \rotatebox{90}{(1,1,0,0)} & \rotatebox{90}{(2,1,0,0)} & \rotatebox{90}{(0,2,0,0)} & \rotatebox{90}{(0,0,1,0)} & \rotatebox{90}{(1,0,1,0)} & \rotatebox{90}{(0,0,0,1)} \\
                \hline
				0  & 1 &    &     &     &      &    &     &     &    &   &     & \\
				1  & 1 &    &     &     &      &    &     &     &    &   &     & \\
				2  & 1 &    &     &     &      &    &     &     &    &   &     & \\
				3  & 1 &    &     &     &      &    &     &     &    &   &     & \\
				4  & 1 & 1  &     &     &      &    &     &     &    &   &     & \\
				5  & 1 & 2  &     &     &      &    &     &     &    &   &     & \\
				6  & 1 & 3  &     &     &      &    &     &     &    &   &     & \\
				7  & 1 & 4  &     &     &      &    &     &     &    &   &     & \\
				8  & 1 & 5  & 4   &     &      & 1  &     &     &    &   &     & \\
				9  & 1 & 6  & 9   &     &      & 2  &     &     &    &   &     & \\
				10 & 1 & 7  & 15  &     &      & 3  &     &     &    &   &     & \\
				11 & 1 & 8  & 22  &     &      & 4  &     &     &    &   &     & \\
				12 & 1 & 9  & 30  & 22  &      & 5  & 12  &     &    & 1 &     &\\
				13 & 1 & 10 & 39  & 52  &      & 6  & 26  &     &    & 2 &     &\\
				14 & 1 & 11 & 49  & 91  &      & 7  & 42  &     &    & 3 &     &\\
				15 & 1 & 12 & 60  & 140 &      & 8  & 60  &     &    & 4 &     &\\
				16 & 1 & 13 & 72  & 200 & 140  & 9  & 80  & 120 & 8  & 5 & 16  & 1 \\
				17 & 1 & 14 & 85  & 272 & 340  & 10 & 102 & 272 & 17 & 6 & 34  & 2 \\
				18 & 1 & 15 & 99  & 357 & 612  & 11 & 126 & 459 & 27 & 7 & 54  & 3 \\
				19 & 1 & 16 & 114 & 456 & 969  & 12 & 152 & 684 & 38 & 8 & 76  & 4 \\
				20 & 1 & 17 & 130 & 570 & 1425 & 13 & 180 & 950 & 50 & 9 & 100 & 5 \\
                \hline
            \end{tabular}
        \end{table}
        
        \begin{definition} \label{def:pmodmonfeasnontrans}
        	For each $m \in \posints$, $k \in \intsge{m}$ and $t \in m\nonnegints$ such that $t \leq k$, we denote by $\pfeasnontrans{m, k, t}$ the number of $\pmodmon{m}{k}$-feasible upper non-transversal patterns that have $t$ of $k$ vertices contained in non-transversal blocks, that is $\pfeasnontrans{m, k, t} = \Sigma_{\substack{t_1, \ldots, t_x \in \nonnegints \\ t_1 + \ldots + xmt_x = t}}\ \pfeasnontrans{m, k, t_1, \ldots, t_x}$ (see Tables \ref{table:pmod2moncards}, \ref{table:pmod3moncards} and \ref{table:pmod4moncards} for example values). 
        \end{definition}

        \begin{theorem} \label{thm:pmodmoncardinality}
            For each $m \in \posints$ and $k \in \intsge{m}$, the cardinality of $\pmodmon{m}{k}$ is given by (see Tables \ref{table:pmod2moncards}, \ref{table:pmod3moncards} and \ref{table:pmod4moncards} for example values),
            \[ \card{\pmodmon{m}{k}} = \sum_{u=0}^{\floor{k/m}}\sum_{l=0}^{\floor{k/m}} \pfeasnontrans{m, k, mu}\pfeasnontrans{m, k, ml}\pfeastrans{m}{k-mu}{k-ml}. \]
            \begin{proof}
				Suppose the number of upper vertices and number of lower vertices contained in non-transversal blocks are $mu$ and $ml$ respectively, where $u, l \in \set{0, \ldots, \floor{\frac{k}{m}}}$. For each of the $\pfeasnontrans{m, k, mu}$ $\pmodmon{m}{k}$-feasible upper non-transversal patterns that have $mu$ of $k$ upper vertices contained in non-transversal blocks and $\pfeasnontrans{m, k, ml}$ $\pmodmon{m}{k}$-feasible lower non-transversal patterns that have $ml$ of $k$ lower vertices contained in non-transversal blocks, there are $\pfeastrans{m}{k-mu}{k-ml}$ planar ways to connect the remaining $k-mu$ upper vertices to the remaining $k-ml$ lower vertices using only transversal blocks $b$ such that $\congmod{\noupverts{b}}{\nolowverts{b}}{m}$.
            \end{proof}
        \end{theorem}
        
        One of the author's anonymous examiners kindly pointed out, in reference to Theorem \ref{thm:pmodmoncardinality}, that \cite{coll:Temperley:NoteOnBaxtersGeneralizationofTLOps} contains an interesting discussion of the $m = 2$ case, with considerable simplification of the formula.
        
        \begin{table}[!ht]
            \caption[ ]{Example values for $\pfeasnontrans{2, k, t}$ and $\card{\pmodmon{2}{k}}$, calculated using the formulas from Definition \ref{def:pmodmonfeasnontrans} and Theorem \ref{thm:pmodmoncardinality}.}
            \label{table:pmod2moncards}
            \centering
            \begin{tabular}{|c|rrrrr|r|}
                \hline
                \diagbox{k}{t} & 0 & 2 & 4 & 6 & 8 & {\normalsize$\card{\pmodmon{2}{k}}$} \\
                \hline
				2 & 1 & 1 &   &   &   & 3 \\
				3 & 1 & 2 &   &   &   & 12 \\
				4 & 1 & 3 & 3 &   &   & 55 \\
				5 & 1 & 4 & 7 &   &   & 273 \\
				6 & 1 & 5 & 12 & 12 &   & 1428 \\
				7 & 1 & 6 & 18 & 30 &   & 7752 \\
				8 & 1 & 7 & 25 & 55 & 55 & 43263 \\
				9 & 1 & 8 & 33 & 88 & 143 & 246675 \\
				10 & 1 & 9 & 42 & 130 & 273 & 1430715 \\
				11 & 1 & 10 & 52 & 182 & 455 & 8414640 \\
				12 & 1 & 11 & 63 & 245 & 700 & 50067108 \\
				13 & 1 & 12 & 75 & 320 & 1020 & 300830572 \\
				14 & 1 & 13 & 88 & 408 & 1428 & 1822766520 \\
				15 & 1 & 14 & 102 & 510 & 1938 & 11124755664 \\
				16 & 1 & 15 & 117 & 627 & 2565 & 68328754959 \\
				17 & 1 & 16 & 133 & 760 & 3325 & 422030545335 \\
				18 & 1 & 17 & 150 & 910 & 4235 & 2619631042665 \\
				19 & 1 & 18 & 168 & 1078 & 5313 & 16332922290300 \\
				20 & 1 & 19 & 187 & 1265 & 6578 & 102240109897695 \\
                \hline
            \end{tabular}
        \end{table}
 
        \begin{table}[!ht]
            \caption[ ]{Example values for $\pfeasnontrans{3, k, t}$ and $\card{\pmodmon{3}{k}}$, calculated using the formulas from Definition \ref{def:pmodmonfeasnontrans} and Theorem \ref{thm:pmodmoncardinality}.}
            \label{table:pmod3moncards}
            \centering
            \begin{tabular}{|c|rrrrr|r|}
                \hline
                \diagbox{k}{t} & 0 & 3 & 6 & 9 & 12 & {\normalsize$\card{\pmodmon{3}{k}}$} \\
                \hline
				3 & 1 & 1 &   &   &   & 5 \\
				4 & 1 & 2 &   &   &   & 16 \\
				5 & 1 & 3 &   &   &   & 54 \\
				6 & 1 & 4 & 4 &   &   & 186 \\
				7 & 1 & 5 & 9 &   &   & 689 \\
				8 & 1 & 6 & 15 &   &   & 2600 \\
				9 & 1 & 7 & 22 & 22 &   & 9856 \\
				10 & 1 & 8 & 30 & 52 &   & 38708 \\
				11 & 1 & 9 & 39 & 91 &   & 153438 \\
				12 & 1 & 10 & 49 & 140 & 140 & 608868 \\
				13 & 1 & 11 & 60 & 200 & 340 & 2467726 \\
				14 & 1 & 12 & 72 & 272 & 612 & 10057082 \\
				15 & 1 & 13 & 85 & 357 & 969 & 40986557 \\
                \hline
            \end{tabular}
        \end{table}
        
        \begin{table}[!ht]
            \caption[ ]{Example values for $\pfeasnontrans{4, k, t}$ and $\card{\pmodmon{4}{k}}$, calculated using the formulas from Definition \ref{def:pmodmonfeasnontrans} and Theorem \ref{thm:pmodmoncardinality}.}
            \label{table:pmod4moncards}
            \centering
            \begin{tabular}{|c|rrrrr|r|}
                \hline
                \diagbox{k}{t} & 0 & 4 & 8 & 12 & 16 & {\normalsize$\card{\pmodmon{4}{k}}$} \\
                \hline
				4 & 1 & 1 &   &   &   & 9 \\
				5 & 1 & 2 &   &   &   & 24 \\
				6 & 1 & 3 &   &   &   & 70 \\
				7 & 1 & 4 &   &   &   & 202 \\
				8 & 1 & 5 & 5 &   &   & 589 \\
				9 & 1 & 6 & 11 &   &   & 1795 \\
				10 & 1 & 7 & 18 &   &   & 5644 \\
				11 & 1 & 8 & 26 &   &   & 17652 \\
				12 & 1 & 9 & 35 & 35 &   & 55335 \\
				13 & 1 & 10 & 45 & 80 &   & 176966 \\
				14 & 1 & 11 & 56 & 136 &   & 575004 \\
				15 & 1 & 12 & 68 & 204 &   & 1862638 \\
				16 & 1 & 13 & 81 & 285 & 285 & 6037799 \\
				17 & 1 & 14 & 95 & 380 & 665 & 19793749 \\
				18 & 1 & 15 & 110 & 490 & 1155 & 65511224 \\
				19 & 1 & 16 & 126 & 616 & 1771 & 216404828 \\
				20 & 1 & 17 & 143 & 759 & 2530 & 715006656 \\
                \hline
            \end{tabular}
        \end{table}
        
        \subsection{The planar modular-$2$ partition monoid, $\card{\pmodmon{2}{k}}$}
        	The calculated values in the right-most column of Table \ref{table:pmod2moncards}, which form the start of the sequence $\set{\card{\pmodmon{2}{k}}: k \in \intsge{2}}$, match the formula $\frac{\binom{3k}{k}}{2k + 1}$ which is listed on the OEIS as sequence $A001764$. Furthermore $\frac{\binom{3k}{k}}{2k + 1}$ enumerates the non-crossing partitions of $\{1, \ldots, 2k\}$ with blocks of even size, which we proceed to use in order to establish that the sequence $\set{\card{\pmodmon{2}{k}}: k \in \intsge{2}}$ matches sequence $A001764$ on the OEIS.
        	
        	\begin{proposition}
	        	The cardinality of the planar mod-$2$ monoid $\pmodmon{2}{k}$ is equal to the cardinality of non-crossing partitions of $\{1, \ldots, 2k\}$ with blocks of even size.
	        	
	        	\begin{proof}
		        	Given a bipartition $\alpha \in \pmodmon{2}{k}$, each block $b \in \alpha$ satisfies $\noupverts{b} + \nolowverts{b} \in 2\ints$, and given a planar partition $\alpha$ of $\{1, \ldots, 2k\}$ with blocks of even size, each block $b \in \alpha$ satisfies $\congmod{\card{b \cap \{1, \ldots, k\}}}{\card{b \cap \{k+1, \ldots, 2k\}}}{2}$.
	        	\end{proof}
        	\end{proposition}
        	
        	\begin{corollary}
	        	For each $k \in \intsge{2}$, $\card{\pmodmon{2}{k}} = \frac{\binom{3k}{k}}{2k + 1}$.
        	\end{corollary}
        	
\begin{figure}[!ht]
				\caption[ ]{The twelve ternary trees with three internal nodes:}
				\label{fig:ternarytrees3internalnodes}
				\vspace{5pt}
				\centering
				\input{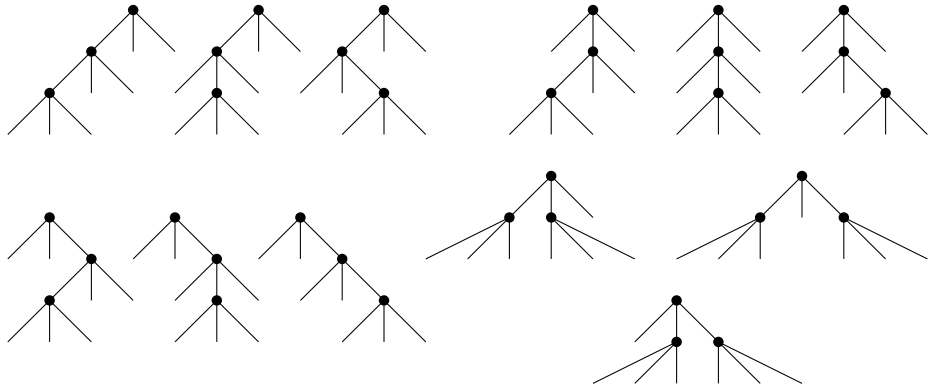}
			\end{figure}

\begin{figure}[!ht]
				\caption[ ]{The twelve rooted plane trees with six edges such that every vertex has an even out degree:}
				\label{fig:eventrees6edges}
				\vspace{5pt}
				\centering
				\input{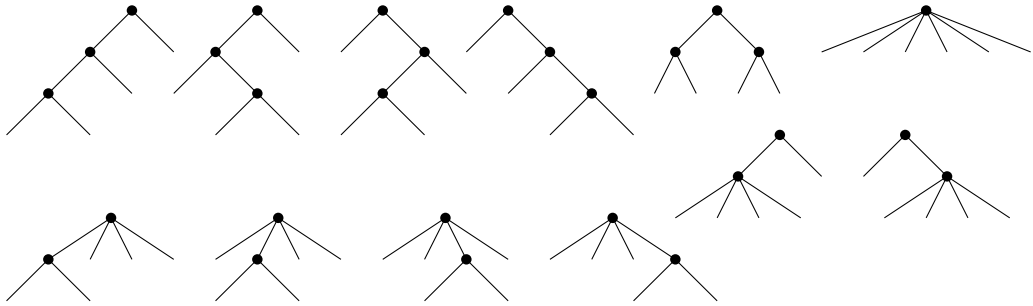}
			\end{figure}

        	It follows from the comments section for sequence $A001764$ on the OEIS that the cardinality of the planar mod-$2$ monoid $\card{\pmodmon{2}{k}}$ is equal to:
        	\begin{enumerate}
        		\item the Pfaff-Fuss-Catalan sequence $C^3_n$;
        		\item the number of lattice paths from $(0,0)$ to $(n,2n)$ that do not step above the line $y=2x$ using precisely $n$ east steps and $2n$ north steps (see Figure \ref{fig:latticepathsfrom(0,0)to(k,2k)} for an example);
        		\item the number of lattice paths from $(0,0)$ to $(2n,0)$ that do not step below the $x$-axis using the step-set $\set{(1,1), (0,-2)}$ (see Figure \ref{fig:latticepathsfrom(0,0)to(2k,0)withdownstep(0,-2)} for an example);
        		\item the number of complete ternary trees with $n$ internal nodes, or $3n$ edges (see Figure \ref{fig:ternarytrees3internalnodes} for an example); and
        		\item the number of rooted plane trees with $2n$ edges such that every vertex has even outdegree, also commonly referred to as \textit{even trees} (see Figure \ref{fig:eventrees6edges} for an example).
        	\end{enumerate}
			
			\begin{figure}[!ht]
				\caption[ ]{The twelve lattice paths from $(0,0)$ to $(3,6)$ that do not step above the line $y=2x$ using precisely $3$ east steps and $6$ north steps:}
				\label{fig:latticepathsfrom(0,0)to(k,2k)}
				\vspace{5pt}
				\centering
				\input{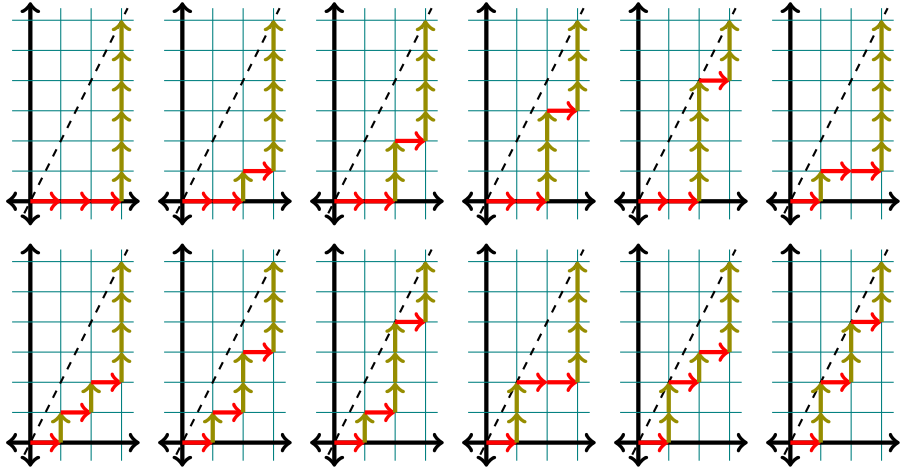}
			\end{figure}
			
			\begin{figure}[!ht]
				\caption[ ]{The twelve lattice paths from $(0,0)$ to $(6,0)$ that do not step below the $x$-axis using the step-set $\set{(1,1), (0,-2)}$:}
				\label{fig:latticepathsfrom(0,0)to(2k,0)withdownstep(0,-2)}
				\vspace{5pt}
				\centering
				\input{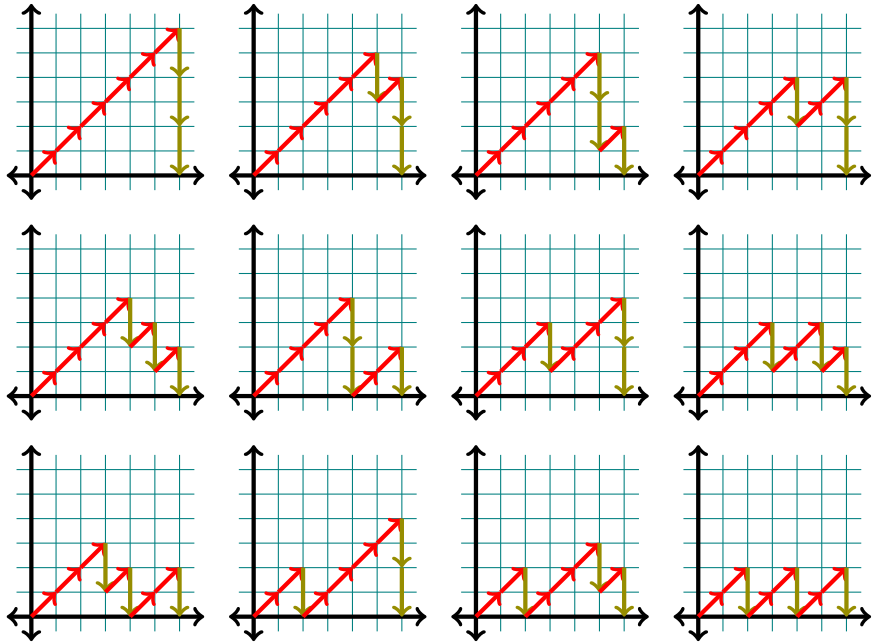}
			\end{figure}
			
    \section{The $m$-apsis generated diagram monoid, $\card{\apsismon{m}{k}}$}
    	Recall from Theorem \ref{thm:apsismon=apsisbound} that the $m$-apsis monoid $\apsismon{m}{k}$ consists of all bipartitions $\alpha \in \pmodmon{m}{k}$ such that either:
		\begin{enumerate}
			\item $\alpha$ is the identity $\id{k}$; or 
			\item $\alpha$ contains at least one upper $m$-apsis and at least one lower $m$-apsis.
		\end{enumerate} 
    
		\begin{definition}
			For each $m \in \intsge{3}$, $k \in \intsge{m}$, $x \in \floor{\frac{k}{m}}$ and $t_1, \ldots, t_x \in \nonnegints$ such that $\Sigma_{i=1}^ximt_i \leq k$, we denote by $\pfeasapsisnontrans{m, k, t_1, \ldots, t_x}$ the number of $\apsismon{m}{k}$-feasible upper non-transversal patterns containing precisely $t_1$ type $(m, 0)$, \ldots, and $t_x$ type $(xm, 0)$ non-transversal blocks.
		\end{definition} 
		
		Note it trivially follows from the $m$-apsis monoid $\apsismon{m}{k}$ being closed under the vertical flip involution $^*$, that is ${\apsismon{m}{k}}^* = \apsismon{m}{k}$, that the number of $\apsismon{m}{k}$-feasible lower non-transversal patterns containing precisely $t_1$ type $(0, m)$, \ldots, and $t_x$ type $(0, xm)$ non-transversal blocks is equal to the number of $\apsismon{m}{k}$-feasible upper non-transversal patterns containing precisely $t_1$ type $(m, 0)$, \ldots, and $t_x$ type $(xm, 0)$ non-transversal blocks.

		Further note $\pfeasmodminusapsisnontrans{m, k, t_1, \ldots, t_m}$ is the number of $\pmodmon{m}{k}$-feasible upper non-transversal patterns containing precisely $t_1$ type $(m, 0)$, \ldots, and $t_x$ type $(xm, 0)$ non-transversal blocks such that none of the $t_1$ type $(m, 0)$ non-transversal blocks is an $m$-apsis.

        \begin{proposition} \label{prop:pfeasapsisnontransrecurrences}
        	For each $m \in \intsge{3}$, $k \in \intsge{m}$, $x \in \floor{\frac{k}{m}}$ and $t_1, \ldots, t_x \in \nonnegints$ such that $\Sigma_{i=1}^ximt_i \leq k$, $\pfeasapsisnontrans{m, k, t_1, \ldots, t_x}$ satisfies the recurrence (see Tables \ref{table:PNB(3,k,t)} and \ref{table:PNB(4,k,t)} for example values):
            \begin{align*}
                \pfeasapsisnontrans{m, k, 0, \ldots, 0} =&\ 1; \\
                \pfeasapsisnontrans{m, k, 0, t_2, \ldots, t_x} =&\ 0; \\
                \pfeasapsisnontrans{m, k, t_1, \ldots, t_x} =&\ \pfeasapsisnontrans{m, k-1, t_1, \ldots, t_x} \\
                     &+ \sum_{\substack{i \in \{1, \ldots, x\} \\ t_i > 0}}\pfeasapsisnontrans{m, k-1, t_1, \ldots, t_i-1, \ldots, t_x} \\
                     &+ \pfeasmodminusapsisnontrans{m, k-m, t_1-1, \ldots, t_x}.
            \end{align*}
				
            \begin{proof}
            	There is only one way to place no non-transversal blocks, giving us $\pfeasapsisnontrans{m, k, 0, \ldots, 0} = 1$. Suppose there is at least one non-transversal block to be placed. If $t_1 = 0$ then we cannot $\apsismon{m}{k}$-feasibly place the $t_2$ type $(2m, 0)$, \ldots, and $t_x$ type $(xm, 0)$ non-transversal blocks, hence $\pfeasapsisnontrans{m, k, 0, t_2, \ldots, t_x} = 0$. If no non-transversal blocks contain the right-most vertex then there are $\pfeasapsisnontrans{k-1, t_1, \ldots, t_x}$ $\apsismon{m}{k}$-feasible upper non-transversal patterns containing precisely $t_1$ type $(m, 0)$, \ldots, and $t_x$ type $(xm, 0)$ non-transversal blocks along the $k-1$ left-most vertices, that is so that no non-transversal blocks contain the right-most vertex $k$.  Let $i \in \{1, \ldots, x\}$. If a type $(im, 0)$ non-transversal block $a$ contains the right-most vertex $k$, then there is at least $\pfeasapsisnontrans{m, k-1, t_1, \ldots, t_i-1, \ldots, t_m}$ $\apsismon{m}{k}$-feasible upper non-transversal patterns containing precisely $t_1$ type $(m, 0)$, \ldots, $t_i-1$ type $(im, 0)$, \ldots, and $t_x$ type $(xm, 0)$ non-transversal blocks along the $k-1$ left-most vertices, after which the type $(im, 0)$ block $a$ is placed along the $im$ right-most unused vertices. Note however if $a$ is a block of type $(m, 0)$ and contains the $m$ right-most vertices then we have only counted the cases with at least two $m$-apses. We must also add the $\pfeasmodminusapsisnontrans{m, k-m, t_1-1, \ldots, t_m}$ ways to $\ppttnmon{k}$-feasibly place the remaining $t_1-1$ type $(m, 0)$, \ldots, and $t_m$ type $(xm, 0)$ upper non-transversal blocks along the $k-m$ left-most vertices such that none of the $t_1-1$ blocks of type $(m, 0)$ are $m$-apses.
            \end{proof}
        \end{proposition}

        \begin{table}[!ht]
            \caption[ ]{Example values for $\pfeasapsisnontrans{3, k, \vec{t}}$ where $\vec{t} = (t_1, t_2, t_3, t_4)$, calculated using the recurrences from Proposition \ref{prop:pfeasapsisnontransrecurrences}.}
            \label{table:PNB(3,k,t)}
            \centering 
            \begin{tabular}{| c | r r r r r r r r r r r r |}
                \hline
                \diagbox{$k$}{$\vec{t}$} & \rotatebox{90}{(0,0,0,0) } & \rotatebox{90}{(1,0,0,0)} & \rotatebox{90}{(2,0,0,0)} & \rotatebox{90}{(3,0,0,0)} & \rotatebox{90}{(4,0,0,0)} & \rotatebox{90}{(0,1,0,0)} & \rotatebox{90}{(1,1,0,0)} & \rotatebox{90}{(2,1,0,0)} & \rotatebox{90}{(0,2,0,0)} & \rotatebox{90}{(0,0,1,0)} & \rotatebox{90}{(1,0,1,0)} & \rotatebox{90}{(0,0,0,1)} \\
                \hline
                0  & 1 &    &    &     &     &   &    &     &   &   &    &   \\
                1  & 1 &    &    &     &     &   &    &     &   &   &    &   \\
                2  & 1 &    &    &     &     &   &    &     &   &   &    &   \\
                3  & 1 & 1  &    &     &     &   &    &     &   &   &    &   \\
                4  & 1 & 2  &    &     &     &   &    &     &   &   &    &   \\
                5  & 1 & 3  &    &     &     &   &    &     &   &   &    &   \\
                6  & 1 & 4  & 3  &     &     & 0 &    &     &   &   &    &   \\
                7  & 1 & 5  & 7  &     &     & 0 &    &     &   &   &    &   \\
                8  & 1 & 6  & 12 &     &     & 0 &    &     &   &   &    &   \\
                9  & 1 & 7  & 18 & 12  &     & 0 & 7  &     &   & 0 &    &   \\
                10 & 1 & 8  & 25 & 30  &     & 0 & 16 &     &   & 0 &    &   \\
                11 & 1 & 9  & 33 & 55  &     & 0 & 27 &     &   & 0 &    &   \\
                12 & 1 & 10 & 42 & 88  & 55  & 0 & 40 & 62  & 0 & 0 & 10 & 0 \\
                13 & 1 & 11 & 52 & 130 & 143 & 0 & 55 & 148 & 0 & 0 & 22 & 0 \\
                14 & 1 & 12 & 63 & 182 & 273 & 0 & 72 & 261 & 0 & 0 & 36 & 0 \\
                15 & 1 & 13 & 75 & 245 & 455 & 0 & 91 & 404 & 0 & 0 & 52 & 0 \\
                \hline
            \end{tabular}
        \end{table}
        
        \begin{table}[!ht]
            \caption[ ]{Example values for $\pfeasapsisnontrans{4, k, \vec{t}}$ where $\vec{t} = (t_1, t_2, t_3, t_4)$, calculated using the recurrences from Proposition \ref{prop:pfeasapsisnontransrecurrences}.}
            \label{table:PNB(4,k,t)}
            \centering
            \begin{tabular}{| c | r r r r r r r r r r r r |}
                \hline
                \diagbox{$k$}{$\vec{t}$} & \rotatebox{90}{(0,0,0,0) } & \rotatebox{90}{(1,0,0,0)} & \rotatebox{90}{(2,0,0,0)} & \rotatebox{90}{(3,0,0,0)} & \rotatebox{90}{(4,0,0,0)} & \rotatebox{90}{(0,1,0,0)} & \rotatebox{90}{(1,1,0,0)} & \rotatebox{90}{(2,1,0,0)} & \rotatebox{90}{(0,2,0,0)} & \rotatebox{90}{(0,0,1,0)} & \rotatebox{90}{(1,0,1,0)} & \rotatebox{90}{(0,0,0,1)} \\
                \hline
				0  & 1 &    &     &     &     &   &     &     &   &   &    &   \\
				1  & 1 &    &     &     &     &   &     &     &   &   &    &   \\
				2  & 1 &    &     &     &     &   &     &     &   &   &    &   \\
				3  & 1 &    &     &     &     &   &     &     &   &   &    &   \\
				4  & 1 & 1  &     &     &     &   &     &     &   &   &    &   \\
				5  & 1 & 2  &     &     &     &   &     &     &   &   &    &   \\
				6  & 1 & 3  &     &     &     &   &     &     &   &   &    &   \\
				7  & 1 & 4  &     &     &     &   &     &     &   &   &    &   \\
				8  & 1 & 5  & 4   &     &     & 0 &     &     &   &   &    &   \\
				9  & 1 & 6  & 9   &     &     & 0 &     &     &   &   &    &   \\
				10 & 1 & 7  & 15  &     &     & 0 &     &     &   &   &    &   \\
				11 & 1 & 8  & 22  &     &     & 0 &     &     &   &   &    &   \\
				12 & 1 & 9  & 30  & 22  &     & 0 & 9   &     &   & 0 &    &   \\
				13 & 1 & 10 & 39  & 52  &     & 0 & 20  &     &   & 0 &    &   \\
				14 & 1 & 11 & 49  & 91  &     & 0 & 33  &     &   & 0 &    &   \\
				15 & 1 & 12 & 60  & 140 &     & 0 & 48  &     &   & 0 &    &   \\
				16 & 1 & 13 & 72  & 200 & 140 & 0 & 65  & 111 & 0 & 0 & 13 & 0 \\
				17 & 1 & 14 & 85  & 272 & 340 & 0 & 84  & 254 & 0 & 0 & 28 & 0 \\
				18 & 1 & 15 & 99  & 357 & 612 & 0 & 105 & 432 & 0 & 0 & 45 & 0 \\
				19 & 1 & 16 & 114 & 456 & 969 & 0 & 128 & 648 & 0 & 0 & 64 & 0 \\
				20 & 1 & 17 & 130 & 570 & 1425& 0 & 153 & 905 & 0 & 0 & 85 & 0 \\
                \hline
            \end{tabular}
        \end{table}
        
        \begin{definition} \label{def:pfeasapsisnontransdef}
        	For each $m \in \intsge{3}$, $k \in \intsge{m}$ and $t \in m\nonnegints$ such that $t \leq k$, we denote by $\pfeasapsisnontrans{m, k, t}$ the number of $\apsismon{m}{k}$-feasible upper non-transversal patterns that have $t$ of $k$ vertices contained in non-transversal blocks, that is  $\pfeasapsisnontrans{m, k, t} = \Sigma_{\substack{t_1 \in \posints, t_2, \ldots, t_x \in \nonnegints \\ t_1 + \ldots + xmt_x = t}} \pfeasapsisnontrans{m, k, t_1, \ldots, t_x}$ (see Tables \ref{table:apsismod3moncards} and \ref{table:apsismod4moncards} for example values). 
        \end{definition}

        \begin{theorem} \label{thm:apsismoncardinalities}
            The cardinality of $\apsismon{m}{k}$ is given by (see Tables \ref{table:apsismod3moncards} and \ref{table:apsismod4moncards} for example values),
            \[ \card{\apsismon{m}{k}} = 1 + \sum_{u=1}^{\floor{k/m}}\sum_{l=1}^{\floor{k/m}} \pfeasapsisnontrans{m, k, mu}\pfeasapsisnontrans{m, k, ml}\pfeastrans{m}{k-mu}{k-ml}. \]
            \begin{proof}
				The identity is the only diagram $\alpha \in \apsismon{m}{k}$ with $\rank{\alpha} = k$. Suppose the number of upper and lower vertices contained in non-transversal blocks are $mu$ and $ml$ respectively, where $u, l \in \set{1, \ldots, \floor{\frac{k}{m}}}$. For each of the $\pfeasapsisnontrans{m, k, mu}$ $\apsismon{m}{k}$-feasible upper non-transversal patterns that have $mu$ of $k$ upper vertices contained in non-transversal blocks and $\pfeasapsisnontrans{m, k, ml}$ $\apsismon{m}{k}$-feasible lower non-transversal patterns that have $ml$ of $k$ lower vertices contained in non-transversal blocks, there are $\pfeastrans{m}{k-mu}{k-ml}$ planar ways to connect the remaining $k-mu$ upper vertices to the remaining $k-ml$ lower vertices using only transversal blocks $b$ such that $\congmod{\noupverts{b}}{\nolowverts{b}}{m}$.
            \end{proof}
        \end{theorem}
 
        \begin{table}[!ht]
            \caption[ ]{Example values for $\pfeasapsisnontrans{3, k, t}$ and $\card{\triapmon{k}}$, calculated using the formulas from Definition \ref{def:pfeasapsisnontransdef} and Theorem \ref{thm:apsismoncardinalities}.}
            \label{table:apsismod3moncards}
            \centering
            \begin{tabular}{|c|rrrrr|r|}
                \hline
                \diagbox{k}{t} & 0 & 3 & 6 & 9 & 12 & $\card{\triapmon{k}}$ \\
                \hline
                3 & 1 & 1 &   &   &   & 2 \\
                4 & 1 & 2 &   &   &   & 5 \\
                5 & 1 & 3 &   &   &   & 19 \\
                6 & 1 & 4 & 3 &   &   & 74 \\
                7 & 1 & 5 & 7 &   &   & 320 \\
                8 & 1 & 6 & 12 &   &   & 1369 \\
                9 & 1 & 7 & 18 & 19 &   & 5732 \\
                10 & 1 & 8 & 25 & 46 &   & 24553 \\
                11 & 1 & 9 & 33 & 82 &   & 104493 \\
                12 & 1 & 10 & 42 & 128 & 127 & 439250 \\
                13 & 1 & 11 & 52 & 185 & 313 & 1871610 \\
                14 & 1 & 12 & 63 & 254 & 570 & 7952611 \\
                15 & 1 & 13 & 75 & 336 & 911 & 33550197 \\
                \hline
            \end{tabular}
        \end{table}
        
        \begin{table}[!ht]
            \caption[ ]{Example values for $\pfeasapsisnontrans{4, k, t}$ and $\card{\apsismon{4}{k}}$, calculated using the formulas from Definition \ref{def:pfeasapsisnontransdef} and Theorem \ref{thm:apsismoncardinalities}.}
            \label{table:apsismod4moncards}
            \centering
            \begin{tabular}{|c|rrrrr|r|}
                \hline
                \diagbox{k}{t} & 0 & 4 & 8 & 12 & 16 & $\card{\triapmon{k}}$ \\
                \hline
				4  & 1 & 1  &     &     &      & 2 \\
				5  & 1 & 2  &     &     &      & 5 \\
				6  & 1 & 3  &     &     &      & 19 \\
				7  & 1 & 4  &     &     &      & 65 \\
				8  & 1 & 5  & 4   &     &      & 217 \\
				9  & 1 & 6  & 9   &     &      & 766 \\
				10 & 1 & 7  & 15  &     &      & 2747 \\
				11 & 1 & 8  & 22  &     &      & 9489 \\
				12 & 1 & 9  & 30  & 31  &      & 32246 \\
				13 & 1 & 10 & 39  & 72  &      & 110817 \\
				14 & 1 & 11 & 49  & 124 &      & 383561 \\
				15 & 1 & 12 & 60  & 188 &      & 1308449 \\
				16 & 1 & 13 & 72  & 265 & 264  & 4430537 \\
				17 & 1 & 14 & 85  & 356 & 622  & 15099409 \\
				18 & 1 & 15 & 99  & 462 & 1089 & 51718033 \\
				19 & 1 & 16 & 114 & 584 & 1681 & 175876825 \\
				20 & 1 & 17 & 130 & 723 & 2415 & 595807899 \\
                \hline
            \end{tabular}
        \end{table}
        
	\section{The modular partition monoid, $\card{\modmon{m}{k}}$}
		Recall from Definition \ref{def:modmon} that the mod-$m$ monoid $\modmon{m}{k}$ consists of all bipartitions $\alpha \in \pttnmon{k}$ such that each block $b \in \alpha$ satisfies $\congmod{\noupverts{b}}{\nolowverts{b}}{m}$.
		
		\begin{definition} \label{def:feastrans}
			For each $m \in \posints$ and $k_1, k_2 \in \nonnegints$, by $\feastrans{m}{k_1}{k_2}$ we denote the number of ways to connect $k_1$ upper vertices to $k_2$ lower vertices using only transversal blocks $b$ such that $\congmod{\noupverts{b}}{\nolowverts{b}}{m}$.
		\end{definition}
		
		\begin{proposition} \label{prop:feastransrecurrences}
			For each $m \in \posints$ and $k_1, k_2 \in \nonnegints$, $\feastrans{m}{k_1}{k_2}$ satisfies the recurrence (see Tables \ref{table:XT(2,k1,k2)}, \ref{table:XT(3,k1,k2)} and \ref{table:XT(4,k1,k2)} for example values):
			\begin{align*}
				\feastrans{m}{0}{0} =&\ 1&\\
				\feastrans{m}{k}{0} =&\ \feastrans{m}{0}{k} = 0 & \text{ if } k > 0;\\
				\feastrans{m}{k_1}{k_2} =&\ 0 & \text{ if } \ncongmod{k_1}{k_2}{m};\\
				\feastrans{m}{k_1}{k_2} =&\ \sum_{\substack{(k_1', k_2') \leq (k_1, k_2)\\ \congmod{k_1'}{k_2'}{m}}} \binom{k_1-1}{k_1'-1}\binom{k_2}{k_2'} \feastrans{m}{k_1-k_1'}{k_2-k_2'} & \text{ if } \congmod{k_1}{k_2}{m}. 
			\end{align*}
			
			\begin{proof}
			 	Clearly we may not connect a positive number of vertices to zero vertices using transversal blocks, nor given $\ncongmod{k_1}{k_2}{m}$ may we connect $k_1$ upper vertices with $k_2$ lower vertices using only transversal blocks $b$ such that $\congmod{\noupverts{b}}{\nolowverts{b}}{m}$. Suppose $\congmod{k_1}{k_2}{m}$ and let $k_1' \in \{1, \ldots, k_1\}$ and $k_2' \in \{1, \ldots, k_2\}$ such that $\congmod{k_1'}{k_2'}{m}$. For each of the $\binom{k_1-1}{k_1'-1}\binom{k_2}{k_2'}$ ways that a transversal block containing the upper-right vertex $k_1$ may connect $k_1'$ upper vertices and $k_2'$ lower vertices, there are $\feastrans{m}{k_1-k_1'}{k_2-k_2'}$ ways in which the remaining $k_1-k_1'$ upper vertices and $k_2-k_2'$ lower vertices may be connected using only transversal blocks $b$ such that $\congmod{\noupverts{b}}{\nolowverts{b}}{m}$. Since there is one way for a single transversal to connect the $k_1$ upper vertices and $k_2$ lower vertices together it is convenient to have $\feastrans{m}{0}{0} = 1$.
			\end{proof}
		\end{proposition}
		
		\begin{table}[!ht]
			\caption[ ]{Example values for $\feastrans{2}{k_1}{k_2}$, calculated using the recurrences in Proposition \ref{prop:feastransrecurrences}.}
			\label{table:XT(2,k1,k2)}
			\centering
			\small
			\begin{tabular}{|c|rrrrrrrrrr|}
				\hline
				\diagbox{$k_1$}{$k_2$} & 0 & 1 & 2 & 3 & 4 & 5 & 6 & 7 & 8 & 9 \\
				\hline
		0 & 1 &   &     &      &       &        &        &          &          &   \\
		1 &   & 1 &     & 1    &       & 1      &        & 1        &          & 1 \\
		2 &   &   & 3   &      & 9     &        & 33     &          & 129      &   \\
		3 &   & 1 &     & 16   &       & 106    &        & 736      &          & 5686 \\
		4 &   &   & 9   &      & 147   &        & 1599   &          & 18027    &   \\
		5 &   & 1 &     & 106  &       & 1996   &        & 31606    &          & 512296 \\
		6 &   &   & 33  &      & 1599  &        & 37383  &          & 791439   &   \\
		7 &   & 1 &     & 736  &       & 31606  &        & 905416   &          & 24489466 \\
		8 &   &   & 129 &      & 18027 &        & 791439 &          & 27370227 &   \\
		9 &   & 1 &     & 5686 &       & 512296 &        & 24489466 &          & 1004077636 \\
				\hline
			\end{tabular}
		\end{table}
		
		\begin{table}[!ht]
			\caption[ ]{Example values for $\feastrans{3}{k_1}{k_2}$, calculated using the recurrences in Proposition \ref{prop:feastransrecurrences}.}
			\label{table:XT(3,k1,k2)}
			\centering
			\begin{tabular}{|c|rrrrrrrrrr|}
				\hline
				\diagbox{$k_1$}{$k_2$} & 0 & 1 & 2 & 3 & 4 & 5 & 6 & 7 & 8 & 9 \\
				\hline
				0 & 1 &   &    &      &      &       &         &        &          &           \\
				1 &   & 1 &    &      & 1    &       &         & 1      &          &           \\
				2 &   &   & 3  &      &      & 11    &         &        & 87       &           \\
				3 &   &   &    & 16   &      &       & 154     &        &          & 2620      \\
				4 &   & 1 &    &      & 131  &       &         & 2647   &          &           \\
				5 &   &   & 11 &      &      & 1521  &         &        & 55731    &           \\
				6 &   &   &    & 154  &      &       & 23562   &        &          & 1419082   \\
				7 &   & 1 &    &      & 2647 &       &         & 469659 &          &           \\
				8 &   &   & 87 &      &      & 55731 &         &        & 11676227 &           \\
				9 &   &   &    & 2620 &      &       & 1419082 &        &          & 353355424 \\
				\hline
			\end{tabular}
		\end{table}
		
		\begin{table}[!ht]
			\caption[ ]{Example values for $\feastrans{4}{k_1}{k_2}$, calculated using the recurrences in Proposition \ref{prop:feastransrecurrences}.}
			\label{table:XT(4,k1,k2)}
			\centering
			\begin{tabular}{|c|rrrrrrrrrr|}
				\hline
				\diagbox{$k_1$}{$k_2$} & 0 & 1 & 2 & 3 & 4 & 5 & 6 & 7 & 8 & 9 \\
				\hline
				0 & 1 &   &    &     &      &       &       &        &          &   \\
				1 &   & 1 &    &     &      & 1     &       &        &          & 1 \\
				2 &   &   & 3  &     &      &       & 13    &        &          &   \\
				3 &   &   &    & 16  &      &       &       & 211    &          &   \\
				4 &   &   &    &     & 131  &       &       &        & 4121     &   \\
				5 &   & 1 &    &     &      & 1496  &       &        &          & 97096 \\
				6 &   &   & 13 &     &      &       & 22518 &        &          &   \\
				7 &   &   &    & 211 &      &       &       & 428891 &          &   \\
				8 &   &   &    &     & 4121 &       &       &        & 10039731 &   \\
				9 &   & 1 &    &     &      & 97096 &       &        &          & 282357916 \\
				\hline
			\end{tabular}
		\end{table}
		
		\begin{definition}
			For each $m \in \posints$, $k \in \intsge{m}$, $x \in \floor{\frac{k}{m}}$ and $t_1, \ldots, t_x \in \nonnegints$ such that $\Sigma_{i=1}^ximt_i \leq k$, we denote by $\feasnontrans{m, k, t_1, \ldots, t_x}$ the number of $\modmon{m}{k}$-feasible upper non-transversal patterns containing precisely $t_1$ type $(m, 0)$, \ldots, and $t_x$ type $(xm, 0)$ non-transversal blocks.
		\end{definition}
        
        Note it trivially follows from the mod-$m$ monoid $\pmodmon{m}{k}$ being closed under the vertical flip involution $^*$, that is ${\modmon{m}{k}}^* = \modmon{m}{k}$, that the number of $\modmon{m}{k}$-feasible lower non-transversal patterns containing precisely $t_1$ type $(0, m)$, \ldots, and $t_x$ type $(0, xm)$ non-transversal blocks is equal to the number of $\modmon{m}{k}$-feasible upper non-transversal patterns containing precisely $t_1$ type $(m, 0)$, \ldots, and $t_x$ type $(xm, 0)$ non-transversal blocks.

        \begin{proposition} \label{prop:feasnontransrecurrences}
        	For each $m \in \posints$, $k \in \intsge{m}$, $x \in \floor{\frac{k}{m}}$ and $t_1, \ldots, t_x \in \nonnegints$ such that $\Sigma_{i=1}^ximt_i \leq k$, $\feasnontrans{m, k, t_1, \ldots, t_x}$ satisfies the recurrence (see Tables \ref{table:XN(2,k,t)}, \ref{table:XN(3,k,t)} and \ref{table:XN(4,k,t)} for example values):
            \begin{align*}
                \feasnontrans{m, k, 0, \ldots, 0} =&\ 1; \\
                \feasnontrans{m, k, t_1, \ldots, t_x} =&\ \feasnontrans{m, k-1, t_1, \ldots, t_x} \\
                &\hspace{-1.75cm}+ \sum_{\substack{i \in \{1, \ldots, x\} \\ t_i > 0}}\binom{k-1}{im-1}\feasnontrans{m, k-im, t_1, \ldots, t_i-1, \ldots, t_x}.
            \end{align*}

            \begin{proof}
                There is only one way to place no non-transversal blocks, giving us $\feasnontrans{m, k, 0, \ldots, 0} = 1$. Suppose there is at least one non-transversal block to be placed. There are $\feasnontrans{m, k-1, t_1, \ldots, t_x}$ $\modmon{m}{k}$-feasible upper non-transversal patterns that contain precisely $t_1$ type $(m, 0)$, \ldots, and $t_x$ type $(xm, 0)$ non-transversal blocks along the $k-1$ left-most vertices, that is so that no non-transversal blocks contain the right-most vertex $k$. Let $i \in \{1, \ldots, x\}$. For each of the $\binom{k-1}{im-1}$ ways a type $(im, 0)$ non-transversal block may contain the right-most vertex $k$ there are $\feasnontrans{m, k-1, t_1, \ldots, t_i-1, \ldots, t_x}$ $\modmon{m}{k}$-feasible upper non-transversal patterns containing precisely $t_1$ type $(m, 0)$, \ldots, $t_i-1$ type $(im, 0)$, \ldots, and $t_x$ type $(xm, 0)$ non-transversal blocks from the remaining $k-im$ vertices.
            \end{proof}
        \end{proposition}
        
        \begin{table}[!ht]
            \caption[ ]{Example values for $\feasnontrans{2, k, \vec{t}}$ where $\vec{t} = (t_1, t_2, t_3)$, calculated using the recurrences in Proposition \ref{prop:feasnontransrecurrences}.}
            \label{table:XN(2,k,t)}
            \centering 
            \begin{tabular}{|c|rrrrrrr|}
                \hline
                \diagbox{$k$}{$\vec{t}$} & \rotatebox{90}{(0,0,0) } & \rotatebox{90}{(1,0,0)} & \rotatebox{90}{(2,0,0)} & \rotatebox{90}{(3,0,0)} & \rotatebox{90}{(0,1,0)} & \rotatebox{90}{(1,1,0)} & \rotatebox{90}{(0,0,1)}  \\
                \hline
				2 & 1 & 1 &   &   &   &   &   \\
				3 & 1 & 3 &   &   &   &   &   \\
				4 & 1 & 6 & 3 &   & 1 &   &   \\
				5 & 1 & 10 & 15 &   & 5 &   &   \\
				6 & 1 & 15 & 45 & 15 & 15 & 15 & 1 \\
				7 & 1 & 21 & 105 & 105 & 35 & 105 & 7 \\
				8 & 1 & 28 & 210 & 420 & 70 & 420 & 28 \\
				9 & 1 & 36 & 378 & 1260 & 126 & 1260 & 84 \\
				10 & 1 & 45 & 630 & 3150 & 210 & 3150 & 210 \\
				11 & 1 & 55 & 990 & 6930 & 330 & 6930 & 462 \\
				12 & 1 & 66 & 1485 & 13860 & 495 & 13860 & 924 \\
				13 & 1 & 78 & 2145 & 25740 & 715 & 25740 & 1716 \\
				14 & 1 & 91 & 3003 & 45045 & 1001 & 45045 & 3003 \\
				15 & 1 & 105 & 4095 & 75075 & 1365 & 75075 & 5005 \\
                \hline
            \end{tabular}
        \end{table}
        
        \begin{table}[!ht]
            \caption[ ]{Example values for $\feasnontrans{3, k, \vec{t}}$ where $\vec{t} = (t_1, t_2, t_3)$, calculated using the recurrences in Proposition \ref{prop:feasnontransrecurrences}.}
            \label{table:XN(3,k,t)}
            \centering 
            \begin{tabular}{|c|rrrrrrr|}
                \hline
                \diagbox{$k$}{$\vec{t}$} & \rotatebox{90}{(0,0,0) } & \rotatebox{90}{(1,0,0)} & \rotatebox{90}{(2,0,0)} & \rotatebox{90}{(3,0,0)} & \rotatebox{90}{(0,1,0)} & \rotatebox{90}{(1,1,0)} & \rotatebox{90}{(0,0,1)}  \\
                \hline
				3 & 1 & 1 &   &   &   &   &   \\
				4 & 1 & 4 &   &   &   &   &   \\
				5 & 1 & 10 &   &   &   &   &   \\
				6 & 1 & 20 & 10 &   & 1 &   &   \\
				7 & 1 & 35 & 70 &   & 7 &   &   \\
				8 & 1 & 56 & 280 &   & 28 &   &   \\
				9 & 1 & 84 & 840 & 280 & 84 & 84 & 1 \\
				10 & 1 & 120 & 2100 & 2800 & 210 & 840 & 10 \\
				11 & 1 & 165 & 4620 & 15400 & 462 & 4620 & 55 \\
				12 & 1 & 220 & 9240 & 61600 & 924 & 18480 & 220 \\
				13 & 1 & 286 & 17160 & 200200 & 1716 & 60060 & 715 \\
				14 & 1 & 364 & 30030 & 560560 & 3003 & 168168 & 2002 \\
				15 & 1 & 455 & 50050 & 1401400 & 5005 & 420420 & 5005 \\
                \hline
            \end{tabular}
        \end{table}
        
        \begin{table}[!ht]
            \caption[ ]{Example values for $\feasnontrans{4, k, \vec{t}}$ where $\vec{t} = (t_1, t_2, t_3)$, calculated using the recurrences in Proposition \ref{prop:feasnontransrecurrences}.}
            \label{table:XN(4,k,t)}
            \centering 
            \begin{tabular}{|c|rrrrrrr|}
                \hline
                \diagbox{$k$}{$\vec{t}$} & \rotatebox{90}{(0,0,0) } & \rotatebox{90}{(1,0,0)} & \rotatebox{90}{(2,0,0)} & \rotatebox{90}{(3,0,0)} & \rotatebox{90}{(0,1,0)} & \rotatebox{90}{(1,1,0)} & \rotatebox{90}{(0,0,1)}  \\
                \hline
				4 & 1 & 1 &   &   &   &   &   \\
				5 & 1 & 5 &   &   &   &   &   \\
				6 & 1 & 15 &   &   &   &   &   \\
				7 & 1 & 35 &   &   &   &   &   \\
				8 & 1 & 70 & 35 &   & 1 &   &   \\
				9 & 1 & 126 & 315 &   & 9 &   &   \\
				10 & 1 & 210 & 1575 &   & 45 &   &   \\
				11 & 1 & 330 & 5775 &   & 165 &   &   \\
				12 & 1 & 495 & 17325 & 5775 & 495 & 495 & 1 \\
				13 & 1 & 715 & 45045 & 75075 & 1287 & 6435 & 13 \\
				14 & 1 & 1001 & 105105 & 525525 & 3003 & 45045 & 91 \\
				15 & 1 & 1365 & 225225 & 2627625 & 6435 & 225225 & 455 \\
                \hline
            \end{tabular}
        \end{table}
        
        \begin{definition} \label{def:feasnontrans}
        	For each $m \in \posints$, $k \in \intsge{m}$ and $t \in m\nonnegints$ such that $t \leq k$, we denote by $\feasnontrans{m, k, t}$ the number of $\modmon{m}{k}$-feasible upper non-transversal patterns that have $t$ of $k$ vertices contained in non-transversal blocks, that is $\feasnontrans{m, k, t} = \Sigma_{\substack{t_1, \ldots, t_x \in \nonnegints \\ t_1 + \ldots + xmt_x = t}}\ \feasnontrans{m, k, t_1, \ldots, t_x}$ (see Tables \ref{table:mod2moncards}, \ref{table:mod3moncards} and \ref{table:mod4moncards} for example values). 
        \end{definition}

        \begin{theorem} \label{thm:feasnontranscardinalities}
            For each $m \in \posints$ and $k \in \intsge{m}$, the cardinality of $\modmon{m}{k}$ is given by (see Tables \ref{table:mod2moncards}, \ref{table:mod3moncards}, \ref{table:mod4moncards} for example values),
            \[ \card{\modmon{m}{k}} = \sum_{u=0}^{\floor{k/m}}\sum_{l=0}^{\floor{k/m}} \feasnontrans{m, k, mu}\feasnontrans{m, k, ml}\feastrans{m}{k-mu}{k-ml}. \]

            \begin{proof}
				Suppose the number of upper vertices and number of lower vertices contained in non-transversal blocks are $mu$ and $ml$ respectively, where $u, l \in \set{0, \ldots, \floor{\frac{k}{m}}}$. For each of the $\feasnontrans{m, k, mu}$ $\modmon{m}{k}$-feasible upper non-transversal patterns that have $mu$ of $k$ upper vertices contained in non-transversal blocks and $\feasnontrans{m, k, ml}$ $\modmon{m}{k}$-feasible lower non-transversal patterns that have $ml$ of $k$ lower vertices contained in non-transversal blocks, there are $\feastrans{m}{k-mu}{k-ml}$ ways to connect the remaining $k-mu$ upper vertices to the remaining $k-ml$ lower vertices using only transversal blocks $b$ such that $\congmod{\noupverts{b}}{\nolowverts{b}}{m}$.
            \end{proof}
        \end{theorem}
        
        Note that the author and James East, upon the author having explained how he established the results from this section, used a very similar counting method to establish recurrence relations for the number of idempotents in the Brauer monoid. James further generalised this result to establish the number of idempotents in the Partition monoid, which tied in with some of what Igor Dolinka was working on at the time and led to the collaborative paper \cite{art:Dolinka:EnumerationOfIdempotentsInDiagramSemigrps}. 
        
        \begin{table}[!ht]
            \caption[ ]{Example values for $\feasnontrans{2, k, t}$ and $\card{\modmon{2}{k}}$, calculated using the formulas in Definition \ref{def:feasnontrans} and Theorem \ref{thm:feasnontranscardinalities}.}
            \label{table:mod2moncards}
            \centering
            \begin{tabular}{|c|rrrr|r|}
                \hline
                \diagbox{k}{t} & 0 & 2 & 4 & 6 & {\normalsize$\card{\modmon{2}{k}}$} \\
                \hline
				2 & 1 & 1 &   &   & 4 \\
				3 & 1 & 3 &   &   & 31 \\
				4 & 1 & 6 & 4 &   & 379 \\
				5 & 1 & 10 & 20 &   & 6556 \\
				6 & 1 & 15 & 60 & 31 & 150349 \\
				7 & 1 & 21 & 140 & 217 & 4373461 \\
				8 & 1 & 28 & 280 & 868 & 156297964 \\
				9 & 1 & 36 & 504 & 2604 & 6698486371 \\
				10 & 1 & 45 & 840 & 6510 & 337789490599 \\
				11 & 1 & 55 & 1320 & 14322 & 19738202807236 \\
				12 & 1 & 66 & 1980 & 28644 & 1319703681935929 \\
				13 & 1 & 78 & 2860 & 53196 & 99896787342523081 \\
				14 & 1 & 91 & 4004 & 93093 & 8484301665702298804 \\
				15 & 1 & 105 & 5460 & 155155 & 802221679220975886631 \\
                \hline
            \end{tabular}
        \end{table}
        
        \begin{table}[!ht]
            \caption[ ]{Example values for $\feasnontrans{3, k, t}$ and $\card{\modmon{3}{k}}$, calculated using the formulas in Definition \ref{def:feasnontrans} and Theorem \ref{thm:feasnontranscardinalities}.}
            \label{table:mod3moncards}
            \centering
            \begin{tabular}{|c|rrrr|r|}
                \hline
                \diagbox{k}{t} & 0 & 3 & 6 & 9 & {\normalsize$\card{\modmon{3}{k}}$} \\
                \hline
				3 & 1 & 1 &   &   & 17 \\
				4 & 1 & 4 &   &   & 155 \\
				5 & 1 & 10 &   &   & 2041 \\
				6 & 1 & 20 & 11 &   & 36243 \\
				7 & 1 & 35 & 77 &   & 826897 \\
				8 & 1 & 56 & 308 &   & 23405595 \\
				9 & 1 & 84 & 924 & 365 & 800555801 \\
				10 & 1 & 120 & 2310 & 3650 & 32417395123 \\
				11 & 1 & 165 & 5082 & 20075 & 1528888375697 \\
				12 & 1 & 220 & 10164 & 80300 & 82865247031515 \\
				13 & 1 & 286 & 18876 & 260975 & 5104104871207161 \\
				14 & 1 & 364 & 33033 & 730730 & 353921927969377043 \\
				15 & 1 & 455 & 55055 & 1826825 & 27403472985911422417 \\
                \hline
            \end{tabular}
        \end{table}
        
        \begin{table}[!ht]
            \caption[ ]{Example values for $\feasnontrans{4, k, t}$ and $\card{\modmon{4}{k}}$, calculated using the formulas in Definition \ref{def:feasnontrans} and Theorem \ref{thm:feasnontranscardinalities}.}
            \label{table:mod4moncards}
            \centering
            \begin{tabular}{|c|rrrr|r|}
                \hline
                \diagbox{k}{t} & 0 & 4 & 8 & 12 & {\normalsize$\card{\modmon{4}{k}}$} \\
                \hline
				4 & 1 & 1 &   &   & 132 \\
				5 & 1 & 5 &   &   & 1531 \\
				6 & 1 & 15 &   &   & 23583 \\
				7 & 1 & 35 &   &   & 463261 \\
				8 & 1 & 70 & 36 &   & 11259867 \\
				9 & 1 & 126 & 324 &   & 330763876 \\
				10 & 1 & 210 & 1620 &   & 11522992578 \\
				11 & 1 & 330 & 5940 &   & 468713029951 \\
				12 & 1 & 495 & 17820 & 6271 & 21971754415317 \\
				13 & 1 & 715 & 46332 & 81523 & 1173833581966501 \\
				14 & 1 & 1001 & 108108 & 570661 & 70790559991302063 \\
				15 & 1 & 1365 & 231660 & 2853305 & 4779273111284582836 \\
                \hline
            \end{tabular}
        \end{table}

        \subsection{The modular-$2$ partition monoid, $\card{\modmon{2}{k}}$}
        	The calculated values in the right-most column of Table \ref{table:mod2moncards}, which form the start of the sequence $\set{\card{\modmon{2}{k}}: k \in \intsge{2}}$, match the number of partitions of $\{1, \ldots, 2k\}$ with blocks of even size, which is listed on the OEIS as sequence $A001764$. We proceed to establish that the sequence $\set{\card{\modmon{2}{k}}: k \in \intsge{2}}$ matches sequence $A001764$ on the OEIS.
        	
        	\begin{proposition}
	        	The cardinality of the mod-$2$ monoid $\modmon{2}{k}$ is equal to the cardinality of partitions of $\{1, \ldots, 2k\}$ with blocks of even size.
	        	
	        	\begin{proof}
		        	Given a bipartition $\alpha \in \modmon{2}{k}$, each block $b \in \alpha$ satisfies $\noupverts{b} + \nolowverts{b} \in 2\ints$. Conversely, given a partition $\alpha$ of $\{1, \ldots, 2k\}$ such that every block has even size, each block $b \in \alpha$ satisfies $\congmod{\card{b \cap \{1, \ldots, k\}}}{\card{b \cap \{k+1, \ldots, 2k\}}}{2}$.
	        	\end{proof}
        	\end{proposition}

    \section{The crossed $m$-apsis generated diagram monoid, $\card{\capsismon{m}{k}}$}
    	Recall from Definition \ref{def:capsismon} that the crossed $m$-apsis monoid $\capsismon{m}{k}$ consists of all bipartitions $\alpha \in \modmon{m}{k}$ such that either:
		\begin{enumerate}
			\item $\alpha$ is a permutation, that is $\alpha \in \symgrp{k}$; or 
			\item $\alpha$ contains at least one type $(m, 0)$ non-transversal block and at least one type $(0, m)$ non-transversal block.
		\end{enumerate} 
		
		\begin{definition}
			For each $m \in \intsge{3}$, $k \in \intsge{m}$, $x \in \floor{\frac{k}{m}}$ and $t_1, \ldots, t_x \in \nonnegints$ such that $\Sigma_{i=1}^ximt_i \leq k$, we denote by $\feasapsisnontrans{m, k, t_1, \ldots, t_x}$ the number of $\capsismon{m}{k}$-feasible upper non-transversal patterns containing precisely $t_1$ type $(m, 0)$, \ldots, and $t_x$ type $(xm, 0)$ non-transversal blocks.
		\end{definition} 
		
		Note it trivially follows from the crossed $m$-apsis monoid $\capsismon{m}{k}$ being closed under the vertical flip involution $^*$, that is ${\capsismon{m}{k}}^* = \capsismon{m}{k}$, that the number of $\capsismon{m}{k}$-feasible lower non-transversal patterns containing precisely $t_1$ type $(0, m)$, \ldots, and $t_x$ type $(0, xm)$ non-transversal blocks is equal to the number of $\capsismon{m}{k}$-feasible upper non-transversal patterns containing precisely $t_1$ type $(m, 0)$, \ldots, and $t_x$ type $(xm, 0)$ non-transversal blocks.

        \begin{proposition} \label{prop:feasapsisnontransrecurrences}
        	For each $m \in \posints$, $k \in \intsge{m}$, $x \in \floor{\frac{k}{m}}$ and $t_1, \ldots, t_x \in \nonnegints$ such that $\Sigma_{i=1}^ximt_i \leq k$, $\feasapsisnontrans{m, k, t_1, \ldots, t_x}$ satisfies the recurrence (see Tables \ref{table:XNB(3,k,t)} and \ref{table:XNB(4,k,t)} for example values):
            \begin{align*}
                \feasapsisnontrans{m, k, 0, \ldots, 0} =&\ 1; \\
                \feasapsisnontrans{m, k, 0, t_2, \ldots, t_x} =&\ 0; \\
                \feasapsisnontrans{m, k, t_1, \ldots, t_x} =&\ \feasapsisnontrans{m, k-1, t_1, \ldots, t_x} \\
                        &\hspace{-1.75cm}+ \sum_{\substack{i \in \{1, \ldots, x\} \\ t_i > 0}}\binom{k-1}{im-1}\feasapsisnontrans{m, k-im, t_1, \ldots, t_i-1, \ldots, t_x}.
            \end{align*}

            \begin{proof}
                There is only one way to place no non-transversal blocks, giving us $\feasapsisnontrans{m, k, 0, \ldots, 0} = 1$. Suppose there is at least one non-transversal block to be placed. If $t_1 = 0$ then we cannot $\capsismon{m}{k}$-feasibly place the $t_2$ type $(2m, 0)$, \ldots, and $t_x$ type $(xm, 0)$ non-transversal blocks, hence $\feasapsisnontrans{m, k, 0, t_2, \ldots, t_x} = 0$. Now there are $\feasapsisnontrans{m, k-1, t_1, \ldots, t_x}$ $\capsismon{m}{k}$-feasible upper non-transversal patterns that contain precisely $t_1$ type $(m, 0)$, \ldots, and $t_x$ type $(xm, 0)$ non-transversal blocks along the $k-1$ left-most vertices, that is so that no non-transversal blocks contain the right-most vertex $k$. Let $i \in \{1, \ldots, x\}$. For each of the $\binom{k-1}{im-1}$ ways a type $(im, 0)$ non-transversal block may contain the right-most vertex $k$ there are $\feasapsisnontrans{m, k-1, t_1, \ldots, t_i-1, \ldots, t_x}$ $\capsismon{m}{k}$-feasible upper non-transversal patterns containing precisely $t_1$ type $(m, 0)$, \ldots, $t_i-1$ type $(im, 0)$, \ldots, and $t_x$ type $(xm, 0)$ non-transversal blocks from the remaining $k-im$ vertices.
            \end{proof}
        \end{proposition}
        
        \begin{table}[!ht]
            \caption[ ]{Example values for $\feasapsisnontrans{3, k, \vec{t}}$ where $\vec{t} = (t_1, t_2, t_3)$, calculated using the recurrences in Proposition \ref{prop:feasapsisnontransrecurrences}.}
            \label{table:XNB(3,k,t)}
            \centering 
            \begin{tabular}{|c|rrrrrrr|}
                \hline
                \diagbox{$k$}{$\vec{t}$} & \rotatebox{90}{(0,0,0) } & \rotatebox{90}{(1,0,0)} & \rotatebox{90}{(2,0,0)} & \rotatebox{90}{(3,0,0)} & \rotatebox{90}{(0,1,0)} & \rotatebox{90}{(1,1,0)} & \rotatebox{90}{(0,0,1)}  \\
                \hline
                0  & 1 &     &      &       &     &      &   \\
                1  & 1 &     &      &       &     &      &   \\
                2  & 1 &     &      &       &     &      &   \\
                3  & 1 & 1   &      &       &     &      &   \\
                4  & 1 & 4   &      &       &     &      &   \\
                5  & 1 & 10  &      &       &     &      &   \\
                6  & 1 & 20  & 10   &       & 0   &      &   \\
                7  & 1 & 35  & 70   &       & 0   &      &   \\
                8  & 1 & 56  & 280  &       & 0   &      &   \\
                9  & 1 & 84  & 840  & 280   & 0   & 84   & 0 \\
                10 & 1 & 120 & 2100 & 2800  & 0   & 840  & 0  \\
                11 & 1 & 165 & 4620 & 15400 & 0   & 4620 & 0  \\
                \hline
            \end{tabular}
        \end{table}
        
        \begin{table}[!ht]
            \caption[ ]{Example values for $\feasapsisnontrans{4, k, \vec{t}}$ where $\vec{t} = (t_1, t_2, t_3)$, calculated using the recurrences in Proposition \ref{prop:feasapsisnontransrecurrences}.}
            \label{table:XNB(4,k,t)}
            \centering 
            \begin{tabular}{|c|rrrrrrr|}
                \hline
                \diagbox{$k$}{$\vec{t}$} & \rotatebox{90}{(0,0,0) } & \rotatebox{90}{(1,0,0)} & \rotatebox{90}{(2,0,0)} & \rotatebox{90}{(3,0,0)} & \rotatebox{90}{(0,1,0)} & \rotatebox{90}{(1,1,0)} & \rotatebox{90}{(0,0,1)}  \\
                \hline
				4  & 1 & 1    &        &         &   &        &   \\
				5  & 1 & 5    &        &         &   &        &   \\
				6  & 1 & 15   &        &         &   &        &   \\
				7  & 1 & 35   &        &         &   &        &   \\
				8  & 1 & 70   & 35     &         & 0 &        &   \\
				9  & 1 & 126  & 315    &         & 0 &        &   \\
				10 & 1 & 210  & 1575   &         & 0 &        &   \\
				11 & 1 & 330  & 5775   &         & 0 &        &   \\
				12 & 1 & 495  & 17325  & 5775    & 0 & 495    & 0 \\
				13 & 1 & 715  & 45045  & 75075   & 0 & 6435   & 0 \\
				14 & 1 & 1001 & 105105 & 525525  & 0 & 45045  & 0 \\
				15 & 1 & 1365 & 225225 & 2627625 & 0 & 225225 & 0 \\
                \hline
            \end{tabular}
        \end{table}
        
        \begin{definition} \label{def:feasapsisnontrans}
        	For each $m \in \posints$, $k \in \intsge{m}$ and $t \in m\nonnegints$ such that $t \leq k$, we denote by $\feasapsisnontrans{m, k, t}$ the number of $\capsismon{m}{k}$-feasible upper non-transversal patterns that have $t$ of $k$ vertices contained in non-transversal blocks, that is $\feasapsisnontrans{m, k, t} = \Sigma_{\substack{t_1 \in \posints, t_2, \ldots, t_x \in \nonnegints \\ t_1 + \ldots + xmt_x = t}}\ \feasapsisnontrans{m, k, t_1, \ldots, t_x}$ (see Tables \ref{table:capsismon3moncards} and \ref{table:capsismon4moncards} for example values).
        \end{definition}

        \begin{theorem} \label{thm:feasapsisnontranscardinalities}
            For each $m \in \posints$ and $k \in \intsge{m}$, the cardinality of $\capsismon{m}{k}$ is given by (see Tables \ref{table:capsismon3moncards} and \ref{table:capsismon4moncards} for example values),
            \[ \card{\capsismon{m}{k}} = k! + \sum_{u=1}^{\floor{k/m}}\sum_{l=1}^{\floor{k/m}} \feasapsisnontrans{m, k, mu}\feasapsisnontrans{m, k, ml}\feastrans{m}{k-mu}{k-ml}. \]

            \begin{proof}
				There are $k!$ elements of rank $k$, more specifically the symmetric group $\symgrp{k}$, otherwise we must have at least one upper and one lower non-transversal block. Suppose the number of upper vertices and number of lower vertices contained in non-transversal blocks are $mu$ and $ml$ respectively, where $u, l \in \set{1, \ldots, \floor{\frac{k}{m}}}$. For each of the $\feasapsisnontrans{m, k, mu}$ $\capsismon{m}{k}$-feasible upper non-transversal patterns that have $mu$ of $k$ upper vertices contained in non-transversal blocks and $\feasapsisnontrans{m, k, ml}$ $\capsismon{m}{k}$-feasible lower non-transversal patterns that have $ml$ of $k$ lower vertices contained in non-transversal blocks, there are $\feastrans{m}{k-mu}{k-ml}$ ways to connect the remaining $k-mu$ upper vertices to the remaining $k-ml$ lower vertices using only transversal blocks $b$ such that $\congmod{\noupverts{b}}{\nolowverts{b}}{m}$.
            \end{proof}
        \end{theorem}
        
        \begin{table}[!ht]
            \caption[ ]{Example values for $\feasapsisnontrans{3, k, t}$ and $\card{\capsismon{3}{k}}$, calculated using the formulas in Definition \ref{def:feasapsisnontrans} and Theorem \ref{thm:feasapsisnontranscardinalities}.}
            \label{table:capsismon3moncards}
            \centering
            \begin{tabular}{|c|rrrr|r|}
                \hline
                \diagbox{k}{t} & 0 & 3 & 6 & 9 & {\normalsize$\card{\capsismon{3}{k}}$} \\
                \hline
                3 & 1 & 1 &   &   & 7 \\
                4 & 1 & 4 &   &   & 40 \\
                5 & 1 & 10 &   &   & 420 \\
                6 & 1 & 20 & 10 &   & 7220 \\
                7 & 1 & 35 & 70 &   & 175315 \\
                8 & 1 & 56 & 280 &   & 5390336 \\
                9 & 1 & 84 & 840 & 364 & 199770928 \\
                10 & 1 & 120 & 2100 & 3640 & 8707927600 \\
                11 & 1 & 165 & 4620 & 20020 & 439169520075 \\
                \hline
            \end{tabular}
        \end{table}
        
        \begin{table}[!ht]
            \caption[ ]{Example values for $\feasapsisnontrans{4, k, t}$ and $\card{\capsismon{4}{k}}$, calculated using the formulas in Definition \ref{def:feasapsisnontrans} and Theorem \ref{thm:feasapsisnontranscardinalities}.}
            \label{table:capsismon4moncards}
            \centering
            \begin{tabular}{|c|rrrr|r|}
                \hline
                \diagbox{k}{t} & 0 & 4 & 8 & 12 & {\normalsize$\card{\capsismon{3}{k}}$} \\
                \hline
				4 & 1 & 1 &   &   & 25 \\
				5 & 1 & 5 &   &   & 145 \\
				6 & 1 & 15 &   &   & 1395 \\
				7 & 1 & 35 &   &   & 24640 \\
				8 & 1 & 70 & 35 &   & 683445 \\
				9 & 1 & 126 & 315 &   & 24291981 \\
				10 & 1 & 210 & 1575 &   & 1012713975 \\
				11 & 1 & 330 & 5775 &   & 48083983200 \\
				12 & 1 & 495 & 17325 & 6270 & 2570506151400 \\
				13 & 1 & 715 & 45045 & 81510 & 153658593860200 \\
				14 & 1 & 1001 & 105105 & 570570 & 10213751655054948 \\
				15 & 1 & 1365 & 225225 & 2852850 & 751055052971960100 \\
                \hline
            \end{tabular}
        \end{table}

    \clearpage{\pagestyle{empty}\cleardoublepage}
\chapter{Green's relations} \label{chap:greensrltns}
	\section{Pattern compatibility $\sim_S$}
	\begin{definition}
		Let $S$ be a submonoid of the partition monoid $\pttnmon{k}$. We say that an $S$-admissible upper pattern $p \in \uppat{S}$ and an $S$-admissible lower pattern $q \in \lowpat{S}$ are \textit{$S$-compatible}, which we denote by $p \sim q$ or less succinctly by $p \sim_{S} q$ when $S$ is contextually ambiguous, if there exists $\alpha \in S$ such that $p$ is the upper pattern of $\alpha$ and $q$ is the lower pattern of $\alpha$, that is $\uppat{\alpha} = p$ and $\lowpat{\alpha} = q$. 
	\end{definition}
	
	\begin{proposition}
		Pattern compatibility is transitive.
		
		\begin{proof}
			Let $p, q, r \in \uppat{S}$. Suppose $p \sim q$ and $q \sim r$. Therefore there exist $\alpha, \beta \in S$ such that $\uppat{\alpha} = p$, $\lowpat{\alpha} = q = \uppat{\beta}$ and $\lowpat{\beta} = r$. It follows from Proposition \ref{prop:productuplowpatwithsamemiddlepat} that $\uppat{\alpha\beta} = p$ and $\lowpat{\alpha\beta} = r$, and hence that $p \sim r$.
		\end{proof}
	\end{proposition}
	
	\begin{proposition}
		Given a submonoid $S$ of the partition monoid $\pttnmon{k}$, if $S$ is closed under the vertical flip involution $^*$, that is $S^* = S$, then pattern compatibility is an equivalence relation.
		
		\begin{proof}
			Given $p \in \uppat{S}$, by definition there exist $\alpha \in S$ such that $\uppat{\alpha} = p$. Since $S$ is closed under $^*$ we have $\alpha^*, \alpha\alpha^* \in S$ and $\uppat{\alpha\alpha^*} = \lowpat{\alpha\alpha^*} = p$, giving us $p \sim p$.
			
			Given $p, q \in \uppat{S}$ satisfying $p \sim q$, by definition there exist $\alpha \in S$ such that $\uppat{\alpha} = p$ and $\lowpat{\alpha} = q$. Since $\alpha^* \in S$, $\uppat{\alpha^*} = \lowpat{\alpha} = q$ and $\lowpat{\alpha^*} = \uppat{\alpha} = p$, $q \sim p$.
		\end{proof}
	\end{proposition}

	\begin{definition}
		Let $S$ be a diagram semigroup that is closed under the vertical flip involution $^*$ and $\alpha, \beta \in S$. We say that \textit{$\alpha$ and $\beta$ have $S$-compatible patterns} if any of the following equivalent conditions are satisfied:
		\begin{enumerate}
			\item $\uppat{\alpha}$ and $\uppat{\beta}$ are $S$-compatible;
			\item $\uppat{\alpha}$ and $\lowpat{\beta}$ are $S$-compatible;
			\item $\lowpat{\alpha}$ and $\uppat{\beta}$ are $S$-compatible; or
			\item $\lowpat{\alpha}$ and $\lowpat{\beta}$ are $S$-compatible.
		\end{enumerate}
	\end{definition}
	
	\section{Green's $\mathcal{D}$ relation on diagram semigroups closed under $^*$}
	\begin{theorem} \label{thm:greensDrltnonregularstardiagsemigrps}
		Let $S$ be a diagram semigroup that is closed under the vertical flip involution $^*$. For each $\alpha, \beta \in S$, $(\alpha, \beta) \in \mathcal{D}$ if and only if $\alpha$ and $\beta$ have compatible patterns.
		
		\begin{proof}
			Suppose $(\alpha, \beta) \in \mathcal{D}$, and hence that there exist $\gamma \in S$ such that $\alpha^*\alpha = \gamma^*\gamma$ and $\gamma\gamma^* = \beta\beta^*$. Then $\uppat{\gamma} = \uppat{\gamma\gamma^*} = \uppat{\beta\beta^*} = \uppat{\beta}$ and $\lowpat{\gamma} = \lowpat{\gamma^*\gamma} = \lowpat{\alpha^*\alpha} = \lowpat{\alpha}$, hence $\alpha$ and $\beta$ have compatible patterns. 
			
			Conversely suppose $\alpha$ and $\beta$ have compatible patterns, and hence that there exist $\gamma \in S$ such that $\uppat{\gamma} = \uppat{\beta}$ and $\lowpat{\gamma} = \lowpat{\alpha}$. It follows from Proposition \ref{prop:productuplowpatwithsamemiddlepat} that $\uppat{\gamma\gamma^*} = \uppat{\gamma} = \uppat{\beta} = \uppat{\beta\beta^*}$ and $\lowpat{\gamma^*\gamma} = \lowpat{\gamma} = \lowpat{\alpha} = \lowpat{\alpha^*\alpha}$. Hence by Proposition \ref{cor:greensrltnsforregulardiagramsemigrps} $\beta\beta^* = \gamma\gamma^*$ and $\gamma^*\gamma = \alpha^*\alpha$, giving us $(\alpha, \beta) \in \mathcal{D}$.
		\end{proof}
	\end{theorem}
	
\begin{figure}[!ht]
		\caption[ ]{Product table for a subsemigroup of $\psyminvmon{k}$ where: pattern compatibility is reflexive but not symmetric; and there exist distinct elements that are not $\mathcal{D}$ related despite their upper and lower patterns being compatible.}
		\label{fig:nonregeg}
		\vspace{5pt}
		\centering
		\input{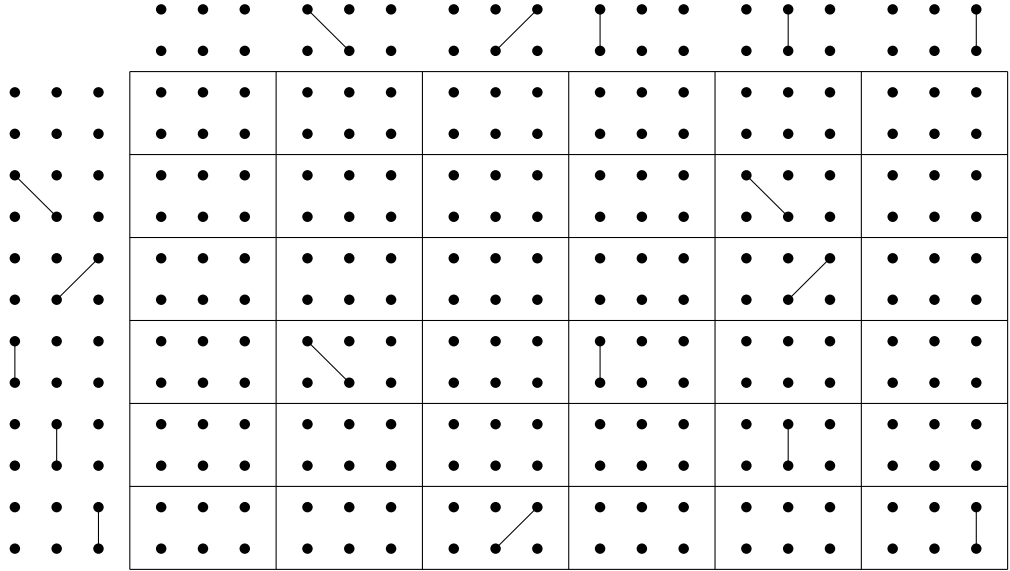}
	\end{figure}

	Figure \ref{fig:nonregeg} establishes that Theorem \ref{thm:greensDrltnonregularstardiagsemigrps} may not be generalised for diagram semigroups that are not closed under the vertical flip involution $^*$. In particular, the upper pattern of the second diagram is compatible with the lower pattern of the third element and vice versa, but the second and third elements are not $\mathcal{D}$ related. 
	
		\section{Green's $\mathcal{D}$ relation on $\pmodmon{m}{k}$ and $\modmon{m}{k}$}			
			Note that the number of $\mathcal{D}$ classes for the planar mod-$m$ monoid was established more thoroughly in \cite{art:Ahmed:OnTheNoOfPrincipalIdealsInDTonalPartMons}. This section will outline how the author reached the same numbers independently.
		
\begin{figure}[!ht]
				\caption[ ]{Dot $\mathcal{D}$ classes for $\pmodmon{2}{4}$ (left) and  $\pmodmon{2}{5}$ (right).}
				\label{fig:dclassespmodmon24and25}
				\vspace{5pt}
				\begin{center}
				\includegraphics[scale = 0.4]{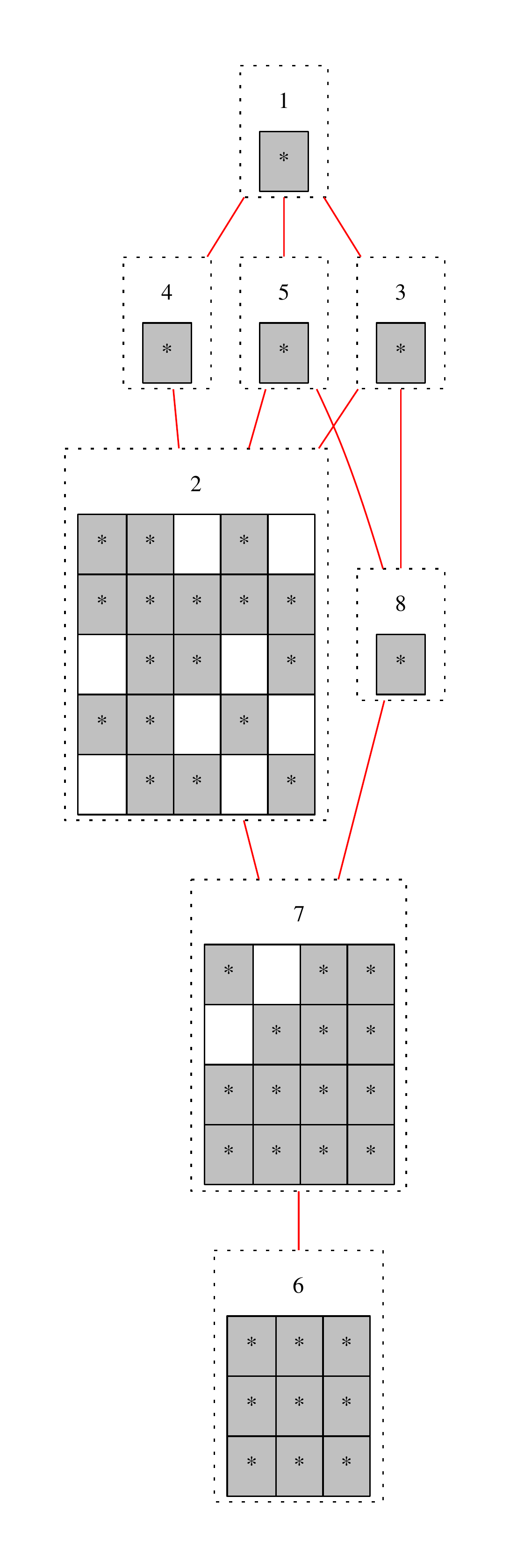}
				\hspace{1cm}
				\includegraphics[scale = 0.3]{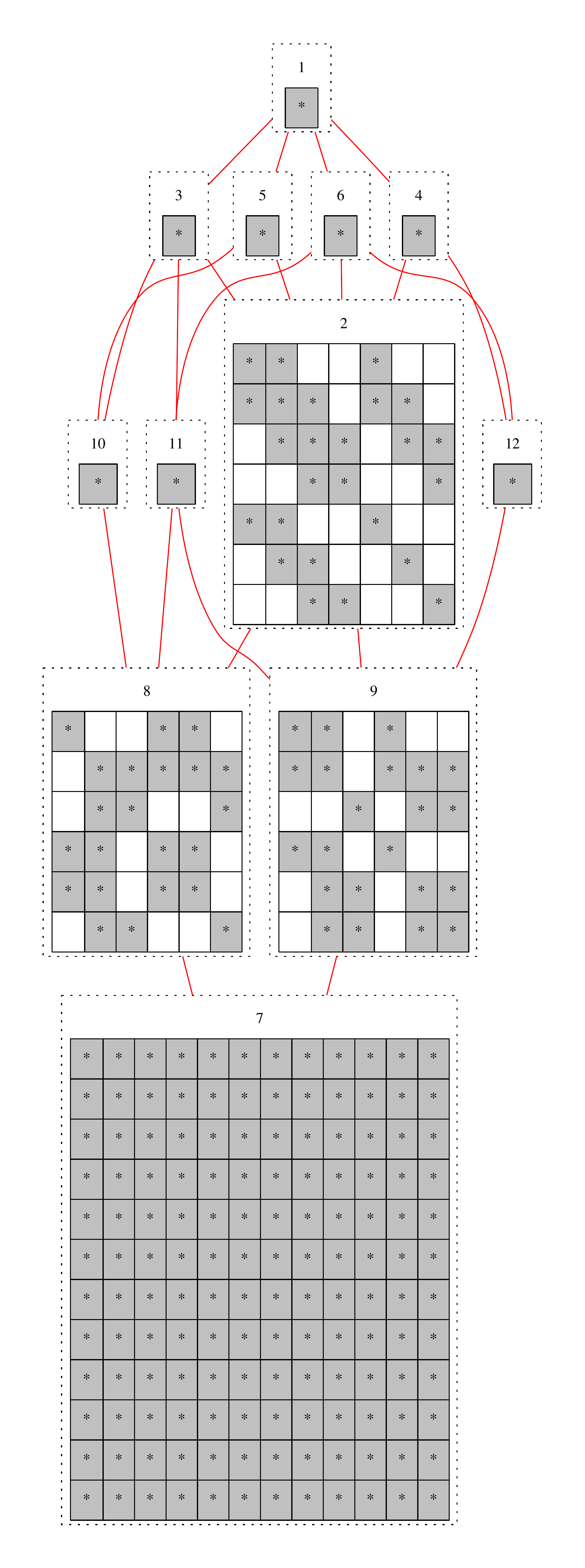}
				\end{center}
			\end{figure}

			It follows directly from Theorem \ref{thm:greensDrltnonregularstardiagsemigrps} that elements of $\pmodmon{m}{k}$ (or $\modmon{m}{k}$) are $\mathcal{D}$ related if and only if they have compatible patterns. Figures \ref{fig:dclassespmodmon22}, \ref{fig:dclassespmodmon23} and \ref{fig:dclassespmodmon24and25} along with Appendices \ref{app:pmoddotD} and \ref{app:moddotD} contain various examples of depictions of Green's $\mathcal{D}$ classes for $\pmodmon{m}{k}$ and $\modmon{m}{k}$ (Note: \textit{Dot D classes} are defined in Subsection 3.8-1 of the manual for the GAP Semigroups package \cite{man:GAPsemigrps}). It further follows that the number of Green's $\mathcal{D}$ classes for either the planar mod-$m$ monoid $\pmodmon{m}{k}$ or $\modmon{m}{k}$, is equal to the number of equivalence classes under pattern compatibility.
			
			\begin{proposition} \label{prop:pmodmoncompatiblepatterns}
				Given $m \in \posints$, $k \in \intsge{m}$ and $p, q \in \uppat{\pmodmon{m}{k}}$, let $r_p = \rank{p}$, $r_q = \rank{q}$ and $s_1, \ldots, s_{r_p}, t_1, \ldots, t_{r_q}$ denote the size of each transversal block, from left to right, contained in $p$ and $q$ respectively.
				\begin{enumerate}
					\item $\congmod{\Sigma^{r_p}_{i=1} s_i, \Sigma^{r_q}_{i=1} t_i}{k}{m}$; and
					\item $p$ and $q$ are $\pmodmon{m}{k}$-compatible if and only if $r_p = r_q$ and for each $i \in \set{1, \ldots, r_p = r_q}$, \newline $\congmod{s_i}{t_i}{m}$.
				\end{enumerate}
				
				\begin{proof}
					It follows from $p$ and $q$ being $\pmodmon{m}{k}$-admissible that the sizes of non-transversal blocks in either $p$ or $q$ are integer multiples of $m$, and consequently that $\congmod{\Sigma^{r_p}_{i=1} s_i, \Sigma^{r_q}_{i=1} t_i}{k}{m}$. Condition (ii) follows trivially from the definition of the planar mod-$m$ monoid.
				\end{proof}
			\end{proposition}
			
			\begin{figure}[!ht]
				\caption[ ]{Given $m=2$ and $k=2$, $\mathcal{D}$ classes for planar mod-$2$ monoid $\pmodmon{2}{2}$ are:}
				\label{fig:dclassespmodmon22}
				\vspace{5pt}
				\centering
				\input{chap_greens_relations/tikz/fig-dclassespmodmon22.tex}
			\end{figure}
	
			\begin{figure}[!ht]
				\caption[ ]{Given $m=2$ and $k=3$, $\mathcal{D}$ classes for planar mod-$2$ monoid $\pmodmon{2}{3}$ are:}
				\label{fig:dclassespmodmon23}
				\vspace{5pt}
				\centering
				\input{chap_greens_relations/tikz/fig-dclassespmodmon23.tex}
			\end{figure}
			
			\begin{proposition} \label{prop:numberofpmodmonDclasses}
				For each $m \in \posints$ and $k \in \nonnegints$, the number of Green's $\mathcal{D}$ classes for the planar mod-$m$ monoid $\pmodmon{m}{k}$ is given by (see Table \ref{table:nopmodmonDclasses} for example values):
				\[\card{\pmodmon{m}{k}/\mathcal{D}} = \begin{cases}2^{\max\set{0, k-1}} & m > k; \\ \card{\pmodmon{m}{k-m}/\mathcal{D}} + \ointparts{m, k} & m \leq k.\end{cases}\]
				
				\begin{proof}
					If $m > k$ then the planar mod-$m$ monoid $\pmodmon{m}{k}$ is equal to the planar uniform block bijections $\puniblockbijmon{k}$, which has $2^{\max\set{0, k-1}}$ $\mathcal{D}$ classes.
					If $m \leq k$ then it follows from Proposition \ref{prop:pmodmoncompatiblepatterns} that Green's $\mathcal{D}$ classes for the planar mod-$m$ monoid $\pmodmon{m}{k}$ may be indexed by ordered partitions of integers from $\set{k' \in \set{0, \ldots, k}: \congmod{k'}{k}{m}}$ into parts of size less than or equal to $m$. Which may be broken up into the sum of the number of ordered partitions of integers from $\set{k' \in \set{0, \ldots, k-m}: \congmod{k'}{k}{m}}$ into parts of size less than or equal to $m$, which is given by $\card{\pmodmon{m}{k-m}/\mathcal{D}}$, and the number of ordered partitions of $k$ into parts of size less than or equal to $m$, which is given by $\ointparts{m, k}$.
				\end{proof}
			\end{proposition}
			
	        \begin{table}[!ht]
	            \caption[ ]{Example values for $\card{\pmodmon{m}{k}/\mathcal{D}}$, calculated using the recurrence in Proposition \ref{prop:numberofpmodmonDclasses}.}
	            \label{table:nopmodmonDclasses}
	            \centering
	            \begin{tabular}{|c|rrrrrrrrrrr|}
	                \hline
	                \diagbox{m}{k} & 0 & 1 & 2 & 3 & 4 & 5 & 6 & 7 & 8 & 9 & 10 \\
	                \hline
					1 & 1 & 2 & 3 & 4 &  5 &  6 &  7 &   8 &   9 &  10 &  11 \\
					2 & 1 & 1 & 3 & 4 & 8 & 12 & 21 & 33 &  55 &  88 & 144 \\
					3 & 1 & 1 & 2 & 5 & 8 & 15 & 29 & 52 &  96 & 178 & 326 \\
					4 & 1 & 1 & 2 & 4 & 9 & 16 & 31 & 60 & 117 & 224 & 432 \\
					5 & 1 & 1 & 2 & 4 & 8 & 17 & 32 & 63 & 124 & 244 & 481 \\
	                \hline
	            \end{tabular}
	        \end{table}
	        
	        Note that the example values computed for $\card{\pmodmon{m}{k}/\mathcal{D}}$, which are given in Table \ref{table:nopmodmonDclasses}, match with the numbers given in \cite{art:Ahmed:OnTheNoOfPrincipalIdealsInDTonalPartMons}. Further note that the numbers in the second row appear to align with sequence $A052952$ on the online encyclopaedia of integer sequences \cite{man:OEIS}.
	        
			\begin{proposition} \label{prop:modmoncompatiblepatterns}
				Given $m \in \posints$, $k \in \intsge{m}$ and $p, q \in \uppat{\modmon{m}{k}}$, let $r_p = \rank{p}$, $r_q = \rank{q}$ and $s_1, \ldots, s_{r_p}, t_1, \ldots, t_{r_q}$ denote the size of each transversal block, from left to right when considering the left-most vertex in each block, contained in $p$ and $q$ respectively.
				\begin{enumerate}
					\item $\congmod{\Sigma^{r_p}_{i=1} s_i, \Sigma^{r_q}_{i=1} t_i}{k}{m}$; and
					\item $p$ and $q$ are $\modmon{m}{k}$-compatible if and only if $r_p = r_q$ and there exist $\pi \in \symgrp{r_p}$ such that for each $i \in \set{1, \ldots, r_p = r_q}$, $\congmod{s_{\pi^{-1}(i)}}{t_{\pi^{-1}(i)}}{m}$.
				\end{enumerate}
				
				\begin{proof}
					It follows from $p$ and $q$ being $\modmon{m}{k}$-admissible that the sizes of non-transversal blocks in either $p$ or $q$ are integer multiples of $m$, and consequently that $\congmod{\Sigma^{r_p}_{i=1} s_i, \Sigma^{r_q}_{i=1} t_i}{k}{m}$. Condition (ii) follows trivially from the definition of the mod-$m$ monoid.
				\end{proof}
			\end{proposition}
	        
			\begin{proposition} \label{prop:numberofmodmonDclasses}
				For each $m \in \posints$ and $k \in \nonnegints$, the number of Green's $\mathcal{D}$ classes for the mod-$m$ monoid $\modmon{m}{k}$ is given by (see Table \ref{table:nomodmonDclasses} for example values):
				\[\card{\modmon{m}{k}/\mathcal{D}} = \begin{cases}\intparts{m, k} & m > k; \\ \card{\modmon{m}{k-m}/\mathcal{D}} + \intparts{m, k} & m \leq k.\end{cases}\]
				
				\begin{proof}
					If $m > k$ then the mod-$m$ monoid $\modmon{m}{k}$ is equal to the uniform block bijections $\uniblockbijmon{k}$, which has $\intparts{m, k}$ $\mathcal{D}$ classes.
					If $m \leq k$ then it follows from Proposition \ref{prop:modmoncompatiblepatterns} that Green's $\mathcal{D}$ classes for the mod-$m$ monoid $\modmon{m}{k}$ may be indexed by partitions of integers from $\set{k' \in \set{0, \ldots, k}: \congmod{k'}{k}{m}}$ into parts of size less than or equal to $m$. Which may be broken up into the sum of the number of partitions of integers from $\set{k' \in \set{0, \ldots, k-m}: \congmod{k'}{k}{m}}$ into parts of size less than or equal to $m$, which is given by $\card{\modmon{m}{k-m}/\mathcal{D}}$, and the number of partitions of $k$ into parts of size less than or equal to $m$, which is given by $\intparts{m, k}$.
				\end{proof}
			\end{proposition}
	        
	        \begin{table}[!ht]
	            \caption[ ]{Example values for $\card{\modmon{m}{k}/\mathcal{D}}$, calculated using the recurrence in Proposition \ref{prop:numberofmodmonDclasses}.}
	            \label{table:nomodmonDclasses}
	            \centering
	            \begin{tabular}{|c|rrrrrrrrrrr|}
	                \hline
	                \diagbox{m}{k} & 0 & 1 & 2 & 3 & 4 & 5 & 6 & 7 & 8 & 9 & 10 \\
	                \hline
					1 & 1 & 2 & 3 & 4 & 5 & 6 &  7 &  8 &  9 & 10 & 11 \\
					2 & 1 & 1 & 3 & 3 & 6 & 6 & 10 & 10 & 15 & 15 & 21 \\
					3 & 1 & 1 & 2 & 4 & 5 & 7 & 11 & 13 & 17 & 23 & 27 \\
					4 & 1 & 1 & 2 & 3 & 6 & 7 & 11 & 14 & 21 & 25 & 34 \\
					5 & 1 & 1 & 2 & 3 & 5 & 8 & 11 & 15 & 21 & 28 & 38 \\
	                \hline
	            \end{tabular}
	        \end{table}
	        
	        Note that the numbers in the second row of Table \ref{table:nomodmonDclasses} appear to align with sequence $A008805$ on the online encyclopaedia of integer sequences \cite{man:OEIS}, and the numbers in the third row appear to align with sequence $A028289$ (also noted by \cite{art:Ahmed:OnTheNoOfPrincipalIdealsInDTonalPartMons}).
	        
	\section{Green's $\mathcal{R}$ and $\mathcal{L}$ relations on $\pmodmon{m}{k}$ and $\modmon{m}{k}$}
			Recall from Corollary \ref{cor:greensrltnsforregulardiagramsemigrps} that for a subsemigroup $S$ of the partition monoid $\pttnmon{k}$ that is closed under the vertical flip involution $^*$ and $\alpha, \beta \in S$:
			\begin{enumerate}
				\item $(\alpha, \beta) \in \mathcal{R}$ if and only if $\uppat{\alpha} = \uppat{\beta}$; and
				\item $(\alpha, \beta) \in \mathcal{L}$ if and only if $\lowpat{\alpha} = \lowpat{\beta}$.
			\end{enumerate}
			
			Hence the number of Green's $\mathcal{R}$ classes for $S$ is equal to the number of $S$-admissible patterns.
			
			\begin{proposition} \label{prop:numberofpmodmonRclasses}
				For each $m \in \posints$ and $k \in \intsge{m}$, the number of Green's $\mathcal{R}$ classes, which is equal to the number of Green's $\mathcal{L}$ classes, for the planar mod-$m$ monoid $\pmodmon{m}{k}$ is given by (see Table \ref{table:nopmodmonRclasses} for example values):
				\[\card{\pmodmon{m}{k}/\mathcal{R}} = \card{\uppat{\pmodmon{m}{k}}} = \begin{cases} 2^{\max\set{0, k-1}} & m > k; \\ \Sigma^{\floor{k/m}}_{u=0} \pfeasnontrans{m, k, mu}2^{\max\set{0, k-mu-1}} & m \leq k. \end{cases}\]
				
				\begin{proof}
					If $m > k$ then the planar mod-$m$ monoid $\pmodmon{m}{k}$ is equal to the planar uniform block bijections $\puniblockbijmon{k}$, which has $2^{\max\set{0, k-1}}$ $\mathcal{R}$ classes. For $m \leq k$, recall from Definition \ref{def:pmodmonfeasnontrans} that for each $u \in \set{0, \ldots, \floor{\frac{k}{m}}}$, $\pfeasnontrans{m, k, mu}$ denotes the number of $\pmodmon{m}{k}$-admissible patterns that have $mu$ of $k$ vertices contained in non-transversal blocks. Since transversal blocks must contain consecutive vertices from the remaining $k-mu$ vertices, the number of ways for the remaining $k-mu$ vertices to be contained in transversal blocks is equal to the number of ordered partitions of $k-mu$, of which there are $2^{\max\set{0, k-mu-1}}$.
				\end{proof}
			\end{proposition}
			
			\begin{table}[!ht]
				\caption[ ]{Example values for $\card{\pmodmon{m}{k}/\mathcal{R}}$, calculated using the recurrence in Proposition \ref{prop:numberofpmodmonRclasses}.}
				\label{table:nopmodmonRclasses}
				\centering
				\begin{tabular}{|c|rrrrrrrrrr|}
			\hline
	\diagbox{$m$}{$k$} & 1 & 2 &  3 &  4 &   5 &   6 &    7 &     8 &     9 &     10 \\ 
			\hline
	        1          & 2 & 6 & 20 & 70 & 252 & 924 & 3432 & 12870 & 48620 & 184756 \\
	        2          & 1 & 3 &  6 & 17 &  39 & 108 &  262 &   717 &  1791 &   4867 \\
	        3          & 1 & 2 &  5 & 10 &  22 &  52 &  113 &   254 &   590 &   1316 \\
	        4          & 1 & 2 &  4 &  9 &  18 &  38 &   80 &   173 &   363 &    772 \\
	        5          & 1 & 2 &  4 &  8 &  17 &  34 &   70 &   144 &   296 &    614 \\ 
	        \hline
				\end{tabular}
			\end{table}
			
			Note that the numbers in the first row of Table \ref{table:nopmodmonRclasses}, which are the number of Green's $\mathcal{R}$ classes for the planar partition monoid $\ppttnmon{k}$, appear to align with the central binomial coefficients (sequence $A000984$ on the OEIS \cite{man:OEIS}).

			\begin{proposition} \label{prop:numberofmodmonRclasses}
				For each $m \in \posints$ and $k \in \intsge{m}$, the number of Green's $\mathcal{R}$ classes, which is equal to the number of Green's $\mathcal{L}$ classes, for the mod-$m$ monoid $\modmon{m}{k}$ is given by (see Table \ref{table:nomodmonRclasses} for example values):
				\[\card{\modmon{m}{k}/\mathcal{R}} = \card{\uppat{\modmon{m}{k}}} = \begin{cases}\bellnos{k} & m > k; \\ \Sigma^{\floor{k/m}}_{u=0} \feasnontrans{m, k, mu}\bellnos{k-mu} & m \leq k. \end{cases}\]
				
				\begin{proof}
					If $m > k$ then the mod-$m$ monoid $\modmon{m}{k}$ is equal to the uniform block bijections $\uniblockbijmon{k}$, which has $\bellnos{k}$ $\mathcal{R}$ classes. For $m \leq k$, recall from Definition \ref{def:pmodmonfeasnontrans} that for each $u \in \set{0, \ldots, \floor{\frac{k}{m}}}$, $\pfeasnontrans{m, k, mu}$ denotes the number of $\pmodmon{m}{k}$-admissible patterns that have $mu$ of $k$ vertices contained in non-transversal blocks. Since each transversal block may contain any of the remaining $k-mu$ vertices, the number of ways for the remaining $k-mu$ vertices to be contained in transversal blocks is equal to the number of set partitions of $k-mu$ vertices, which is given by the Bell number $\bellnos{k-mu}$.
				\end{proof}
			\end{proposition}
			
			\begin{table}[!ht]
				\caption[ ]{Example values for $\card{\modmon{m}{k}/\mathcal{R}}$, calculated using the recurrence in Proposition \ref{prop:numberofmodmonRclasses}.}
				\label{table:nomodmonRclasses}
				\centering
				\begin{tabular}{|c|rrrrrrrrrr|}
			\hline
			\diagbox{$m$}{$k$} & 1 & 2 & 3 & 4 & 5 & 6 & 7 & 8 & 9 & 10 \\
			\hline
			1 & 2 & 6 & 22 & 94 & 454 & 2430 & 14214 & 89918 & 610182 & 4412798 \\
			2 & 1 & 3 &  8 & 31 & 122 &  579 &  2886 & 16139 &  95358 &  611111 \\
			3 & 1 & 2 &  6 & 19 &  72 &  314 &  1479 &  7668 &  43184 &  259515 \\
			4 & 1 & 2 &  5 & 16 &  57 &  233 &  1052 &  5226 &  28023 &  161845 \\
			5 & 1 & 2 &  5 & 15 &  53 &  209 &   919 &  4420 &  23037 &  129206 \\
	        \hline
				\end{tabular}
			\end{table}
			
			Note that the numbers in the first row of Table \ref{table:nomodmonRclasses}, which are the number of Green's $\mathcal{R}$ classes for the partition monoid $\pttnmon{k}$, appear to align with sequence $A001861$ on the OEIS \cite{man:OEIS}.
    \clearpage{\pagestyle{empty}\cleardoublepage}
\chapter{Presentations} \label{chap:presentations}
	\section{Planar modular-$2$ partition monoid $\pmodmon{2}{k}$}
		In this section we conjecture relations that appear to give, when combined with the appropriate diapsis and $(2,2)$-transapsis generators, a presentation of the planar mod-$2$ monoid $\pmodmon{2}{k}$ for all $k \in \intsge{2}$. 
		
		\begin{conjecture} \label{conj:pmod2presentation}
			For each $k \in \intsge{2}$, the planar mod-$2$ monoid $\pmodmon{2}{k}$ is characterised by the generators $\set{\transapgen{i}, \hookgen{i}: i = 1, \ldots, k-1}$ along with the relations:
			\begin{enumerate}
				\item $\transapgen{i}^2 = \transapgen{i}$;
				\item $\transapgen{j}\transapgen{i} = \transapgen{i}\transapgen{j}$ for all $|j - i| \geq 1$;
				\item $\hookgen{i}^2 = \hookgen{i}$;
				\item $\hookgen{i}\hookgen{j}\hookgen{i} = \hookgen{i}$ for all $|j-i| = 1$;
				\item $\hookgen{j}\hookgen{i} = \hookgen{i}\hookgen{j}$ for all $|j-i| \geq 2$;
				\item $\transapgen{i}\hookgen{i} = \hookgen{i}\transapgen{i} = \hookgen{i}$;
				\item $\transapgen{i}\hookgen{j}\transapgen{i} = \transapgen{i}\transapgen{j}$ for all $j - i = 1$; and
				\item $\transapgen{j}\hookgen{i} = \hookgen{i}\transapgen{j}$ for all $ |j-i| \geq 2$.
			\end{enumerate}
		\end{conjecture}
		
		The author was able to verify Conjecture \ref{conj:pmod2presentation} up to $k=7$ using GAP (see Appendix \ref{app:pmod2presentations} for code). Since the generators of the planar mod-$2$ monoid $\pmodmon{2}{k}$ all commute for differences in indices greater than or equal to two, the author would be extremely surprised if Conjecture \ref{conj:pmod2presentation} does not hold.
		
		While the author has not been able to establish the result in Conjecture \ref{conj:pmod2presentation} for $k \in \intsge{8}$, we will proceed to establish a number of results which may be of use to the eventual establishment for all $k \in \intsge{2}$. We begin by noting a number of further relations that are implied by the relations in Conjecture \ref{conj:pmod2presentation}, then establish a non-equivalent upper bound on reduced $\pmodmon{2}{k}$-words in normal form, with a number of further conjectures for enumerating the number of reduced $\pmodmon{2}{k}$-words in normal form contained within the established upper bound.
		
		\begin{proposition} \label{prop:morepmod2monrltns}
			The following relations are implied by Relations (i) - (vii) from Conjecture \ref{conj:pmod2presentation}:
			\begin{enumerate}
				\item $\hookgen{i}\transapgen{i+1}\hookgen{i} = \hookgen{i}$;
				\item $\hookgen{i+1}\hookgen{i}\transapgen{i+1} = \hookgen{i+1}\transapgen{i}$;
				\item $\transapgen{i+1}\hookgen{i}\hookgen{i+1} = \transapgen{i}\hookgen{i+1}$;
				\item $\transapgen{i+1}\hookgen{i}\transapgen{i+1} = \transapgen{i}\transapgen{i+1}$;
				\item $\transapgen{i}\hookgen{i+1}\hookgen{i} = \transapgen{i+1}\hookgen{i}$;
				\item $\hookgen{i}\hookgen{i+1}\transapgen{i} = \hookgen{i}\transapgen{i+1}$; and
				\item $\hookgen{i+1}\transapgen{i}\hookgen{i+1} = \hookgen{i+1}$.
			\end{enumerate}
			
			\begin{proof}
				For each $i, j \in \set{1, \ldots, k}$ such that $j - i = 1$:
				\begin{enumerate}
					\item $(\hookgen{i})\transapgen{j}\hookgen{i} = \hookgen{i}(\transapgen{i}\transapgen{j})\hookgen{i} = (\hookgen{i}\transapgen{i})\hookgen{j}(\transapgen{i}\hookgen{i}) = \hookgen{i}\hookgen{j}\hookgen{i} = \hookgen{i}$;
					\item $\hookgen{j}(\hookgen{i})\transapgen{j} = \hookgen{j}\hookgen{i}(\transapgen{i}\transapgen{j}) = \hookgen{j}(\hookgen{i}\transapgen{i})\hookgen{i}\transapgen{i} = (\hookgen{j}\hookgen{i}\hookgen{j})\transapgen{i} = \hookgen{j}\transapgen{i}$;
					\item $\transapgen{j}(\hookgen{i})\hookgen{j} = (\transapgen{j}\transapgen{i})\hookgen{i}\hookgen{j} = (\transapgen{i}\transapgen{j})\hookgen{i}\hookgen{j} = \transapgen{i}\hookgen{j}(\transapgen{i}\hookgen{i})\hookgen{j} = \transapgen{i}(\hookgen{j}\hookgen{i}\hookgen{j}) = \transapgen{i}\hookgen{j}$;
					\item $\transapgen{j}(\hookgen{i})\transapgen{j} = (\transapgen{j}\hookgen{i}\hookgen{j})\hookgen{i}\transapgen{j} = \transapgen{i}(\hookgen{j}\hookgen{i}\transapgen{j}) = (\transapgen{i}\hookgen{j}\transapgen{i}) = \transapgen{i}\transapgen{j}$;
					\item $\transapgen{i}(\hookgen{j})\hookgen{i} = (\transapgen{i}\transapgen{j})\hookgen{j}\hookgen{i} = \transapgen{j}\hookgen{i}(\transapgen{j}\hookgen{j})\hookgen{i} = \transapgen{j}\hookgen{i}\hookgen{j}\hookgen{i} = \transapgen{j}\hookgen{i}$;
					\item $\hookgen{i}(\hookgen{j})\transapgen{i} = \hookgen{i}\hookgen{j}(\transapgen{j}\transapgen{i}) = \hookgen{i}\hookgen{j}(\transapgen{i}\transapgen{j}) = \hookgen{i}(\hookgen{j}\transapgen{j})\hookgen{i}\transapgen{j} = (\hookgen{i}\hookgen{j}\hookgen{i})\transapgen{j} = \hookgen{i}\transapgen{j}$; and
					\item $\hookgen{j}\transapgen{i}(\hookgen{j}) = \hookgen{j}(\transapgen{i}\transapgen{j})\hookgen{j} = (\hookgen{j}\transapgen{j})\hookgen{i}(\transapgen{j}\hookgen{j}) = (\hookgen{j}\hookgen{i}\hookgen{j}) = \hookgen{j}$.
				\end{enumerate}
			\end{proof}
		\end{proposition}
		
		\subsection{A bound on reduced $\pmodmon{2}{k}$-words} \label{subsec:pmod2boundednormals}
			\begin{definition}
				We refer to elements of the free semigroup of the planar mod-$2$ monoid $\freesemigrp{(\pmodmon{2}{k})}$ as \textit{$\pmodmon{2}{k}$-words}, and say that a $\pmodmon{2}{k}$-word is \textit{reduced} if it may not be written with fewer generators using the relations from Conjecture \ref{conj:pmod2presentation}.
			\end{definition}

			\begin{definition}
				For each $i \in \set{1, \ldots, k-1}$ we denote by $G_i$ the set $\set{\hookgen{i}, \transapgen{i}}$ containing the $i$th diapsis generator and the $i$th $(2,2)$-transapsis generator.
			\end{definition}
			
			\begin{corollary}\label{cor:GiGjGiisreducible}
				For each $k \in \intsge{2}$ and $i, j \in \{1, ..., k-1\}$ such that $|j-i| = 1$, every $\pmodmon{2}{k}$-word in $G_iG_jG_i$ is reducible.
				
				\begin{proof}
					Follows from Relation (iv) in Conjecture \ref{conj:pmod2presentation} along with the relations established in Proposition \ref{prop:morepmod2monrltns}.
				\end{proof}
			\end{corollary}
		
			Note every $\pmodmon{2}{k}$-word may be written in the form $g_{1, i_1}\ldots g_{n, i_n}$ where for each $j \in \set{1, \ldots, n}$, $i_j \in \set{1, \ldots, k-1}$ and $g_{j, i_j} \in G_{i_j}$. For example given the $\pmodmon{2}{4}$-word $\hookgen{3}\transapgen{2}\hookgen{1}$, letting $i_1 = 3$, $i_2 = 2$, $i_3 = 1$ $g_{1, 3} = \hookgen{3}$, $g_{2, 2} = \transapgen{2}$ and $g_{3, 1} = \hookgen{1}$, we have $\hookgen{3}\transapgen{2}\hookgen{1} = g_{1, i_1}g_{2, i_2}g_{3, i_3}$ where for each $j \in \set{1, 2, 3}$, $i_j \in \set{1, 2, 3}$ and $g_{j, i_j} \in G_{i_j}$.
			
			\begin{proposition}
				In any reduced $\pmodmon{2}{k}$-word $g_{1, i_1}\ldots g_{n, i_n}$, the maximal index $m = \max\big\{i_j: j \in \{1, \ldots, n\}\big\}$ occurs precisely once.
				\begin{proof}
					Similar to Lemma 2.2 in \cite{art:Ridout:StandardModules} we use induction on the maximal index of our reduced word. When $m=1$, if any $g_{j, 1} = \hookgen{1}$ then $g_{1, 1}\ldots g_{n, 1} = \hookgen{1}$, otherwise $g_{1, 1}\ldots g_{n, 1} = \transapgen{1}$. Let $m > 1$ with our statement holding for all $m' \in \set{1, \ldots, m-1}$. Suppose there exist $j \neq j'$ such that $i_j = m = i_{j'}$, and hence that our reduced word has the form $\ldots\ g_{j, m}Wg_{j', m}\ldots$ where $W$ is either the empty product (when $j' = j+1$) or a reduced $\pmodmon{2}{k}$-word whose maximal index $m'$ satisfies $m' < m$. We proceed to argue that the existence of $j, j'$ forms a contradiction. If $W$ is the empty product then $\ldots\ g_{j, m}g_{j+1, m}\ldots$ may trivially be further reduced, which contradicts the existence of $j, j'$. If $m' < m-1$ then we may trivially commute $W$ with $g_{j', m}$, allowing $\ldots\ g_{j, m}Wg_{j', m}\ldots = \ldots\ g_{j, m}g_{j', m}W\ldots$, which may further be reduced again contradicting the existence of $j, j'$. Finally, if $m' = m-1$, then there exists $q \in \set{j+1, \ldots, j'-1}$ such that $W = W_Lg_{q, m'}W_R$ where the maximal indices $m_L$ of $W_L$ and $m_R$ of $W_R$ satisfy $m_L, m_R < m-1$. Hence we may write $\ldots\ g_{j, m}Wg_{j', m}\ldots = \ldots\ g_{j, m}W_Lg_{q, m-1}W_Rg_{j', m}\ldots = \ldots\ W_Lg_{j, m}g_{q, m-1}g_{j', m}W_R\ldots$. Yet again the existence of $j, j'$ is contradicted since $g_{j, m}g_{q, m-1}g_{j', m} \in G_mG_{m-1}G_m$ and, as established in Corollary \ref{cor:GiGjGiisreducible}, every element of $G_mG_{m-1}G_m$ is reducible. Therefore no such $j, j'$ may exist, by induction the maximal index must occur precisely once.
				\end{proof}
			\end{proposition}
			
			\begin{proposition} \label{prop:firststepnormalisingreducedpmod2word}
				If $W$ is a reduced $\pmodmon{2}{k}$-word with maximal index $m$ then $W$ may be rewritten as $W = W'g_{m}\ldots g_{l}$ where $W'$ is a reduced $\pmodmon{2}{k}$-word with maximal index less than $m$, $l \in \set{1, \ldots, m}$ and for each $q \in \set{l, \ldots, m}$, $g_{l} \in G_{l}$.
				\begin{proof}
					If there were a gap in the sequence of indices
					following $W'$, then the elements to the right of the gap all commute with the generators to the left of the gap (up to $g_{m}$), so may be relocated to the left of $g_{m}$ then absorbed in to $W'$. If the sequence of indices following $W'$ increased then $m$ would not be the maximal index, furthermore since the maximal index occurs precisely once, the sequence of indices following $W'$ may not remain constant.
				\end{proof}
			\end{proposition}
			
			Note that the same process may be repeated on $W'$ in Proposition \ref{prop:firststepnormalisingreducedpmod2word}.
			
			\begin{proposition} \label{prop:pmod2normalbound}
				For each $k \in \intsge{2}$, any reduced $\pmodmon{2}{k}$-word $W$ may be rewritten as $W = r_{j_1, i_1}\ldots r_{j_n, i_n}$ where:
				\begin{enumerate}
					\item $n \in \posints$ and $i_1, \ldots, i_n, j_1, \ldots, j_n \in \set{1, \ldots, k-1}$ such that: \begin{enumerate}
						\item for each $l \in \set{1, \ldots, n}$, $j_l \geq i_l$;
						\item for each $l \in \set{1, \ldots, n-1}$, $i_l < i_{l+1}$ and $j_l < j_{l+1}$, and
					\end{enumerate}
					\item for each $l \in \set{1, \ldots, n}$, $r_{j_l, i_l} \in G_{j_l}\ldots G_{i_l}$.
				\end{enumerate}
				
				\begin{proof}
					Repeatedly applying the process in Proposition \ref{prop:firststepnormalisingreducedpmod2word} gives us the increasing nature of the $j_l$, as they are the maximal indices at each step. Suppose there exists $l \in \set{1, \ldots, n-1}$ such that $i_l \geq i_{l+1}$. Let $r_{j_l, i_l} = g_{j_l}\ldots g_{i_l}$ and $r_{j_{l+1}, i_{l+1}} = h_{j_{l+1}}\ldots h_{i_{l+1}}$, since $i_l \leq j_l < j_{l+1}$ we have $r_{j_l, i_l}r_{j_{l+1}, i_{l+1}} = g_{j_l}\ldots g_{i_l}h_{j_{l+1}}\ldots h_{i_{l}}\ldots h_{i_{l+1}} = g_{j_l}\ldots g_{i_l+1}h_{j_{l+1}}\ldots (g_{i_l}h_{i_{l}+1}h_{i_{l}})\ldots h_{i_{l+1}}$. That $g_{i_l}h_{i_{l}+1}h_{i_{l}}$ may be reduced contradicts our word being reduced. Thus we must also have $i_l < i_{l+1}$ for all $l \in \set{1, \ldots, n-1}$.
				\end{proof}
			\end{proposition}
			
			Hence a $\pmodmon{2}{k}$-word is written in the form required for Proposition \ref{prop:pmod2normalbound} when the replacement of any $(2,2)$-transapsis generators with the same index diapsis generators forms a $\jonesmon{k}$-word in normal form (from Definition \ref{def:Jonesnormalform}). If we think of diapsis and $(2,2)$-transapsis generators as being lower and upper case analogues of each other, then this may be thought of as replacing any instances of upper case letters with the appropriate lower case letters.
			
			Unfortunately distinct reduced $\pmodmon{2}{k}$-words in the form required for Proposition \ref{prop:pmod2normalbound} may produce the same bipartition once the product has been performed. For example $\transapgen{2}\transapgen{1} = \transapgen{1}\transapgen{2}$ and are distinct reduced $\pmodmon{2}{k}$-words in the form required for Proposition \ref{prop:pmod2normalbound}.
			
			Ideally we want to establish an exhaustive set of unique reduced $\pmodmon{2}{k}$-words, to be identified as $\pmodmon{2}{k}$-words in normal form, much like $\jonesmon{k}$-words in normal form from Definition \ref{def:Jonesnormalform}. One option is to use the lexicographical ordering of $\pmodmon{2}{k}$-words induced by a total ordering of the generators, then assign to each element $\alpha \in \pmodmon{2}{k}$ the lexicographically lowest-ordered reduced $\pmodmon{2}{k}$-word corresponding to $\alpha$. 
			
			When the analogous procedure is done on $\jonesmon{k}$-words, in order for the selected reduced $\jonesmon{k}$-words to be in the normal form identified in Definition \ref{def:Jonesnormalform}, the generators must be ordered by $\hookgen{i} \leq \hookgen{j}$ if and only if $i \leq j$. If, for example, we used the inverse order then $\hookgen{3}\hookgen{1}$ would be assigned over $\hookgen{1}\hookgen{3}$, while only the latter is a $\jonesmon{k}$-word in normal form. Similarly, for the selected reduced $\pmodmon{2}{k}$-word to be in the form required for Proposition \ref{prop:pmod2normalbound} then we need lower indexed generators to be ordered lower than higher indexed generators, however the order of generators at each fixed index is not as important.
	
			Appendix \ref{app:pmod2words} contains candidates for $\pmodmon{2}{k}$-words in normal form using the above procedure both when diapsis generators are ordered lower than $(2,2)$-transapsis generators at each fixed index and vice-versa, and where diapsis and $(2,2)$-transapsis generators have been replaced by lower and upper case letters respectively from the English alphabet. Note that the candidate words have been ordered based on the run of decreasing indices that they end with, similar to how we considered the run of decreasing indices at the end of $\jonesmon{k}$-words in normal form in Definition \ref{def:noJknormalformsendingwithrun}.
			
			The first thing the author noticed with the candidates for $\pmodmon{2}{k}$-words in normal form from Appendix \ref{app:pmod2words} was that for each decreasing run of indices, no neighbouring indices both belong to $(2,2)$-transapsis generators.
		
			\begin{definition}
				For each $i, j \in \set{1, \ldots, k-1}$ such that $j > i$, denote by $\Run{j}{i}$ the set of all $g_j\ldots g_i \in G_j\ldots G_i$ such that there does not exist $l \in \set{i, \ldots, j-1}$ that satisfies $g_l = \transapgen{l}$ and $g_{l+1} = \transapgen{l+1}$.
			\end{definition}
			
			\begin{example} \label{eg:Runji}
				Examples of $\Run{j}{i}$ are:
				\begin{enumerate}
					\item $\Run{i}{i} = \{\hookgen{i},\ \transapgen{i}\}$;
					\item $\Run{i+1}{i} = \{\hookgen{i+1}\hookgen{i},\ \hookgen{i+1}\transapgen{i},\ \transapgen{i+1}\hookgen{i}\}$; and
					\item $\Run{i+2}{i} = \{\hookgen{i+2}\hookgen{i+1}\hookgen{i},\ \hookgen{i+2}\hookgen{i+1}\transapgen{i},\ \hookgen{i+2}\transapgen{i+1}\hookgen{i},\ \transapgen{i+2}\hookgen{i+1}\hookgen{i},\ \transapgen{i+2}\hookgen{i+1}\transapgen{i}\}$.
				\end{enumerate}
			\end{example}
			
			\begin{proposition}
				For each $i, j \in \posints$ such that $i \leq j$, $\card{\Run{j}{i}} = \fibno{j-i+3}$, that is the $(j-i+3)$th Fibonacci number.
				
				\begin{proof}
					Example \ref{eg:Runji} verifies the cases where $0 \leq j-i \leq 3$. The remaining cases follow inductively when noting that $\{\Run{j-1}{i}\hookgen{j}, \Run{j-2}{i}\hookgen{j-1}\transapgen{i}\}$ trivially partitions $\Run{j}{i}$.
				\end{proof}
			\end{proposition}
			
			\begin{conjecture} \label{conj:pmod2normalbound}
				For each $k \in \intsge{2}$, any reduced $\pmodmon{2}{k}$-word may be written in the form $r_{j_1, i_1}\ldots r_{j_n, i_n}$ where:
				\begin{enumerate}
					\item $n \in \posints$ and $i_1, \ldots, i_n, j_1, \ldots, j_n \in \set{1, \ldots, k-1}$ such that: \begin{enumerate}
						\item for each $l \in \set{1, \ldots, n}$, $j_l \geq i_l$;
						\item for each $l \in \set{1, \ldots, n-1}$, $i_l < i_{l+1}$ and $j_l < j_{l+1}$, and
					\end{enumerate}
					\item for each $l \in \set{1, \ldots, n}$, $r_{j_l, i_l} \in \Run{j_l}{i_l}$.
				\end{enumerate}
			\end{conjecture}
			
			Both sets of candidates for $\pmodmon{2}{k}$-words in normal form from Appendix \ref{app:pmod2words} satisfy Conjecture \ref{conj:pmod2normalbound}, as do both sets of candidates when $k=7$ (which was not included in Appendix \ref{app:pmod2words} due to the number of words requiring an unreasonable amount of space).
			
			Further note, even if Conjecture \ref{conj:pmod2normalbound} holds, there still exist words in the form required for Conjecture \ref{conj:pmod2normalbound} that are reducible (see Figure \ref{fig:boundednonreducedpmod2words} for an example), and distinct reduced $\pmodmon{2}{k}$-words in the required form whose products are equal (see Figure \ref{fig:boundedequivreducedpmod2words} for an example). Table \ref{table:candidateRjivalues} contains the number of candidate $\pmodmon{2}{k}$-words in normal form based on which $\Run{j}{i}$ the candidate word ends with a run from, once reducible and equivalent $\pmodmon{2}{k}$-words have been removed. Note that these numbers agree both when diapsis generators are ordered lower than $(2,2)$-transapsis generators at each fixed index and vice-versa from Appendix \ref{app:pmod2words}.
			
			\begin{table}[!ht]
				\caption[ ]{Number of planar mod-$2$ normal form words ending with a run from $\Run{j}{i}$.}
				\label{table:candidateRjivalues}
				\centering
				\begin{tabular}{| c | r r r r r r |}
					\hline
					\diagbox{$k_1$}{$k_2$} & 1 & 2 & 3 & 4 & 5 & 6 \\
					\hline
					1  & 2  &     &     &      &      &      \\
					2  & 3  & 6   &     &      &      &      \\
					3  & 5  & 14  & 24  &      &      &      \\
					4  & 8  & 32  & 68  & 110  &      &      \\
					5  & 13 & 65  & 183 & 348  & 546  &      \\
					6  & 21 & 128 & 428 & 1036 & 1855 & 2856 \\
					\hline
				\end{tabular}
			\end{table}

		 	\begin{figure}[!ht]
		  		\caption[ ]{$\pmodmon{2}{k}$-word in the form required for Conjecture \ref{conj:pmod2normalbound} that is reducible:}
		  		\label{fig:boundednonreducedpmod2words}
				\vspace{5pt}
				\centering
		  		\input{chap_presentations/tikz/fig-boundednonreducedpmod2words.tex}
			\end{figure}
			
			\begin{figure}[!ht]
		  		\caption[ ]{reduced $\pmodmon{2}{k}$-words in the form required for Conjecture \ref{conj:pmod2normalbound} that form the same product:}
		  		\label{fig:boundedequivreducedpmod2words}
				\vspace{5pt}
				\centering
		  		\input{chap_presentations/tikz/fig-boundedequivreducedpmod2words.tex}
			\end{figure}
			
			\begin{conjecture}
				The number of required $\pmodmon{2}{k}$-words in normal form ending with a run from $\Run{j}{1}$ is equal to the $(j+2)$th Fibonacci number (sequence $A000045$ on the OEIS \cite{man:OEIS}).
			\end{conjecture}
			
			\begin{conjecture}
				The number of required $\pmodmon{2}{k}$-words in normal form ending with a run from $\Run{j}{2}$ is equal to $(2j-1)\fibno{j+1} - 2^{j-2} + 1$ where $\fibno{j+1}$ is the $(j+1)$th Fibonacci number.
			\end{conjecture}
 
 			\begin{conjecture} \label{conj:wordsendingfromRunjj}
				The number of required $\pmodmon{2}{k}$-words in normal form ending with a run from $\Run{j}{j}$ forms the sequence $A046646$ on the OEIS \cite{man:OEIS} (that is, the number of certain rooted planar maps).
			\end{conjecture}
			
			It was noted on the OEIS \cite{man:OEIS} page for $A046646$ by Emeric Deutsch that for each $k \geq 2$, the $k$th number in sequence $A046646$ is double the $(k-1)$th number in the sequence $A001764$, that is double the number of $\pmodmon{2}{k-1}$-words. Hence Conjecture \ref{conj:wordsendingfromRunjj} would establish that every reduced $\pmodmon{2}{k-1}$-word (including the identity) remains reduced when either the $(k-1)$th diapsis generator or the $(k-1)$th $(2,2)$-transapsis generator is appended on the right.

	\nocite{man:GAP4, man:GAPsemigrps, man:OEIS}
    \bibliographystyle{plainnat}
    \clearpage{\pagestyle{empty}\cleardoublepage}
    \bibliography{references}
   	\addcontentsline{toc}{chapter}{References} 
   	
   	\appendix
   	    \pagestyle{fancy}
   	    \fancyhf{}
   	    \renewcommand{\headrulewidth}{0.5pt}
   	    \renewcommand{\sectionmark}[1]{\markright{\textbf{\thesection:}\ #1}} 
   	    \renewcommand{\chaptermark}[1]{\markboth{\textbf{\chaptername\ \thechapter:}\ #1}{}}
   	    \fancyhead[RO,LE]{}
   	    \fancyhead[RO]{\small\rightmark}
   	    \fancyhead[LE]{\small\leftmark}
   	    \fancyfoot[RO,LE]{\vspace{0.25cm}\textbf{\thepage}}
   	    
   	    \titleformat{\section}[display]
    	{\normalfont\Large\normalfont\bfseries}
    	{}
    	{-3pc}
    	{\Large\normalfont\bfseries$\blacksquare$\ \thesection\ }
    	[\vspace{-1ex}]
   	\setcounter{chapter}{1}
   	\setcounter{section}{0}
   	\clearpage{\pagestyle{plain}\cleardoublepage}
   	\chapter*{Appendices}\addcontentsline{toc}{chapter}{Appendices}\vspace{-1cm}
   	\markboth{APPENDICES}{}
	\section{Dot $\mathcal{D}$ classes for $\pmodmon{m}{k}$} \label{app:pmoddotD}
	Appendix \ref{app:pmoddotD} contains dot $\mathcal{D}$ classes for various values of $m, k \in \posints$ on the planar mod-$m$ monoid $\pmodmon{m}{k}$, which were generated using the Semigroups package \cite{man:GAPsemigrps} for GAP \cite{man:GAP4}.

	\begin{figure}[!ht]
		\caption[ ]{Dot $\mathcal{D}$ classes for $\pmodmon{2}{3}$ (left) and  $\pmodmon{2}{4}$ (right).}
		\label{fig:dclassespmodmon23and24}
		\vspace{-25pt}
		\begin{center}
		\includegraphics[scale = 0.65]{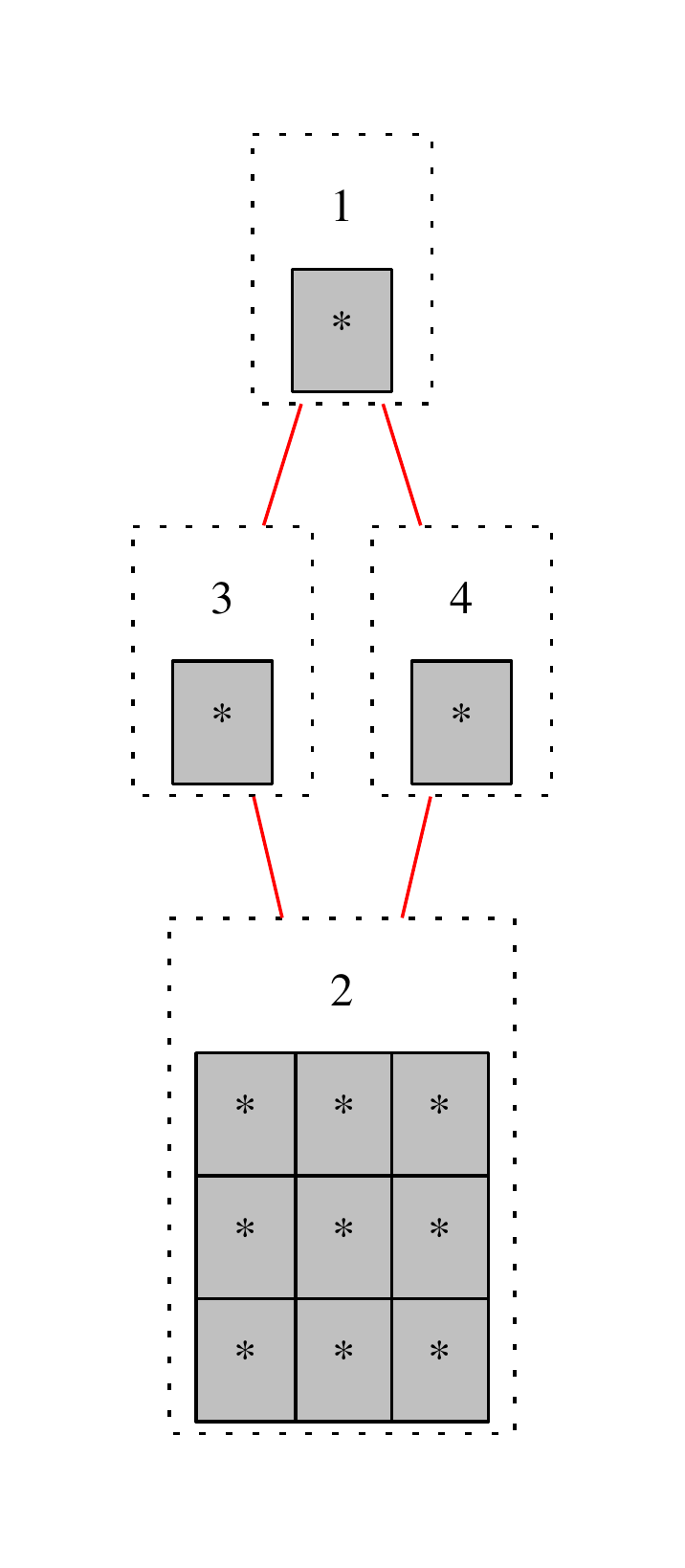}
		\hspace{0.5cm}
		\includegraphics[scale = 0.4]{chap_greens_relations/dotd/pmod/pmod24dotd.pdf}
		\end{center}
	\end{figure}

	\begin{figure}[!ht]
		\caption[ ]{Dot $\mathcal{D}$ classes for $\pmodmon{2}{5}$.}
		\label{fig:dclassespmodmon25}
		\vspace{-0.25cm}
		\centering
		\includegraphics[scale = 0.45]{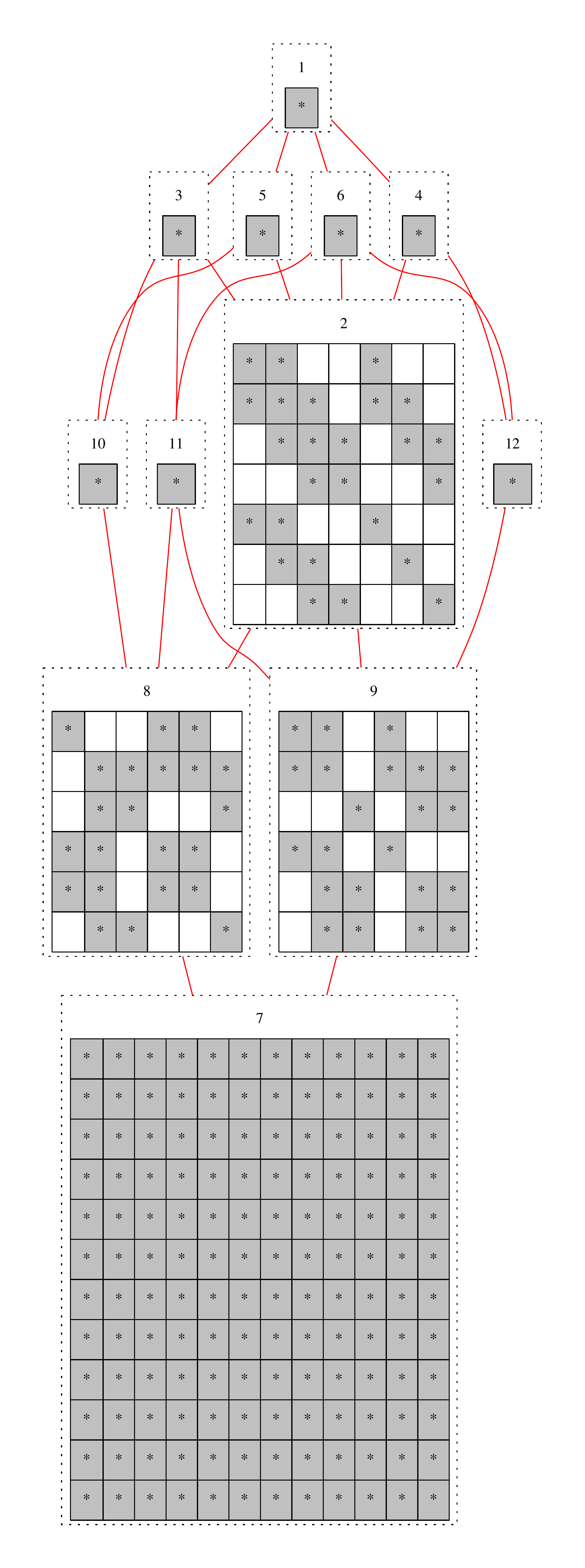}
	\end{figure}
	
	\begin{figure}[!ht]
		\caption[ ]{Dot $\mathcal{D}$ classes for $\pmodmon{2}{6}$.}
		\label{fig:dclassespmodmon26}
		\vspace{-15pt}
		\centering
		\includegraphics[scale = 0.2]{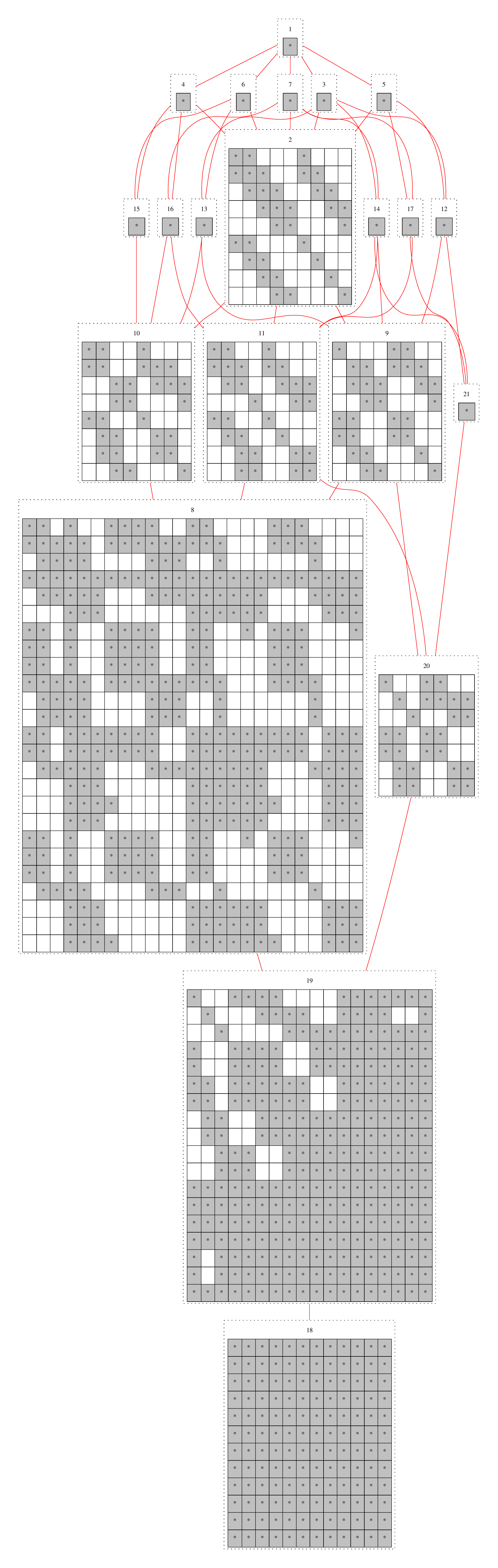}
	\end{figure}
	
	\begin{figure}[!ht]
		\caption[ ]{Dot $\mathcal{D}$ classes for $\pmodmon{3}{3}$ (left) and  $\pmodmon{3}{4}$ (right).}
		\label{fig:dclassespmodmon33}
		\vspace{5pt}
		\begin{center}
		\includegraphics[scale = 0.7]{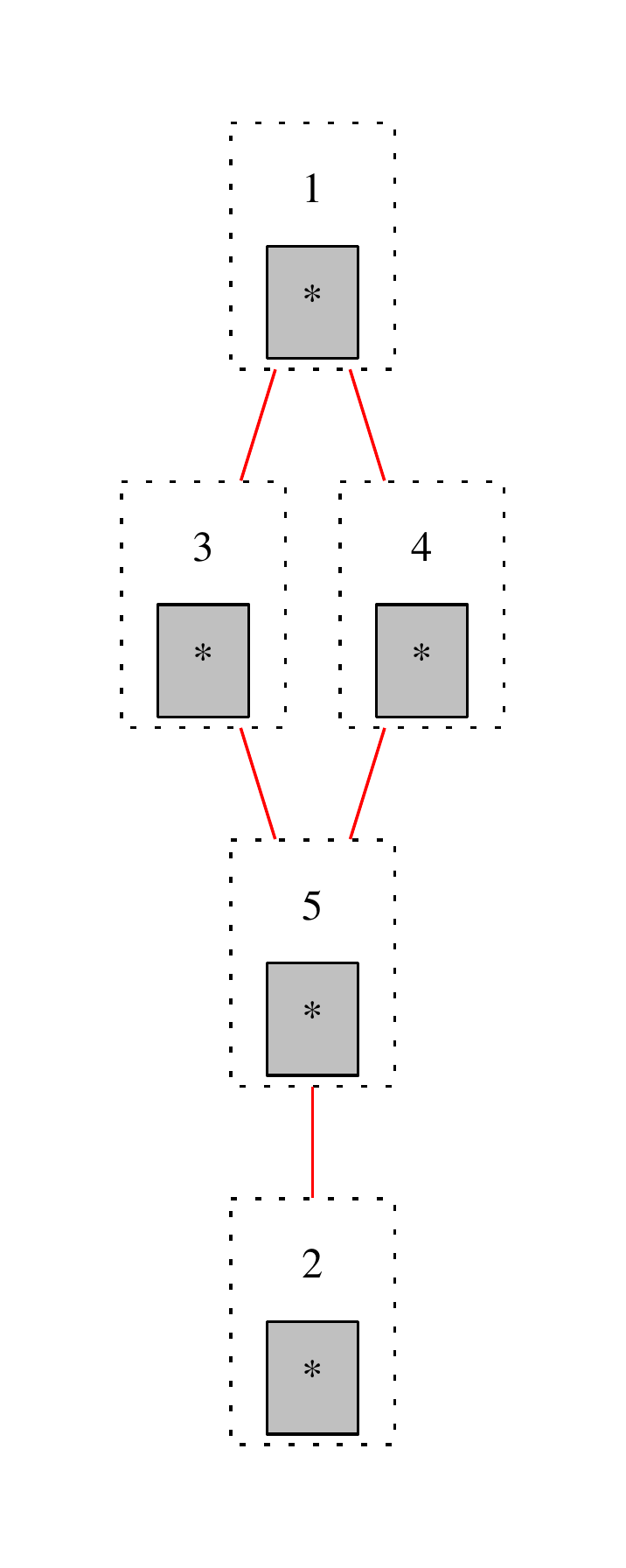}
		\hspace{0.5cm}
		\includegraphics[scale = 0.7]{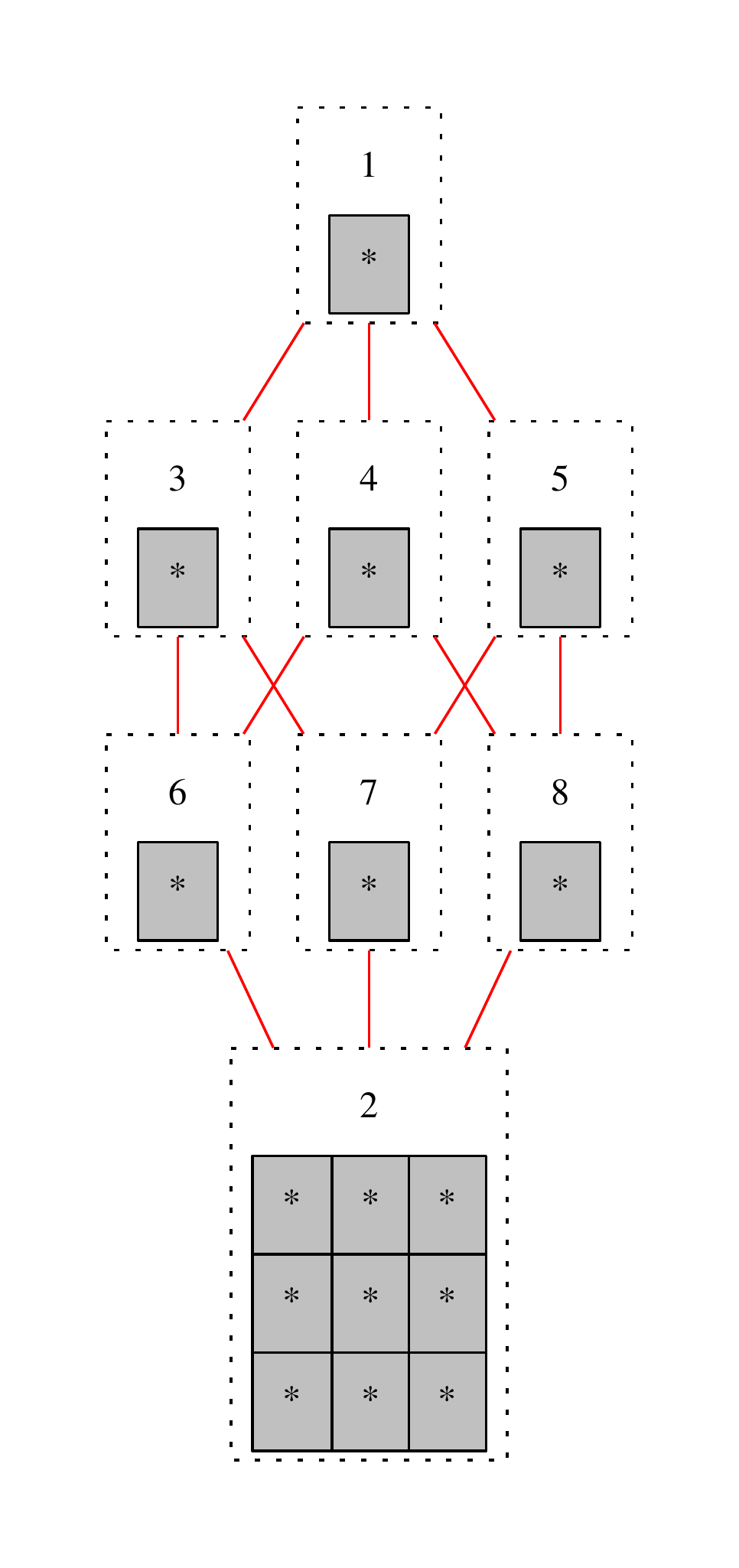}
		\end{center}
	\end{figure}
	
	\begin{figure}[!ht]
		\caption[ ]{Dot $\mathcal{D}$ classes for $\pmodmon{3}{5}$.}
		\label{fig:dclassespmodmon35}
		\vspace{5pt}
		\centering
		\includegraphics[scale = 0.7]{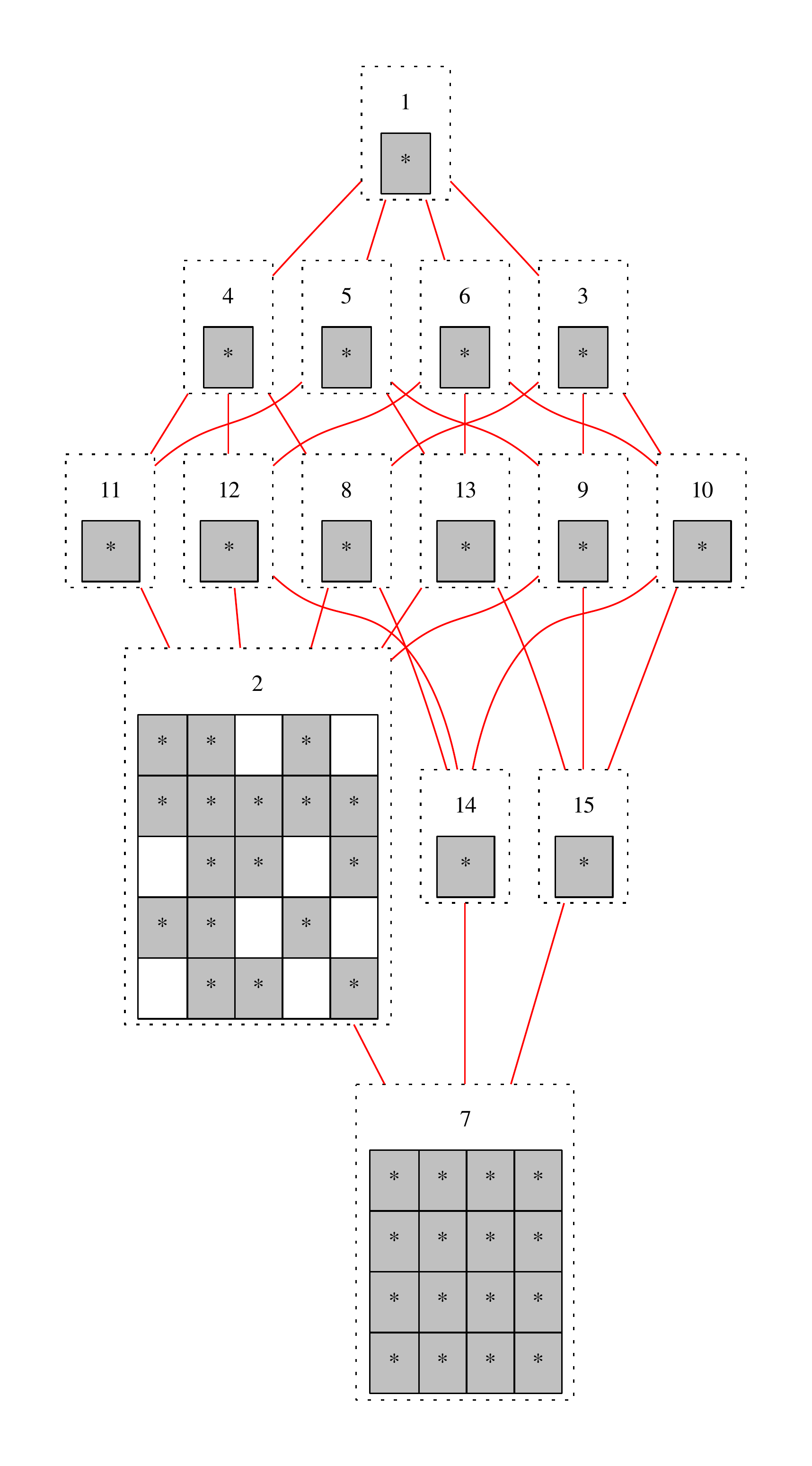}
	\end{figure}
	
	\begin{figure}[!ht]
		\caption[ ]{Dot $\mathcal{D}$ classes for $\pmodmon{3}{6}$.}
		\label{fig:dclassespmodmon36}
		\vspace{-15pt}
		\centering
		\includegraphics[scale = 0.43]{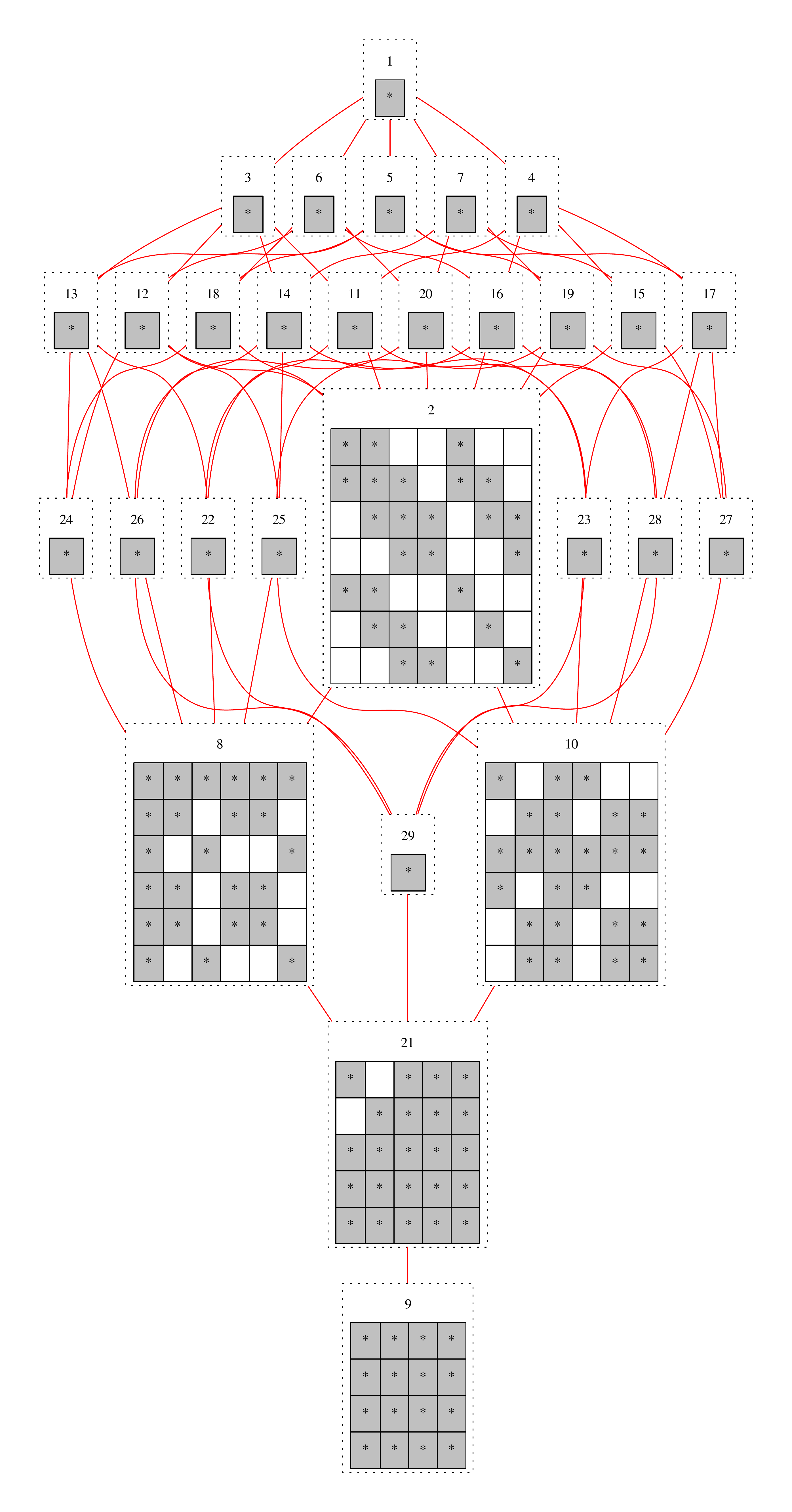}
	\end{figure}
	
	\begin{figure}[!ht]
		\caption[ ]{Dot $\mathcal{D}$ classes for $\pmodmon{4}{4}$.}
		\label{fig:dclassespmodmon44}
		\vspace{5pt}
		\centering
		\includegraphics[scale = 0.8]{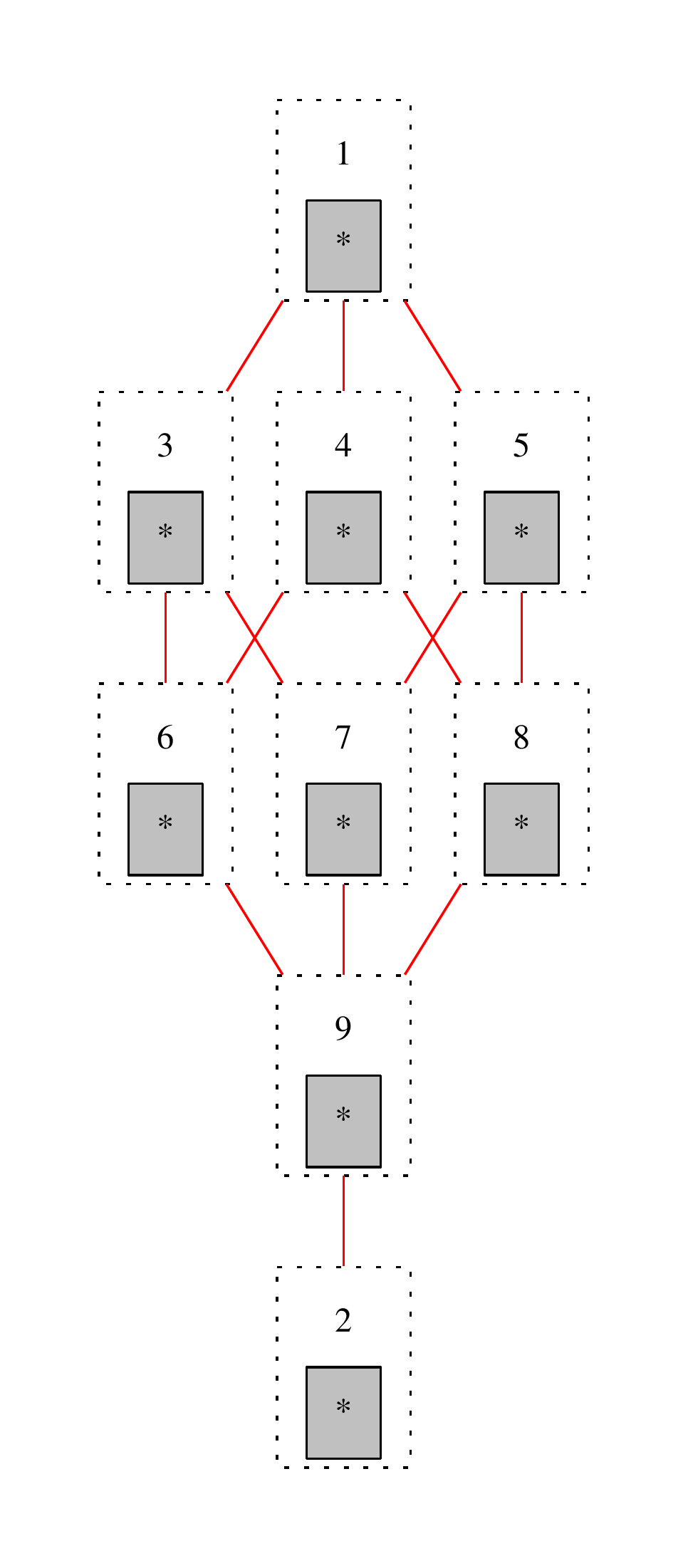}
	\end{figure}
	
	\begin{figure}[!ht]
		\caption[ ]{Dot $\mathcal{D}$ classes for $\pmodmon{4}{5}$.}
		\label{fig:dclassespmodmon45}
		\vspace{5pt}
		\centering
		\includegraphics[scale = 0.8]{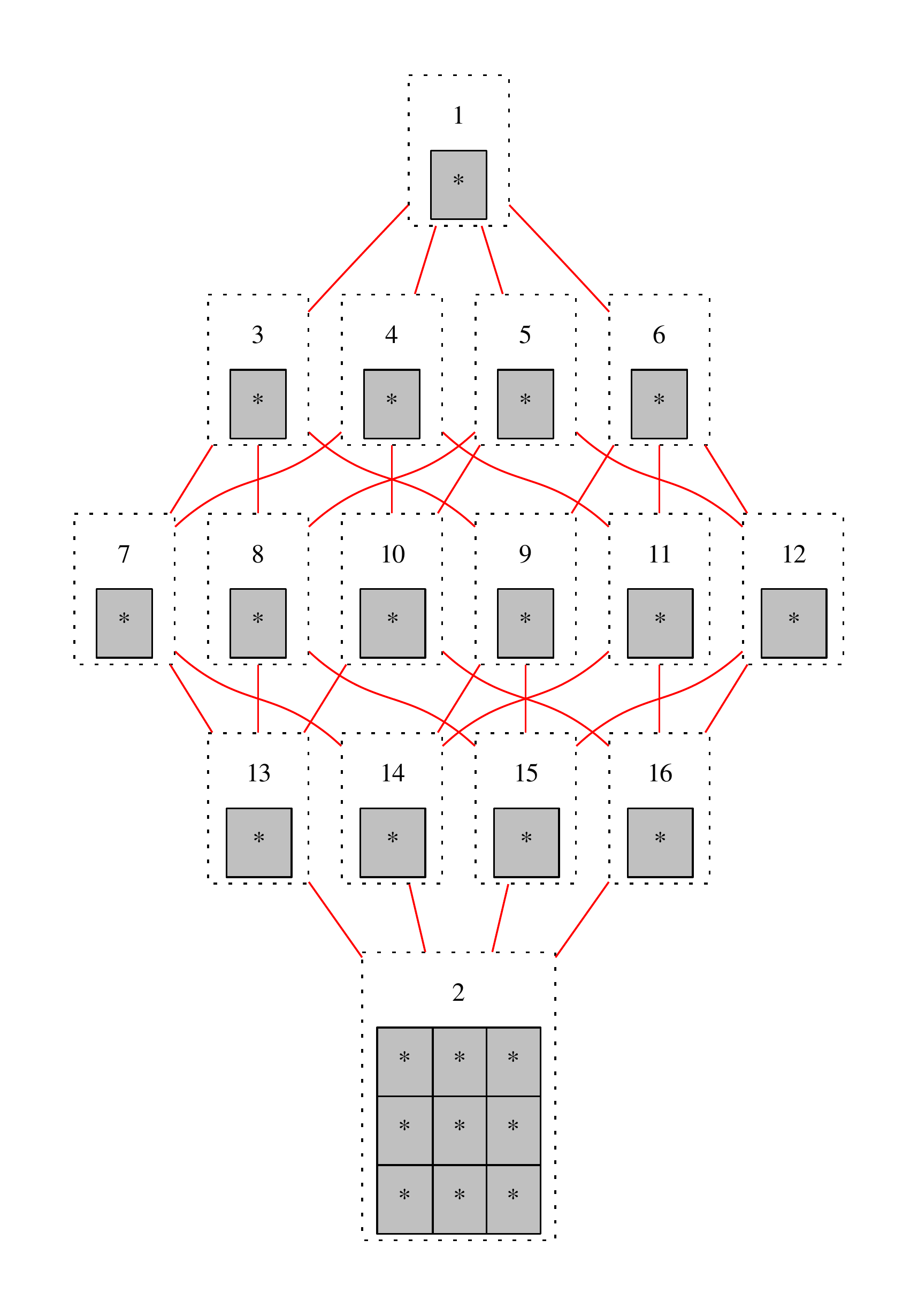}
	\end{figure}
	
	\begin{figure}[!ht]
		\caption[ ]{Dot $\mathcal{D}$ classes for $\pmodmon{4}{6}$.}
		\label{fig:dclassespmodmon46}
		\vspace{5pt}
		\centering
		\includegraphics[scale = 0.52]{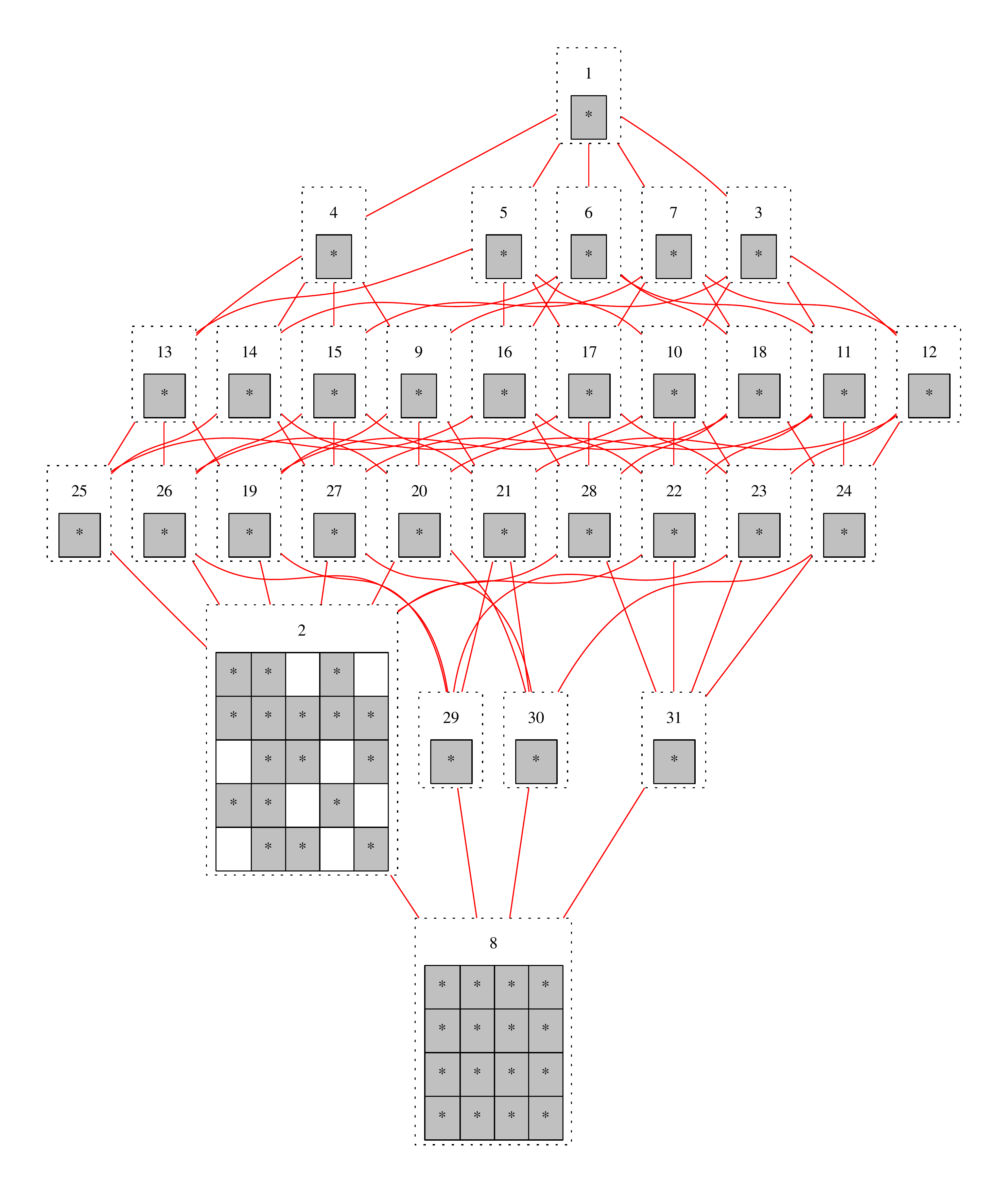}
	\end{figure}
	
	\begin{figure}[!ht]
		\caption[ ]{Dot $\mathcal{D}$ classes for $\pmodmon{5}{5}$.}
		\label{fig:dclassespmodmon55}
		\vspace{5pt}
		\centering
		\includegraphics[scale = 0.75]{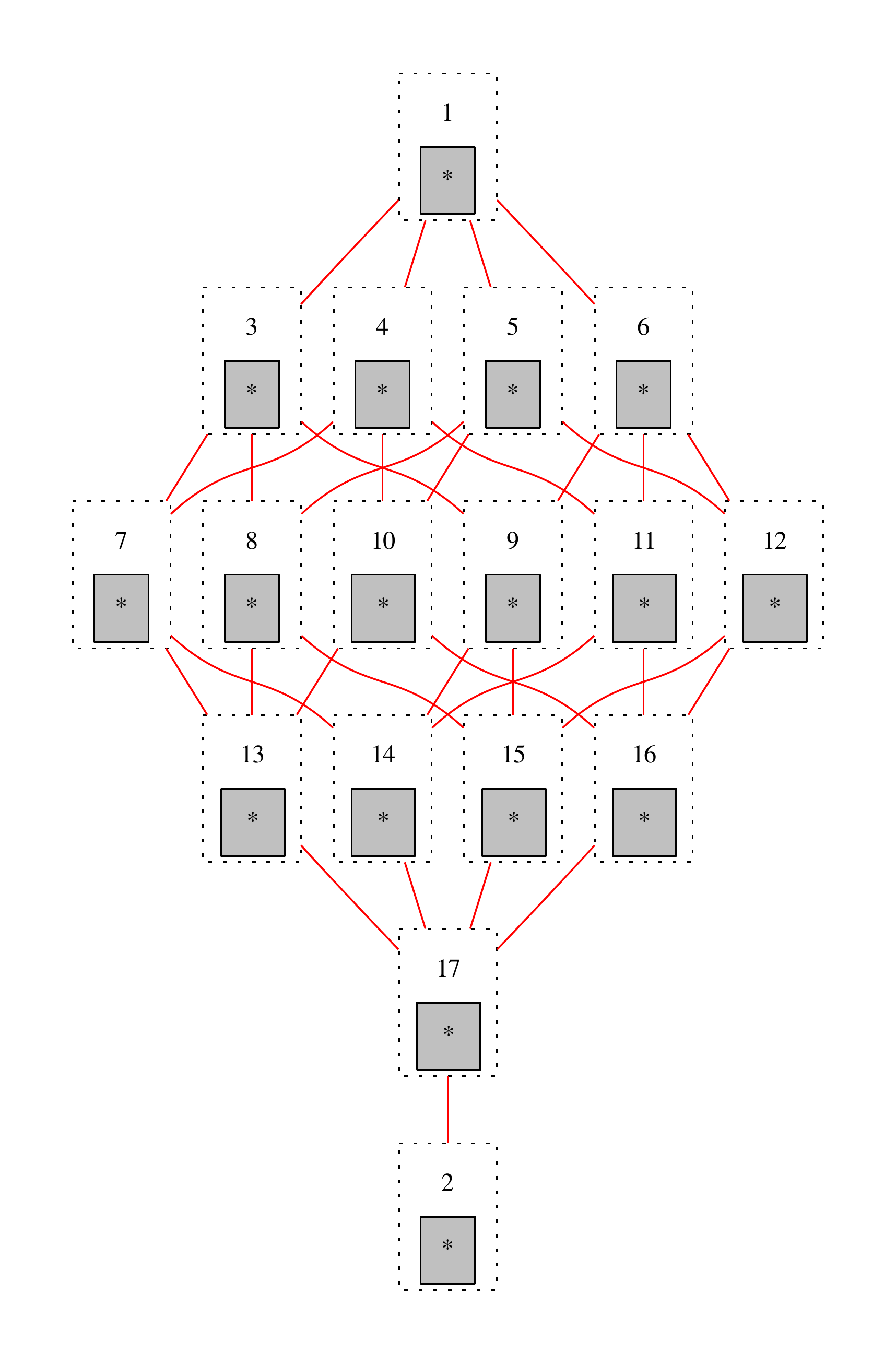}
	\end{figure}
	
	\begin{figure}[!ht]
		\caption[ ]{Dot $\mathcal{D}$ classes for $\pmodmon{5}{6}$.}
		\label{fig:dclassespmodmon56}
		\vspace{5pt}
		\centering
		\includegraphics[scale = 0.55]{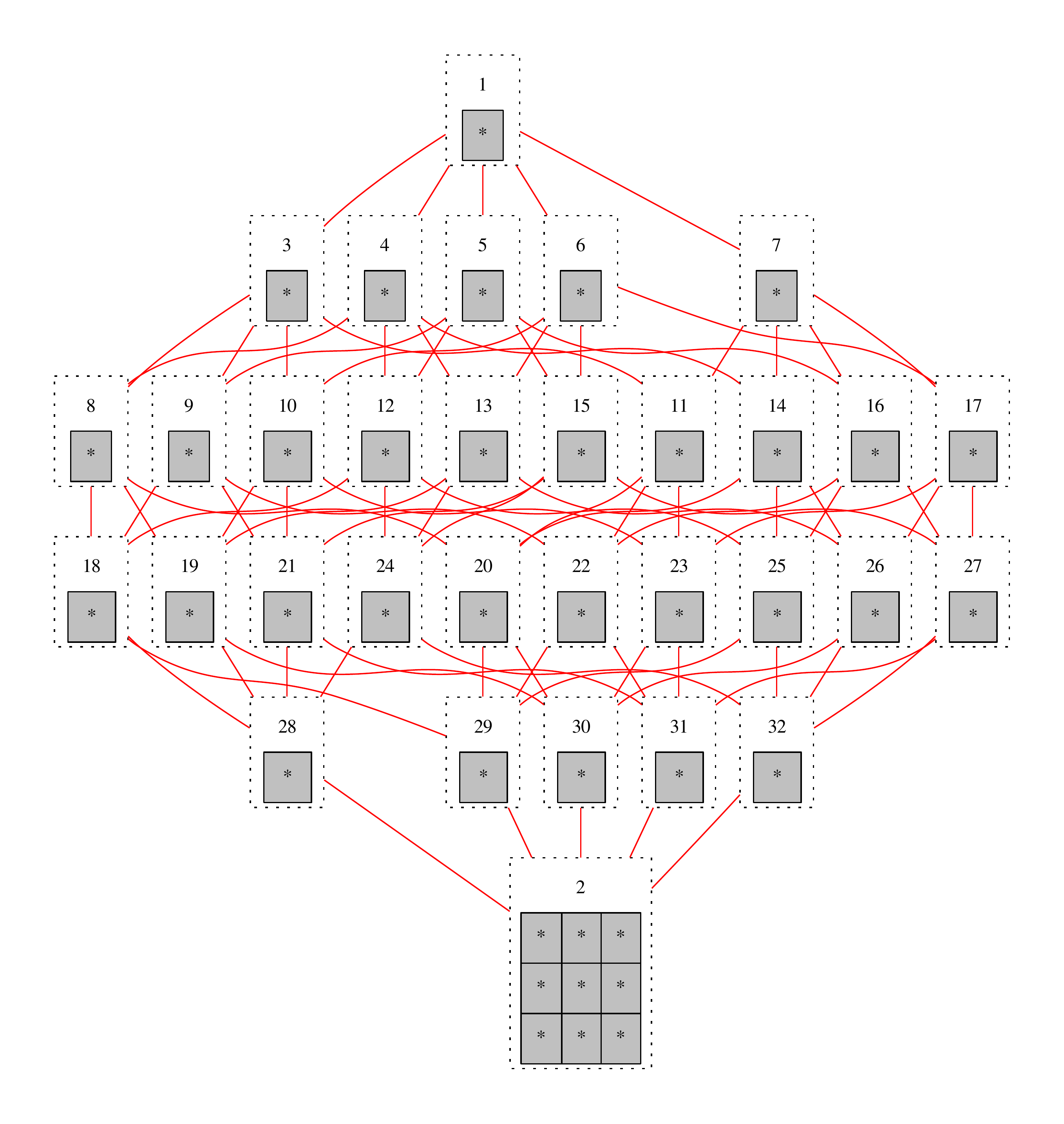}
	\end{figure}

	\clearpage{\pagestyle{plain}\cleardoublepage}
\section{Dot $\mathcal{D}$ classes for $\modmon{m}{k}$} \label{app:moddotD}
	Appendix \ref{app:moddotD} contains dot $\mathcal{D}$ classes for various values of $m, k \in \posints$ on the mod-$m$ monoid $\modmon{m}{k}$, which were generated using the Semigroups package \cite{man:GAPsemigrps} for GAP \cite{man:GAP4}.
	
	\begin{figure}[!ht]
		\caption[ ]{Dot $\mathcal{D}$ classes for $\modmon{2}{2}$ (left) and  $\modmon{2}{3}$ (right).}
		\label{fig:dclassesmodmon22and23}
		\vspace{5pt}
		\begin{center}
		\includegraphics[scale = 0.7]{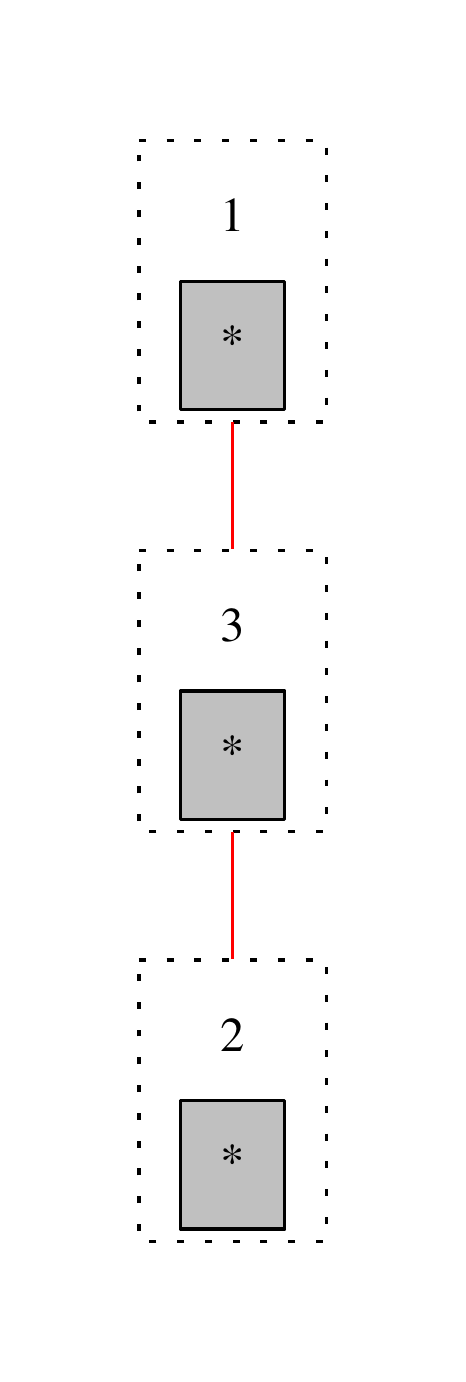}
		\hspace{2cm}
		\includegraphics[scale = 0.7]{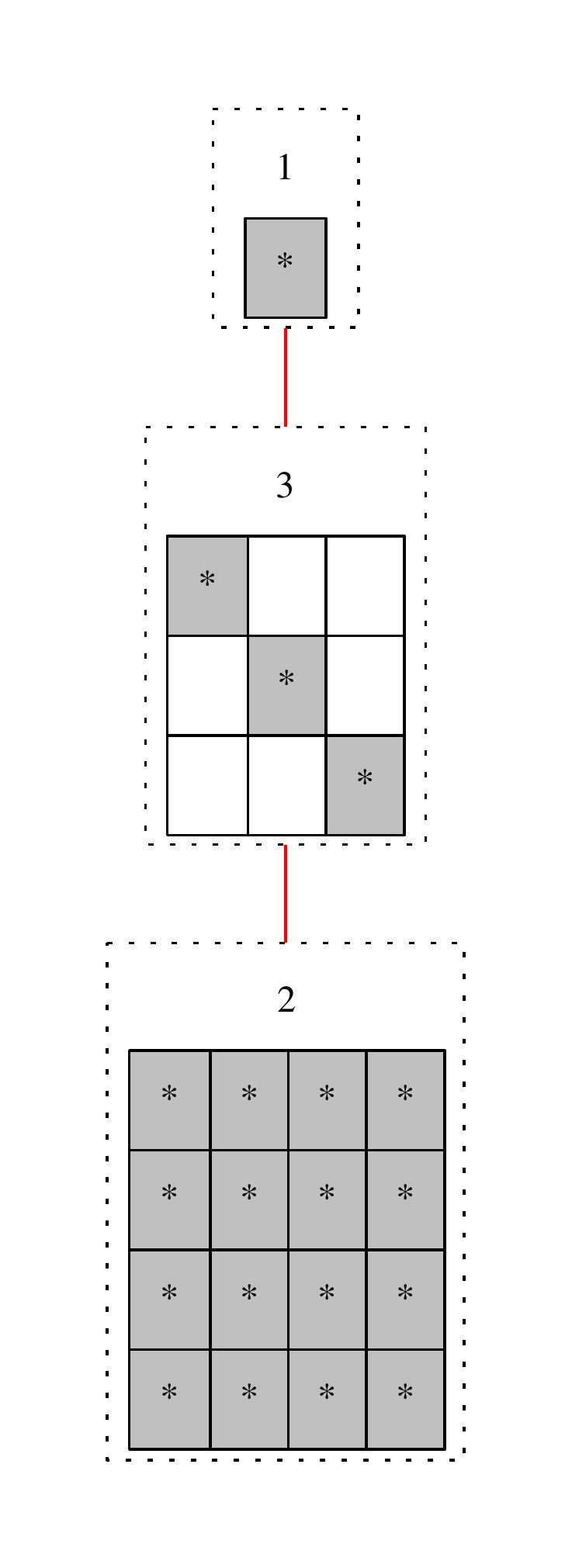}
		\end{center}
	\end{figure}
	
	\begin{figure}[!ht]
		\caption[ ]{Dot $\mathcal{D}$ classes for $\modmon{2}{4}$.}
		\label{fig:dclassesmodmon24}
		\vspace{-25pt}
		\begin{center}
			\includegraphics[scale = 0.44]{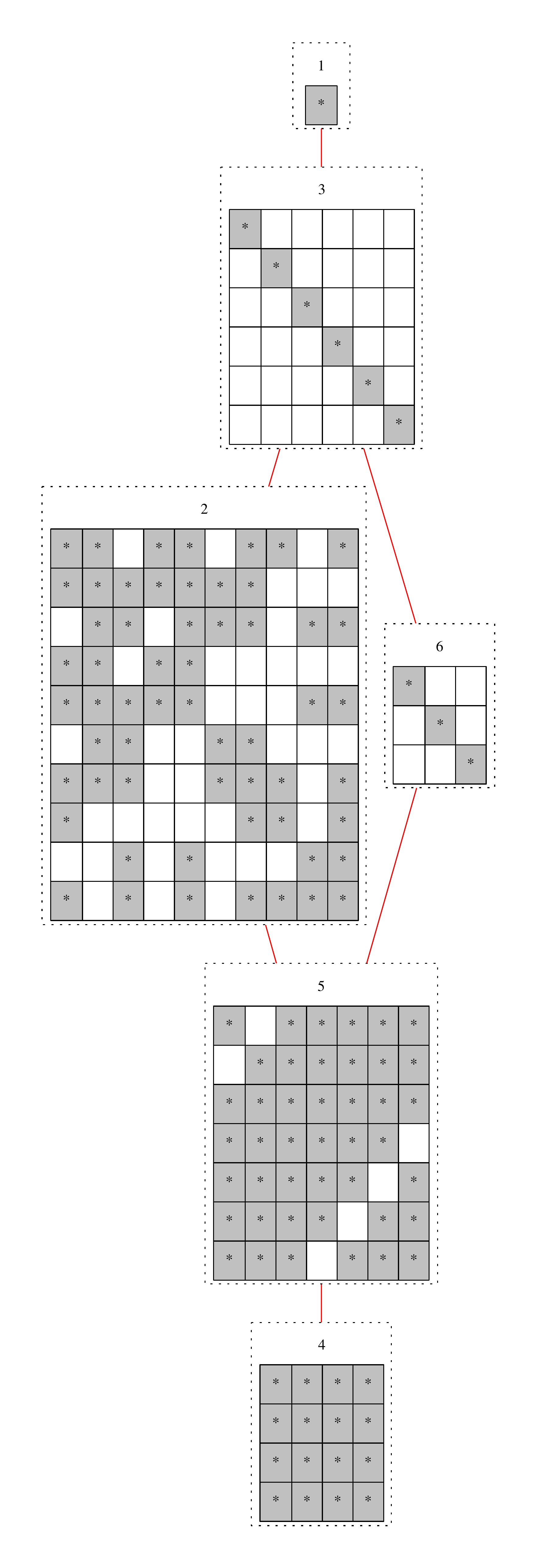}
		\end{center}
	\end{figure}
	
	\begin{figure}[!ht]
		\caption[ ]{Dot $\mathcal{D}$ classes for $\modmon{2}{5}$.}
		\label{fig:dclassesmodmon25}
		\vspace{-15pt}
		\begin{center}
			\includegraphics[scale = 0.145]{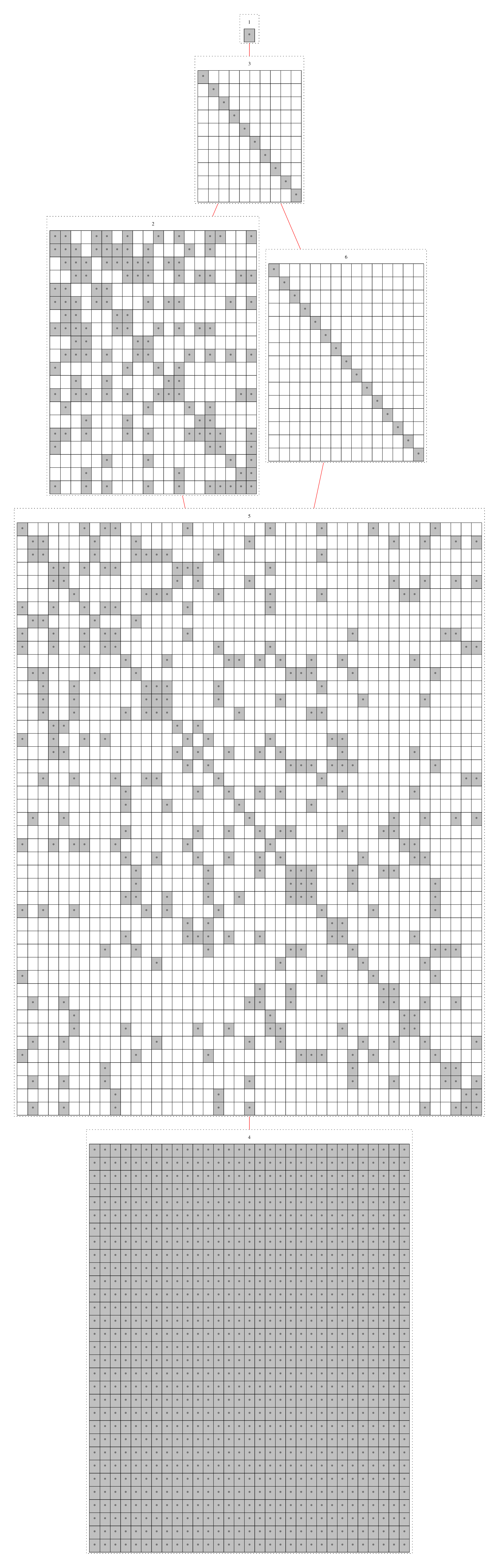}
		\end{center}
	\end{figure}
	
	\begin{figure}[!ht]
		\caption[ ]{Dot $\mathcal{D}$ classes for $\modmon{3}{3}$ (left) and  $\modmon{3}{4}$ (right).}
		\label{fig:dclassesmodmon33and34}
		\vspace{5pt}
		\begin{center}
		\includegraphics[scale = 0.7]{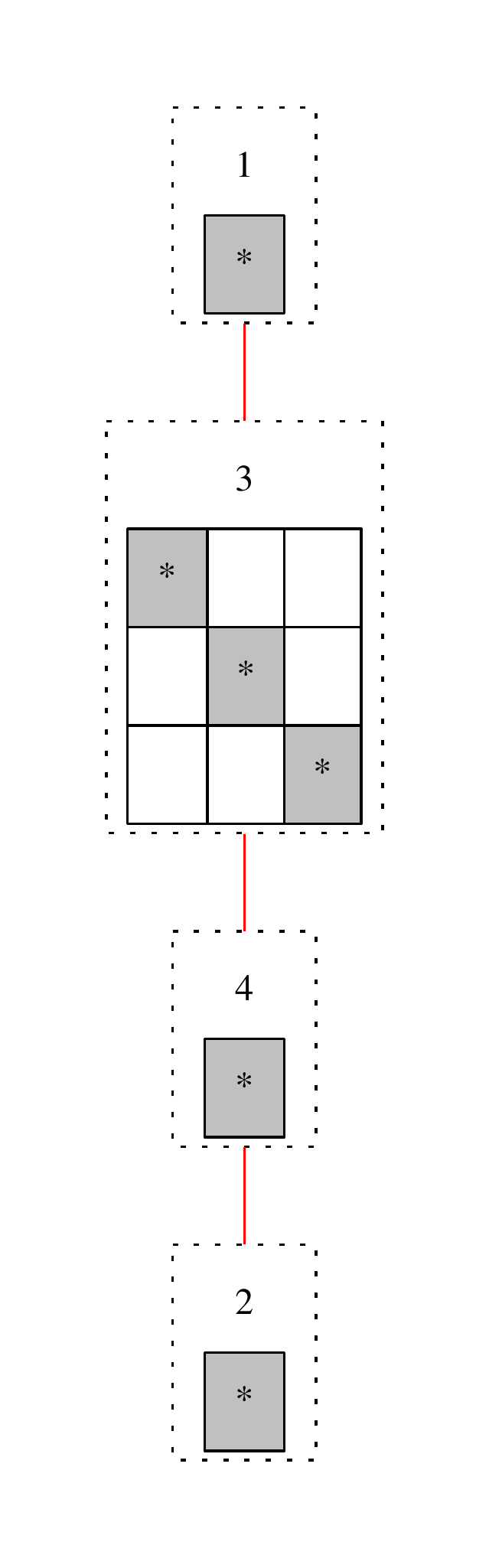}
		\hspace{0.5cm}
		\includegraphics[scale = 0.6]{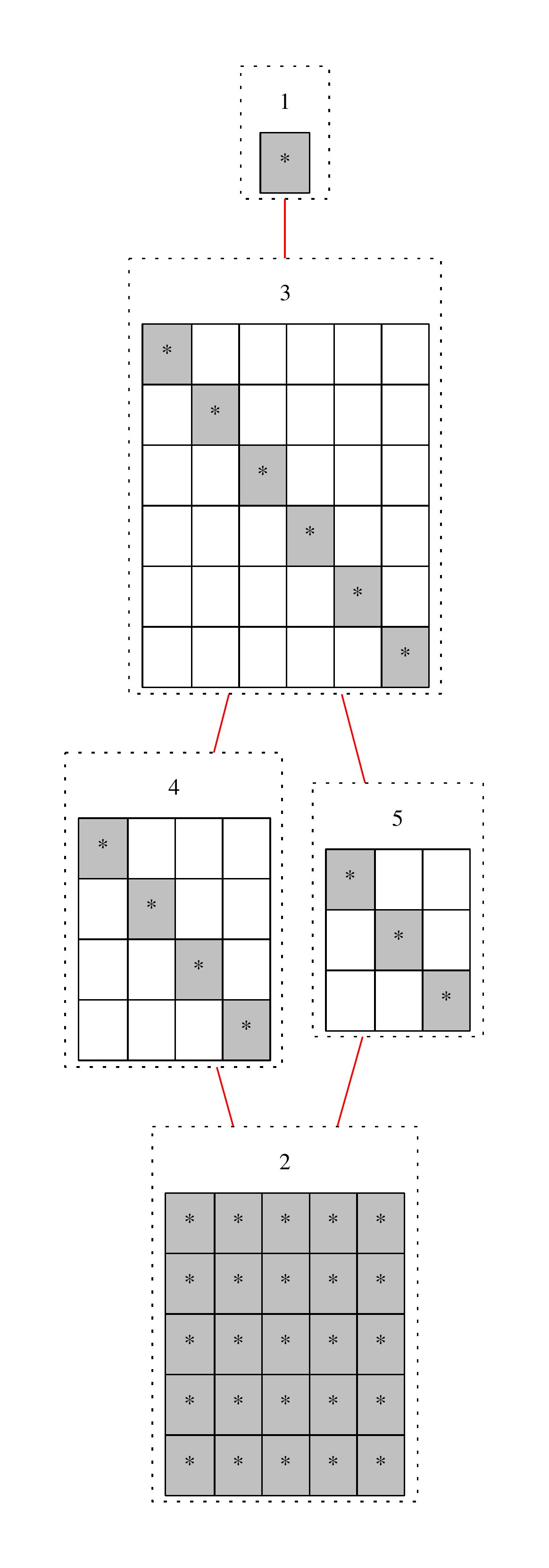}
		\end{center}
	\end{figure}
	
	\begin{figure}[!ht]
		\caption[ ]{Dot $\mathcal{D}$ classes for $\modmon{3}{5}$.}
		\label{fig:dclassesmodmon35}
		\vspace{-25pt}
		\begin{center}
		\includegraphics[scale = 0.28]{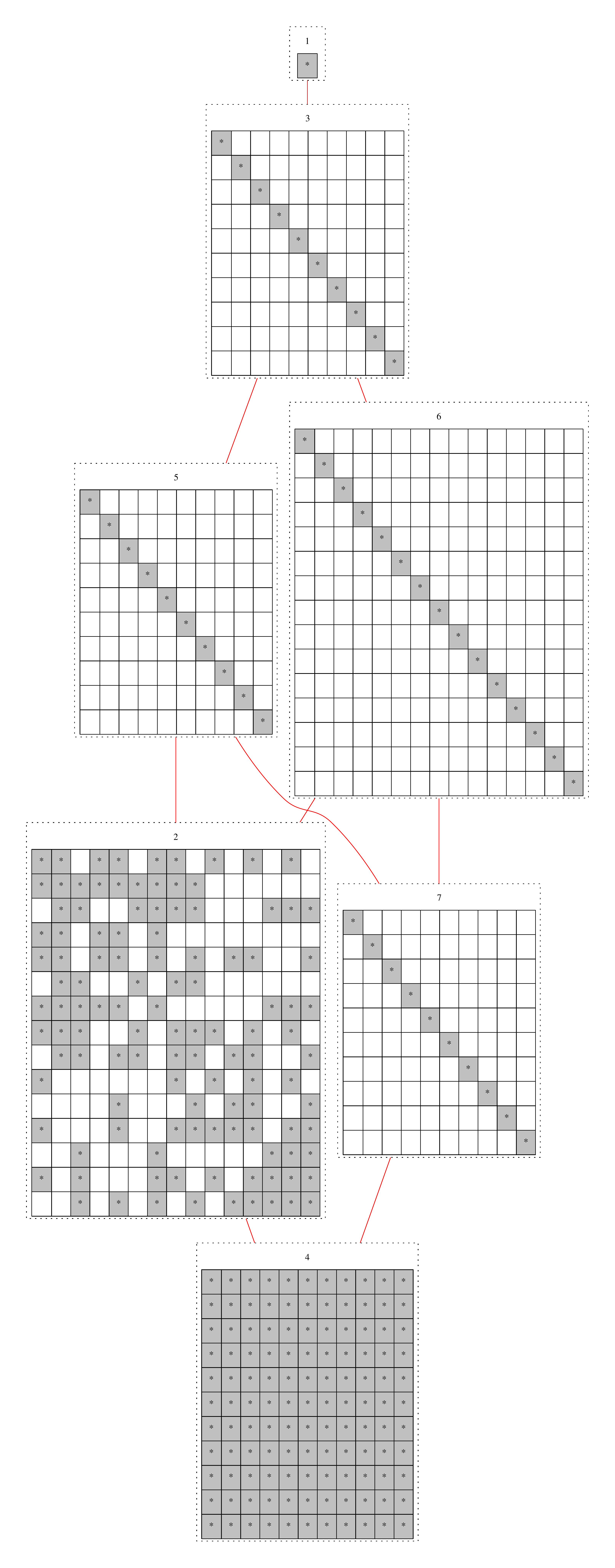}
		\end{center}
	\end{figure}
	
	\begin{figure}[!ht]
		\caption[ ]{Dot $\mathcal{D}$ classes for $\modmon{4}{4}$.}
		\label{fig:dclassesmodmon44}
		\vspace{-35pt}
		\begin{center}
		\includegraphics[scale = 0.7]{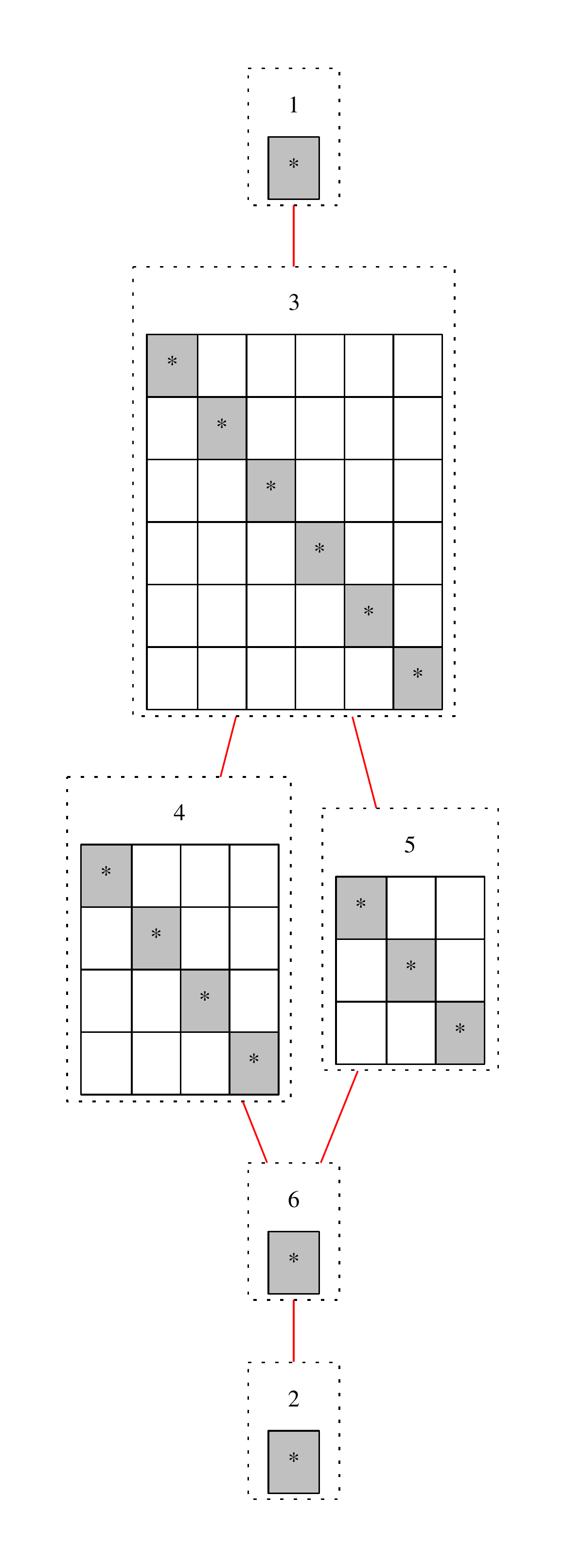}
		\end{center}
	\end{figure}
	
	\begin{figure}[!ht]
		\caption[ ]{Dot $\mathcal{D}$ classes for $\modmon{4}{5}$.}
		\label{fig:dclassesmodmon45}
		\vspace{-25pt}
		\begin{center}
		\includegraphics[scale = 0.33]{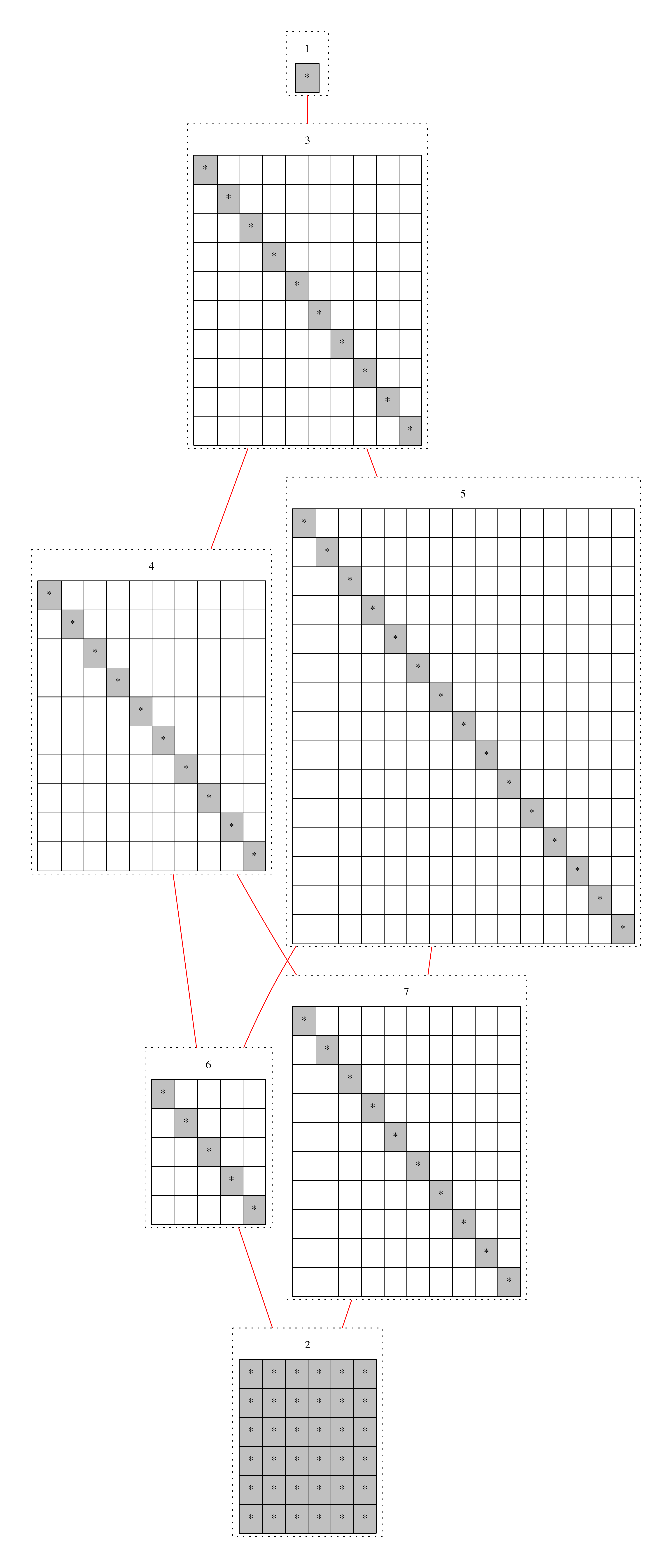}
		\end{center}
	\end{figure}
	
	\clearpage{\pagestyle{plain}\cleardoublepage}
\section{GAP code for $\pmodmon{2}{k}$ presentation up to $k=7$} \label{app:pmod2presentations}
	Appendix \ref{app:pmod2presentations} contains GAP \cite{man:GAP4} code for presentations of the planar mod-$2$ monoid $\pmodmon{2}{k}$ for all $k \in \set{2, 3, 4, 5, 6, 7}$ (using the relations from Conjecture \ref{conj:pmod2presentation}).
	
	Note that we have relabelled the diapsis generators using $\hookgen{1} \mapsto a$, $\hookgen{2} \mapsto b$, $\hookgen{3} \mapsto c$, $\hookgen{4} \mapsto d$, $\hookgen{5} \mapsto f$, $\hookgen{6} \mapsto g$, and relabelled the $(2,2)$-transapsis generators using $\transapgen{1} \mapsto A$, $\transapgen{2} \mapsto B$, $\transapgen{3} \mapsto C$, $\transapgen{4} \mapsto D$, $\transapgen{5} \mapsto F$, $\transapgen{6} \mapsto G$
\singlespacing
\begin{lstlisting}
freemon2 := FreeMonoid("a", "A");
a := GeneratorsOfMonoid(freemon2)[1];
A := GeneratorsOfMonoid(freemon2)[2];
PMod2 := freemon2/[
	[a^2, a],


	[A^2, A],


	[a*A, a],
	[A*a, a],
];
Size(PMod2);
Elements(PMod2);
\end{lstlisting}

\begin{lstlisting}
freemon3 := FreeMonoid("a", "A", "b", "B");
a := GeneratorsOfMonoid(freemon3)[1];
A := GeneratorsOfMonoid(freemon3)[2];
b := GeneratorsOfMonoid(freemon3)[3];
B := GeneratorsOfMonoid(freemon3)[4];
PMod3 := freemon3/[
	[a^2, a],  [b^2, b],

	[a*b*a, a],
	[b*a*b, b],


	[A^2, A], [B^2, B], 

	[B*A, A*B],


	[a*A, a], [b*B, b],
	[A*a, a], [B*b, b],

	[A*b*A, A*B],
];
Size(PMod3);
Elements(PMod3);
\end{lstlisting}

\newpage

\begin{lstlisting}
freemon4 := FreeMonoid("a", "A", "b", "B", "c", "C");
a := GeneratorsOfMonoid(freemon4)[1];
A := GeneratorsOfMonoid(freemon4)[2];
b := GeneratorsOfMonoid(freemon4)[3];
B := GeneratorsOfMonoid(freemon4)[4];
c := GeneratorsOfMonoid(freemon4)[5];
C := GeneratorsOfMonoid(freemon4)[6];
PMod4 := freemon4/[
	[a^2, a],  [b^2, b], [c^2, c],

	[a*b*a, a], [b*c*b, b],
	[b*a*b, b], [c*b*c, c],

	[c*a, a*c],


	[A^2, A], [B^2, B], [C^2, C],

	[B*A, A*B], [C*B, B*C],
	[C*A, A*C],


	[a*A, a], [b*B, b], [c*C, c],
	[A*a, a], [B*b, b], [C*c, c],

	[A*b*A, A*B], [B*c*B, B*C],

	[c*A, A*c],

	[C*a, a*C],
];
Size(PMod4);
Elements(PMod4);
\end{lstlisting}

\newpage

\begin{lstlisting}
freemon5 := FreeMonoid("a", "A", "b", "B", "c", "C", "d", "D");
a := GeneratorsOfMonoid(freemon5)[1];
A := GeneratorsOfMonoid(freemon5)[2];
b := GeneratorsOfMonoid(freemon5)[3];
B := GeneratorsOfMonoid(freemon5)[4];
c := GeneratorsOfMonoid(freemon5)[5];
C := GeneratorsOfMonoid(freemon5)[6];
d := GeneratorsOfMonoid(freemon5)[7];
D := GeneratorsOfMonoid(freemon5)[8];
PMod5 := freemon5/[
	[a^2, a],  [b^2, b], [c^2, c], [d^2, d],

	[a*b*a, a], [b*c*b, b], [c*d*c, c],
	[b*a*b, b], [c*b*c, c], [d*c*d, d],

	[c*a, a*c], [d*b, b*d],
	[d*a, a*d],


	[A^2, A], [B^2, B], [C^2, C], [D^2, D], 

	[B*A, A*B], [C*B, B*C], [D*C, C*D],
	[C*A, A*C], [D*B, B*D],
	[D*A, A*D],


	[a*A, a], [b*B, b], [c*C, c], [d*D, d],
	[A*a, a], [B*b, b], [C*c, c], [D*d, d],

	[A*b*A, A*B], [B*c*B, B*C], [C*d*C, C*D],

	[c*A, A*c], [d*B, B*d],
	[d*A, A*d],

	[C*a, a*C], [D*b, b*D],
	[D*a, a*D],
];
Size(PMod5);
Elements(PMod5);
\end{lstlisting}

\newpage

\begin{lstlisting}
freemon6 := FreeMonoid("a", "A", "b", "B", "c", "C", "d", "D", "f", "F");
a := GeneratorsOfMonoid(freemon6)[1];
A := GeneratorsOfMonoid(freemon6)[2];
b := GeneratorsOfMonoid(freemon6)[3];
B := GeneratorsOfMonoid(freemon6)[4];
c := GeneratorsOfMonoid(freemon6)[5];
C := GeneratorsOfMonoid(freemon6)[6];
d := GeneratorsOfMonoid(freemon6)[7];
D := GeneratorsOfMonoid(freemon6)[8];
f := GeneratorsOfMonoid(freemon6)[9];
F := GeneratorsOfMonoid(freemon6)[10];
PMod6 := freemon6/[
	[a^2, a],  [b^2, b], [c^2, c], [d^2, d], [f^2, f],

	[a*b*a, a], [b*c*b, b], [c*d*c, c], [d*f*d, d],
	[b*a*b, b], [c*b*c, c], [d*c*d, d], [f*d*f, f],

	[c*a, a*c], [d*b, b*d], [f*c, c*f],
	[d*a, a*d], [f*b, b*f],
	[f*a, a*f],


	[A^2, A], [B^2, B], [C^2, C], [D^2, D], [F^2, F], 

	[B*A, A*B], [C*B, B*C], [D*C, C*D], [F*D, D*F],
	[C*A, A*C], [D*B, B*D], [F*C, C*F],
	[D*A, A*D], [F*B, B*F],
	[F*A, A*F],


	[a*A, a], [b*B, b], [c*C, c], [d*D, d], [f*F, f],
	[A*a, a], [B*b, b], [C*c, c], [D*d, d], [F*f, f],

	[A*b*A, A*B], [B*c*B, B*C], [C*d*C, C*D], [D*f*D, D*F],

	[c*A, A*c], [d*B, B*d], [f*C, C*f],
	[d*A, A*d], [f*B, B*f],
	[f*A, A*f],

	[C*a, a*C], [D*b, b*D], [F*c, c*F],
	[D*a, a*D], [F*b, b*F],
	[F*a, a*F],
];
Size(PMod6);
Elements(PMod6);
\end{lstlisting}

\newpage

\begin{lstlisting}
freemon7 := FreeMonoid("a", "A", "b", "B", "c", "C", "d", "D", "f", "F", "g", "G");
a := GeneratorsOfMonoid(freemon7)[1];
A := GeneratorsOfMonoid(freemon7)[2];
b := GeneratorsOfMonoid(freemon7)[3];
B := GeneratorsOfMonoid(freemon7)[4];
c := GeneratorsOfMonoid(freemon7)[5];
C := GeneratorsOfMonoid(freemon7)[6];
d := GeneratorsOfMonoid(freemon7)[7];
D := GeneratorsOfMonoid(freemon7)[8];
f := GeneratorsOfMonoid(freemon7)[9];
F := GeneratorsOfMonoid(freemon7)[10];
g := GeneratorsOfMonoid(freemon7)[11];
G := GeneratorsOfMonoid(freemon7)[12];
PMod7 := freemon7/[
	[a^2, a],  [b^2, b], [c^2, c], [d^2, d], [f^2, f], [g^2, g],

	[a*b*a, a], [b*c*b, b], [c*d*c, c], [d*f*d, d], [f*g*f, f],
	[b*a*b, b], [c*b*c, c], [d*c*d, d], [f*d*f, f], [g*f*g, g],

	[c*a, a*c], [d*b, b*d], [f*c, c*f], [g*d, d*g],
	[d*a, a*d], [f*b, b*f], [g*c, c*g],
	[f*a, a*f], [g*b, b*g],
	[g*a, a*g],


	[A^2, A], [B^2, B], [C^2, C], [D^2, D], [F^2, F], [G^2, G],

	[B*A, A*B], [C*B, B*C], [D*C, C*D], [F*D, D*F], [G*F, F*G],
	[C*A, A*C], [D*B, B*D], [F*C, C*F], [G*D, D*G],
	[D*A, A*D], [F*B, B*F], [G*C, C*G],
	[F*A, A*F], [G*B, B*G],
	[G*A, A*G],


	[a*A, a], [b*B, b], [c*C, c], [d*D, d], [f*F, f], [g*G, g],
	[A*a, a], [B*b, b], [C*c, c], [D*d, d], [F*f, f], [G*g, g],

	[A*b*A, A*B], [B*c*B, B*C], [C*d*C, C*D], [D*f*D, D*F], [F*g*F, F*G],

	[c*A, A*c], [d*B, B*d], [f*C, C*f], [g*D, D*g],
	[d*A, A*d], [f*B, B*f], [g*C, C*g],
	[f*A, A*f], [g*B, B*g],
	[g*A, A*g],

	[C*a, a*C], [D*b, b*D], [F*c, c*F], [G*d, d*G],
	[D*a, a*D], [F*b, b*F], [G*c, c*G],
	[F*a, a*F], [G*b, b*G],
	[G*a, a*G],
];
Size(PMod7);
\end{lstlisting}
\doublespacing
	\clearpage{\pagestyle{plain}\cleardoublepage}
\section{Reduced $\pmodmon{2}{k}$-words in normal form} \label{app:pmod2words}
	Appendix \ref{app:pmod2words} contains $\pmodmon{2}{k}$-words for various values of $k$. The words were generated by running (with the specified values of $k$) `\texttt{Elements(PModk);}' proceeding the code for each presentation in Appendix \ref{app:pmod2presentations}, then reordered based on the end run as outlined in Subsection \ref{subsec:pmod2boundednormals}. Note that the presentation code in Appendix \ref{app:pmod2presentations} only covers one ordering of the generators, while in this appendix we will consider two orderings of the generators.

\singlespacing
\begin{lstlisting} [language=Bash]	
planar mod-2 words when k=3, 
with generators ordered: a,A,b,B
========================================================================
a, A, 
------------------------------------------------------------------------
b*a, b*A, B*a, 
========================================================================
b, B, 
a*b, a*B, 
A*b, A*B, 
\end{lstlisting}

\begin{lstlisting}[language=Bash]
planar mod-2 words when k=3, 
with generators ordered: A,a,B,b
========================================================================
A, a,
------------------------------------------------------------------------
B*a, b*A, b*a
========================================================================
B, b,
A*B, A*b, 
a*B, a*b,
\end{lstlisting}
	
\newpage
	
\begin{lstlisting}[language=Bash]
planar mod-2 words when k=4, 
with generators ordered: a,A,b,B,c,C
========================================================================
a, A, 
------------------------------------------------------------------------
b*a, b*A, B*a, 
------------------------------------------------------------------------
c*b*a, c*b*A, c*B*a, C*b*a, C*b*A, 
========================================================================
b, B, 
a*b, a*B, 
A*b, A*B, 
------------------------------------------------------------------------
    c*b,     c*B,     C*b, 
  a*c*b,   a*c*B,   a*C*b, 
  A*c*b,   A*c*B,   A*C*b, 
b*a*c*b, b*a*c*B, b*a*C*b, 
B*a*c*b, B*a*c*B, 
========================================================================
    c,     C, 
  a*c,   a*C, 
  A*c,   A*C, 
  b*c,   b*C, 
  B*c,   B*C, 
a*b*c, a*b*C, 
a*B*c, a*B*C, 
A*b*c, A*b*C, 
A*B*c, A*B*C, 
b*a*c, b*a*C, 
b*A*c, b*A*C, 
B*a*c, B*a*C,
\end{lstlisting}
	
\newpage
	
\begin{lstlisting}[language=Bash]
planar mod-2 words when k=4, 
with generators ordered: A,a,B,b,C,c
========================================================================
A, a, 
------------------------------------------------------------------------
B*a, b*A, b*a, 
------------------------------------------------------------------------
C*b*A, C*b*a, c*B*a, c*b*A, c*b*a, 
========================================================================
  B,   b, 
A*B, A*b, 
a*B, a*b, 
------------------------------------------------------------------------
    C*b,     c*B,     c*b, 
  A*C*b,   A*c*B,   A*c*b, 
  a*C*b,   a*c*B,   a*c*b, 
         B*a*c*B, B*a*c*b, 
b*A*C*b, 
         b*a*c*B, b*a*c*b, 
========================================================================
    C,     c, 
  A*C,   A*c, 
  a*C,   a*c, 
  B*C,   B*c, 
  b*C,   b*c, 
A*B*C, A*B*c, 
A*b*C, A*b*c, 
a*B*C, a*B*c, 
a*b*C, a*b*c, 
B*a*C, B*a*c, 
b*A*C, b*A*c, 
b*a*C, b*a*c, 
\end{lstlisting}
	
\newpage
	
\begin{lstlisting}[language=Bash]
planar mod-2 words when k=5, 
with generators ordered: a,A,b,B,c,C,d,D
========================================================================
a, A, 
------------------------------------------------------------------------
b*a, b*A, B*a, 
------------------------------------------------------------------------
c*b*a, c*b*A, c*B*a, C*b*a, C*b*A, 
------------------------------------------------------------------------
d*c*b*a, d*c*b*A, d*c*B*a, d*C*b*a, d*C*b*A, D*c*b*a, D*c*b*A, D*c*B*a, 
========================================================================
  b,   B, 

a*b, a*B, 
A*b, A*B, 
------------------------------------------------------------------------
    c*b,     c*B,     C*b, 

  a*c*b,   a*c*B,   a*C*b, 
  A*c*b,   A*c*B,   A*C*b, 

b*a*c*b, b*a*c*B, b*a*C*b, 
B*a*c*b, B*a*c*B,
------------------------------------------------------------------------  
      d*c*b,       d*c*B,       d*C*b,       D*c*b,       D*c*B, 

    a*d*c*b,     a*d*c*B,     a*d*C*b,     a*D*c*b,     a*D*c*B, 
    A*d*c*b,     A*d*c*B,     A*d*C*b,     A*D*c*b,     A*D*c*B, 

  b*a*d*c*b,   b*a*d*c*B,   b*a*d*C*b,   b*a*D*c*b,   b*a*D*c*B, 
  B*a*d*c*b,   B*a*d*c*B,                B*a*D*c*b,   B*a*D*c*B, 

c*b*a*d*c*b, c*b*a*d*c*B, c*b*a*d*C*b, c*b*a*D*c*b, c*b*a*D*c*B, 
C*b*a*d*c*b, C*b*a*d*c*B, C*b*a*d*C*b, 
========================================================================
    c,     C, 
 
  a*c,   a*C, 
  A*c,   A*C, 

  b*c,   b*C, 
  B*c,   B*C, 

a*b*c, a*b*C, 
a*B*c, a*B*C, 
A*b*c, A*b*C, 
A*B*c, A*B*C, 
 
b*a*c, b*a*C, 
b*A*c, b*A*C, 
B*a*c, B*a*C,
------------------------------------------------------------------------






        d*c,         d*C,         D*c, 

      a*d*c,       a*d*C,       a*D*c, 
      A*d*c,       A*d*C,       A*D*c, 

      b*d*c,       b*d*C,       b*D*c, 
      B*d*c,       B*d*C,       B*D*c, 

    a*b*d*c,     a*b*d*C,     a*b*D*c, 
    a*B*d*c,     a*B*d*C,     a*B*D*c, 
    A*b*d*c,     A*b*d*C,     A*b*D*c, 
    A*B*d*c,     A*B*d*C,     A*B*D*c, 

    b*a*d*c,     b*a*d*C,     b*a*D*c, 
    b*A*d*c,     b*A*d*C,     b*A*D*c, 
    B*a*d*c,     B*a*d*C,     B*a*D*c, 

    c*b*d*c,     c*b*d*C,     c*b*D*c, 
    C*b*d*c,     C*b*d*C, 

  a*c*b*d*c,   a*c*b*d*C,   a*c*b*D*c, 
  a*C*b*d*c,   a*C*b*d*C, 
  A*c*b*d*c,   A*c*b*d*C,   A*c*b*D*c, 
  A*C*b*d*c,   A*C*b*d*C, 

  c*b*a*d*c,   c*b*a*d*C,   c*b*a*D*c, 
  c*b*A*d*c,   c*b*A*d*C, 
  C*b*a*d*c,   C*b*a*d*C, 
  C*b*A*d*c,   C*b*A*d*C, 
 
b*a*c*b*d*c, b*a*c*b*d*C, b*a*c*b*D*c, 
b*a*C*b*d*c, b*a*C*b*d*C, 
B*a*c*b*d*c, B*a*c*b*d*C, B*a*c*b*D*c, 
========================================================================
d, D, 

a*d, a*D, 
A*d, A*D, 

b*d, b*D, 
B*d, B*D, 

c*d, c*D, 
C*d, C*D, 

a*b*d, a*b*D, 
a*B*d, a*B*D, 
A*b*d, A*b*D, 
A*B*d, A*B*D, 

a*c*d, a*c*D, 
a*C*d, a*C*D, 
A*c*d, A*c*D, 
A*C*d, A*C*D, 
 
b*a*d, b*a*D, 
b*A*d, b*A*D, 
B*a*d, B*a*D, 

b*c*d, b*c*D, 
b*C*d, b*C*D, 
B*c*d, B*c*D, 
B*C*d, B*C*D, 

c*b*d, c*b*D, 
c*B*d, c*B*D, 
C*b*d, C*b*D, 

a*b*c*d, a*b*c*D, 
a*b*C*d, a*b*C*D, 
a*B*c*d, a*B*c*D, 
a*B*C*d, a*B*C*D, 
A*b*c*d, A*b*c*D, 
A*b*C*d, A*b*C*D, 
A*B*c*d, A*B*c*D, 
A*B*C*d, A*B*C*D, 

a*c*b*d, a*c*b*D, 
a*c*B*d, a*c*B*D, 
a*C*b*d, a*C*b*D, 
A*c*b*d, A*c*b*D, 
A*c*B*d, A*c*B*D, 
A*C*b*d, A*C*b*D, 

b*a*c*d, b*a*c*D, 
b*a*C*d, b*a*C*D, 
b*A*c*d, b*A*c*D, 
b*A*C*d, b*A*C*D, 
B*a*c*d, B*a*c*D, 
B*a*C*d, B*a*C*D, 

c*b*a*d, c*b*a*D, 
c*b*A*d, c*b*A*D, 
c*B*a*d, c*B*a*D, 
C*b*a*d, C*b*a*D, 
C*b*A*d, C*b*A*D, 

b*a*c*b*d, b*a*c*b*D, 
b*a*c*B*d, b*a*c*B*D, 
b*a*C*b*d, b*a*C*b*D, 
B*a*c*b*d, B*a*c*b*D, 
B*a*c*B*d, B*a*c*B*D, 
\end{lstlisting}
	
\newpage
	
\begin{lstlisting}[language=Bash]
planar mod-2 words when k=5, 
with generators ordered: A,a,B,b,C,c,D,d
========================================================================
A, a, 
------------------------------------------------------------------------
B*a, b*A, b*a, 
------------------------------------------------------------------------
C*b*A, C*b*a, c*B*a, c*b*A, c*b*a, 
------------------------------------------------------------------------
D*c*B*a, D*c*b*A, D*c*b*a, d*C*b*A, d*C*b*a, d*c*B*a, d*c*b*A, d*c*b*a, 
========================================================================
  B,   b, 

A*B, A*b, 
a*B, a*b, 
------------------------------------------------------------------------
    C*b,     c*B,     c*b, 

  A*C*b,   A*c*B,   A*c*b, 
  a*C*b,   a*c*B,   a*c*b, 

         B*a*c*B, B*a*c*b, 
b*A*C*b, 
         b*a*c*B, b*a*c*b, 
------------------------------------------------------------------------
      D*c*B,       D*c*b,       d*C*b,       d*c*B,       d*c*b, 

    A*D*c*B,     A*D*c*b,     A*d*C*b,     A*d*c*B,     A*d*c*b, 
    a*D*c*B,     a*D*c*b,     a*d*C*b,     a*d*c*B,     a*d*c*b, 

  B*a*D*c*B,   B*a*D*c*b,                B*a*d*c*B,   B*a*d*c*b, 
                            b*A*d*C*b, 
  b*a*D*c*B,   b*a*D*c*b,                b*a*d*c*B,   b*a*d*c*b, 

                          C*b*A*d*C*b, 
                                       C*b*a*d*c*B, C*b*a*d*c*b, 
c*B*a*D*c*B, c*B*a*D*c*b, 
                          c*b*A*d*C*b, 
                                       c*b*a*d*c*B, c*b*a*d*c*b, 
========================================================================
    C,     c, 

  A*C,   A*c, 
  a*C,   a*c, 

  B*C,   B*c, 
  b*C,   b*c, 

A*B*C, A*B*c, 
A*b*C, A*b*c, 
a*B*C, a*B*c, 
a*b*C, a*b*c, 

B*a*C, B*a*c, 
b*A*C, b*A*c, 
b*a*C, b*a*c, 
------------------------------------------------------------------------

        D*c,         d*C,         d*c, 

      A*D*c,       A*d*C,       A*d*c, 
      a*D*c,       a*d*C,       a*d*c, 

      B*D*c,       B*d*C,       B*d*c, 
      b*D*c,       b*d*C,       b*d*c, 

    A*B*D*c,     A*B*d*C,     A*B*d*c, 
    A*b*D*c,     A*b*d*C,     A*b*d*c, 
    a*B*D*c,     a*B*d*C,     a*B*d*c, 
    a*b*D*c,     a*b*d*C,     a*b*d*c, 

    B*a*D*c,     B*a*d*C,     B*a*d*c, 
    b*A*D*c,     b*A*d*C,     b*A*d*c, 
    b*a*D*c,     b*a*d*C,     b*a*d*c, 

                 C*b*d*C,     C*b*d*c, 
    c*B*D*c, 
                 c*b*d*C,     c*b*d*c, 

               A*C*b*d*C,   A*C*b*d*c, 
  A*c*B*D*c, 
               A*c*b*d*C,   A*c*b*d*c, 
               a*C*b*d*C,   a*C*b*d*c, 
  a*c*B*D*c, 
               a*c*b*d*C,   a*c*b*d*c, 

               C*b*A*d*C,   C*b*A*d*c, 
               C*b*a*d*C,   C*b*a*d*c, 
  c*B*a*D*c, 
               c*b*A*d*C,   c*b*A*d*c, 
               c*b*a*d*C,   c*b*a*d*c, 

B*a*c*B*D*c, 
             B*a*c*b*d*C, B*a*c*b*d*c, 
             b*A*C*b*d*C, b*A*C*b*d*c, 
b*a*c*B*D*c, 
             b*a*c*b*d*C, b*a*c*b*d*c, 
========================================================================
        D,         d, 

      A*D,       A*d, 
      a*D,       a*d, 

      B*D,       B*d, 
      b*D,       b*d, 

      C*D,       C*d, 
      c*D,       c*d, 

    A*B*D,     A*B*d, 
    A*b*D,     A*b*d, 
    a*B*D,     a*B*d, 
    a*b*D,     a*b*d, 

    A*C*D,     A*C*d, 
    A*c*D,     A*c*d, 
    a*C*D,     a*C*d, 
    a*c*D,     a*c*d, 

    B*a*D,     B*a*d, 
    b*A*D,     b*A*d, 
    b*a*D,     b*a*d, 

    B*C*D,     B*C*d, 
    B*c*D,     B*c*d, 
    b*C*D,     b*C*d, 
    b*c*D,     b*c*d, 

    C*b*D,     C*b*d, 
    c*B*D,     c*B*d, 
    c*b*D,     c*b*d, 

  A*B*C*D,   A*B*C*d, 
  A*B*c*D,   A*B*c*d, 
  A*b*C*D,   A*b*C*d, 
  A*b*c*D,   A*b*c*d,  
  a*B*C*D,   a*B*C*d, 
  a*B*c*D,   a*B*c*d, 
  a*b*C*D,   a*b*C*d, 
  a*b*c*D,   a*b*c*d, 

  A*C*b*D,   A*C*b*d, 
  A*c*B*D,   A*c*B*d, 
  A*c*b*D,   A*c*b*d, 
  a*C*b*D,   a*C*b*d, 
  a*c*B*D,   a*c*B*d, 
  a*c*b*D,   a*c*b*d, 

  B*a*C*D,   B*a*C*d, 
  B*a*c*D,   B*a*c*d, 
  b*A*C*D,   b*A*C*d, 
  b*A*c*D,   b*A*c*d, 
  b*a*C*D,   b*a*C*d, 
  b*a*c*D,   b*a*c*d, 

  C*b*A*D,   C*b*A*d, 
  C*b*a*D,   C*b*a*d, 
  c*B*a*D,   c*B*a*d, 
  c*b*A*D,   c*b*A*d, 
  c*b*a*D,   c*b*a*d, 

B*a*c*B*D, B*a*c*B*d, 
B*a*c*b*D, B*a*c*b*d, 
b*A*C*b*D, b*A*C*b*d, 
b*a*c*B*D, b*a*c*B*d, 
b*a*c*b*D, b*a*c*b*d, 
\end{lstlisting}

\newpage	
\begin{landscape}

\end{landscape}
\doublespacing

    \chapter*{Notation}
\addcontentsline{toc}{chapter}{Notation}
\markboth{NOTATION}{}
\begin{longtable}{ c | l}
	$\nullset$							& null set \\
	$\setdiff{X}{Y}$					& $\set{x \in X: x \notin Y}$ \\
	$\powerset{X}$						& power set of $X$ \\
	\hline
	$\mathbb{R}$						& real numbers \\
	$\ints$								& integers \\
	$\posints$							& positive integers \\
	$\nonnegints$						& non-negative integers \\
	$\intsg{x}$							& $\{z \in \ints: z > x\}$\\
\end{longtable}

    \singlespacing
    \newpage
    \renewcommand\listfigurename{\sffamily Figure Index}
    \renewcommand{\cftfigfont}{\bfseries Figure }
    \addcontentsline{toc}{chapter}{Figure Index}
    \begin{multicols}{2}
    	\listoffigures
    \end{multicols}
    \newpage
    \renewcommand\listtablename{\sffamily Table Index}
    \renewcommand{\cfttabfont}{\bfseries Table }
    \addcontentsline{toc}{chapter}{Table Index}
    \begin{multicols}{2}
    \listoftables
    \end{multicols}
\end{document}